\documentclass[11pt,fullpage]{article}

\usepackage{amsmath}
\usepackage{amssymb}
\usepackage{amsfonts}
\usepackage{xr}

\numberwithin{equation}{section}
\newtheorem{thm}{Theorem}[section] 
\newtheorem{prp}[thm]{Proposition}
\newtheorem{lmm}[thm]{Lemma}   
\newtheorem{crl}[thm]{Corollary} 
\newtheorem{dfn}[thm]{Definition}

\def\under#1{\underline{#1}}
\def\ov#1{\overline{#1}}
\def\lra{\longrightarrow}
\def\Lra{\Longrightarrow}
\def\Llra{\Longleftrightarrow}

\def\l{\left}
\def\r{\right}
\def\lan{\langle}
\def\ran{\rangle}
\def\llan{\lan\lan}
\def\rran{\ran\ran}

\def\e_ref#1{(\ref{#1})}
\def\wh#1{\widehat{#1}}
\def\lr#1{\lan{#1}\ran}

\def\al{\alpha}
\def\be{\beta}
\def\de{\delta}
\def\ep{\epsilon}
\def\ga{\gamma}
\def\io{\iota}
\def\la{\lambda}
\def\na{\nabla}
\def\om{\omega}
\def\si{\sigma}
\def\th{\theta}
\def\ups{\upsilon}
\def\ve{\varepsilon}
\def\vp{\varpi}
\def\ze{\zeta}

\def\Ga{\Gamma}
\def\La{\Lambda}
\def\Om{\Omega}
\def\Si{\Sigma}
\def\Th{\Theta}

\def\eset{\emptyset}
\def\i{\infty}
\def\t{\top}
\def\lap{\bigtriangleup}
\def\P{\Bbb{P}^n}
\def\PP{\Bbb{P}^2}
\def\PPP{\Bbb{P}^3}

\def\inj{\text{inj}}   
\def\ev{\text{ev}} 
\def\coker{\text{coker}}
\def\ind{\text{ind}}
\def\codim{\text{codim}}
\def\RT{\text{RT}}

\begin{document}

\title{Enumerative vs. Symplectic Invariants and Obstruction Bundles}
\author{Aleksey Zinger
\thanks{Partially supported by NSF Graduate Research Fellowship
and NSF grant DMS-9803166}}
\date{\today}
\maketitle

\thispagestyle{empty}

\tableofcontents

\section{Introduction}

\subsection{Background and Results}

\noindent
Suppose $(\Si,j)$ is a nonsingular Riemann surface of genus $g\ge 2$
and $(V,J,\om)$ is a Kahler manifold of complex dimension~$n$.
If $\la\!\in\! H_2(V;\Bbb{Z})$, denote by ${\cal H}_{\Si,\la}(V)$ 
the set of simple $(J,j)$-holomorphic maps~$u$ 
from $\Si$ to $V$ such that  \hbox{$u_*[\Si]=\la$}.
Let 
$\mu=(\mu_1,\ldots,\mu_N)$
be an $N$-tuple of proper oriented submanifolds of $V$ such~that 
\begin{equation}\label{c_codim}
\codim~\mu\equiv\sum_{l=1}^{l=N}\hbox{codim}~\mu_l
=2\big(\lan c_1(V,J),\la\ran-n(g-1)+N\big).
\end{equation}
For many Kahler manifolds $(V,J,\om)$ and choice of constraints~$\mu$,
the cardinality of the set
\begin{equation}
\label{intro_1}
{\cal H}_{\Si,\la}(\mu)=
\big\{(\Si;y_1,\ldots,y_N;u)\!: u\!\in\!{\cal H}_{\Si,\la}(V);~
y_l\!\in\!\Si,~u(y_l)\!\in\!\mu_l~\forall l=1,\ldots,N\big\}
\end{equation}
is finite and depends only on the homology classes of $\mu_1,\ldots,\mu_N$.
The cardinality $|{\cal H}_{\Si,\la}(\mu)|$ of the set 
${\cal H}_{\Si,\la}(\mu)$ is then an {\it enumerative}
invariant of the complex manifold~$(V,J)$.
Such numbers for algebraic manifolds $(V,J)$, 
e.g.~the complex projective spaces $\P$,
have been of great interest in algebraic geometry for a long~time.\\

\noindent 
If $(V,\om,J)$ is a semipositive symplectic manifold, 
the {\it symplectic} invariant of~$(V,\om)$, 
$$\RT_{g,\la}(;\mu)\equiv \RT_{g,\la}(;\mu_1,\ldots,\mu_N)$$
of~\cite{RT}, is a well-defined integer.
Due to the two composition laws of~\cite{RT}, 
this symplectic invariant is often more readily computable 
than the enumerative invariant~$|{\cal H}_{\Si,\la}(\mu)|$.
In fact, all such symplectic invariants of~$\P$ are easily computable.
It is also shown in~\cite{RT} that the appropriately defined
genus-zero invariants of $\P$ agree with 
the corresponding symplectic invariants.
On the other hand, even for $\PP$ and for genus~one,
the two invariants are no longer equal.
In~\cite{I}, the difference
$$\RT_{1,\la}(\mu_1;\mu_2,\ldots,\mu_N)-
|{\cal H}_{\Si,\la}(\mu)|$$
is commuted for genus-one surfaces $\Si$ and
all projective spaces using an obstruction-bundle approach, 
first introduced by~\cite{T} in a very different setting.
In~\cite{Z2}, the difference
$$\RT_{2,\la}(;\mu)-|{\cal H}_{\Si,\la}(\mu)|$$
is commuted for genus-two surfaces~$\Si$ 
for $\PP$ and~$\PPP$ using a similar approach.\\

\noindent
The main results of this paper are {\it detailed} descriptions
of gluing pseudoholomorphic maps in symplectic geometry.
Theorem~\ref{si_str} concerns the adoption of
the obstruction-bundle idea of~\cite{T} into the symplectic setting.
Theorem~\ref{str_global} describes the local structure
of moduli spaces of holomorphic rational maps under certain 
regularity conditions.
In addition to providing the necessary analytic justification
for the computations in~\cite{I} and~\cite{Z2},
the two theorems lay foundations for 
comparing enumerative and symplectic invariants 
in many Kahler manifolds, in particular,
the complex homogeneous manifolds.
Furthermore, in~\cite{Z2}, the second theorem is used
to enumerate cuspidal rational curves in $\PP$ and 
two-component tacnodal rational curves in~$\PPP$.
The same arguments can be used to count curves with 
higher-order degeneracies.\\

\noindent
The author is grateful to T.~Mrowka for pointing out the paper~\cite{I},
encouraging the author to work out all of the analytic
issues arising in~\cite{I},
and sharing some of his expertise in applications of global analysis
over countless hours of conversations.
The author also thanks G.~Tian  for first introducing him
to Gromov's symplectic invariants and helping him understand~\cite{LT}.

\subsection{Summary}
\label{summary}

\noindent
Let $\La^{0,1}\pi_{\Si}^*T^*\Si\otimes\pi_V^*TV\!\lra\!\Si\!\times\! V$
denote the bundle of $(J,j)$-antilinear homomorphisms from
$\pi_{\Si}^*T^*\Si$ to~$\pi_V^*TV$.
If 
$$\nu\in\Ga\big(\Si\times V;\La^{0,1}\pi_{\Si}^*T^*\Si\otimes\pi_V^*TV
\big),$$ 
denote by ${\cal M}_{\Si,\nu,\la}$ the set of all smooth maps~$u$ 
from $\Si$ to $\P$ such that $u_*[\Si]\!=\!\la$ and 
$\bar{\partial}u|_z\!=\!\nu|_{(z,u(z))}$ for all $z\!\in\!\Si$.
If $\mu$ is an $N$-tuple of constraints as above,~put
$${\cal M}_{\Si,\nu,\la}(\mu)=
\big\{(\Si;y_1,\ldots,y_N;u)\!: u\!\in\!{\cal M}_{\Si,\nu,\la};~
y_i\!\in\!\Si,~u(y_i)\!\in\!\mu_i ~\forall i=1,\ldots,N\big\}.$$
If $(V,\om,J)$ is semipositive,  for generic choices of $\nu$ and $\mu$, 
${\cal M}_{\Si,\nu,\la}$ is a smooth finite-dimensional oriented manifold, 
and ${\cal M}_{\Si,\nu,\la}(\mu)$ 
is a zero-dimensional finite submanifold of 
${\cal M}_{\Si,\nu,\la}\times\Si^N$, 
whose cardinality (with sign) depends only the homology
classes of $\la_1,\ldots,\la_N$; see~\cite{RT}.
The symplectic invariant $\RT_{g,\la}(;\mu)$ is the signed
cardinality of the set~${\cal M}_{\Si,\nu,\la}(\mu)$.\\

\noindent
If $\|\nu_i\|_{C^0}\lra 0$ and 
   $(\Si;\under{y}_i;u_i)\!\in\! {\cal M}_{\Si,\nu_i,\la}(\mu)$, 
then a subsequence of $\{(\Si;\under{y}_i,u_i)\}_{i=1}^{\i}$ must converge 
in the Gromov topology to one of the following:\\
(1) an element of ${\cal H}_{\Si,\la}(\mu)$;\\
(2) $(\Si_T;\under{y};u)$, where $\Si_T$ is a bubble tree of 
$S^2$'s attached to $\Si$ with marked points $y_1,\ldots,y_N$,
and \hbox{$u\!:\Si_T\lra V$} is a holomorphic map such that
$u(y_i)\!\in\!\mu_i$ for $i=1,\ldots,N$, and\\
(2a) $u|\Si_T$ is simple and the tree contains at least one $S^2$;\\
(2b) $u|\Si_T$ is multi-covered;\\
(2c) $u|\Si_T$ is constant and the tree contains at least one $S^2$.\\
This convergence statement says that for all $t$ sufficiently small
all the elements of ${\cal M}_{\Si,t\nu,\la}(\mu)$
lie near one of the spaces described by (1)-(2c).
In many practical applications, it is easy to show that
there is a bijection between
the elements of ${\cal H}_{\Si,\la}(\mu)$
and the nearby elements of ${\cal M}_{\Si,t\nu,\la}(\mu)$;
see Proposition~\ref{holom_prp}.
This is the case for all projective spaces,
provided $\la$ is a sufficiently high multiple of the~line
(depending on the genus~$g$).
In~\cite{I} and~\cite{Z2}, Cases (2a) and (2b) do not occur,
but they may have to be considered when dealing
with higher-dimensional projective spaces
or higher genera.
Thus, the convergence statement implies that
if the signed cardinality of ${\cal M}_{\Si,t\nu,\la}(\mu)$
is $\RT_{g,\la}(;\mu)$ for all $t\!>\!0$ sufficiently small,
the number of elements of ${\cal M}_{\Si,t\nu,\la}(\mu)$
that lie near the spaces described by (2) is 
exactly \hbox{$\RT_{2,\la}(;\mu)\!-\!|{\cal H}_{\Si,\la}(\mu)|$}.
The goal of this paper is to describe this difference in terms
of the spaces of holomorphic maps themselves,
which can be viewed as an {\it enumerative} object,
rather than a {\it symplectic} one.
We do need to assume that certain spaces of holomorphic maps
are smooth, but they do not need to have the expected dimension.\\

\noindent
While there is a very good understanding of what constitutes 
a stable map, 
there is little in a way of commonly accepted notation  
for stable maps and various spaces of stable maps.
In Section~\ref{bubble_sect}, we recall the definition of bubble 
or stable maps as well as set up analytically convenient notation.
Our notation for bubble maps evolved from that
of D.~McDuff's lectures at Harvard.
We also define various spaces of bubble maps 
and restate the definition of the Gromov topology 
on the set of all bubble maps in our notation.\\

\noindent
In Section~\ref{gluing}, we describe a gluing construction
and an obstruction bundle setup.
One difficulty in the gluing procedure is that a certain
operator has eigenvalues that tend to zero as the gluing parameter
tends to zero, but then disappear as the gluing parameter hits zero.
This is not really dealt with in \cite{I}, but there are now
several standard approaches to this problem.
We use the modified Sobolev norms of \cite{LT}, redefined  
in the notation of  Section~\ref{bubble_sect}.
The main goal of Section~\ref{gluing} is to describe a correspondence 
between the elements of ${\cal M}_{\Si,t\nu,\la}(\mu)$ lying near the maps
of type (2) and the zero set of a map between two bundles;
see Theorem~\ref{si_str}.
The domain and target bundles 
will be the bundle of gluing parameters
and the obstruction bundle, respectively.
We also give a local analytic description 
of spaces of stable rational maps into $V$ under certain assumptions.
\hbox{If $V\!=\!\P$}, these spaces are very familiar in algebraic geometry,
but the description of Theorem~\ref{str_global} is used for a transition 
from the analysis of the gluing problem to 
the topology of moduli spaces of stable rational maps in~\cite{Z2}.
\\

\noindent
Section~\ref{tech_sec} contains proofs of continuity,
injectivity, and surjectivity of the gluing maps.
These are usually omitted in the literature, 
but in the given case 
one has to choose the obstruction bundle setup carefully
to ensure that these properties of the gluing map actually hold.
In particular, Section~\ref{tech_sec} contains
what~\cite{LT} may mean by ``asymptotic analysis near the nodes,''
which they omit.
The Appendix deals with even more technical details of the analysis.

\subsection{Fundamental Notation}
\label{notation}

\noindent
In this subsection, we collect the most frequently used
combinatorial and analytic notation.

\begin{dfn}
\label{rooted_tree}
A finite partially ordered set $I$ is a \under{linearly ordered set}
if for all \hbox{$i_1,i_2,h\!\in\! I$} such that $i_1,i_2\!<\!h$, 
either $i_1\!\le\! i_2$ \hbox{or $i_2\!\le\! i_1$.}\\
A linearly ordered set $I$ is a \under{rooted tree} if
$I$ has a unique minimal element, i.e.
there exists \hbox{$\hat{0}\!\in\! I$} such that $\hat{0}\!\le\! h$ 
for {all $h\!\in\! I$.}\\
If $I$ and $I'$ are linearly ordered sets, 
a bijection $\phi\!: I\!\lra\! I'$ is an 
\under{isomorphism of linearly ordered sets} if
for all $h,i\!\in\! I$, $i\!<\!h$ if and only if $\phi(i)\!<\!\phi(h)$.
\end{dfn}

\noindent
Let $I$ be a linearly ordered set.
We denote the subset of the non-minimal elements of $I$ by 
$\hat{I}$,~i.e.
$$\hat{I}=\l\{h\!\in\! I\!: i<h\hbox{~for some~}i\!\in\! I\r\}.$$
For every $h\!\in\!\hat{I}$,  the set $\{i\!\in\! I\!: i<h\}$ 
has a unique maximal element~$\io_h$, i.e.
$$\io_h<h \quad\hbox{and}\quad
i\le\io_h\hbox{~for all~}i\!\in\! I\hbox{~s.t.~}i<h.$$
For reasons made clear in Subsection~\ref{bubble_trees},
$\io\!:\hat{I}\lra I$ will be called the attaching map of~$I$.
It~is clear from Definition~\ref{rooted_tree} that 
$I$ has a unique splitting $I=\bigsqcup\limits_{k\in K}I_k$
such that $I_k\subset I$ is a~rooted tree.
The attaching map of~$I$ restricts to 
the attaching map of each $I_k$, which will still be denoted by~$\io$.
\\

\noindent
Let $I$ be a rooted tree.
We denote the unique minimal element of $I$ by $\hat{0}_I$, 
or simply by $\hat{0}$ if there is no ambiguity.
If $I^*$, $I_*$, and $I*$ are rooted trees, 
we will write $\hat{I}^*$, $\hat{I}_*$, and $\hat{I}*$
for  $\wh{I^*}$, $\wh{I_*}$, and $\wh{I*}$, respectively;
here $*$ denotes any string of symbols.
If $i\!\in\! I$, let 
$$D_iI=\{h\!\in\! I\!: h>i\},\qquad
\bar{D}_iI=D_iI\cup\{i\}.$$
Every rooted tree $I$ has a number of subsets that are rooted trees;
the subsets $\bar{D}_iI$ are one example.
If $H$ is a subset of $I$, the set
$$I^{(H)}\equiv\l\{i\!\in\! I\!: i\not> h~\forall h\in H\r\}$$
is also a rooted tree.
If $i\!\in\! I$, denote $I^{(\{i\})}$ by~$I^{(i)}$.
If $H$ is a subset of $\hat{I}$, let
$$I^H=\l\{i\!\in\! I\!: i\not\ge h~\forall h\in H\r\},\quad
I(H)=H\cup\{\hat{0}_I\}.$$
If $h\!\in\!\hat{I}$,  denote $I^{\{h\}}$ by~$I^h$.\\

\noindent
If $M_1$ and $M_2$ are two sets, let $M_1+M_2$ be the disjoint union
of $M_1$ and~$M_2$.
Finally, if $N$ is a nonnegative integer, let \hbox{$[N]=\{1,\ldots,N\}$}.
\\

\noindent
We now introduce some analytic notation.
Let $\be\!:\Bbb{R}\lra[0,1]$ be a smooth function such~that
\begin{equation}\label{cutoff_fun}
\be(t)=\begin{cases}
0,&\hbox{if~}t\le 1;\\
1,&\hbox{if~}t\ge 2,
\end{cases}
\qquad\hbox{and}\qquad
\be'(t)>0~~\hbox{if~}t\!\in\!(1,2).
\end{equation}
If $r>0$, let $\be_r\!\in\! C^{\i}(\Bbb{R};\Bbb{R})$ be given by 
$\be_r(t)=\be(r^{-\frac{1}{2}}t)$.
Note that
\begin{equation}\label{cutoff_fun2}
\hbox{supp}(\be_r)=[r^{\frac{1}{2}},2r^{\frac{1}{2}}],~~
\|\be_r'\|_{C^0}\le C_{\be}r^{-\frac{1}{2}},\hbox{~~and~~}
\|\be_r''\|_{C^0}\le C_{\be}r^{-1}.
\end{equation}
Throughout the paper, $\be$ and $\be_r$ will refer to these 
smooth bump functions.\\

\noindent
Let $q_N,q_S\!: \Bbb{C}\lra S^2\subset\Bbb{R}^3$ be the stereographic
projections mapping the origin in $\Bbb{C}$ to the north and 
south poles, respectively. Explicitly,
\begin{equation}
\label{stereo_pr}
q_N(z)=\l(\frac{2z}{1+|z|^2},\frac{1-|z|^2}{1+|z|^2}\r)
\in\Bbb{C}\times\Bbb{R},~~~
q_S(z)=\l(\frac{2z}{1+|z|^2},\frac{-1+|z|^2}{1+|z|^2}\r).
\end{equation}
We denote the south pole of $S^2$, i.e. the point $(0,0,-1)\!\in\!\Bbb{R}^3$,
by $\i$. Let
\begin{equation}
\label{special_vect}
e_{\i}=(0,0,1)=dq_S\Big|_0\l(\frac{\partial}{\partial s}\r)\in T_{\i}S^2,
\end{equation}
where we write $z=s+it\in\Bbb{C}$.
We identify $\Bbb{C}$ with  $S^2-\{\i\}$ via the map $q_N$.
If \hbox{$x,z\!\in\! S^2\!-\!\{\i\}$}, define $\phi_xz\!\in\!\Bbb{C}$ by
\begin{equation}
\label{exp_def1}
\phi_xz=z-x\equiv q_N^{-1}(z)-q_N^{-1}(x).
\end{equation}
Note that the map $\phi_x\!: S^2-\{\i\}\lra\Bbb{C}$ is biholomorphic.
If $g$ is a Riemannian metric on a Riemann surface~$(\Si,j)$
of positive genus, $x\!\in\!\Si$ and $v\!\in\! T_x\Si$,
we write $\exp_{g,x}v\!\in\!\Si$ for the exponential of~$v$
defined with respect to the Levi-Civita connection of~$g$.
Let $\inj_gx$ denote the corresponding injectivity radius at~$x$
and $d_g$ the distance function.
\hbox{If $x,z\!\in\!\Si$} are such that \hbox{$d_g(x,z)<\inj_gx$}, 
define \hbox{$\phi_{g,x}z\!\in\! T_x\Si$ by}
\begin{equation}
\label{exp_def2}
\exp_{g,x}\phi_{g,x}z=z,~~~|\phi_{g,x}z|_{g,x}<\inj_gx.
\end{equation}
Note that if $g$ is flat on a neighborhood $U$ of $x$ in 
$\Si$, then $\phi_{g,x}|U$ is holomorphic.\\

\noindent
Let $g_V$ be the Kahler metric of $(V,J,\om)$.
Denote the corresponding Levi-Civita connection,
exponential map, and distance function by $\na^V$,
$\exp_V$ and  $d_V$, respectively.
For every \hbox{$\la\!\in\! H_2(V;\Bbb{Z})$}, let
\hbox{$|\la|=\lan \om,\la\ran$}. 
The number $|\la|$ is the $g_V$-energy of 
any element of~${\cal H}_{\Si,\la}$; see~\cite{MS}.
By rescaling $\om$, it can be assumed that $|\la|\ge 1$,
whenever $\la\neq0$ \hbox{and ${\cal H}_{S^2,\la}\neq\eset$}.
If $g$ is any Kahler metric on $(V,J)$, denote the corresponding 
Levi-Civita connection, exponential map, distance function, 
injectivity radius,  and
the parallel transport along the geodesic \hbox{for $X\!\in\! TV$}
by $\na^g$, $\exp_g$, $d_g$, $\inj_g$, and $\Pi_{g,X}$, respectively.
If $q\!\in\! V$ and $\de\!\in\!\Bbb{R}$, let
\hbox{$B_g(q,\de)=\big\{q'\!\in\!V\!: d_g(q,q')\le\de\big\}$}.
In our construction, we allow $g$ vary in a smooth family.
Without causing any additional difficulty in the gluing construction,
consideration of such families simplifies computations in
specific cases such as in~\cite{Z2}.
If $(S,j)$ is a smooth Riemann surface and $u\in C^{\i}(S;V)$, put
\begin{equation*}\begin{split}
\Ga(u)=\Ga(S;u^*TV),&\quad\Ga^1(u)=\Ga(S;T^*S\otimes u^*TV);\\
\Ga^{0,1}(u)=\Ga(S;\La^{0,1}T^*S\otimes u^*TV),&\quad
\bar{\partial}u=\frac{1}{2}\big(du+J\circ du\circ j\big)\in\Ga^{0,1}(u).
\end{split}\end{equation*}
We denote by $D_V$ and $D_g$ the linearization of
$\bar{\partial}$-operator with respect to the metrics $g_V$ and $g$
on~$V$, respectively.
Since both metrics are Kahler, $D_V$ and $D_g$ commute with $J$
and have no zeroth-order term; see~\cite{Z1}.\\

\noindent
It should be mentioned that it is not essential for 
the main gluing construction  described in this paper that 
$(V,J,g)$ is Kahler or even symplectic.
If $(V,J,g)$ is not Kahler, we would need to choose an orientation
on certain spaces of holomorphic maps and
take the induced orientation on the cokernel bundle;
see Subsection~\ref{str_sub1}.
Dropping the Kahler assumption would have almost no effect on the analysis,
but would slightly complicate the notation.

\section{Spaces of Bubble Maps}
\label{bubble_sect}

\subsection{Bubble Trees}
\label{bubble_trees}

\noindent
Let $S$ be either the Riemann sphere $S^2$ or
a smooth Riemann surface $\Si$ of genus at least~$2$.
Allowing the genus-one case would lead to somewhat more complicated notation,
but would have no effect on the analysis done in Section~\ref{gluing}.
Denote by $S^*$ the open subset \hbox{$S-\{\i\}$} if $S=S^2$ 
and $S$ itself \hbox{if $S=\Si$}.

\begin{dfn}
\label{bubble_tree_dfn}
A \under{bubble tree based on $S$} is a tuple 
$\t=\big(S,I;x\big)$, where\\
(1) $I$ is a rooted tree and $x\!: \hat{I}\!\lra\!S\cup S^2$ is a map;\\
(2) $x_h\!\in\! S^*$ if $\io_h\!=\!\hat{0}$ and 
$x_h\!\in\! S^2-\{\i\}$ otherwise;\\
(3) if $h_1\neq h_2$ and $\io_{h_1}=\io_{h_2}$, $x_{h_1}\neq x_{h_2}$.
\end{dfn}

\noindent
Given a bubble tree $\t$ as above, let $\Si_{\t}$ be the complex curve
$$\Si_{\t}=\Big((\{\hat{0}\}\times S)\sqcup
\bigsqcup_{h\in\hat{I}}(\{h\}\times S^2)\Big)\Big/\!\sim,$$
where $(h,\i)\sim (\io_h,x_h)$ for $h\!\in\!\hat{I}$.
The subset $\Si_{\t,\hat{0}}\equiv\{\hat{0}\}\times S$ of $\Si_{\t}$
will be called the {\it principal component} of $\t$ or $\Si_{\t}$.
If $h\!\in\!\hat{I}$, let $\Si_{\t,h}=\{h\}\times S^2$.
For every $i\!\in\! I$, $\Si_{\t,i}$ will be called 
the $i${\it{th} bubble component} of $\t$ or $\Si_{\t}$
or simply a {\it bubble component}.
Let $\Si_{\t,i}^*$ and $\Si_{\t}^*$ denote the open subsets of 
smooth points of $\Si_{\t,i}$ and $\Si_{\t}$, respectively,  i.e.
\begin{equation*}
\label{bubble_trees_e1}
\Si_{\t,i}^*=
\begin{cases}
\Si_{\t,i}-\{(i,\i)\}-\{(i,x_h): \io_h=i\},&\hbox{if~}i\in\hat{I};\\
S-\{(\hat{0},x_h): \io_h=\hat{0}\},&\hbox{if~}i=\hat{0};
\end{cases}~~~~~~~~
\Si_{\t}^*=\bigcup_{i\in I}\Si_{\t,i}^*.
\end{equation*}
The complement of $\Si_{\t,i}^*$ in $\Si_{\t}$ are the {\it singular
points} of~$\Si_{\t}$. \\
 
\noindent
For every $i\!\in\! I$, let 
$\t^{(i)}\!=\!\big(S,I^{(i)};x|\hat{I}^{(i)}\big)$.
Similarly, for each $h\!\in\!\hat{I}$, let 
\hbox{$\t^h\!=\!\big(S,I^h;x|\hat{I}^h\big)$.}
These tuples are again bubble trees based on~$S$.
The complex curve~$\Si_{\t^{(i)}}$ is
obtained from $\Si_{\t}$ by dropping
all bubble components descendent from the $i$th bubble component.
The curve $\Si_{\t^h}$ is obtained by dropping
the $h$th bubble component along with 
all bubble components descendent from~it.\\

\noindent
If $S\!=\!S^2$ and $h\!\in\!\hat{I}$, we denote the map 
$\phi_{x_h}$ defined in \e_ref{exp_def1} by $\phi_{\t,h}$.
If $z\!\in\!\Si_{\t,i}$, put
\begin{equation}
\label{special_dst1}
r_{\t,h}(z)=\begin{cases}
\l|\phi_{\t,h}z\r|,&\hbox{if~}i=\io_h\hbox{~and~}z\neq\i;\\
100,&\hbox{otherwise}.
\end{cases}
\end{equation}
If $\de>0$, let $B_{\t,h}(\de)=\{z\!\in\!\Si_{\t}\!: r_{\t,h}(z)<\de\}$.
Put
\begin{equation}
\label{inj_radius1}
r_{\t}=\min_{h\in\hat{I}}
\Big(\l|q_S^{-1}(x_h)\r|,\min\{r_{\t,h}(\io_l,x_l): l\neq h\}\Big).
\end{equation}
If $S\!=\!\Si$ and $h\!\in\!\hat{I}$ is such that $\io_h\!\in\!\hat{I}$,
we again let $\phi_{\t,h}$ denote the function $\phi_{x_h}$ of 
\e_ref{exp_def1} and define $r_{\t,h}$ and $B_{\t,h}(\de)$ as above.
If $g$ is a Riemannian metric on $\Si$, $\io_h=\hat{0}$, and 
$z\!\in\!\Si_{\t,i}$, put
\begin{equation}
\label{special_dst2}
r_{\t,g,h}(z)=\begin{cases}
d_g(x_h,z),&\hbox{if~}i=\hat{0};\\
100,&\hbox{otherwise}.
\end{cases}
\end{equation}
We denote by $\phi_{\t,g,h}$ the function $\phi_{g,x_h}$ of \e_ref{exp_def2}
and by $B_{\t,g,h}(\de)$ the ball $B_g(x_h,\de)$.
Put
\begin{equation}
\label{inj_radius2}
r_{\t}g=\min\Big(
\min_{\io_h=\hat{0}}\{r_{\t,g,h}(\io_l,x_l): l\neq h\},
\min_{\io_h\neq\hat{0}}\l(\l|q_S^{-1}(x_h)\r|,
        \min\{r_{\t,h}(\io_l,x_l): l\neq h\}\r)\Big).
\end{equation}
We say $g$ is a $\t$-{\it{admissible Riemannian metric}} on $\Si$
if there exists $\de>0$ such that for all $h\!\in\!\hat{I}$ with 
$\io_h=\hat{0}$
the metric $g$ is flat on~$B_{\t,g,h}(\de)$.

\subsection{The Basic Gluing Construction}
\label{basic_gluing}

\noindent
In this subsection, we describe a gluing construction on $\t$ 
which is the basis of all the other gluing constructions
in this paper.
Lemma~\ref{basic_gluing_lmm} 
plays a very important role in the next section
and in the explicit computations of~\cite{Z2}.\\

\noindent
Let $\t=\big(S,I;x\big)$  be a bubble tree.
If $h\!\in\!\hat{I}$, put
\begin{equation}\label{glue_data1}
F_{h,\t}^{(0)}=\begin{cases}
\Bbb{C},&\hbox{if~}x_h\in S^2;\\
T_{x_h}\Si,&\hbox{if~}x_h\in\Si,
\end{cases}\qquad
F_{\t}^{(0)}=\bigoplus\limits_{h\in\hat{I}}F_{h,\t}^{(0)}.
\end{equation}
If $S\!=\!S^2$, for any $\de>0$, put
$$F_{\t,\de}^{(0)}=\Big\{\ups=(\t,v_{\hat{I}})\!: 
v_{\hat{I}}\!\in\! F_{\t}^{(0)},
|\ups|\equiv\sum_{h\in\hat{I}}|v_h|<\de\Big\}.$$
Let $\de_{\t}\!\in\!(0,1)$ be such that $8\de_{\t}^{\frac{1}{2}}<r_{\t}$.
If $S=\Si$ and $g$ is an admissible metric on $\Si$, put
$$F_{\t,g,\de}^{(0)}=\Big\{\ups=(\t,v_{\hat{I}})\!: 
v_{\hat{I}}\!\in\! F_{\t}^{(0)},
|\ups|_g\equiv\sum_{i_h=\hat{0}}|v_h|_g+
\sum_{i_h\neq\hat{0}}|v_h|<\de\Big\},$$
where $|v_h|_g=|v_h|_{g,x_h}$.
Let $\de_{\t}g\!\in\!(0,1)$ be such that 
$8(\de_{\t}g)^{\frac{1}{2}}\!<\!r_{\t}g$  
and the metric $g$ is flat on $B_g\big(x_h,4(\de_{\t}g)^{\frac{1}{2}}\big)$ 
for all $h\!\in\!\hat{I}$ with $\io_h\!=\!\hat{0}$.
We now construct a family of bubble trees parameterized by
$F_{\t,\de}^{(0)}$ with  $\de\!=\!\de_{\t}$ if $S\!=\!S^2$ and
by $F_{\t,g,\de}^{(0)}$ with  $\de\!=\!\de_{\t}g$ \hbox{if $S\!=\!\Si$.}\\

\noindent
First,  for every $h\!\in\!\hat{I}$ and  $v_h\!\in\! F_{\t,h}^{(0)}$ with 
$$|v_h|\!\in\!(0,\de)\hbox{~~if~}x_h\!\in\! S^2
\quad\hbox{and}\quad 
|v_h|_g\!\in\!(0,\de)\hbox{~~if~}x_h\!\in\!\Si,$$
we define local stretching maps 
$$q_{h,(x_h,v_h)}\!:\Si_{\t^{(\io_h)}}\lra\Si_{\t^{(h)}}
\hbox{~~if~}x_h\!\in\! S^2
\quad\hbox{and}\quad
q_{g,h,(x_h,v_h)}:\Si_{\t^{(\io_h)}}\lra\Si_{\t^{(h)}}
\hbox{~~if~}x_h\!\in\!\Si$$
as follows.
If $x_h\!\in\! S^2$, let 
$p_{h,(x_h,v_h)}\!:
B_{\t,h}(\de^{\frac{1}{2}})\lra\Bbb{C}\cup{\i}$ be given~by
\begin{equation}
\label{basic_gluing_e2}
p_{h,(x_h,v_h)}(z)=\l(1-\be_{|v_h|}(2|\phi_{\t,h}z|)\r)
\ov{\l(\frac{v_h}{\phi_{\t,h}z}\r)}.
\end{equation}
Define $q_{h,(x_h,v_h)}\!:\Si_{\t^h}\lra\Si_{\t^{(h)}}$ by
\begin{equation}
\label{basic_gluing_e3}
q_{h,(x_h,v_h)}(z)=\begin{cases}
\big(h,q_S(p_{h,(x_h,v_h)}(z))\big),&
              \hbox{if~}
\frac{1}{2} \le |v_h|^{-\frac{1}{2}}r_{\t,h}(z)\le 1;\\
\l(\io_h,\phi_{\t,h}^{-1}\l(\be_{|v_h|}(|\phi_{\t,h}z|)(\phi_{\t,h}z)\r)\r),
    &\hbox{if~}1\le |v_h|^{-\frac{1}{2}}r_{\t,h}(z)\le 2;\\
z,& \hbox{otherwise}.
\end{cases}
\end{equation}
Note that $q_{h,(x_h,v_h)}$ is smooth and is a diffeomorphism, except on 
the circle $r_{\t,h}(z)=|v_h|^{\frac{1}{2}}$ in~$\Si_{\t,\io_h}$.
The map $q_{h,(x_h,v_h)}$ stretches $B_{\t,h}(|v_h|^{\frac{1}{2}})$ around
the sphere $\Si_{\t,h}$.
If $x_h\!\in\!\Si$, similarly to the above, let 
   $p_{g,h,(x_h,v_h)}\!:B_{\t,g,h}(\de^{\frac{1}{2}})\lra\Bbb{C}\cup{\i}$ 
be given by
\begin{equation}
\label{basic_gluing_e2b}
p_{g,h,(x_h,v_h)}(z)=\l(1-\be_{|v_h|_g}(2|\phi_{\t,g,h}z|_g)\r)
\ov{\l(\frac{v_h}{\phi_{\t,g,h}z}\r)}.
\end{equation}
Note that the ratio $v_h/\phi_{\t,g,h}z$   is well-defined as an extended
complex number, since $T_{x_h}\Si$ is one-dimensional and $v_h\neq0$.
Define $q_{g,h,(x_h,v_h)}\!:\Si_{\t^h}\lra\Si_{\t^{(h)}}$ by
\begin{equation}
\label{basic_gluing_e3b}
q_{g,h,(x_h,v_h)}(z)=\begin{cases}
\l(h,q_S(p_{g,h,(x_h,v_h)}(z))\r),&
              \hbox{if~}
      \frac{1}{2}\le|v_h|_g^{-\frac{1}{2}}r_{\t,g,h}(z)\le 1;\\
\l(\io_h,\phi_{g,\t,h}^{-1}\l(\be_{|v_h|_g}(|\phi_{\t,g,h}z|_g)
(\phi_{\t,g,h}z)\r)\r),
       &\hbox{if~}1\le |v_h|_g^{-\frac{1}{2}}r_{\t,g,h}(z)\le2;\\
z,& \hbox{otherwise}.
\end{cases}
\end{equation}
Similarly to the case $x_h\!\in\! S^2$, $q_{g,h,(x_h,v_h)}$ is smooth and
is a diffeomorphism, except on the circle 
\hbox{$r_{\t,g,h}(z)=|v_h|_g^{\frac{1}{2}}$} in $\Si_{\t,i_h}$.\\

\noindent
If $S\!=\! S^2$, for every $h\!\in\! I$ and $\ups\!\in\! F_{\t,\de}^{(0)}$, 
we now define a bubble tree $\t_h(\ups)$
and a smooth map $q_{\ups,h}\!:\Si_{\t_h(\ups)}\lra\Si_{\t^{(h)}}$.
Choose an ordering of $I$ consistent with its partial ordering.
If \hbox{$h\!=\!\hat{0}$}, we take $I_h(\ups)=\{\hat{0}\}$,
$\t_h(\ups)\!=\!\big(S,I_h(\ups);\big)$, and $q_{\ups,h}=Id_S$.
Suppose $h\!\neq\!\hat{0}$ and 
$$\t_{h-1}(\ups)=\big(S,\hat{I}_{h-1}(\ups);x_h(\ups)\big)$$
with $I_{h-1}(\ups)\subset I$.
If $v_h=0$, put 
$$I_h(\ups)=I_{h-1}(\ups)\cup\{h\},\hspace{3mm}
\big(\io_{h,l}(\ups),x_{h,l}(\ups)\big)=
\begin{cases}
\big(\io_{h-1,l}(\ups),x_{h-1,l}(\ups)\big),&\hbox{if~} l\!\in\! I_{h-1}(\ups);\\
q_{\ups,\io_h}^{-1}(\io_h,x_h),&\hbox{otherwise}.
\end{cases}$$
Let $q_{\ups,h}|\Si_{\t_{h-1}(\ups)}=q_{\ups,h-1}$
and $q_{\ups,h}(h,z)=(h,z)$.
If $v_h\neq 0$, let 
$$I_h(\ups)=I_{h-1}(\ups),\quad
\big(\io_{h,l}(\ups),x_{h,l}(\ups)\big)=
\big(\io_{h-1,l}(\ups),x_{h-1,l}(\ups)\big).$$
Take $q_{\ups,h}=q_{h,(x_h,v_h)}\circ q_{\ups,h-1}$.
Inductively this procedure defines a bubble tree 
$\t(\ups)=\t_{h^*}(\ups)$ based on $S$ and a smooth map 
$q_{\ups}=q_{\ups,h^*}\!: \Si_{\t(\ups)}\lra \Si_{\t}$, 
which is a diffeomorphism outside of \hbox{$|I-I(\ups)|$} disjoint circles,
where $h^*$ is the largest element of~$I$.
The resulting bubble tree and map are independent of
the choice of the extension of the partial ordering.
While the domains of the maps $q_{\ups,h}$ do depend on such a choice,
whenever we make use of the maps $q_{\ups,h}$ below,
the result will also be independent of the choice.
If $S\!=\!\Si$, for every $h\!\in\! I$ and $\ups\!\in\! F_{\t,g,\de}^{(0)}$, 
we define bubble tree $\t_{g,h}(\ups)$ and maps
$q_{g,\ups,h}\!:\Si_{\t_{g,h}(\ups)}\lra\Si_{\t^{(h)}}$
similarly to the above, but replacing $q_{h,(x_h,v_h)}$
by $q_{g,h,(x_h,v_h)}$ whenever $\io_h\!=\!\hat{0}$.
We let $\t_g(\ups)\!=\!\t_{g,h^*}(\ups)$ 
and \hbox{$q_{g,\ups}\!=\!q_{g,\ups,h^*}$}.
As before, $q_{g,\ups}$ is smooth and a diffeomorphism outside
of $|I-I(\ups)|$ disjoint circles.\\

\noindent
If $S\!=\!S^2$ and $v_h\!\neq\! 0$, put
\begin{equation}\label{ann_def1}\begin{split}
&A_{\ups,h}^+=q_{\ups,\io_h}^{-1}\Big(\big\{z\!\in\!\Si_{\t,\io_h}\!: 
|v_h|^{\frac{1}{2}}\le r_{\t,h}(z)\le 2|v_h|^{\frac{1}{2}}\big\}\Big);\\
&A_{\ups,h}^-=q_{\ups,\io_h}^{-1}\Big(\big\{z\!\in\!\Si_{\t,\io_h}\!: 
\frac{1}{2}|v_h|^{\frac{1}{2}}\le r_{\t,h}(z)\le |v_h|^{\frac{1}{2}}
\big\}\Big).
\end{split}\end{equation}
Note that $A_{\ups,h}^{\pm}\subset
  \Si_{\t(\ups),i^*_h(\ups)}$, where 
$$i^*_h(\ups)=\min\{i\!\in\! I\!: i<h\hbox{~and~}
v_{h'}\neq0\hbox{~if~}i<h'<h\}
=\max\{i\!\in\! I(\ups)\!: i<h\}.$$
If $S\!=\!\Si$ and $v_h\!\neq\! 0$, we similarly define
\begin{equation}\label{ann_def2}\begin{split}
&A_{g,\ups,h}^+=q_{g,\ups,\io_h}^{-1}\Big(\big\{z\!\in\!\Si_{\t,\io_h}\!: 
|v_h|_g^{\frac{1}{2}}\le r_{\t,g,h}(z)\le 2|v_h|_g^{\frac{1}{2}}
\big\}\Big);\\
&A_{g,\ups,h}^-=q_{g,\ups,\io_h}^{-1}\Big(\big\{z\!\in\!\Si_{\t,\io_h}\!: 
\frac{1}{2}|v_h|_g^{\frac{1}{2}}\le r_{\t,g,h}(z)\le 
|v_h|_g^{\frac{1}{2}}\big\}\Big),
\end{split}\end{equation}
where $|v_h|_g$ and $r_{\t,g,h}$  denote $|v_h|$ and $r_{\t,h}$
if $\io_h\!\in\!\hat{I}$.

\begin{lmm}
\label{basic_gluing_lmm}
If $S\!=\!S^2$, the map $q_{\ups}$ is holomorphic outside 
of the annuli $A_{\ups,h}^{\pm}$ with $v_h\!\neq\!0$.
\hbox{For such~$h$},
\begin{gather*}
\|dq_{h,(x_h,v_h)}\|_{C^0(q_{\ups,\io_h}(A_{\ups,h}^{\pm}))}\le C;\\
\bar{\partial}\big(q_{\ups}\circ q_{\ups,\io_h}^{-1}\big)
\big|_z=-2
|v_h|^{-\frac{1}{2}}\l(\frac{v_h}{\phi_{\t,h}z} \r)
 dq_S\big|_{p_{h,(x_h,v_h)}z}\circ
\partial\be\big|_{2|v_h|^{-\frac{1}{2}}\phi_{\t,h}z}
\circ d\phi_{\t,h}\big|_z
\quad\forall z\!\in\! q_{\ups,\io_h}(A_{\ups,h}^-),
\end{gather*}
where the norm is computed with respect to the standard metric on $S^2$, 
and $\be$ is a viewed as a function on $\Bbb{C}$ 
via the standard norm on $\Bbb{C}$.
If $S\!=\!\Si$, the map $q_{g,\ups}$ is holomorphic outside 
of the annuli $A_{g,\ups,h}^{\pm}$ with $v_h\!\neq\!0$.
\hbox{For such~$h$},
\begin{gather*}
\|dq_{g,h,(x_h,v_h)}\|_{C^0(q_{g,\ups,\io_h}(A_{g,\ups,h}^{\pm}))}\le C_g
\hbox{~~if~}\io_h=\hat{0};\quad
\|dq_{h,(x_h,v_h)}\|_{C^0(q_{g,\ups,\io_h}(A_{g,\ups,h}^{\pm}))}\le C
\hbox{~~if~~}\io_h\neq\hat{0};\\
\bar{\partial}\big(q_{g,\ups}\circ q_{g,\ups,\io_h}^{-1}\big)
\big|_z=-2
|v_h|^{-\frac{1}{2}}\l(\frac{v_h}{\phi_{\t,h}z} \r)
 dq_S\big|_{p_{h,(x_h,v_h)}z}\circ
\partial\be\big|_{2|v_h|^{-\frac{1}{2}}\phi_{\t,h}z}
\circ d\phi_{\t,h}\big|_z
\quad\forall z\!\in\! q_{g,\ups,\io_h}(A_{\ups,h}^-),
\end{gather*}
where we regard $\be$ as a function on $T_{x_h}\Si$ via 
the metric~$g$ and denote $\phi_{g,\t,h}$ by $\phi_{\t,h}$
if $\io_h\!=\!\hat{0}$. 
\end{lmm}

\noindent
{\it Proof:} The first statement in each of the two cases
is immediate from the construction.
The estimates on the differential of 
$q_{h,(x_h,v_h)}$ and $q_{g,h,(x_h,v_h)}$
follow from~\e_ref{cutoff_fun2}.
Suppose $S\!=\!\Si$, \hbox{$\io_h\!=\!\hat{0}$}, 
\hbox{$v_h\neq0$} and $z\!\in\! A_{g,\ups,h}^-$.
Since $q_{g,\ups}\!=\!q_{g,h,(x_h,v_h)}$ on $A_{g,\ups,h}^-$ and 
$q_S$ is anti-holomorphic, 
from \e_ref{basic_gluing_e2b} and~\e_ref{basic_gluing_e3b} we obtain
\begin{equation*}\begin{split}
\bar{\partial}q_{g,\ups}\big|_z&=
dq_S\big|_{p_{g,h,(x_h,v_h)}z}\circ \partial p_{g,h,(x_h,v_h)}\big|_z\\
&=-2|v_h|^{-\frac{1}{2}}_g\l(\frac{v_h}{\phi_{\t,g,h}z} \r)
dq_S\big|_{p_{g,h,(x_h,v_h)}z}\circ
\partial\be \big|_{2|v_h|^{-\frac{1}{2}}_g\phi_{\t,g,h}z}
\circ d\phi_{\t,g,h}\big|_z.
\end{split}\end{equation*}
The other cases are proved similarly, since 
$q_{g,\ups}\circ q_{g,\ups,\io_h}^{-1}=q_{h,(x_h,v_h)}$ on 
$q_{g,\ups,\io_h}(A_{g,\ups,h}^-)$
and a similar statement holds in the case $S=S^2$.

\subsection{Curves with Marked Points}
\label{curves}

\begin{dfn}
\label{curves_dfn}
If $M$ is a finite set, a 
\under{curve with $M$-marked points based on $S$} is a tuple
$${\cal C}=
\big(S,M,I;x,(j,y)\big),\qquad\hbox{where}~~$$
(1) $\t_{\cal C}\!\equiv\!\big(S,I;x\big)$ 
is a bubble tree \hbox{based on~$S$},
and $j\!:M\!\lra\! I$ and $y\!:M\!\lra\! S\cup S^2$
are maps;\\
(2) $j_l\!\in\! I$, $(j_l,y_l)\!\in\!\Si_{\t_{\cal C},j_l}^*$, 
and $y_l\neq\i$ for all $l\in M$;\\
(3) for any $l_1,l_2\!\in\! M$ with $l_1\neq l_2$ and $j_{l_1}=j_{l_2}$,
$y_{l_1}\neq y_{l_2}$.\\
The curve ${\cal C}$ is \under{stable} 
if $\big|\{h\!:\io_h=i\}\big|+
\big|\{l\!:j_l=i\}\big|\ge 2$ for all $i\!\in\!\hat{I}$
if $S\!=\! \Si$ and all \hbox{$i\!\in\! I$ if $S\!=\!S^2$}.
\end{dfn}

\noindent
Via the construction in Subsection~\ref{bubble_trees}, 
such a tuple ${\cal C}$ corresponds to a complex curve 
$\Si_{\cal C}\equiv\Si_{\t_{\cal C}}$ with marked points  
$\{(j_l,y_l)\}_{l\in M}$.
For each $i\!\in\! I$, we denote by 
$\Si_{{\cal C},i}$ and $\Si_{{\cal C},i}^*$ the surfaces 
$\Si_{\t_{\cal C},i}$ and $\Si_{\t_{\cal C},i}^*$, respectively.
\\

\noindent
With notation as above, for every $h\!\in\!\hat{I}$, let 
$F_{h,{\cal C}}^{(0)}$ and $F_{\cal C}^{(0)}$ denote the spaces 
$F_{h,\t_{\cal C}}^{(0)}$ and $F_{\t_{\cal C}}^{(0)}$, respectively.
If $S\!=\!S^2$, put
\begin{equation}\label{inj_curve1}
r_{\cal C}=\min\Big(r_{\t_{\cal C}},
\min_{l\in M}\big(\l|q_S^{-1}y_l\r|,
\min\{r_{\t_{\cal C},h}(j_l,y_l)\!:h\!\in\!\hat{I}\},
\min\{|\phi_{y_l}^{-1}y_h|\!:h\neq l,j_h=j_l\}\big)\Big).
\end{equation}
Let $\de_{\cal C}\!\in\!(0,1)$
be such that $16(|I|\!+\!|M|)\de_{\cal C}^{\frac{1}{2}}<r_{\cal C}$. 
If $\ups\!=\!({\cal C},v_{\hat{I}})$ with
$v_{\hat{I}}\!\in\! F_{\cal C}^{(0)}$ and $|\ups|<\de_{\cal C}$,
we now construct a curve ${\cal C}(\ups)$ with $M$-marked points as follows. 
Let 
$$\t(\ups)=\big(S,I(\ups);x(\ups)\big)
\quad\hbox{and}\quad q_{\ups}\!:\Si_{\t(\ups)}\lra\Si_{\cal C}$$
be the bubble tree and the smooth map defined in 
Subsection~\ref{basic_gluing}.
Then we take
$${\cal C}(\ups)=\big(S,M,I(\ups);x(\ups),(j(\ups),y(\ups))\big),$$
where $(j_l(\ups),y_l(\ups))\!\in\!\Si_{\t(\ups),j_l(\ups)}$ is defined by
$$q_{\ups}(j_l(\ups),y_l(\ups))=(j_l,y_l).$$
Similarly, if $S\!=\!\Si$ and $g$ is an admissible 
Riemannian metric on $\Si$, put
\begin{equation}\label{inj_curve2}\begin{split}
r_{\cal C}g=\min\Big( r_{\t_{\cal C}}g,
 \min_{j_l=\hat{0}}\big( \min_{\io_h=\hat{0}}\{r_{\t_{\cal C},g,h}(j_l,y_l)\},
 \min\{|\phi_{g,y_l}^{-1}y_h|_g: h\neq l,j_h=\hat{0}\}\big),~~~~~~~~~&\\
\min_{j_l\neq\hat{0}}\big( \l|q_S^{-1}y_l\r|,
\min_{\io_h\neq{0}}\{r_{\t_{\cal C},h}(j_l,y_l)\},
\min\{|\phi_{y_l}^{-1}y_h|:h\neq l,j_h=j_l \}\big)\Big).&
\end{split}\end{equation}
Let $\de_{\cal C}g\!\in\!(0,1)$ be such that 
$16(|I|\!+\!|M|)(\de_{\cal C}g)^{\frac{1}{2}}\!<\!r_{\cal C}g$ 
and  $g$ is flat in $B_g\big(x_h,8(\de_{\cal C}g)^{\frac{1}{2}}\big)$ 
for all $h\!\in\!\hat{I}$ with $\io_h\!=\!\hat{0}$.
If $\ups\!\in\! F_{\cal C}^{(0)}$ and $|\ups|_g<\de_{\cal C}g$,
we construct a curve ${\cal C}_g(\ups)$ with $M$-marked points
in the same way as above, but replacing $q_{\ups}$ and $\t(\ups)$
by $q_{g,\ups}$ \hbox{and $\t_g(\ups)$.}

\begin{dfn}
\label{curves_iso}
An~ \under{isomorphism of curves with $M$-marked points}~
${\cal C}=\big(S,M,I;x,(j,y)\big)$ and 
\hbox{${\cal C}'=\big(S,M,I';x',(j',y')\big)$},
is  a tuple of maps,
$$\phi_0: I\lra I',\hspace{5mm}
\phi_{1,\hat{0}}:S\lra S,\hspace{5mm}
\phi_{1,h}:S^2\lra S^2\hspace{2mm}\hbox{~for~}h\!\in\! I,
\qquad\hbox{where}$$
(1) $\phi_0$ is an isomorphism of the linearly ordered
sets $I$ and $I'$ and $\phi_0(j_l)=j_l'$ for all $l\!\in\! M$;\\
(2) $\phi_{1,i}$ is a biholomorphic map for all $i\!\in\! I$
and $\phi_{1,\hat{0}}$ is the identity map if $S=\Si$;\\
(3) $\phi_{1,i}(\i)=\i$ for all $i\!\in\! I$ if $S= S^2$ and 
for all $i\!\in\!\hat{I}$ if $S=\Si$;\\
(4) $\phi_{1,\io_h}(x_h)=x_{\phi_0(h)}'$ for $h\!\in\!\hat{I}$ and
$\phi_{1,j_l}(y_l)=y_l'$ for all $l\!\in\! M$.
\end{dfn}

\noindent
Such a set of maps corresponds to a continuous map
$\Si_{\cal C}\lra\Si_{{\cal C}'}$ that maps the $l$th marked point 
$(j_l,y_l)$ on $\Si_{\cal C}$ to the $l$th marked point $(j_l',y_l')$ on 
$\Si_{{\cal C}'}$ and is biholomorphic on each component of~$\Si_{\cal C}$.
Note that if ${\cal C}$ is stable,
${\cal C}$ has no nontrivial automorphisms.
Let $[{\cal C}]$ denote the equivalence class
of~${\cal C}$ in the set of all curves based on $S$
with marked points.
Denote by $\ov{\cal M}_{S,M}$ the set of all equivalence classes
of stable curves based on $S$ with $M$-marked points.
If $S=S^2$, $\ov{\cal M}_{S,M}$ can be identified with 
the moduli space $\ov{\cal M}_{0,|M|+1}$ of all stable rational curves
with $|M|+1$ marked points, or more canonically 
with the space $\ov{\cal M}_{0,M+\{\hat{0}\}}$ of all stable rational curves
with the marked points labeled by the set $M+\{\hat{0}\}$.
If $S=\Si$ has genus bigger than two and is generic,
$\ov{\cal M}_{S,M}$ is the closed subset of $\ov{\cal M}_{g,M}$
consisting of all stable curves of genus~$g$ with $M$-marked points
that have a fixed complex structure on the principal component.
If $S$ has genus two, $\ov{\cal M}_{S,M}$ is a double cover
of the corresponding set for $g=2$, since any smooth genus-two curve
has a holomorphic automorphism of order two; see~\cite[p254]{GH}.
The reason we require \hbox{$\phi_{1,\hat{0}}=Id_{\Si}$} is that
the symplectic invariant of~\cite{RT} 
disregards the automorphisms of~$\Si$.

\subsection{Bubble Maps}
\label{bubble_maps}

\begin{dfn}
\label{bubble_dfn}
A \under{$V$-valued bubble map} is a tuple
$b=\big(S,M,I;x,(j,y),u\big)$,
where\\
(1) $I$ is a linearly ordered set, which is a rooted tree if 
$S\!=\!\Si$ and $u\!:I\!\lra\!C^{\i}(S;V)\cup C^{\i}(S^2;V)$
is a map;\\
(2) if $I\!=\!\bigsqcup\limits_{k\in K} I_k$ is the splitting of $I$ 
into rooted trees, then $M\!=\!\bigsqcup\limits_{k\in K}M_k$
for some subsets $M_k$ of~$M$ such that
${\cal C}_k\!=\!\big(S,M_k,I_k;x|\hat{I}_k,(j,y)|M_k\big)$ 
is an $M_k$-marked curve based on~$S$;\\
(3) $u_h\!:S\lra V$ if $h\!\in\! I-\hat{I}$ 
and $u_h\!:S^2\lra V$ if $h\!\in\!\hat{I}$ is a smooth map
such that \hbox{$u_h(\i)\!=\!u_{\io_h}(x_h)$} for all $h\!\in\!\hat{I}$;\\
(4) for all $i\!\in\!\hat{I}$ if $S=\Si$ and $i\!\in\! I$ if $S=S^2$,
$$|\{h\!\in\!\hat{I}\!: \io_h=i\}|+|\{l\!\in\! M\!:j_l=i\}|<2
\Lra u_{i*}[S^2]\neq 0\in H_2(V;\Bbb{Z}).$$
The bubble map $b$ is \under{simple} if $I$ is a rooted tree;
$b$ is \under{holomorphic} if  $\bar{\partial}u_i=0$ for 
\hbox{all $i\!\in\! I$}. 
\end{dfn}

\noindent
With notation as in Definition~\ref{bubble_dfn},
every bubble map~$b$ corresponds to a continuous map
$$u_b\!: \Si_b\equiv\bigsqcup_{k\in K}\Si_{{\cal C}_k}\lra V,$$ 
which is smooth on the components of $\Si_{{\cal C}_k}$.
If $h\!\in\!\hat{I}_k$, let $F_{h,b}^{(0)}=F_{h,{\cal C}_k}^{(0)}$.
Similarly, let
$$F_b^{(0)}=\bigoplus\limits_{k\in K}F_{{\cal C}_k}^{(0)},\quad
\Si_b^*=\bigcup_{k\in K}\Si_{{\cal C}_k}^*\subset\Si_b,\quad
\Si_{b,i}^*=\Si_{{\cal C}_k,i}^*\subset
       \Si_{b,i}\equiv\Si_{{\cal C}_k,i},$$
whenever $i\!\in\! I_k\!\subset\! I$.
If $b$ is simple, denote by $\t_b$ the bubble tree 
$\t_{{\cal C}_k}$ for the unique \hbox{element $k\!\in\! K$.}

\begin{dfn}
\label{bubble_isom}
~An~~ \under{isomorphism of~ $V$-valued bubble maps}~~~
$b=\big(S,M,I;x,(j,y),u\big)$~ and\\ 
\hbox{$b'=\big(S,M,I';x',(j',y'),u'\big)$
is a tuple of maps}
$$\phi_0\!: I\lra I',\quad 
\phi_{1,i}\!:S\lra S~~\hbox{for~}i\!\in\! I\!-\!\hat{I},\quad
\phi_{1,i}\!:S^2\lra S^2~~\hbox{for~}i\!\in\!\hat{I},
\qquad\hbox{where}$$
(1) $\phi_0$ is an isomorphism of the linearly ordered sets $I$
and $I'$ with $\phi_0(j_l)\!=\!j_l'$ for all $l\!\in\! M$;\\
(2) $\phi_{1,i}$ is  a biholomorphic map for all $i\!\in\! I$
and is the identity map if $S=\Si$ and $i\!\not\in\!\hat{I}$;\\
(3) $\phi_{1,i}(\i)=\i$ for all $i\!\in\! I$ if $S=S^2$ and  
for all $i\!\in\!\hat{I}$ if $S=\Si$;\\
(4) $\phi_{1,\io_h}(x_h)=x_{\phi_0(h)}'$ for all $h\!\in\!\hat{I}$ and
$\phi_{1,j_l}(y_l)=y_l'$ for all $l\!\in\! M$;\\
(5) $u_{\phi_0(i)}'\circ \phi_{1,i}=u_i$ for all $i\!\in\! I$.
\end{dfn}

\noindent
Such a set of maps corresponds to a continuous map $\Si_b\!\lra\!\Si_{b'}$ 
that  maps the marked points of~$b$ to the marked points of $b'$,
intertwines the maps $u_b\!:\Si_b\lra V$ and $u_{b'}\!:\Si_{b'}\lra V$, 
and is biholomorphic on each component $\Si_{b,i}$ of $\Si_b$.
Let $G_b$ denote the group of automorphisms of the bubble map $b$.
This group is necessarily finite by the stability condition~(4)
of Definition~\ref{bubble_dfn}.
If \hbox{$\la\!\in\! H_2(V;\Bbb{Z})$, let}
\begin{gather*}
\bar{C}^{\i}_{(\la;M)}(S;V)=
\big\{ b\!=\!\big(S,M,I;x,(j,y),u\big)
\hbox{~is $V$-valued bubble map}\!:
\sum_{i\in I}u_{i*}[\Si_{b,i}]=\la\big\}\big/\sim;\\
C^{\i}_{(\la;M)}(S;V)=\big\{ 
b\!=\!\big(S,\{\hat{0}\};,(\hat{0},y),u_{\hat{0}}\big) 
\hbox{~is $V$-valued bubble map}\!:
u_{\hat{0}*}[S]\!=\!\la\big\}\big/\sim,
\end{gather*}
where the equivalence relation is given by isomorphisms
of $V$-valued bubble maps.
If $\mu\!=\!\mu_{\tilde{M}}$ 
is an $\tilde{M}$-tuple of submanifolds of~$V$, let
\begin{gather*}
\bar{C}^{\i}_{(\la;M)}(S;\mu)=
\big\{b\!=\!\big[S,M,I;x,(j,y),u\big]\!\in\! 
\bar{C}^{\i}_{(\la;M)}(S;V)\!:
u_{j_l}(y_l)\!\in\!\mu_l~\forall l\!\in\! M\cap\tilde{M}\big\},\\
C^{\i}_{(\la;M)}(S;\mu)=
\big\{b=\big[S,M,\{\hat{0}\};,(\hat{0},y),u_{\hat{0}}\big]
\!\in\! C^{\i}_{(\la;M)}(S;V)\!:
u_{\hat{0}}(y_l)\!\in\!\mu_l~\forall l\!\in\! M\cap\tilde{M}\big\}.
\end{gather*}
A topology on $\bar{C}^{\i}_{(\la;M)}(S;V)$ and its subsets
$C^{\i}_{(\la;M)}(S;V)$, $\bar{C}^{\i}_{(\la;M)}(S;\mu)$,
$C^{\i}_{(\la;M)}(S;\mu)$ is defined below.

\begin{dfn}
\label{gr_top}
Let $b^*\!=\!\big(S,M,I^*;x^*,(j^*,y^*),u^*\big)$ and 
$b_k\!=\!\big(S,M,I_k;x_k,(j_k,y_k),u_k\big)$ 
be simple bubble maps.
If $S\!=\!S^2$, the sequence $\{b_k\}$ \under{converges} to $b^*$ if
for all $k$ sufficiently large one can~choose\\
(a) $M$-marked curves 
${\cal C}_k=\big(S,M,I^*;x'_k,(j^*,y^*)\big)$, and\\
(b) vectors $(v_k)_{\hat{I}^*}\in F_{{\cal C}_k}^{(0)}$
with $16|v_k|<r_{{\cal C}_k}^2$,\\
such that with $\ups_k=\l({\cal C}_k,(v_k)_{\hat{I}^*}\r)$, \\
(1) $\lim\limits_{k\lra\i}x_{k,h}'=x_h^*$ for all $h\!\in\!\hat{I}$,
and $\lim\limits_{k\lra\i}|\ups_k|=0$;\\
2) ${\cal C}(\ups_k)=\big(S,M,I_k;x_k,(j_k,y(\ups_k))\big)$, 
$$\lim_{k\lra\i}q_{\ups_k}(j_{k,l},y_{k,l})=(j_l^*,y_l^*)
~~\forall l\!\in\! M,\quad\hbox{and}\quad
\lim_{k\lra\i}\sup_{z\in\Si_{{\cal C}(\ups_k)}}
   d_V(u_{b^*}(q_{\ups_k}(z)),u_{b_k}(z))=0.$$
If $S\!=\!\Si$, convergence is defined in the same~way, but
$|v_k|$ and ${\cal C}(\ups_k)$ are replaced by  
$|v_k|_g$ and  ${\cal C}_g(\ups_k)$, respectively,
for a $\t_{b^*}$-admissible metric \hbox{$g$ on~$\Si$}.
\end{dfn}

\noindent
This notion of convergence is independent of 
the choice of an admissible metric on~$\Si$.
Definition~\ref{gr_top} induces a topology on the space 
$\bar{C}^{\i}_{(\la;M)}(S;V)$,
which will be referred to as the Gromov topology.

\subsection{Stratums of Bubble Maps}
\label{stratification}

\noindent
In this subsection, we introduce the notion of a bubble type.
We then define various spaces of holomorphic bubble maps indexed
by bubble types and vector bundles over~them.

\begin{dfn}
\label{strata_dfn}
A \under{bubble type} is a tuple
${\cal T}=\big(S,M,I;j,\la\big)$,
such that\\
(1) $I$ is a linearly ordered set,
and $j\!:M\!\lra\!I$ and $\la\!:I\!\lra\!H_2(S;\Bbb{Z})$ are maps;\\
(2) for all $i\!\in\!\hat{I}$ if $S=\Si$ and all $i\!\in\! I$ if $S=S^2$,
$\la_i\neq0$ if $|\{h\!: \io_h=i\}|+|\{l\!: j_l=i\}|<2$.\\
Bubble type ${\cal T}$ is \under{simple} if $I$ is a rooted tree;
${\cal T}$ is \under{basic} if $\hat{I}=\eset$.\\
Two bubble types
${\cal T}\!=\!\big(S,M,I;j,\la\big)$ and
${\cal T}'\!=\!\big(S,M,I';j',\la'\big)$
are \under{equivalent} if there exists an isomorphism of
linearly ordered sets $\phi_0\!: I\lra I'$ such that
$\phi_0(j_l)=j_l'$ for all $l\!\in\! M$ and
\hbox{$\la_{\phi_0(i)}'=\la_i$} for all $i\!\in\! I$.\\
If  ${\cal T}^*\!=\!\big(S,M,I^*;j^*,\la^*\big)$ and
${\cal T}\!=\!\big(S,M,I;j_M,\la\big)$ are two bubble types,
\under{${\cal T}^*\!<\!{\cal T}$} if $I\!\subset\! I^*$, 
$$j_l=\max\{i\!\in\! I\!: i\le j_l^*\}~~\hbox{for all~} l\!\in\! M
\quad\hbox{and}\quad
\la_i=\sum_{i=\max\{i'\le h: i'\in I\}}\!\!  \la_h^*.$$
If ${\cal T}\!=\!\big(S,M,I;j,\la\big)$ is a bubble type,
a \under{${\cal T}$-bubble map} is a bubble map
\hbox{$b\!=\!\big(S,M,I;x,(j,y),u\big)$}
such that $u_{i*}[\Si_{b,i}]\!=\!\la_i\in H_2(V;\Bbb{Z})$ for all
$i\!\in\! I$.
\end{dfn}

\noindent
The splitting of $I$ into rooted trees $I_k$ induces
a splitting of ${\cal T}$ into simple bubble types
$${\cal T}_k=\big(S,M_k,I_k;j_k,\la_k\big),$$
where $j_k$ and $\la_k$ are the restrictions
of $j$ and $\la$ to  $M_k$ and $I_k$, respectively.
Similarly, each ${\cal T}$-bubble map $b$ corresponds 
to a $K$-tuple of bubble maps $b_K=(b_k)_{k\in K}$, where $b_k$
is a ${\cal T}_k$-bubble~map.\\

\noindent
We denote the equivalence class of the bubble type ${\cal T}$
by~$[{\cal T}]$ and the groups of automorphisms 
of~${\cal T}$ by $\hbox{Aut}({\cal T})$.
The partial ordering on the set of bubble types induces
a partial ordering on the set of their equivalence classes.
If $b$ and $b'$ are ${\cal T}$- and ${\cal T}'$-bubble maps,
respectively, such that $[b]\!=\![b']$, then $[{\cal T}]\!=\![{\cal T}']$.
Furthermore, if $\{b_k\}$ is a~sequence of ${\cal T}$-bubble maps,
$b^*$ is ${\cal T}^*$-bubble map, and $[b_k]$ converges to $[b^*]$
with respect to the Gromov topology, then
$[{\cal T}^*]\le[{\cal T}]$.\\

\noindent
Let ${\cal T}\!=\!\big(S,M,I;j,\la\big)$ be a bubble type.
We denote by $\lr{\cal T}$ the basic bubble type such 
\hbox{that $\lr{\cal T}\!\ge\!{\cal T}$}.
It can be described explicitly as follows.
Let $I\!=\!\!\bigsqcup\limits_{k\in K}I_k$ be the splitting
of~$I$~into rooted trees and $M\!=\!\!\bigsqcup\limits_{k\in K}M_k$
the corresponding splitting of~$M$;
see Definition~\ref{bubble_dfn}.
It can be assumed that $K\!=\! I\!-\!\hat{I}$ and $k$ is the unique
minimum element of~$I_k$.
For every $k\!\in\! K$ and \hbox{$l\!\in\! M_k$, let}
$$\la_k'=\sum_{i\in I_k}\la_i,\qquad
j_l'=k.$$
Then $\lr{\cal T}=\big(S,M,K;j',\la'\big)$.\\

\noindent
If ${\cal T}\!=\!\big(S,M,I;j,\la\big)$ is a simple bubble type
and $i\!\in\! I$, let
$D_i{\cal T}=D_iI$ and $\bar{D}_i{\cal T}=\bar{D}_iI$.
If~$H$ is a subset of $\hat{I}$, let
${\cal T}(H)\!=\!\big(S,M,H\cup\hat{0};j',\la'\big)$,
where 
$$j_l'=\max\{i\!\in\! H\cup\hat{0}\!: i\le j_l\}
\qquad\hbox{and}\quad
\la_i'=\sum_{i_H^*\le h\le i}\la_h
\hbox{~~with~~}
i_H^*=\max\{i^*\!\in\! H\cup\hat{0}: i^*\le i\}.$$
Then ${\cal T}(H)$ is again a bubble type.
The bubble type ${\cal T}(H)$ is the bubble type obtained 
by gluing ${\cal T}$-bubble maps with the parameter $v_{\hat{I}}$
such that $v_h=0$ if and only if $h\!\in\! H$; see the next section.
Finally, we will denote by $(\la;M)$ the equivalence class
of the basic simple bubble type~$(S^2,M,\{\hat{0}\};\hat{0},\la)$.
\\

\noindent
Given a bubble type ${\cal T}\!=\!(S,M,I;j,\la)$,
let $d({\cal T})\!: I\!\lra\!\Bbb{R}$ be given by\\
\begin{equation}\label{weights_dfn}
d_i({\cal T})=|\la_i|+|\{l\!\in\! M\!: j_l=i\}|+
        \sum_{\io_h=i}d_h({\cal T}) \quad\forall i\!\in\! I.
\end{equation}
Since $I$ is a linearly ordered set, the numbers $d_i({\cal T})$ are
uniquely defined by~\e_ref{weights_dfn}.
If 
$$b=\big(S,M,I;x,(j,y),u\big)$$ 
is a ${\cal T}$-bubble map, $b$ is  ${\cal T}${\it{-balanced}}
if for all~$i\!\in\!\hat{I}$ \\
(B1) $\int_{\Bbb{C}}|du_i\circ q_N|^2z+
\sum\limits_{\io_h=i}d_h({\cal T})x_h+
\sum\limits_{j_l=i}y_l=0$;\\
(B2) $\int_{\Bbb{C}}|du_i\circ q_N|^2\be(|z|)+
\sum\limits_{\io_h=i}d_h({\cal T})\be(|x_h|)+
\sum\limits_{j_l=i}\be(|y_l|)=\frac{1}{2}$.\\
The integrals above are computed with respect to 
the metric $g_V$ on~$V$.
 Recall that we consider $\Bbb{C}$ to be a subset of $S^2$ via the map $q_N$.
Thus, $x_h$ and $y_l$ can be viewed as complex numbers,
as done above.
If $S\!=\!S^2$ and $b$ is as above, $b$ is
{\it completely} ${\cal T}${\it{-balanced}} (or {\it{cb}})
if (B1) and (B2) hold for \hbox{all $i\!\in\! I$.}\\

\noindent
Denote by ${\cal H}_{\cal T}$ the set of all holomorphic 
${\cal T}$-bubble maps.
Let
$$PSL_2^{(0)}=\{g\!\in\! PSL_2\!: g(\i)=\i\},\qquad
{\cal G}_{\cal T}=\prod_{h\in\hat{I}}PSL_2^{(0)}.$$
The group ${\cal G}_{\cal T}$ acts on ${\cal H}_{\cal T}$ as follows. If 
$$b=\big(S,M,I;x,(j,y),u\big)\in{\cal H}_{\cal T}
\quad\hbox{and}\quad g=g_{\hat{I}}\in{\cal G}_{\cal T},$$
define $gb=\big(S,M,I;gx,(j,gy),(gu)\big)$ by
$$(gx)_h=\begin{cases}
g_{\io_h}x_h,&\hbox{if~}\io_h\!\in\!\hat{I};\\
x_h,&\hbox{if~}\io_h\!\not\in\!\hat{I};
\end{cases}\quad
(gy)_l=\begin{cases}
g_{j_l}y_l,&\hbox{if~}j_l\!\in\!\hat{I};\\
y_l,&\hbox{if~}j_l\!\not\in\!\hat{I};
\end{cases}\hspace{3mm}
(gu)_i=\begin{cases}
g_i\cdot u_i,& \hbox{if~}i\!\in\!\hat{I};\\
u_i,& \hbox{if~}i\!\not\in\!\hat{I},
\end{cases}$$
where for any map $f\!: S^2\lra V$ and $g\!\in\! PSL_2$, we define
$$g\cdot f\!:S^2\lra V  \quad\hbox{by}\quad
\{g\cdot f\}(z)=f(g^{-1}z).$$
Let ${\cal M}^{(0)}_{\cal T}\!\subset\!{\cal H}_{\cal T}$ 
denote the subset of ${\cal T}$-balanced holomorphic maps.
The group $G_{\cal T}\!\times\!\hbox{Aut}({\cal T})$, where
$$G_{\cal T}\equiv\prod_{h\in\hat{I}}S^1\subset{\cal G}_{\cal T},$$
acts on ${\cal M}^{(0)}_{\cal T}$ and all the stabilizers are finite.
Denote the quotient by ${\cal M}_{\cal T}$, and let
$$\bar{\cal M}_{\cal T}=\bigcup_{{\cal T}'\le{\cal T}}
{\cal M}_{{\cal T}'}.$$
If $\hbox{Aut}({\cal T})\!=\!\{1\}$,
corresponding to the quotient 
${\cal M}_{\cal T}={\cal M}^{(0)}_{\cal T}/G_{\cal T}$, 
we obtain $|\hat{I}|$ line (orbi)-bundles 
$$\{L_h{\cal T}\lra{\cal M}_{\cal T}\!: h\!\in\!\hat{I}\},$$
which carry natural norms:
$$|[b,c_h]|=|c_h|\qquad\hbox{if}\quad b\in{\cal M}^{(0)}_{\cal T}
\quad\hbox{and}\quad c_h\in\Bbb{C}.$$
If $\hbox{Aut}({\cal T})\!\neq\!\{1\}$,
the fiber products and connect sums of the above line bundles taken
over each orbit of $\hbox{Aut}({\cal T})$ are well-defined. 
Let $F^{(0)}_h{\cal T}\!\lra\!{\cal M}^{(0)}_{\cal T}$ be 
the bundle with the fiber $F_{h,b}^{(0)}$ 
at~$b\!\in\! {\cal M}^{(0)}_{\cal T}$\!,~i.e.
$$F^{(0)}_h{\cal T}=
\begin{cases}
{\cal M}^{(0)}_{\cal T}\times\Bbb{C},&\hbox{if~}x_h\!\in\! S^2;\\
\pi_h^*T\Si,&\hbox{if~}x_h\!\in\!\Si,
\end{cases}
\qquad\hbox{where~~}\pi_h(b)=x_h,$$
with notation as above.
The action of $G_{\cal T}$ on ${\cal M}^{(0)}_{\cal T}$
lifts to an action on  each bundle 
$F_h^{(0)}{\cal T}$~by
$$g\cdot(b,v_h)=
\begin{cases}
\big(g\cdot b,g_{\io_h}g_h^{-1}v_h\big),&\hbox{if~}\io_h\!\in\!\hat{I};\\
\big(g\cdot b,g_h^{-1}v_h\big),&\hbox{if~}\io_h\!\not\in\!\hat{I}.
\end{cases}$$
Here and in the rest of the paper,
we identify $S^1$ with the unit complex
numbers in the usual way.
Let $F_h{\cal T}$ be the line orbi-bundle over ${\cal M}_{\cal T}$ given by
$$F_h{\cal T}= F^{(0)}_h{\cal T}\big/
\big(G_{\cal T}\times\hbox{Aut}({\cal T})\big).$$
This bundle has a natural norm unless $\io_h\!=\!\hat{0}$ and $S\!=\!\Si$.
In such a case, any metric $g$ on $\Si$ induces a norm on 
$F_h{\cal T}$.
Let 
$$F^{(0)}{\cal T}=\bigoplus\limits_{h\in\hat{I}}F^{(0)}_h{\cal T},\quad
F^{(0)}_b{\cal T}=F^{(0)}{\cal T}\big|b;\quad
F{\cal T}=\bigoplus\limits_{h\in\hat{I}}F_h{\cal T},\quad
F^{(0)}_{[b]}{\cal T}\!=\!F^{(0)}{\cal T}\big|[b].$$
Note that if ${\cal T}^*\!<\!{\cal T}$, 
$G_{\cal T}$ is naturally a subgroup of $G_{{\cal T}^*}$
and thus acts on ${\cal M}_{{\cal T}^*}^{(0)}$ and 
the line bundles~$F^{(0)}_h{\cal T}^*$.
Furthermore, there is a natural decomposition
$G_{{\cal T}^*}\!=\!G_{\cal T}\!\times\! G$ for a certain group~$G$.
In particular, the $G_{\cal T}$-action on ${\cal M}_{\cal T}^{(0)}$
and $F^{(0)}_h{\cal T}$ induces an action of~$G_{{\cal T}^*}$.\\

\noindent
If $S\!=\!S^2$, let 
$${\cal B}_{\cal T}=
\big\{b\!=\!\big(S,M,I;x,(j,y),u\big)\!\in\! {\cal H}_{\cal T}\!:
b\hbox{~is cb};~
u_{i_1}(\i)\!=\!u_{i_2}(\i)~\forall i_1,i_2\!\in\! I\!-\!\hat{I}\big\}.$$
Denote by ${\cal U}_{\cal T}^{(0)}\!\subset\!{\cal M}_{\cal T}$
the quotient ${\cal B}_{\cal T}/(G_{\cal T}\!\times\!\hbox{Aut}({\cal T}))$.
The group 
$$G_{\cal T}^*\equiv\prod_{i\in I-\hat{I}}S^1$$
acts on ${\cal U}_{\cal T}^{(0)}$ 
and ${\cal M}_{\cal T}$ as follows.
If 
$$[b]=\big[(S^2,M,I;x,(j,y),u)\big]
\in{\cal M}_{\cal T}
\quad\hbox{and}\quad g=(g_i)_{i\in I-\hat{I}}\!\in\! G_{\cal T}^*,$$
define $g[b]\!=\!\big[(S^2,M,I;gx,(j,gy),gu)\big]$ by
$$(gx)_h=\begin{cases}
x_h,&\hbox{if~}\io_h\!\in\!\hat{I};\\
g_{\io_h}x_h,&\hbox{if~}\io_h\!\not\in\!\hat{I};
\end{cases}\hspace{3mm}
(gy)_l=\begin{cases}
y_l,&\hbox{if~}j_l\!\in\!{I};\\
g_{j_l}y_l,&\hbox{if~}j_l\!\not\in\!\hat{I};
\end{cases}\hspace{3mm}
(gu)_i=\begin{cases}
u_i,& \hbox{if~}i\!\in\!\hat{I};\\
g_i\cdot u_i,& \hbox{if~}i\!\not\in\!\hat{I}.
\end{cases}$$
As in the previous paragraph, all the stabilizers are finite.
Furthermore, this $G_{{\cal T}^*}$-action on ${\cal M}_{\cal T}$
naturally lifts to an action on ${\cal M}_{\cal T}^{(0)}$
and along with the $G_{\cal T}$-action on ${\cal M}_{\cal T}^{(0)}$
induces an action of 
$\tilde{G}_{\cal T}\equiv G_{\cal T}^*\times G_{\cal T}$
on ${\cal M}_{\cal T}^{(0)}$ as well as on $F_h^{(0)}{\cal T}$ by
$$(g^*,g)\cdot(b,v_h)=
\begin{cases}
\big((g^*,g)\cdot b,g_{\io_h}g_h^{-1}v_h\big),
&\hbox{if~}\io_h\!\in\!\hat{I};\\
\big((g^*,g)\cdot b,g_{\io_h}^*g_h^{-1}v_h),
&\hbox{if~}\io_h\!\not\in\!\hat{I}.
\end{cases}$$
Note that $G_{{\cal T}'}=G_{\cal T}$ whenever 
${\cal T}'\le{\cal T}$.
Let 
$${\cal U}_{\cal T}={\cal U}_{\cal T}^{(0)}/G_{\cal T}^*,\qquad
\ov{\cal U}_{\cal T}^{(0)}=
 \bigcup_{{\cal T}'\le{\cal T}}{\cal U}_{{\cal T}'}^{(0)}, \qquad
\ov{\cal U}_{\cal T}=\bigcup_{{\cal T}'\le{\cal T}}{\cal U}_{{\cal T'}}.$$
With respect to the Gromov topology,
the space $\ov{\cal U}_{\cal T}^{(0)}$ is  Hausdorff
and compact if $(V,\om,J)$ is semipositive; see~\cite{RT}.
Furthermore, $G_{\cal T}^*$ acts continuously on 
${\cal U}_{{\cal T}'}^{(0)}$ as can be easily seen 
from Definition~\ref{gr_top}.
If follows that $\ov{\cal U}_{\cal T}$ is also Hausdorff and
compact if $(V,\om,J)$ is semipositive in the quotient topology.
Denote by $\{L_i{\cal T}\lra\ov{\cal U}_{\cal T}\!:i\!\in\! I-\hat{I}\}$ 
the line orbi-bundles corresponding to the quotient
$\ov{\cal U}_{\cal T}=\ov{\cal U}_{\cal T}^{(0)}/G_{\cal T}^*$.
Let 
$${\cal F}_h{\cal T}=
 \Big(F_h{\cal T}^{(0)}\big|{\cal B}_{\cal T}\Big)/\tilde{G}_{\cal T}
\lra{\cal U}_{\cal T},
\quad
{\cal F}_{h,[b]}{\cal T}={\cal F}_h{\cal T}\big|[b];\quad
{\cal F}{\cal T}=\bigoplus_{h\in\hat{I}}{\cal F}_h{\cal T},\quad
{\cal F}_{[b]}{\cal T}={\cal F}{\cal T}\big|[b].$$
The line bundles ${\cal F}_h{\cal T}$ have natural norms,
defined as in the previous paragraph.\\

\noindent
If ${\cal T}\!=\!(S,M,I;j,\la)$ is a bubble type
 and $b\!=\!\big(S,M,I;x,(j,y),u\big)$
is a ${\cal T}$-bubble map, for any $l\!\in\! M$ let 
$\ev_l\!:{\cal H}_{\cal T}\!\lra\! V$ be the map given by 
$$\ev_l\big((S,M,I;x,(j,y),u)\big)=
u_{j_l}(y_l).$$
This map descends to the quotients defined above and
induces continuous maps on the spaces
$\ov{\cal M}_{\cal T}$, $\ov{\cal U}_{\cal T}^{(0)}$, and
$\ov{\cal U}_{\cal T}$.
If $\mu\!=\!\mu_{\tilde{M}}$ is an $\tilde{M}$-tuple of submanifolds in $V$,
where $\tilde{M}$ is possibly different from~$M$,~put
$${\cal H}_{\cal T}(\mu)=
\{b\!\in\!{\cal H}_{\cal T}\!: \ev_l(b)\!\in\!\mu_l~\forall 
 l\!\in\! M\cap\tilde{M}\}.$$
Define spaces ${\cal M}_{\cal T}^{(0)}(\mu)$, 
${\cal M}_{\cal T}(\mu)$, $\ov{\cal M}_{\cal T}(\mu)$,  etc.~similarly.
If $S\!=\!S^2$, we define another evaluation~map,
$$\ev\!: {\cal B}_{\cal T}\lra V\quad\hbox{by}\quad
\ev\big((S^2,M,I;x,(j,y),u)\big)=u_{\hat{0}}(\i),$$
where $\hat{0}$ is any minimal element of $I$.
This map induces continuous maps on the spaces
$\bar{\cal U}_{\cal T}^{(0)}$ \hbox{and $\bar{\cal U}_{\cal T}$}.

\section{The Gluing Construction and the Obstruction Bundle}
\label{gluing}

\subsection{Summary and Notation}
\label{summary3}

\noindent
We now present a gluing construction on the spaces ${\cal M}_{\cal T}(\mu)$ 
such that ${\cal H}_{\cal T}$ is a smooth manifold with 
the tangent bundle isomorphic to the kernel of the linearization
of the $\bar{\partial}$-operator, as defined below.
This is well-known to be the case if the linearization 
of the $\bar{\partial}$-operator is surjective; see~\cite{MS}.
However, surjectivity of the  linearization is not a necessary 
condition; see~\cite{Z2} for examples.
In fact, there are two main cases of primarily interest to~us.
The first is when 
\hbox{${\cal T}\!=\!\big(S^2,M;j_M,\la)$} and 
the linearization of the $\bar{\partial}$-operator is indeed surjective.
In this case, we give an analytic description of
a neighborhood of ${\cal U}_{\cal T}(\mu)$
in~$\bar{\cal U}_{\lr{\cal T}}(\mu)$ for a generic set 
of constraints~$\mu$.
The second case is when 
\hbox{${\cal T}\!=\!\big(\Si,M,I;j,\la)$} and 
the cokernels of the linearization of the $\bar{\partial}$-operator
form a vector bundle over~${\cal H}_{\cal T}$,
which will  the analogue of
Taubes's obstruction bundle of~\cite{T}
in the gluing construction below.
Using the same analysis as in the first case,
we describe any sufficiently nice element of
 $C^{\i}_{(\la^*;M)}(\Si;\mu)$  lying near~${\cal M}_{\cal T}(\mu)$,
\hbox{where $\la^*\!=\!\sum\la_i$}.
The elements of ${\cal M}_{\Si,t\nu,\la^*}(\mu)$ lying near 
${\cal M}_{\cal T}(\mu)$ will correspond to the zero set 
of a certain section of the obstruction bundle.\\

\noindent
For our gluing construction, we fix a smooth family  
\hbox{$\{g_{V,b}\!: b\!\in\!{\cal M}_{\cal T}\}$} 
of Kahler metrics on~$(V,J)$.
We assume that this family is
$\big(G_{\cal T}\!\times\!\hbox{Aut}({\cal T})\big)$-invariant 
if $S\!=\!\Si$ and 
$\big(\tilde{G}_{\cal T}\!\times\!\hbox{Aut}({\cal T})\big)$-invariant 
if $S\!=\!S^2$.
If $b\!\in\!{\cal M}_{\cal T}$, $X,Y\!\in\! T_qV$, and 
$u\!:(D,j)\!\lra\! V$ is a smooth map from a one-dimensional complex 
manifold,~let
$$\exp_{b,q}X= \exp_{g_b,q}X,\quad
\na^b=\na^{g_b},\quad
\Pi_{b,X}Y=\Pi_{g_b,X}Y,\quad
D_{b,u}=D_{g_b,u};$$
see Subsection~\ref{notation} for more details.
If $S\!=\!\Si$, we also choose a smooth family 
$$\{g_{{\cal T},x}\!: x\!=\!(x)_{\{h:\io_h=\hat{0}\}};~ x_h\!\in\!\Si;~ 
x_{h_1}\!\neq\! x_{h_2}\hbox{~if~}h_1\!\neq\! h_2\}$$
of Riemannian metrics on $\Si$ such that each metric 
$g_{{\cal T},x}$ is flat on a neighborhood of $x_h$ in $\Si$
for all $h\!\in\!\hat{I}$ with $\io_h\!=\!\hat{0}$.
Existence of such a family of metrics is shown in \cite{FO}.
It can be assumed that all these metrics have the same volume
as the standard metric on $S^2$, i.e. $4\pi$.
If
$$b=\big(\Si,M,I;x,(j,y),u\big)\!\in\!{\cal H}_{\cal T},$$
let $g_{b,\hat{0}}$ denote the metric 
$g_{{\cal T},(x)_{\{h:\io_h=\hat{0}\}}}$ on $\Si$.
If $i\!\in\!\hat{I}$, we write $g_{b,i}$ for the standard metric on~$S^2$.
Similarly, if $S=S^2$, for all $i\!\in\! I$, we write $g_{b,i}$
for the standard metric on $S^2$.\\

\noindent
If $b\!=
\!\big(S,M,I;x,(j,y),u\big)\!\in\!{\cal H}_{\cal T}$,
let
\begin{gather*}
\Ga'(b)=\bigoplus_{i\in I}\Ga(u_i); \quad
\Ga(b)=\Ga(u_b)=\big\{\xi_I\!\in\! \Ga'(b)\!:
 \xi_h(\i)\!=\!\xi_{\io_h}(x_h)~\forall h\!\in\!\hat{I}\big\}; \\
\Ga^1(b)=\Ga^1(u_b)=\bigoplus_{i\in I}\Ga^1(u_i); \qquad
\Ga^{0,1}(b)=\Ga^{0,1}(u_b)=\bigoplus_{i\in I}\Ga^{0,1}(u_i).
\end{gather*}
Define $D_b\!: \Ga(b)\!\lra\!\Ga^{0,1}(b)$ by
$$\big(D_b\xi_I\big)_i=D_{b,u_i}\xi_i\qquad\forall i\!\in\! I.$$
We denote the kernel of operator $D_b$ on $\Ga(b)$
by~$\Ga_-(b)$.
If $\xi\!\in\!\Ga(u_i)$ or $\xi\!\in\!\Ga^1(u_i)$,  let $\|\xi\|_{b,C^k}$
and $\|\xi\|_{b,2}$ denote the $C^k$- and $L^2$-norms of $\xi$
computed with respect to the metrics $g_{V,b}$ on $V$ and
$g_{b,i}$ on $\Si_{b,i}$.
If $\xi\!=\!\xi_I\!\in\!\Ga'(b)$ or $\xi\!\in\!\Ga^1(b)$, put
$$\|\xi\|_{b,C^k}=\sum_{i\in I}\|\xi_i\|_{b,C^k},~~~
\|\xi\|_{b,2}=\sum_{i\in I}\|\xi_i\|_{b,2}.$$
Let $\pi_{b,-}\!: \Ga(b)\lra\Ga_-(b)$ be the $(L^2,b)$-orthogonal
projection map.\\

\noindent
The space ${\cal P}_b{\cal T}$ of perturbations of bubble 
map $b$ is
the collection of tuples \hbox{$\si=(\xi_{\hat{I}};w_{\hat{I}+M})$}, where
$$\xi_i\in\Ga(u_i)~\forall i\!\in\! I,\quad
w_h\in F_{h,b}^{(0)}~\forall h\!\in\!\hat{I},\quad
w_l\in\begin{cases}
\Bbb{C},&\hbox{if~}l\!\in\! M~\&~\Si_{b,j_l}=S^2;\\
T_{y_l}\Si,&\hbox{if~}l\!\in\! M~\&~\Si_{b,j_l}=\Si.\end{cases}$$
If $\si$ is sufficiently small, we define
$\exp_b\si\!=\!\big(S,M,I;x(\si),(j,y(\si)),u_{\si}\big)$
by
$$x_h(\si)=\begin{cases}
x_h+w_h,&\hbox{if~}\Si_{b,i_h}=S^2;\\
\exp_{g_{b,\hat{0}},x_h}w_h,&\hbox{if~}\Si_{b,i_h}=\Si;
\end{cases}~~~
y_l(\si)=\begin{cases}
y_l+w_l,&\hbox{if~}\Si_{b,j_l}=S^2;\\
\exp_{g_{b,\hat{0}},y_l}w_l,&\hbox{if~}\Si_{b,j_l}=\Si;
\end{cases}$$
and  $u_{\si,i}=\exp_{b,u_i}\xi_i$.
If $z\!\in\!\Si$, let $|v|_b=|v|_{g_{b,\hat{0}},x}$.
For consistency, if $v\!\in\!\Bbb{C}$, let $|v|_b=|v|$.
Along with the $(L^2,b)$-norm on the vector fields defined above,
we obtain an inner-product on the space of tuples $\si$ as above.\\

\noindent
In order to get a good description of the spaces
${\cal M}_{\cal T}^{(0)}$ as submanifolds of ${\cal H}_{\cal T}$,
we describe an action of an open subset of $0$
in $\l(\Bbb{C}\oplus\Bbb{R}\r)^{\hat{I}}$
on bubble maps and distinguished elements 
$\si_{(b,i)}^{(k)}\!\in\! {\cal P}_b{\cal T}$
that correspond to this~action.
If $(c,r,\th)\!=\!(c,r,\th)_{\hat{I}}\!\in\!
         \big(\Bbb{C}\!\times\!\Bbb{R}\!\times\!\Bbb{R}\big)^{\hat{I}}$ 
and $b$ is a bubble map as above, we~define 
$$(c,r,\th)\cdot b=
      \big(S,M,I;(c,r,\th)x,(j,(c,r,\th)y),(c,r,\th)u\big)$$
by setting
\begin{gather*}
\big((c,r,\th)x\big)_h=e^{i\th_{\io_h}}(1+r_{\io_h})(x_h+c_{\io_h}),\qquad
\big((c,r,\th)y\big)_l=e^{i\th_{j_l}}(1+r_{j_l})(y_l+c_{j_l}),\\
\big((c,r,\th)u\big)_i\big(q_N(z)\big)=
u_i\big(q_N\big((1+r_i)^{-1}e^{-i\th_i}z-c_i\big)\big).
\end{gather*}
If $(c,r,\th)$ is sufficiently small, $(c,r,\th)\cdot b$ 
                                          is again a bubble map, 
i.e.~the maps into $V$ still agree at the nodes,
and the nodes and the marked points are still all distinct.
In fact, the values of the maps at the nodes or the marked points 
do not change,~i.e.
$$\big((c,r,\th)u\big)_{\io_h}\!\big(((c,r,\th)x)_h\big)
\!=\!u_{\io_h}\!(x_h),~
\big((c,r,\th)u\big)_h(\i)\!=\!u_h\!(\i),~
\big((c,r,\th)u\big)_{j_l}\!\big(((c,r,\th)y)_l\big)\!=\!u_{j_l}\!(y_l).$$
Furthermore, if $b\!\in\!{\cal H}_{\cal T}$, 
$(c,r,\th)\cdot b\!\in\!{\cal H}_{\cal T}$.
If $b$ is of type ${\cal T}$,
the above describes the action of a neighborhood of the identity in 
${\cal G}_{\cal T}$ on the space of stable maps of type ${\cal T}$.
The action by $\Bbb{C}$ corresponds to the translations  
of $\Bbb{C}$, by the first $\Bbb{R}$-component to dilations about 
the origin, and by the last $\Bbb{R}$-component to rotations
about the~origin.
In addition, if $S\!=\!S^2$, ${\cal T}$ is simple, and
$(c,r,\th)\!\in\!\Bbb{C}\!\times\!\Bbb{R}\!\times\!\Bbb{R}$
is sufficiently small,
we can define $(c,r,\th)\cdot b$ similarly
by viewing $(c,r,\th)$ as $\big(c_{\hat{0}},r_{\hat{0}},\th_{\hat{0}}\big)$
in the above definition of the action.\\

\noindent
If $u\!\in\! C^{\i}(S^2;V)$, define 
$\xi_u^{(1)},\ldots,\xi_u^{(4)}\in\Ga(u)$ by:
\begin{gather*}
\xi_u^{(1)}({q_N(z)})
=-d(u\circ q_N)\Big|_z\frac{\partial}{\partial s},\qquad
\xi_u^{(2)}({q_N(z)})
=-d(u\circ q_N)\Big|_z\frac{\partial}{\partial t}\\
\xi_u^{(3)}({q_N(z)})=-d(u\circ q_N)\Big|_z
\Big(s\frac{\partial}{\partial s}+t\frac{\partial}{\partial t}\Big)
=-rd(u\circ q_N)\Big|_z\frac{\partial}{\partial r},\\
\xi_u^{(3)}({q_N(z)})=d(u\circ q_N)\Big|_z
\Big(t\frac{\partial}{\partial s}-s\frac{\partial}{\partial t}\Big)
=-d(u\circ q_N)\Big|_z\frac{\partial}{\partial\th}.
\end{gather*}
where we write $z=s+it\in\Bbb{C}$ and $r=\sqrt{s^2+t^2}$.
These vector fields extend smoothly by zero over the south pole.
For any $x\!\in\! S^2-\{\i\}$, let 
$w_x^{(1)},\ldots,w_x^{(4)}\!\in\!\Bbb{C}$ be given~by 
$$w_x^{(1)}=1,\quad w_x^{(2)}=i,\quad w_x^{(3)}=x,\quad w_x^{(4)}=ix.$$
If $b$ is a bubble map as above, $k=1,\ldots,4$,
$i^*\!\in\!\hat{I}$ if $S=\Si$ and $i^*\!\in\!\hat{I}$ if $S=S^2$, 
let 
$$\si_{(b,i^*)}^{(k)}=
\big((\xi_{(b,i^*)}^{(k)})_I,(w_{(b,i^*)}^{(k)})_{\hat{I}+M}\big)$$
be given by
$$\xi_{(b,i^*),i}^{(k)}=
\begin{cases}
\xi_{u_i}^{(k)},&\hbox{if~}i\!=\!i^*;\\
0,&\!\hbox{if~}i\!\neq\! i^*;
\end{cases}\quad
w_{(b,i^*),h}^{(k)}=
\begin{cases}
w_{x_h}^{(k)},&\io_h\!=\!i^*;\\
0,&\io_h\!\neq\! i^*;
\end{cases}\quad
w_{(b,i^*),l}^{(k)}=
\begin{cases}
w_{y_l}^{(k)},&j_l\!=\!i^*;\\
0,&j_l\!\neq\! i^*.
\end{cases}$$
The tuples $\si_{(b,i^*)}^{(k)}$ correspond to the infinitesimal
action of ${\cal G}_{\cal T}$ on the space of stable maps of 
type~${\cal T}$.
\\

\noindent
Finally, if $X$ is any space, $F\!\lra\! X$ a normed vector bundle,
and $\de\!: X\!\lra\!\Bbb{R}$ is any function, let
$$F_{\de}=\big\{(b,v)\!\in\! F\!: |v|_b<\de(b)\big\}.$$
Similarly, if $\Om$ is a subset of $F$, let 
$\Om_{\de}\!=\!F_{\de}\cap\Om$.
If $\ups\!=\!(b,v)\!\in\! F$, denote by $b_{\ups}$ the image of $\ups$
under the bundle projection map, i.e. $b$ in this case.

\subsection{The Basic Setup}
\label{str_sub1}

\noindent
In this subsection, we describe 
our assumptions on the smooth structure of~${\cal H}_{\cal T}$
and state some of their implications.

\begin{dfn}
\label{reg_dfn}
Bubble type ${\cal T}\!=\!(S^2,M,I;j_M,\la)$ 
is \under{$(V,J)$-regular} if for all
$$b=\big(S,M,I;x,(j,y),u\big)
\in{\cal H}_{\cal T}$$
(1) $D_{b,u_i}\!:\Ga(u_i)\lra\Ga^{0,1}(u_i)$ 
is onto for all $i\!\in\! I$;\\
(2)  $\ker D_{b,u_i}\lra T_{u_i(\i)}V$, $\xi\lra\xi(\i)$, 
is onto for all $i\!\in\!I$.
\end{dfn}

\begin{dfn}
\label{semireg_dfn}
Simple bubble type ${\cal T}\!=\!(S,M,I;j,\la)$ 
is \under{$(V,J)$-semiregular} if \\
(1) the space ${\cal H}_{(S,\eset,\{\hat{0}\};,\la_{\hat{0}})}$ is a complex 
manifold of the same dimension as $\ker D_{g_V,b}$ for all
\hbox{$b\!\in\! {\cal H}_{(S,\eset,\{\hat{0}\};,\la_{\hat{0}})}$}, and 
there exist
$\de,C\!\in\! C^{\i}({\cal H}_{(S,\eset,\{\hat{0}\};,\la_{\hat{0}})};
\Bbb{R}^+)$ and for each 
$b\!=\!(S,\eset,\{\hat{0}\};,u_{\hat{0}})
\in {\cal H}_{(S,\eset,\{\hat{0}\};,\la_{\hat{0}})}$
a~map
$$h_{(S,\eset,\{\hat{0}\};,\la_{\hat{0}}),b}\!: 
\big\{\xi\!\in\!\ker D_{g_V,b}\!: 
 \|\xi\|_{g_V,C^0}<\de(b)\big\} \lra \Ga(u_{\hat{0}})$$
such that
$$\|h_{(S,\eset,\{\hat{0}\};,\la_{\hat{0}}),b}(\xi)\|_{g_V,b}
            \!\le\! C(b)\|\xi\|_{g_V,b}^2,
\quad
\big\|h_{(S,\eset,\{\hat{0}\};,\la_{\hat{0}}),b}(\xi)\!-\!
h_{(S,\eset,\{\hat{0}\};,\la_{\hat{0}}),b}(\xi')\big\|_{g_V,C^0}
         \!\le\! C(b)\big\|\xi\!-\!\xi'\big\|_{g_V,C^0},$$
for all \hbox{$\xi,\xi'\!\in\!\ker D_{g_V,b}$} 
with $\|\xi\|_{g_V,C^0},\|\xi\|_{g_V,C^0}<\de(b)$ and
the map
$$H_{(S,\eset,\{\hat{0}\};,\la_{\hat{0}}),b}\!:
\big\{\xi\!\in\!\ker D_{g_V,b}\!: 
\|\xi\|_{g_V,C^0}\!<\!\de(b)\big\}
\!\lra\! {\cal H}_{(S,\eset,\{\hat{0}\};,\la_{\hat{0}})},~~
\xi\!\lra\!\exp_{g_V,u_{\hat{0}}}
 \big(\xi\!+\!h_{(S,\eset,\{\hat{0}\};,\la_{\hat{0}}),b}(\xi)\big)$$
is an orientation-preserving
 diffeomorphism onto an open neighborhood of $b$
in ${\cal H}_{(S,\eset,\{\hat{0}\};,\la_{\hat{0}})}$.
Furthermore, the family of maps 
$\{H_{(S,\eset,\{\hat{0}\};,\la_{\hat{0}}),b}\!:
b\in{\cal H}_{(S,\eset,\{\hat{0}\};,\la_{\hat{0}})}\}$
is continuous.\\
(2)  for all
$b\!=\!\big(S,M,I;x,(j,y),u\big)\!\in\!{\cal H}_{\cal T}$\\
(2a) $D_{b,u_h}\!:\Ga(u_h)\!\lra\!\Ga^{0,1}(u_h)$ 
is onto for all $h\!\in\!\hat{I}$;\\
(2b)  $\ker D_{b,u_h}\lra T_{u_h(\i)}V$, $\xi\lra\xi(\i)$, 
is onto for all $h\!\in\!\hat{I}$.
\end{dfn}

\noindent
{\it Remarks:}
(1) All conditions in both definitions above
are independent of the choice of metric~on~$V$.\\
(2) Condition (1) of Definition~\ref{semireg_dfn} says
that ${\cal H}_{(S,\eset,\{\hat{0}\};,\la_{\hat{0}})}$ is a smooth manifold
modeled on $\ker D_b$ for 
\hbox{$b\!\in\!{\cal H}_{(S,\eset,\{\hat{0}\};,\la_{\hat{0}})}$},
as would be the case if $D_b\!:\Ga(u_b)\lra\Ga^{0,1}(u_b)$
were surjective.\\
(3) The conditions of Definitions~\ref{reg_dfn}
and~\ref{semireg_dfn} insure that ${\cal H}_{\cal T}$
is a smooth manifold; see Proposition~\ref{balanced_case} below.
However, (2) of Definition~\ref{reg_dfn} and 
(2b) of Definition~\ref{semireg_dfn} are somewhat stronger than
necessary to show that ${\cal H}_{\cal T}$ is smooth.
They allow us to obtain the second part of (1) in 
Proposition~\ref{balanced_case}, which is used in the proof
of surjectivity of the gluing map; see Subsection~\ref{surject_b}.
These two conditions hold for all complex homogeneous manifolds;
see~\cite{RT}.\\

\noindent
Note that if~${\cal T}$ is semiregular, the homotopy  
invariance of the index implies that the vector spaces
$$\Ga_-(b)\equiv \coker D_b\approx\ker D_b^*
\subset\Ga^{0,1}(b),
\quad b\in{\cal H}_{\cal T},$$
form a vector bundle over ${\cal H}_{\cal T}$.
Here $D_b^*$ denotes the formal adjoint of $D_b$
with respect to a metric $g$ on~$S$;
it is a $J$-linear operator.
The space $\ker D_b^*$ is independent of 
a conformal choice of the metric~$g$.
The bundle \hbox{$\Ga_-\!\lra\!{\cal H}_{\cal T}$} will be called
the ${\cal T}$-{\it{cokernel bundle}}.
It is $G_{\cal T}$-equivarent,
and thus descends to a bundle $\Ga_-\!\lra\!{\cal M}_{\cal T}$,
which will be the analogue of Taubes's obstruction
in our gluing setting.\\

\noindent
Let ${\cal T}\!=\!\big(S,M,I;j,\la\big)$ be a bubble type.
If $b\!=\!\big(S,M,I;x,(j,y),u\big)\!\in\!{\cal H}_{\cal T}$,~put
\begin{equation*}
\begin{split}
{\cal K}_b{\cal T}=\Big\{
\si\!=\!(\xi,w_{\hat{I}+M})\!\in\!{\cal P}_b{\cal T}\!: 
\xi_i\!\in\!\ker(D_{b,u_i})~\forall i\!\in\! I;~
\lan\si,\si_{(b,h)}^{(k)}\ran=0~\forall h\in\hat{I},k\in[4];~
\qquad\quad&\\
\xi_h(\i)=\xi_{\io_h}(x_h)+du_{\io_h}\big|_{x_h}w_h~\forall h\!\in\!\hat{I}
\Big\}.&
\end{split}\
\end{equation*}
If $\si\!=\!(\xi,w_{\hat{I}+M})\in{\cal K}_b{\cal T}$,~let
$$\|\si\|_{b,C^k}=\|\xi\|_{b,C^k}+\sum_{h\in\hat{I}}|w_h|_b+
\sum_{l\in M}|w_l|_b.$$
We take the default norm on ${\cal K}_b{\cal T}$ be given by
$\|\cdot\|_{b,C^0}$.
If $b$ is as above, \hbox{$b\!=\!\big(S,M,I;x',(j,y'),u\big)$},
and $\de\!>\!0$, we say $d(b,b')\!<\!\de$ if there exists 
$\si\!\in\!{\cal P}_b{\cal T}$ such that $\exp_b\si\!=\!b'$
and $\|\si\|_{b,C^0}\!\le\!\de$.

\begin{prp}
\label{balanced_case} 
(1) If ${\cal T}\!=\!(S,M,I;j,\la)$ 
is a regular or semiregular bubble type,
${\cal H}_{\cal T}$ is a complex manifold and 
there exists $\ep_{\cal T},C_{\cal T}\!\in\! 
 C^{\i}({\cal M}_{\cal T}^{(0)};\Bbb{R}^+)$ with
the following property.
\hbox{If $b^*\!\in\!{\cal H}_{\cal T}$ and}
$$b=\big(S,M,I;x,(j,y),u\big)
\quad\hbox{is s.t.}\quad
d(b^*,b)<\ep_{\cal T}(b^*)\hbox{~and~}
\bar{\partial}u_i=0~\forall i\!\in\! I,$$
there exist $\xi_i\!\in\!\Ga(u_i)$ for $i\in\hat{I}$
such
$$\|\xi_i\|_{g_V,C^0}\le C_{\cal T}(b^*)\sum_{h\in\hat{I}} 
d_V\big(u_{\io_h}(x_h),u_h(\i)\big)
\quad\hbox{and}\quad
b'=(S,M,I;x,(j,y),u')\in {\cal H}_{\cal T},$$ 
where $u_{\hat{0}}'=u_{\hat{0}}$ and 
$u_i'=\exp_{g_V,u_i}\xi_i$ if $i\!\in\!\hat{I}$.\\
(2) The space ${\cal M}_{\cal T}^{(0)}$ is a smooth oriented  manifold
on which the group~$G_{\cal T}$ acts smoothly.
The maps
\begin{alignat*}{2}
&\ev\!: {\cal M}_{\cal T}^{(0)}\lra V,&\qquad
&\ev\big(S,M,I;x,(j,y),u\big)=u_{\hat{0}}(\i),\\
&\ev_l\!: {\cal M}_{\cal T}^{(0)}\lra V,&\qquad
&\ev_l\big(S,M,I;x,(j,y),u\big)=u_{j_l}(y_l),\\
&du_i|_z: {\cal M}_{\cal T}^{(0)}\lra T^*\Si_{{\cal T},i}\otimes u_i^*TV,
&\qquad 
&du_i|_z\big(S,M,I;x,(j,y),u\big)=du_i|_z,
\end{alignat*}
are smooth. In particular, $u_i\lra\|du_i\|_{b,C^0}$
defines a continuous function on ${\cal M}_{\cal T}^{(0)}$.\\
(3) There exist $\de_{\cal T},C_{\cal T}\!\in \!
 C^{\i}({\cal M}_{\cal T}^{(0)};\Bbb{R}^+)$ and smooth maps 
$$h_{{\cal T},b}=h_{{\cal T},b}^{(1)}\oplus h_{{\cal T},b}^{(2)}:
{\cal K}_b{\cal T}_{\de_{\cal T}(b)}\lra\Ga'(b)\oplus 
 \l(\Bbb{C}\oplus\Bbb{R}\r)^{\hat{I}} ,$$
such that
$\big\|h_{{\cal T},b}(\si)\big\|_{b,C^0}\le C_{\cal T}(b)\|\si\|_{b,C^0}^2$, 
\begin{equation*}
\big\|h_{{\cal T},b}(\si)-h_{{\cal T},b}(\si')\big\|_{b,C^0}
\le C_{\cal T}(b)
\l(\|\si\|_{b,C^0}+\|\si'\|_{b,C^0}\r)\|\si-\si'\|_{b,C^0},
\end{equation*}
and each map
\begin{gather*}
H_{{\cal T},b}^{(0)}\!: 
\big\{(\si,\th)\!\in\!
{\cal K}_b{\cal T}_{\de_{\cal T}(b)}\times\Bbb{R}^{\hat{I}}\!: |\th|<\pi\big\}
\lra{\cal M}_{\cal T}^{(0)},\\  
H_{{\cal T},b}^{(0)}(b,\si,\th)=
\big(h_{{\cal T},b}^{(2)}(\si),\th\big)\cdot
              \exp_b\big(\si+h_{{\cal T},b}^{(1)}(\si)\big),
\end{gather*}
is  orientation-preserving diffeomorphism
onto an open neighborhood of~$b$ in~${\cal M}_{\cal T}^{(0)}$.
\end{prp}

\noindent
{\it Proof:} (1) Let 
${\cal T}_i\!=\!(\Si_{{\cal T},i},\{l\!:j_l\!=\!i\}\!+\!\{h\!:\io_h\!=\!i\},
\{\hat{0}\};\hat{0},\la_i)$.
By (1) of Definition~\ref{reg_dfn}, and (1) and (2a) 
of Definition~\ref{semireg_dfn}, 
${\cal H}_{{\cal T}_i}$ is a complex manifold
for \hbox{all $i\!\in\! I$}. Let 
$$\lap_V^{\hat{I}}=
\Big\{ (q,q)_{\hat{I}}\!\in \prod_{\hat{I}}(V\!\times\! V)\!:
q_h\!\in\! V\Big\}.$$
The submanifold $\lap_V^{\hat{I}}$ is the $\hat{I}$-product
of the diagonal \hbox{in $V\!\times\! V$}.
Since $V$ is oriented, so is the normal bundle of~$\lap_V^{\hat{I}}$.
Claim (1) of the proposition follows by applying 
the Implicit Function Theorem,
(2) of Definition~\ref{reg_dfn} and 
(2b) of Definition~\ref{semireg_dfn} to the smooth map
$$\ev_{\hat{I}}\!:
\prod_{i\in I}{\cal H}_{{\cal T}_i}\lra
            \prod_{\hat{I}}(V\!\times\! V),\quad
\ev_h\big((S,M,I;x,(j,y),u)\big)=
\big(u_h(\i),u_{\io_h}(x_h)\big).$$
Note that ${\cal H}_{\cal T}=
             \ev_{\hat{I}}^{-1}\big(\lap_V^{\hat{I}}\big)$.\\
(2) For any $u\!\in\! C^{\i}(S^2;V)$, define 
$\tilde{\Psi} u\!\in\!\Bbb{C}$, $\Psi^{(3)}u\!\in\!\Bbb{R}$,
and $\Psi u\!\in\!\Bbb{C}\!\times\!\Bbb{R}$~by
\begin{equation*}\begin{split}
\Psi u=\big(\tilde{\Psi} u,\Psi^{(3)}u\big)=
\Big(\int_{\Bbb{C}}|du\circ q_N|^2z,
                \int_{\Bbb{C}}|du\circ q_N|^2\be(|z|)-\frac{1}{2} \Big),
\end{split}\end{equation*}
where the integrals are computed using the metric~$g_V$.
For $i^*\!\in\!\hat{I}$ if $S\!=\!\Si$ and $i^*\!\in\! I$ if $S\!=\!S^2$, 
we define~maps 
\begin{gather*}
\Psi_{{\cal T},i^*}\!:
 \prod_{i\in I}{\cal H}_{{\cal T}_i}\lra\Bbb{C}\!\times\!\Bbb{R}
\qquad\hbox{by}\\
\Psi_{{\cal T},i^*}\!\!\big(S,M,I;x,(j,y),u\big)\!=\!
\Big(\! \tilde{\Psi}u_{i^*}\!\!+
\!\!\!\sum_{\io_h=i^*}\!\!d_h({\cal T})x_h\!+
\!\!\!\sum_{j_l=i^*}\! y_l,
\Psi^{(3)}u_{i^*}\!\!+
\!\!\!\sum_{\io_h=i^*}\!\!d_h({\cal T})\be(|x_h|)\!+
\!\!\!\sum_{j_l=i^*}\!\be(|y_l|)\!\Big).
\end{gather*}
These maps $\Psi_{{\cal T},i^*}$ are smooth, since the smooth structure
on all ${\cal H}_{{\cal T}_i}$ is described similarly
to (1) of Definition~\ref{semireg_dfn}.
Furthermore, if $b\!\in\!{\cal M}_{\cal T}^{(0)}$,  
$i^*\!\in\!\hat{I}$, and $k^*\!=\!1,2,3$,
since \hbox{$\Psi_{{\cal T},i}(b)\!=\!0$} for all~$i$
and $\be'$ does not change sign, by Lemma~\ref{balancing_lmm},
$$d\Psi_{{\cal T},i^*}^{(k^*)}\Big|_b\si_{(b,i)}^{(k)}~~
\begin{cases}
=0,&\hbox{if~}i\neq i^*;\\
\neq0,&\hbox{if~}i=i^*,k=k^*;\\
=0,&\hbox{if~}k\neq k^* \neq 3
\end{cases}$$
where $k=1,2,3$.
By (2) of Definition~\ref{reg_dfn} and (2b) of Definition~\ref{semireg_dfn},
it follows that the map
$$\prod_{i\in I}{\cal H}_{{\cal T}_i}\lra
\l(\Bbb{C}\times\Bbb{R}\r)^{\hat{I}}\times
\prod_{\hat{I}}(V\!\times\! V),\quad
b\lra
\Big(\big(\Psi_{{\cal T},i}(b)\big)_{i\in\hat{I}},\ev_{\hat{I}}(b)\Big)$$
is transversal to the  submanifold  $\{0\}\!\times\!\lap_V^{\hat{I}}$.
The preimage of  this submanifold
is precisely the space~${\cal M}_{\cal T}^{(0)}$.
Thus, ${\cal M}_{\cal T}^{(0)}$ is 
a smooth oriented manifold by the Implicit Function Theorem.

\begin{lmm}
\label{balancing_lmm}
For any $k\!\in\![4]$ and $u\in C^{\i}(S^2;V)$, $\xi^{(k)}(\i)=0$.
Furthermore,
\begin{gather}
\label{balancing_lmm_e1}
\tilde{\Psi}\l((c,r,\th)\cdot u\r)
   =(1+r)\big(\tilde{\Psi}u+c\|du\|_2^2\big)\quad
\forall (c,r)\in\Bbb{C}\times\Bbb{R};\\
\label{balancing_lmm_e2}
\frac{d}{dr}\Psi^{(3)}\l((0,r,\th)\cdot u\r)\Big|_{r=0}
=\int_{\Bbb{C}}|d(u\circ q_N)|^2\be'(|z|)|z|,
\end{gather}
where $(c,r)\cdot u$ is defined as in Section~\ref{summary3}.
Finally, $D_u\xi_u^{(k)}=0$ if~$\bar{\partial}u=0$.
\end{lmm}

\noindent
{\it Proof:} The first and last statements are immediate.
To prove \e_ref{balancing_lmm_e1}, we use the change of\\
\hbox{$z\lra(1+r)^{-1}z-c$.}
\begin{equation*}\begin{split}\label{balanced_str_e1}
\int_{\Bbb{C}}\big|d\big(((c,r)\cdot u)\circ q_N\big)\big|^2z
&=\int_{\Bbb{C}}(1+r)^{-2}\big|d(u\circ q_N)\big|_{(1+r)^{-1}z-c}^2z\\
&=(1+r)\int_{\Bbb{C}}\big|d\big(u\circ q_N)\big|^2(z+c)
=(1+r)\big(\tilde{\Psi}u+c\|du\|_2^2\big),
\end{split}\end{equation*}
Similarly,
\begin{equation*}\begin{split}
\frac{d}{dr}
\int_{\Bbb{C}}\big|d\big((r\cdot u)\circ q_N\big)\big|^2\be(|z|)\Big|_{r=0}
& =\frac{d}{dr}\int_{\Bbb{C}}
 \big|d(u\circ q_N)\big|^2\be\big((1+r)|z|\big)\Big|_{r=0} \\
&=\int_{\Bbb{C}}\big|d(u\circ q_N)\big|^2\be'\big(|z|\big)|z|.
\end{split}\end{equation*}
The lemma is now proved, since the action by the $\th$-component
does not change~$\tilde{\Psi}$.\\

\noindent
If ${\cal T}\!=\!(S^2,M,I;j,\la)$ is a regular 
bubble type, with notation as above, let
$$\tilde{\cal K}_b{\cal T}=\big\{
\si=(\xi_I,w_{M+\hat{I}})\!\in\! {\cal K}_b{\cal T}\!: 
\lr{\si,\si_{(b,\hat{0})}^k}=0~\forall k\!\in\![4],
~\xi_{i_1}(\i)=\xi_{i_1}(\i)~\forall i_1,i_2\!\in\! I\!-\!\hat{I}\big\}.$$
By (3) of Definition~\ref{reg_dfn} and the same argument
as in the proof of Proposition~\ref{balanced_case},
we can construct smooth maps 
$h_{{\cal T},b}^{(1)}\times h_{{\cal T},b}^{(2)}:
 \tilde{\cal K}_b{\cal T}_{\de(b)}\lra 
 \Ga'(b)\times(\Bbb{C}\times\Bbb{R})^I$
such that each map
\begin{gather*}
H_{{\cal T},b}^{(0)}: 
 \big\{(b,\si,\th)\in {\cal K}_b\tilde{\cal T}_{\de(b)}\times\Bbb{R}^I:
 |\th|<\pi\big\}\lra{\cal B}_{\cal T},\\
H_{{\cal T},b}^{(0)}(\si,\th)=(h_{{\cal T},b}^{(2)}(\si),\th)\cdot
              \exp_b\big(\si+h_{{\cal T},b}^{(1)}(\si)\big),
\end{gather*}
is orientation-preserving diffeomorphism onto
an open neighborhood of~$b$ in ${\cal B}_{\cal T}$.

\subsection{Construction of Nearly Holomorphic Bubble Maps}
\label{approx_maps_sct1}
 
\noindent
Let ${\cal T}\!=\!(S,M,I;j,\la)$ be a  simple bubble type.
In this subsection, for all $b\!\in\!{\cal M}_{\cal T}^{(0)}$
and \hbox{$\ups\!=\!(b,v_{\hat{I}})$} with 
$v_{\hat{I}}\!\in\! F_b^{(0)}{\cal T}$
sufficiently small, we construct a bubble map $b(\ups)$
with domain~$\Si_{\ups}$, where $\Si_{\ups}$ is as
Subsection~\ref{basic_gluing}.
The map $u_{b(\ups)}$ will be just the composite $u_b\circ q_{\ups}$.
We then define a Riemannian metric $g_{\ups,i}$ 
and  a nonnegative function $\rho_{\ups,i}$
on each component $\Si_{\ups,i}$ of~$\Si_{\ups}$.
The metrics will be such that the $C^0$-norm of the differential 
of $q_{\ups}$ is bounded independently of~$v_{\hat{I}}$.
The nonnegative functions are used to modify the Sobolev norms,
in such that a way that the norm of the inverse of the operator $D_{b(\ups)}$
on certain subspaces of $\Ga(b(\ups))$ is bounded independently 
of~$v_{\hat{I}}$.\\

\noindent
By Proposition~\ref{balanced_case}, ${\cal M}_{\cal T}^{(0)}$ is
a smooth manifold.
If $S\!=\!S^2$, let 
  $\de_{\cal T}\!\in\! C^{\i}({\cal M}_{\cal T}^{(0)};\Bbb{R}^+)$
be a \hbox{$\tilde{G}_{\cal T}$-invariant} function 
such that $\de_{\cal T}(b)<r_{\t_b}$ for all 
$b\!\in\!{\cal M}_{\cal T}^{(0)}$.
If $S\!=\!\Si$, let  
     \hbox{$\de_{\cal T}\!\in\! C^{\i}({\cal M}_{\cal T}^{(0)};\Bbb{R}^+)$}
be a $G_{\cal T}$-invariant function such that for~all 
 $$b=\big(\Si,M,I;x,(j,y),u\big)
\in{\cal M}_{\cal T}^{(0)}$$
(A1) $4\de_{\cal T}$ is smaller than the function~$\de$
of Lemma~\ref{si_metrics_l1a};\\
(A2) $4\de_{\cal T}(b)<r_{{\cal C}_b}g_{b,\hat{0}}$.\\
In both cases, it can be assumed that $\de_{\cal T}$
does not exceed~$\frac{1}{4}$.\\

\noindent
If $H$ is a subset of $\hat{I}$, put
\begin{equation*}\begin{split}
F^{(H)}{\cal T}&=\big\{\ups\!=\!(b,v_{\hat{I}})\!\in \!
F^{(0)}{\cal T}\!: v_h\!=\!0\hbox{~if and if~}h\!\in\! H\big\},\\
F^{H}{\cal T}&=\big\{\ups\!=\![b,v_{\hat{I}}]\!\in\! F{\cal T}\!: 
v_h\!=\!0\hbox{~if and if~}h\!\in\! H\big\}.
\end{split}\end{equation*}
For any $\ups\!=\!(b,v_{\hat{I}})\!\in\! F^{(0)}{\cal T}$, 
let $|\ups|$ denote $|\ups|_{g_b}$ if $S\!=\!\Si$.
From now on, we assume that 
\hbox{$\de\!\in\! C^{\i}({\cal M}_{\cal T}^{(0)};\Bbb{R}^+)$} is
a $G_{\cal T}$-equivarent function if $S\!=\!\Si$
and a $\tilde{G}_{\cal T}$-equivarent function if $S\!=\!S^2$
such that $8\de^{\frac{1}{2}}\le\de_{\cal T}$. 
If 
$$\ups=\big(b_{\ups},v_{\hat{I}}\big)=
\big((S,M,I;x,(j,y),u),v_{\hat{I}}\big)
\in F^{(0)}{\cal T}_{\de},$$
let $q_{\ups}: \Si_{\ups}\lra\Si_{b_{\ups}}$ be the smooth map defined in 
Subsection~\ref{basic_gluing} for
$$\ups=\big({\cal C},v_{\hat{I}}\big)=
\big((S,M,I;x,(j,y)),v_{\hat{I}}\big),$$
using the metric $g_{b_{\ups},\hat{0}}$ on $\Si$ if $S=\Si$.
Let $u_{\ups}=u_{b_{\ups}}\circ q_{\ups}$ and 
$b({\ups})=\big({\cal C}(\ups),u_{\ups}\big)$.\\

\noindent
We now define a Riemannian metric  $g_{\ups,i}$ on $\Si_{\ups,i}$
for each $i\!\in\! I(\ups)\!\subset\! I$.
Along the way, we construct a metric $g_{\ups,i}$ on $\Si_{b_{\ups},i}$
for each $i\!\in\! I$.
Suppose $i\!\in\! I$ and for all $h\!\in\!\hat{I}$ such that $\io_h\!=\!i$,
we have constructed a metric $g_{\ups,h}$ on $\Si_{b_{\ups},h}$.
For each $h\!\in\!\hat{I}$ such that $\io_h\!=\!i$ and $v_h\!\neq\!0$,
let $\tilde{g}_{\ups,i,h}$ denote the metric on
$B_{b_{\ups},h}\big(2\de(b_{\ups})^{\frac{1}{2}}\big)$
which is the pullback of the metric $g_{\ups,h}$ by the map
$$z\lra q_N\Big(\frac{\phi_{b_{\ups},h}}{v_h}\Big),
\qquad\hbox{where}\quad
\phi_{b_{\ups},h}=\begin{cases}
\phi_{\t_{b_{\ups}},h},&\hbox{if~}x_h\in S^2;\\
\phi_{\t_{b_{\ups}},g_b,h},&\hbox{if~}x_h\in\Si.
\end{cases}$$
This metric is conformal with the original metric $g_{b_{\ups},i}$ on 
$\Si_{b_{\ups},i}$,
because the maps $\phi_{b,h}$ are holomorphic on
the set $\{r_{b,h}\le\de_{\cal T}(b)\}$
and the metric $g_{\ups,h}$ 
is conformal with the standard metric on $\Bbb{C}$.
Thus, there exists a smooth positive function 
$\la_{\ups,i,h}$ such that
$\tilde{g}_{\ups,i,h}=\la^2_{\ups,i,h}g_{b_{\ups},i}$.
Let \hbox{$\la_{\ups,i}\!\in\! C^{\i}(\Si_{b_{\ups},i};\Bbb{R}^+)$} be given by
$$\la_{\ups,i}(z)=\begin{cases}
\la_{\ups,i,h}(z)+\be_{|v_h|}\l(r_{b_{\ups},h}(z)\r)
\l(1-\la_{\ups,i,h}(z)\r),& 
\hbox{if~} \io_h=i\hbox{~and~}r_{b_{\ups},h}(z)\le 2|v_h|^{\frac{1}{2}};\\
1,& \hbox{if~}r_{b_{\ups},h}(z)\ge 2|v_h|^{\frac{1}{2}}~\forall h\in\hat{I}.
\end{cases}$$
Since $I$ is a rooted tree, this procedure defines
metrics $g_{\ups,i}$ for each $i\!\in\! I(\ups)$.\\

\noindent
In addition, we define a smooth nonnegative function $\rho_{\ups,i}$ on 
$\Si_{\ups,i}$ for each $i\!\in\! I(\ups)$.
As in the previous paragraph, along the way 
we define a function $\rho_{\ups,i}$ for each $i\!\in\! I$.
Suppose $i\!\in\! I$ and for all $h\!\in\!\hat{I}$ such that $\io_h=i$,
we have constructed a smooth function $\rho_{\ups,h}$ on $\Si_{b_{\ups},h}$.
For $h\!\in\!\hat{I}$ with $\io_h=i$ and $z\!\in\!\Si_{b_{\ups},i}$
with $|z|_h\equiv r_{b_{\ups},h}(z)\le 2\de_{\cal T}(b_{\ups})$, put 
$$\rho_{\ups,i}(z)=
\begin{cases}
\rho_{\ups,h}(q_{h,v_h}z)+
\be\Big(\frac{\de_{\cal T}(b_{\ups})|z|_h}{|v_h|}\Big)
\Big\{\Big(|z|_h^2+\frac{|v_h|^2}{|z|_h^2}\Big)-
\tilde{\rho}_{\ups,h}(q_{h,v_h}z)\Big\},&
\hbox{if~}|z|_h\le\de_{\cal T}(b_{\ups});\\
\Big(|z|_h^2+\frac{|v_h|^2}{|z|_h^2}\Big)+
\be\Big(\frac{|z|_h}{\de_{\cal T}(b_{\ups})}\Big)
\Big\{1-\Big(|z|_h^2+\frac{|v_h|^2}{|z|_h^2}\Big)\Big\},&
\hbox{if~}|z|_h\ge\de_{\cal T}(b_{\ups}),
\end{cases}$$
if $v_h\!\neq\!0$,
where $q_{h,v_h}$ is defined as in Section~\ref{curves},
using the metric $g_{b_{\ups},\hat{0}}$ on $\Si$ if $S=\Si$.
If $v_h=0$ and $z$ is as above, let 
$$\rho_{\ups,i}(z)=
|z|_h^2+\be\Big(\frac{|z|_h}{\de_{\cal T}(b_{\ups})}\Big)
\big\{1-|z|_h^2\big\}.$$
If $|z|_h\ge 2\de_{\cal T}(b_{\ups})$ for all $h\in\hat{I}$ with $\io_h=i$
and $v_i\neq0$ if $i>0$, set $\rho_{\ups,i}(z)=1$.
Otherwise, let 
$$\rho_{\ups,i}(z)=
|q_S^{-1}(z)|^2+\be\big(\de_{\cal T}(b_{\ups})|q_S^{-1}(z)|\big)
\big\{1-|q_S^{-1}(z)|^2\big\}.$$
This construction defines nonnegative functions $\rho_{\ups,i}$
on $\Si_{\ups,i}$ for all $i\!\in\! I(\ups)$.\\

\noindent
We finally define norms on the spaces $\Ga(u_{\ups})$ and $\Ga^1(u_{\ups})$.
If $\eta_i\!\in\!\Ga^1(u_{\ups,i})$, put
\begin{equation}
\label{lp_norm}
2\|\eta_i\|_{\ups,p;i}=\Big(\int_{\Si_{\ups,i}}|\eta_i|^p\Big)^{\frac{1}{p}}+
\Big(\int_{\Si_{\ups,i}}\rho_{\ups,i}^{-\frac{p-2}{p}}|\eta_i|^2
 \Big)^{\frac{1}{2}},
\end{equation}
where $|\eta_i|$ and the integrals are computed with respect to the metric
$g_{\ups,i}$ on $\Si_{\ups,i}$ and $g_{V,b_{\ups}}$ on~$V$.
Denote by $\|\eta_i\|_{\ups,C^0;i}$ the $C^0$-norm of $\eta_i$ with
respect to these metrics.
If $\eta\!=\!\eta_{I(\ups)}\!\in\!\Ga^1(u_{\ups})$, let
$$\|\eta\|_{\ups,p}=\sum_{i\in I(\ups)}\|\eta_i\|_{\ups,p;i},\qquad
\|\eta\|_{\ups,C^0}=\sum_{i\in I(\ups)}\|\eta_i\|_{\ups,C^0;i}.$$
Similarly, for any $\xi_i\!\in\!\Ga(u_{\ups,i})$, put
\begin{equation}
\label{lp1_norm}
2\|\xi_i\|_{\ups,p;i}=\Big(\int_{\Si_{\ups,i}}|\xi_i|^p\Big)^{\frac{1}{p}}+
\Big(\int_{\Si_{\ups,i}}
          \rho_{\ups,i}^{-\frac{p-2}{p}}|\xi_i|^2\Big)^{\frac{1}{2}};\quad
\|\xi\|_{\ups,p,1;i}=\|\xi_i\|_{\ups,p;i}+\|\na\xi_i\|_{\ups,p;i},
\end{equation}
where we again use the metrics $g_{\ups,i}$ on $\Si_{\ups,i}$ and 
$g_{V,b_{\ups}}$ on $V$  as in \e_ref{lp_norm}.
Denote by $\|\xi_i\|_{\ups,C^0;i}$ the \hbox{$C^0$-norm} of $\xi_i$ with
respect to the metric $g_{V,b_{\ups}}$ on~$V$.
If $\xi\!=\!\xi_{I(\ups)}\!\in\!\Ga(u_{\ups})$, let
$$\|\xi\|_{\ups,p}=\sum_{i\in I(\ups)}\|\xi_i\|_{\ups,p;i},\qquad
\|\xi\|_{\ups,p,1}=\sum_{i\in I(\ups)}\|\xi_i\|_{\ups,p,1;i},\qquad
\|\xi\|_{\ups,C^0}=\sum_{i\in I(\ups)}\|\xi_i\|_{\ups,C^0;i}.$$
Note that even though the functions $\rho_{\ups,i}^{-\frac{p-2}{p}}$
have poles at the singular points of $\Si_{\ups}$,
all smooth one-forms and vector fields have finite norms defined by
\e_ref{lp_norm} and \e_ref{lp1_norm}, respectively,
\hbox{since $\frac{p-2}{p}\!<\!1$}.
We denote by $L^p_1(\ups)$ the completion of $\Ga(u_{\ups})$
with respect to the $(\ups,p,1)$-norm and by 
  $L^p(\ups)$ the completion of $\Ga^{0,1}(u_{\ups})$
with respect to the  $(\ups,p)$-norm.
Finally, 
let 
$$D_{\ups}\!: \Ga(u_{\ups})\lra\Ga^{0,1}(u_{\ups})$$ 
denote the linearization of the $\bar{\partial}$-operator  at $u_{\ups}$
with respect to the metric $g_{V,b_{\ups}}$ on~$V$.

\begin{lmm}
\label{approx_maps}
If ${\cal T}$ is a simple bubble type and $p\!>\!2$,    
there exist $\de,C\!\in\! C^{\i}({\cal M}_{\cal T}^{(0)};\Bbb{R}^+)$ 
such that for all $\ups\!\in\! F^{(0)}{\cal T}_{\de}$,\\
(1) $\|du_{\ups}\|_{\ups,C^0}\le C(b_{\ups})$ and 
$\|\bar{\partial}u_{\ups}\|_{\ups,p}\le C(b_{\ups})|\ups|^{\frac{1}{p}}$;\\
(2) $\|D_{\ups}\xi\|_{\ups,p}\le C(b_{\ups})\|\xi\|_{\ups,p,1}$
for all $\xi\!\in\!\Ga(u_{\ups})$;\\
(3) $\|\xi\|_{\ups,C^0}\le C(b_{\ups})\|\xi\|_{\ups,p,1}$
for all $\xi\!\in\!\Ga(u_{\ups})$;\\
(4) $\|\xi\|_{\ups,p,1}\le  C(b_{\ups})
 \big(\|D_{\ups}\xi\|_{\ups,p}+\|\xi\|_{\ups,p}\big)$
for all $\xi\!\in\!\Ga(u_{\ups})$.
\end{lmm}

\noindent
{\it Proof:} If $h\!\in\! I\!-\!I(\ups)$ and $S=S^2$, 
let $A_{\ups,h}^{\pm}$ be the annulus as in Subsection~\ref{basic_gluing}.
If $S\!=\!\Si$, let $A_{\ups,h}^{\pm}$ denote $A_{g_{b_{\ups}},\ups,h}^{\pm}$.
By definition of the norms, $q_{\ups}$ is an isometry outside 
of such annuli, and  by Lemma~\ref{basic_gluing_lmm}
the $C^0$-norm of $dq_{\ups}$ is bounded on such annuli
independently of~$v_{\hat{I}}$.
Thus, the first part of (1) follows from (2) of 
Proposition~\ref{balanced_case}.
Since $\rho_{\ups}\ge|v_h|$ on $A_{\ups,h}$, 
the second part of (1) follows from Lemma~\ref{basic_gluing_lmm}.
Statement (2) of the lemma is immediate from the definition
of the norms.
The last two claims are proved in the Appendix;
see Proposition~\ref{c0_bound} and~\ref{elliptic_bound}.
In fact, the $C^0$-norm of~$\xi$ is bounded by the usual
$L^p_1$-norm of~$\xi$.

\subsection{Scale of Variations}
\label{scale_sec}

\noindent
In Subsection~\ref{gluing_map}, we consider perturbations of 
the bubble maps~$\{b(\ups)\}$ in directions ``away''
from the space of such bubble maps. 
More precisely, we look at replacing $u_{\ups}$ by 
$\exp_{b_{\ups},u_{\ups}}\xi$ with $\xi$ lying in 
a certain subspace of $L^p_1(\ups)$ complementary
to ``the tangent space'' of the space of maps~$\{b(\ups)\}$.
If~${\cal T}$ is regular, one obvious candidate for such a subspace is the 
$(L^2,\ups)$-orthogonal complement of the kernel of~$D_{\ups}$.
While the construction in Subsection~\ref{gluing_map} would go through,
we would run into significant difficulty showing injectivity
and surjectivity of the gluing map; see Subsections~\ref{local_inject}
and~\ref{surj_c}.
In this subsection, we start by describing a choice of 
the complementary subspace which will work
for the purposes of Subsections~\ref{gluing_map},
\ref{local_inject}, and~\ref{surj_c}.
We then describe norms on the tangent spaces to~$F{\cal T}$
and the properties of our setup that are sufficient
to show  injectivity and surjectivity of the gluing map.\\

\noindent
Suppose $\ups\!=\!\big((S,M,I;x,(j,y),u),v\big)\!\in\!
F^{(0)}{\cal T}_{\de}$, where ${\cal T}$ is a simple bubble type
as before.
For any $\xi\!\in\!\Ga(b_{\ups})$,
define \hbox{$R_{\ups}\xi\!\in\! L^p_1(\ups)$~by}
$$\{R_{\ups}\xi\}(z)=\xi(q_{\ups}(z)).$$
Note that $R_{\ups}\xi$ is smooth outside of the $|I\!-\!I(\ups)|$ circles 
mapped by $q_{\ups}$ to the nodes of $\Si_{\ups_b}$ and is continuous
everywhere, since $\Ga(b_{\ups})$ is the set of smooth vector fields on
the components of $\Si_{b_{\ups}}$ that agree at the nodes.
It follows that $R_{\ups}\xi$ is indeed of class~$L^p_1$.
Let $\Ga_-(\ups)$ be the image of $\ker(D_{b_{\ups}})$
 under the map~$R_{\ups}$. 
This space models the ``tangent bundle'' to the space
of maps~$\{b(\ups)\}$.
Denote by $\Ga_+(\ups)$ its $(L^2,g_{\ups})$-orthogonal complement in 
$L^p_1(\ups)$.
Let $\pi_{\ups,-}$ and $\pi_{\ups,+}$ be the 
$(L^2,g_{\ups})$-orthogonal projections onto $\Ga_-(\ups)$ and $\Ga_+(\ups)$, 
respectively.\\

\noindent
With $H\!\subset\!\hat{I}$ and  $\ups\!\in\! F^{(H)}{\cal T}_{\de}$, let
\begin{gather*}
T_{\ups}F^H{\cal T}=\big\{
\vp\!=\!(\xi,w_{\hat{I}+M},\th_{\hat{I}},r_{\hat{I}-H})\!:
(\xi,w_{\hat{I}+M})\!\in\!{\cal K}_{b_{\ups}}{\cal T};~
\th_h,r_h\!\in\!\Bbb{R}\big\};\\
\tilde{T}_{\ups}F^H{\cal T}=\big\{
 (\xi,w_{\hat{I}+M},\th_{\hat{I}},r_{\hat{I}-H})\!\in\!  T_{\ups}F^H{\cal T}\!:
   w_h=0~\forall h\!\in\!H\big\}.  
\end{gather*}
Given $\vp$ as above, put
$$\|\vp\|=\|\xi\|_{b_{\ups},C^0}+
\sum_{h\in\hat{I}}|w_h|_{b_{\ups}}+
\sum_{l\in M}|w_l|_{b_{\ups}}
+\sum_{h\in\hat{I}}|\th_h|+\sum_{h\in\hat{I}- H}\!\!|r_h|.$$
If $\de_{\cal T}$ and $H_{{\cal T},b_{\ups}}^{(0)}$ are as in 
Proposition~\ref{balanced_case} and
$\|\vp\|<\de_{\cal T}(b_{\ups})$, put
\begin{gather*}
b_{\vp}\equiv\big(S,M,I;x(\vp),(j,y(\vp)),u(\vp)\big)
=H_{{\cal T},b_{\ups}}^{(0)}\big(\xi,w_{\hat{I}+M};\th_{\hat{I}}\big)
\in{\cal M}_{\cal T}^{(0)},\\
v_h(\vp)=\begin{cases}
(1+r_h)
\begin{cases}
v_h,&\hbox{if~}x_h\in S^2;\\
d\phi_{b_{\ups},h}^{-1}|_{\phi_{b_{\ups},h}x_h(\vp)}v_h,&\hbox{if~}x_h\in\Si;
\end{cases}&
\hbox{if~}h\not\in H;\\
0,& \hbox{if~}h\in H;
\end{cases}
\qquad
\ups(\vp)\equiv \big(b_{\vp},(v(\vp))_{\hat{I}}\big).
\end{gather*}
Then $\ups(\vp)\!\in\! F^{(H)}{\cal T}_{2\de}$ if 
$\|\vp\|<\de(b_{\ups})$ for some
$\de\!\in\! C^{\i}({\cal M}_{\cal T}^{(0)};\Bbb{R}^+)$
sufficiently small.
If $H\!=\!\eset$, 
\hbox{$T_{\ups}F^H{\cal T}\!=\!\tilde{T}_{\ups}F^H{\cal T}$} models 
the tangent space of $[\ups]$ in~$F^H{\cal T}$.
If $H\!\neq\!\eset$, the bundle~$F^H{\cal T}$ and
the construction in the previous subsection lift 
to a bundle ${^HF}{\cal T}$ over 
$${\cal M}_{\cal T}^H\equiv
{\cal M}_{\cal T}^{(0)}\big/
\big\{g_{\hat{I}}\!\in\! G_{\cal T}\!: g_h=1~\forall g\!\in\! H\}.$$
Then $T_{\ups}F^H{\cal T}$ models the tangent space of $[\ups]$ 
in~${^HF}{\cal T}$.
On the other hand, 
$\tilde{T}_{\ups}F^H{\cal T}$ models the tangent space of $[\ups]$
in the restriction of~${^HF}{\cal T}$ to the subspace
$$\big\{
\big[b'\!=\!(S,M,I;x',(j,y'),u)\big]\!\in\!{\cal M}_{\cal T}^H\!:
x_h'\!=\!x_h~\forall h\!\in\! H\big\}.$$
The reason for defining subspaces $\tilde{T}_{\ups}F^H{\cal T}$
is that if $x_h'\!\neq\! x_h$ for some $h'\!\in\! H$,
$b(\ups)$ and $b(\ups')$ do not have the same singular points
for all $\ups\!\in\! F^{(H)}_b{\cal T}$ and
$\ups\!\in\! F^{(H)}_{b'}{\cal T}$.
Since the perturbation construction of Subsection~\ref{gluing_map}
does not change the singular points of $b(\ups)$ and~$b(\ups')$, 
the resulting bubble maps
$\tilde{b}(\ups)$ and $\tilde{b}(\ups')$ will necessarily be different.
\\

\noindent
We now define norms on $T_{\ups}F^H{\cal T}$, which make 
the estimates in Lemma~\ref{inj_l1} dependent only on~$b_{\ups}$.
\hbox{If $h\!\in\!\hat{I}\!-\!H$,~let} 
$$w_h'=\phi_{b_{\ups},h}q_{\ups(\vp),\io_h}
\big(q_{\ups,\io_h}^{-1}(\io_h,x_h)\big)
    \in F_{h,b_{\ups}}^{(0)}
\qquad\hbox{if~~}
q_{\ups(\vp),\io_h}\l(q_{\ups,\io_h}^{-1}(\io_h,x_h)\r)\in\Si_{b_{\ups},\io_h}.$$
In such a case, let $\|\vp\|_{\ups,h}\!=\!\l|\frac{w_h'+w_h}{v_h}\r|$.
Otherwise, put \hbox{$\|\vp\|_{\ups,h}\!=\!1$}.
Let $\|\vp\|_{\ups}\!=\!\|\vp\|\!+
              \!\sum\limits_{h\in\hat{I}-H}\|\vp\|_{\ups,h}$.\\

\noindent
In order to simplify notation, we replace $\ups(\vp)$
by ${\vp}$ whenever there is no ambiguity.
If $\|\vp\|_{\ups}$ is sufficiently small, define
$\ze_{\vp}\!\in\!\Ga'(u_{\ups})$ by
$$\exp_{b_{\ups},u_{\ups}}\ze_{\vp}=u_{\vp},\quad
\|\ze_{\vp}\|_{b_{\ups},C^0}<\inj~g_{V,b_{\ups}}.$$
Similarly, for $h\!\in\! H$ and $l\!\in\! M$, define 
$w_h(\vp)\!\in\! T_{x_h(\ups)}\Si_{\ups,\io_h(\ups)}$ and 
$w_l(\vp)\!\in\! T_{y_l(\ups)}\Si_{\ups,j_l(\ups)}$ by
\begin{gather*}
\exp_{g_{\ups},x_h(\ups)}w_h(\vp)=x_h\l(\ups(\vp)\r),\quad
|w_h(\vp)|\equiv|w_h(\vp)|_{g_{\ups}}<\inj_{x_h(\ups)}g_{\ups};\\
\exp_{g_{\ups},y_l(\ups)}w_l(\vp)=y_l\l(\ups(\vp)\r),\quad
|w_l(\vp)|\equiv|w_l(\vp)|_{g_{\ups}}<\inj_{y_l(\ups)}g_{\ups}.
\end{gather*}
If $\vp\!\in\!\tilde{T}F_{\ups}^H{\cal T}$ and 
$\xi\!\in\!\Ga(u_{\ups})$, let $R_{\vp}\xi\!\in\!\Ga(u_{\vp})$
be the vector field given by
$$R_{\vp}\xi(z)=\Pi_{b_{\ups},\ze_{\vp}(z)}\xi(z).$$
Note that since $b(\ups)$ and $b(\vp)$ have the same singular
points whenever $\vp\!\in\!\tilde{T}F_{\ups}^H{\cal T}$,
$\Pi_{b_{\ups},\ze_{\vp}}$ does indeed map $\Ga(u_{\ups})$
to~$\Ga(u_{\vp})$.
If $\eta\!\in\!\Ga^1(u_{\ups})$, we define 
$R_{\vp}\eta\!\in\!\Ga^1(u_{\vp})$ similarly.
Let $S_{\vp}$ denote the inverse of~$R_{\vp}$.

\begin{lmm}
\label{inj_l1}
There exist 
$\de,C\!\in\! C^{\i}({\cal M}_{\cal T}^{(0)};\Bbb{R}^+)$
such that for  all $\ups\!\in\! F^{(H)}{\cal T}_{\de}$
and \hbox{$\vp\!\in\!\tilde{T}_{\ups}F^H{\cal T}_{\de}$},\\ 
(1) $C(b_{\ups})^{-1}\|\vp\|_{\ups}\le
\|\ze_{\vp}\|_{\ups,p,1}+\sum\limits_{h\in H}\!|w_h(\vp)|_{g_{\ups}}+
\sum\limits_{l\in M}\!|w_l(\vp)|_{g_{\ups}}
\le C(b_{\ups})\|\vp\|_{\ups}$;\\
(2) $\big\|\frac{g_{V,b_{\vp}}}{g_{V,b_{\ups}}}-1\big\|_{C^3}\le 
C(b_{\ups})\|\vp\|$, 
$\big\|\frac{g_{\vp}}{g_{\ups}}-1\big\|_{C^0}\le 
C(b_{\ups})\|\vp\|_{\ups}$ and
$\big\|\frac{\rho_{\vp}}{\rho_{\ups}}-1\big\|_{C^0}\le 
 C(b_{\ups})\|\vp\|_{\ups}$;\\
(3) $\big\|S_{\vp}du_{\vp}-du_{\ups}\big\|_{\ups,p}\le 
C(b_{\ups})\|\vp\|_{\ups}$ and
$\big\|S_{\vp}\bar{\partial}u_{\vp}-\bar{\partial}u_{\ups}
\big\|_{\ups,p}\le C(b_{\ups})|\ups|^{\frac{1}{p}}\|\vp\|_{\ups}$;\\
(4) 
$\big\|S_{\vp}\nu-\nu\big\|_{\ups,p}\le C(b_{\ups})\big\|\vp\|_{\ups}$;\\
(5) $\big\|S_{\vp}D_{\vp}R_{\vp}\xi-D_{\ups}\xi\big\|_{\ups,p}
\le C(b_{\ups})\|\vp\|_{\ups}\|\xi\|_{\ups,p,1}$  
and
$\big\|S_{\vp}\pi_{\vp,\pm}R_{\vp}\xi-
 \pi_{\ups,\pm}\xi\big\|_{\ups,p,1}\le 
C(b_{\ups})\|\vp\|_{\ups}\|\xi\|_{\ups,p,1}$ for all 
$\xi\!\in\!\Ga(u_{\ups})$.
\end{lmm}

\noindent
{\it Proof:} The first statement of (2) is clear.
Proofs of (1), the last two claims of~(2),
(3), and the last claim of~(4)
are direct, though lengthy, computations, all of the same nature.
The statement of~(4) is immediate from~(1).
The first claim of (5) follows from~(2)
and basic Riemannian geometry estimates
as~in~\cite{Z1}.\\

\noindent
{\it Remark:} The second of Claim (5) above is proved by choosing 
an orthonormal basis $\{\xi_{b,i}\}$ for the kernel of $D_b$ 
for $b$ lying near $b_{\ups}$ in ${\cal M}^{(0)}_{\cal T}$,
so that each $\xi_{b,i}$ varies smoothly with~$b$.
Then the claim follows immediately from an estimate on
\hbox{$S_{\vp}R_{\ups(\vp)}\xi_{b_{\vp},i}-R_{\ups}\xi_{b,i}$},
since the projection maps can be expressed in terms of inner-products
with~$\xi_{b,i}$.
Note that if we had defined $\Ga_-(\ups)$ to be the kernel 
of $D_{\ups}$ in the case ${\cal T}$ is regular,
this claim, if true, would have been much harder to prove
because of the presence of small eigenvalues of 
$D_{\ups}^*D_{\ups}$;
see Subsection~\ref{gluing_map} for more details.\\

\noindent
If $\ups\!\in\! F^{(0)}{\cal T}_{\de}$,
$(\Si_{\ups},g_{\ups})$ can be viewed as a connected sum
of the surfaces $\{(\Si_{{\cal T},i},g_{b_{\ups},i})\}$ with very thin necks.
If $\vp\!\in\!{\cal K}_{b_{\ups}}{\cal T}\!\subset
    \! T_{\ups}F^{\eset}{\cal T}$ 
is as above and
$|w_h|\ge 2|v_h|^{\frac{1}{2}}$, the maps
$u_{\ups}\!:\Si\lra V$ and  $u_{\vp}\!:\Si\lra V$ 
are very far apart in the $C^0$-norm
even if $\|\vp\|$ is small.
However, we can still compare the two maps 
and the various objects of Lemma~\ref{inj_l1}, appropriately defined, 
on  the corresponding direct summands.
If the gluing map of Subsection~\ref{gluing_map} is defined only 
on~$F^{(\eset)}{\cal T}_{\de}$, and not on~$F{\cal T}_{\de}$,
we need to be able to do such comparisons in order
to adjust the gluing map in the presence of constraints~$\mu$; 
see Subsection~\ref{orient2_sec}.\\

\noindent
In order to state an analogue  of Lemma~\ref{inj_l1}
with $\|\vp\|_{\ups}$ for 
\hbox{$\vp\!\in\!\tilde{T}_{\vp}F^H{\cal T}$}
replaced by $\|\vp\|$ for
$$\vp\in{\cal K}_{b_{\ups}}{\cal T}
\subset\tilde{T}_{\vp}F^{\eset}{\cal T},$$
for each $\vp\!\in\! {\cal K}_{b_{\ups}}{\cal T}_{\de(b)}$,
with $\de$ sufficiently small,
we construct a smooth map
\hbox{$\tilde{q}_{\vp}\!:(\Si_{\ups},g_{\ups})\lra(\Si_{\vp},g_{\vp})$},
which is almost an isometry.
The map will depend only on the elements \hbox{$w_h\!\in\! F_{b,h}^{(0)}$}.
The structure of the construction is similar to
the construction of the map~$q_{\ups}$ in Subsection~\ref{basic_gluing}.
For each $h\!\in\!\hat{I}$ with $\io_h\!=\!\hat{0}$,
let 
$\tilde{p}_{h,\vp}\!: 
B_{b_{\ups},h}\big(4\de_{\cal T}(b_{\ups})\big)\!\lra\!\Si$
be the (holomorphic) 
$(g_{b_{\ups},\hat{0}},g_{b_{\vp},\hat{0}},)$-isometry
provided by Lemma~\ref{si_metrics_l1a}.
Define 
\hbox{$\tilde{q}_{h,\vp}\!:\Si\lra\Si$ by}
$$\tilde{q}_{h,\vp}(z)\!=\!
\begin{cases}
\phi_{b_{\ups},h}^{-1}\Big\{
\phi_{b_{\ups},h}\tilde{p}_{h,\vp}(z)
\!+\!\be_{\de_{\cal T}^2(b_{\ups})}\big(r_{b_{\ups},h}(z)\big)
\big(\phi_{b_{\ups},h}(z)\!-\!
\phi_{b_{\ups},h}\tilde{p}_{h,\vp}(z)\big)\Big\},\!
&\hbox{if~}r_{b_{\ups},h}(z)\le 2\de_{\cal T}(b_{\ups});\\
z,\! &\hbox{if~}r_{b_{\ups},h}(z)\ge 2\de_{\cal T}(b_{\ups}).
\end{cases}$$
If $h\!\in\!\hat{I}$ and $\io_h\!\neq\!\hat{0}$, we similarly define
\hbox{$\tilde{q}_{h,(x_h,w_h)}\!:\Si_{b,\io_h}\lra\Si_{b,\io_h}$ by}
$$
\tilde{q}_{h,\vp}(z)=
\begin{cases}
\phi_{b_{\ups},h}^{-1}\Big\{
\phi_{b_{\ups},h}(z)+w_h
-\be_{\de_{\cal T}^2(b_{\ups})}\big(r_{b_{\ups},h}(z)\big)w_h\Big\},
&\hbox{if~}r_{b_{\ups},h}(z)\le 2\de_{\cal T}(b_{\ups});\\
z,&\hbox{if~}r_{b_{\ups},h}(z)\ge 2\de_{\cal T}(b_{\ups}).
\end{cases}$$
Let $\tilde{q}_{\vp,\hat{0}}=Id_{\Si}$.
If $h\!\in\!\hat{I}$ and 
$\tilde{q}_{\vp,\io_h}\!:\Si\lra\Si$ has been constructed,
let 
$$\tilde{q}_{\vp,h}(z)=
\begin{cases}     
 q_{\vp,\io_ h}^{-1}\big(\tilde{p}_{h,\vp}(z)\big(
     q_{\vp,\io_h}(\tilde{q}_{\vp,\io_h}(z))\big)\big),
&\hbox{if~}
r_{b_{\ups,h},h}\big(q_{\vp,\io_h}(z)\big)\le 2\de_{\cal T}(b_{\ups});\\ 
\tilde{q}_{\vp,\io_h}(z),&\hbox{if~}
r_{b_{\ups,h},h}\big(q_{\vp,\io_h}(z)\big)\ge 2\de_{\cal T}(b_{\ups}).
\end{cases}$$
Going through all of $I$, we obtain a map 
$\tilde{q}_{\vp}\!:\Si\lra\Si$,
which shifts the connect-summands of $(\Si,g_{\ups})$
to the connect-summands of $(\Si,g_{\vp})$.
The important properties of such maps $\tilde{q}_{\vp}$
as summarized below.

\begin{lmm}
\label{shifts_lmm}
There exist $\de,C\!\in\! C^{\i}({\cal M}_{\cal T}^{(0)};\Bbb{R})$
and a smooth family  of maps
$$\big\{ \tilde{q}_{\vp}\!:\Si\lra\Si~|~
\vp\!\in \!
{\cal K}_{b_{\ups}}{\cal T}_{\de(b_{\ups})}\!\subset\!
T_{\ups}F_{\ups}^{(\eset)}{\cal T},~
\ups\!\in\! F_{b_{\ups}}^{(\eset)}{\cal T}_{\de(b_{\ups})}\big\},
\qquad \hbox{such that}$$
(1) $\tilde{q}_{0}=Id_{\Si}$ and 
$q_{\ups}=q_{\vp}\circ\tilde{q}_{\vp}\in 
\Si_{b_{\ups},i}^*=\Si_{b_{\vp},i}^*$ 
outside of the annuli
$$A_{\vp}=\tilde{q}_{\vp,\io_h}^{-1}
q_{\vp,\io_h}^{-1}
\Big( \big\{z\!\in\!\Si_{b_{\ups},\io_h}\!: 
\de_{\cal T}(b_{\ups})\le r_{b,h}(z)\le 2\de_{\cal T}(b_{\ups})\}
\Big),$$
which contain no marked points of $b(\ups)$ or $b(\vp)$.\\
(2) $\Big|
\frac{\tilde{q}_{\vp}^*g_{\vp}}{g_{\ups}}-
\frac{\tilde{q}_{\vp'}^*g_{\vp'}}{g_{\ups}}\Big|
\le C(b_{\ups})\|\vp-\vp'\|$
for all $\vp,\vp'\in {\cal K}_{b_{\ups}}{\cal T}_{\de(b_{\ups})}$.
\end{lmm}

\noindent
These maps $q_{\vp}$ allow us to compare 
operators on vector fields and one-forms
on $(\Si_{\ups},u_{\ups})$  and $(\Si_{\vp},u_{\vp})$
whenever $\|\vp\|$ is sufficiently small.
Define $\ze_{\vp}'\!\in\!\Ga(u_{\ups})$ by
$$\exp_{b_{\ups},u_{\ups}}\ze_{\vp}'=
u_{\vp}\circ \tilde{q}_{\vp},
\quad\|\ze_{\vp}'\|_{b_{\ups},C^0}\le\inj~g_{b_{\ups}}.$$
For $\xi\!\in\!\Ga(u_{\ups})$, let 
$R_{\vp}'\xi\!\in\!\Ga(u_{\vp})$ be given by 
$$\{R_{\vp}'\xi\}(z)
=\Pi_{b_{\ups},\ze_{\vp}'(\tilde{q}_{\vp}^{-1}(z))}
\xi\big(\tilde{q}_{\vp}^{-1}(z)\big).$$
Similarly, of $\eta\!\in\!\Ga^{0,1}(u_{\ups})$, let 
$R_{\vp}'\eta\in\Ga^{0,1}(u_{\vp})$ be given by
$$\{R_{\vp}'\eta\}\big|_z
=\Pi_{b_{\ups},\ze_{\vp}'(\tilde{q}_{\vp}^{-1}(z))}
\circ
\eta\big|_{\tilde{q}_{\vp}^{-1}(z)}
\circ\partial\tilde{q}_{\vp}^{-1}\big|_z.$$
Denote by $S_{\vp}'$ the inverse of $R_{\vp}'$.
Similarly to Lemma~\ref{inj_l1}, we have

\begin{lmm}
\label{inj_l1c}
There exist 
$\de,C\!\in\! C^{\i}({\cal M}_{\cal T}^{(0)};\Bbb{R}^+)$
such that for  all $\ups\!\in\! F^{(\eset)}{\cal T}_{\de}$
and \hbox{$\vp\!\in\!{\cal K}_{b_{\ups}}{\cal T}
\!\subset\!\tilde{T}_{\ups}F^{\eset}{\cal T}_{\de}$},\\ 
(1) $C(b_{\ups})^{-1}\|\vp\|\le
\|\ze_{\vp}'\|_{\ups,p,1}+
\sum\limits_{l\in M}\!|w_l(\vp)|_{g_{\ups}}
\le C(b_{\ups})\|\vp\|$;\\
(2) $\big\|S_{\vp}'du_{\vp}-du_{\ups}\big\|_{\ups,p}\le 
C(b_{\ups})\|\vp\|$ and
$\big\|S_{\vp}'\bar{\partial}u_{\vp}-\bar{\partial}u_{\ups}
\big\|_{\ups,p}\le C(b_{\ups})|\ups|^{\frac{1}{p}}\|\vp\|$;\\
(3) 
$\big\|S_{\vp}'\nu-\nu\big\|_{\ups,p}\le C(b_{\ups})\big\|\vp\|$;\\
(4) $\big\|S_{\vp}'D_{\vp}R_{\vp}'\xi-D_{\ups}\xi\big\|_{\ups,p}
\le C(b_{\ups})\|\vp\|\|\xi\|_{\ups,p,1}$  
and
$\big\|S_{\vp}'\pi_{\vp,\pm}R_{\vp}'\xi-
 \pi_{\ups,\pm}\xi\big\|_{\ups,p,1}\le 
C(b_{\ups})\|\vp\|\|\xi\|_{\ups,p,1}$ for all 
$\xi\!\in\!\Ga(u_{\ups})$.

\end{lmm}

\subsection{Obstruction Bundle Setup}
\label{obs_setup}

\noindent
In the next subsection, we look for solution of
the equation \hbox{$\bar{\partial}\exp_{b_{\ups},u_{\ups}}\xi=t\nu$}
with $\xi$ lying in a fixed complement of $\Ga_-(\ups)$.
If $t$ is sufficiently small, we are able to solve this equation
up to an element of a vector bundle
of the same rank as the dimension of $\Ga_-(b_{\ups})$,
called obstruction bundle.
This element is the obstruction to solving the equation.
There are choices to be made for this obstruction bundle
as well as for the subspace complementary to~$\Ga_-(\ups)$.
We describe in this subsection what conditions these choices
must satisfy for the gluing construction to work properly.\\

\noindent
If $b^*\!=\!\big(S,M,I;x^*,(j,y^*),u^*\big)
   \!\in\!{\cal M}_{\cal T}^{(0)}$
and $b\!=\!\big(S,M,I;x,(j,y),u_I\big)\!=\!H_{{\cal T},b^*}(\si,\th)$
for some \hbox{$\si\!\in\!{\cal K}_{b^*}{\cal T}$} and
$\th\!\in\!\Bbb{R}^{\hat{I}}$, let 
$\xi_{b^*,b}\!=\!\xi_{b^*,b,I}\!\in\!\Ga'(b)$ be given~by
$$\exp_{b^*,u_i^*}\xi_{b^*,b,i}=u_i,\quad
\big\|\xi_{b^*,b,i}\big\|_{C^0}<\inj~g_{V,b^*}.$$
Let $\Pi_{b^*,b}=\Pi_{b^*,\xi_{b^*,b}}$.

\begin{dfn}
\label{C0conv_dfn}
Suppose 
$b^*\!=\!\big(S,M,I;x^*,(j,y^*),u^*\big)$,
$b_k\!=\!\big(S,M,I;x_k,(j,y_k),u_k\big)
                     \!\in\!{\cal M}_{\cal T}^{(0)}\!$,
and $\ups_k\!=\!(b_k,v_k)\!\in\! F^{(0)}{\cal T}$
are such that the sequences $\{b_k\}$ and $\{|\ups_k|_{b_k}\}$
converge to $b^*\!\in\!{\cal M}_{\cal T}^{(0)}$
and $0\!\in\!\Bbb{R}$, respectively.\\
(1) The sequence $\{\xi_k\!\in\! L^p_1(u_{\ups_k})\}$
\under{$C^0$-converges} to $\xi^*\!\in\!\Ga'(b^*)$ if \\
(1a) the sequence 
$\big\{\Pi_{b^*,b_k}^{-1}(\xi_k\circ q_{\ups_k}^{-1})\big\}$
$C^0$-converges to $\xi^*$ on compact subsets of~$\Si_{b^*}^*$;\\
(1b) there exists $C>0$ such that 
$\|\xi_k\|_{\ups_k,p,1}<C$ for all~$k$.\\
(2) The sequence of subspaces $\{V_k\!\subset\!\Ga(u_{\ups_k})\}$
\under{$C^0$-converges} to  subspace $V^*\!\subset\!\Ga(b^*)$
if there exists a sequence of bases
$\big\{\{\xi_{k,i}\}_{i=1}^{i=N}\!\subset\! V_k\big\}$
such that\\
(2a) for each $i$ fixed, 
the sequence $\{\xi_{k,i}\}$ $C^0$-converges to 
some $\xi_i^*\!\in\! V^*$;\\
(2b) the set $\{\xi_i^*\}$ has cardinality $N$ and is basis for $V^*$.
\end{dfn}

\begin{lmm}
\label{C0L2}
If the sequence $\{\ups_k\}\subset F^{(0)}{\cal T}$ converges
to $b^*\!\in\!{\cal M}_{\cal T}^{(0)}$
and the sequences $\{\xi_k\!\in\! L^p_1(\ups_k)\}$
and $\{\tilde{\xi}_k\!\in\! L^p_1(\ups_k)\}$
converge to $\xi^*\!\in\!\Ga'(b^*)$ and $\tilde{\xi}^*\!\in\!\Ga'(b^*)$,
respectively,
$$\lim_{k\lra\i}\llan \xi_k,\tilde{\xi}_k\rran_{\ups_k,2}
=\llan \xi^*,\tilde{\xi}^*\rran_{\ups^*,2}.$$
\end{lmm}

\noindent
{\it Proof:}
If $\ups_k\!\lra\! b^*$, the metrics $g_{V,b_{\ups_k}}$ on $V$
and $g_{b_{\ups_k},i}$ on $\Si_{{\cal T},i}$
$C^0$-converge to $g_{V,b^*}$ and $g_{b^*,i}$, respectively.
On the other hand, by (1b) of Definition~\ref{C0conv_dfn}
and (2) of Lemma~\ref{approx_maps}, there exists $C\!>\!0$
such~that
$$\|\xi_k\|_{\ups_k,C^0},
\|\tilde{\xi}_k\|_{\ups_k,C^0}<C\qquad\forall k.$$
Thus, the claim follows from (1a) of  Definition~\ref{C0conv_dfn}.

\begin{dfn}
\label{obs_setup_dfn1}
Suppose $\Om$ is an open subset of $F^{(\eset)}{\cal T}$
such that $b(\ups)$ is defined for all $\ups\!\in\!\Om$.
A $\big(G_{\cal T}\!\times\!\hbox{Aut}({\cal T})\big)$-invariant 
smooth complex subbundle 
$\tilde{\Ga}_-\!\lra\!\Om$ of the Banach bundle $L^p_1\!\lra\!\Om$
is a \under{tangent-space model} over $\Om$ if\\
(1) for every sequence $\{\ups_k\}\!\subset\!\Om$
converging to $b^*\!\in\!{\cal M}_{\cal T}^{(0)}$,
a subsequence of $\{\tilde{\Ga}_-(\ups_k)\}$ 
$C^0$-converges to a subspace $V^*\subset\Ga(b)$
such that $\pi_{b,-}\!:V^*\lra\Ga_-(b^*)$ 
is an isomorphism;\\
(2) if $\bar{\pi}_{\ups,-}\!:L^p_1(\ups)\!\lra\!\tilde{\Ga}_-(\ups)$
is the $(L^2,\ups)$-orthogonal projection,
there exist $\de,C\!\in\! C^{\i}({\cal M}_{\cal T}^{(0)};\Bbb{R}^+)$
such that for all $\ups\!\in\!\Om_{\de}$
and all $\xi\!\in\!\Ga(u_{\ups})$,\\
(2a)
$\big\|S_{\vp}\bar{\pi}_{\vp,-}R_{\vp}\xi-
\bar{\pi}_{\ups,-}\xi\big\|_{\ups,2}\le
C(b_{\ups})\|\vp\|_{\ups}\|\xi\|_{\ups,p,1}$
for all $\vp\!\in\! T_{\ups}F^{(\eset)}{\cal T}_{\de(b_{\nu})}$;\\
(2b) $\big\|S_{\vp}'\bar{\pi}_{\vp,-}R_{\vp}'\xi-
\bar{\pi}_{\ups,-}\xi\big\|_{\ups,2}\le
C(b_{\ups})\|\vp\|\|\xi\|_{\ups,p,1}$
for all $\vp\!\in\! {\cal K}_{b_{\ups}}{\cal T}\!\subset\! 
             T_{\ups}F^{(\eset)}{\cal T}_{\de(b_{\nu})}$.
\end{dfn}

\noindent
One example of a tangent-space model is 
$\{\Ga_-(\ups)\!:\ups\!\in\! F^{(\eset)}{\cal T}_{\de}\}$.
In such a case, the limit $V^*$ in (1) of Definition~\ref{obs_setup_dfn1}
is $\Ga_-(b^*)$ and thus depends only on $b^*$, and 
not on the sequence~$\{\ups_k\}$.
However, for computational reasons, it is sometimes advantageous
to work with other choices.
With the choices in~\cite{Z2}, 
the limit $V^*$ in (1) of Definition~\ref{obs_setup_dfn1}
in fact depends on the sequence \hbox{unless $|\hat{I}|=1$}.
The following lemma collects some of the implications 
of (1) of Definition~\ref{obs_setup_dfn1}.
Condition (2) is needed in Subsections~\ref{local_inject} 
and~\ref{surj_c}.
For any tangent space model over $\Om$ and $\ups\!\in\!\Om$, 
we denote the $(L^2,\ups)$-orthogonal complement of
$\tilde{\Ga}_-(\ups)$ by~$\tilde{\Ga}_+(\ups)$.
Write $\tilde{\Ga}_+^{0,1}(\ups)$ for the image 
of $\tilde{\Ga}_+(\ups)$ under the operator~$\tilde{D}_{\ups}$.

\begin{lmm}
\label{tangent_l1}
Let $\tilde{\Ga}_-\lra\Om$ be a tangent-space model.
Then there exist $C,\de\!\in\! C^{\i}({\cal M}_{\cal T}^{(0)};\Bbb{R})$
such that for all $\ups\!\in\!\Om_{\de}$\\
(1a) $\|\xi\|_{\ups,p,1}\le C(b_{\ups})\|\xi\|_{\ups,2}$
for all $\xi\!\in\!\tilde{\Ga}_-(\ups)$;\\
(1b) $\|\bar{\pi}_{\ups,-}\xi\|_{\ups,p,1}\le
C(b_{\ups})\|\xi\|_{\ups,p,1}$ for all $\xi\!\in\!\Ga(u_{\ups})$;\\
(2a) $L^p_1(\ups)=\Ga_-(\ups)\oplus\tilde{\Ga}_+(\ups)$;\\
(2b) if $\tilde{\pi}_-$ and $\tilde{\pi}_+$ are the projection maps
corresponding to the above decomposition,
$$\|\tilde{\pi}_{\ups,\pm}\xi\|_{\ups,p,1}\le
     C(b_{\ups})\|\xi\|_{\ups,p,1}\qquad
\forall\xi\!\in\!\Ga(u_{\ups}).$$
\end{lmm}

\noindent
{\it Proof:}
(1) Suppose there exists a sequence $\{\ups_k\!\in\!\Om\}$
converging to $b^*\!\in\!{\cal M}_{\cal T}^{(0)}$
and a sequence \hbox{$\{\xi_k\!\in\!\tilde{\Ga}_-(\ups_k)\}$}
such that $\|\xi_k\|_{\ups_k,p,1}\!=\!1$, while 
$\|\xi\|_{\ups_k,2}\!\lra\!0$.
Since $\|\xi_k\|_{\ups_k,p,1}\!=\!1$,
by (2) of Lemma~\ref{approx_maps2} and (1) of 
Definition~\ref{obs_setup_dfn1}, 
a subsequence of $\{\xi_k\}$
\hbox{$C^0$-converges} to some nonzero \hbox{$\xi^*\!\in\!\Ga(b^*)$}.
However, since $\|\xi_k\|_{\ups_k,2}\lra0$,
$\|\xi^*\|_{b^*,2}=0$ by Lemma~\ref{C0L2}.
This is a contradiction, and thus (1a) holds.
Claim (1b) is an immediate consequence of (1a)
and (2) of Lemma~\ref{approx_maps2}.\\
(2)  Claim (2a) is equivalent to saying that no nonzero element
of $\tilde{\Ga}_-(\ups)$ is orthogonal to~$\Ga_-(\ups)$.
So, suppose $\ups_k\!\lra\! b^*\!\in\!{\cal M}_{\cal T}^{(0)}$
and $\{\xi_k\!\in\!\tilde{\Ga}_-(\ups_k)\}$ is such that
 $\xi_k$ is orthogonal to~$\Ga_-(\ups)$ and 
\hbox{$\|\xi_k\|_{\ups_k,p,1}\!=\!1$}.
Since $\xi_k\!\in\!\tilde{\Ga}_-(\ups_k)$ and $\|\xi_k\|_{\ups_k,p,1}=1$,
by (1) of Definition~\ref{obs_setup_dfn1},
a subsequence of $\{\xi_k\}$ converges to some nonzero  
$\xi^*\!\in\!\Ga(b^*)$.
By Lemma~\ref{C0L2},
$\xi^*$ is orthogonal to~$\Ga_-(\ups)$.
However, this contradicts the second part of (1) of 
Definition~\ref{obs_setup_dfn1}.\\
(3) Due to (1b),
Claim (2b) is equivalent to saying that there exist
$C,\de\!\in\! C^{\i}({\cal M}_{\cal T};\Bbb{R})$
such~that 
$$\|\xi\|_{\ups,p,1}\le C(b_{\ups})
\|\bar{\pi}_{\ups,-}\xi\|_{\ups,p,1}
\qquad\forall\ups\!\in\!\Om_{\de}\hbox{~~and~~}
\xi\!\in\!\Ga_-(\ups).$$
Suppose there exists a sequence $\{\ups_k\}\!\subset\!\Om$
converging to some  $b^*\!\in\!{\cal M}_{\cal T}^{(0)}$
and a sequence \hbox{$\{\xi_k\!\in\!\Ga_-(\ups_k)\}$}
such that $\|\bar{\pi}_{\ups_k,-}\xi_k\|_{\ups_k,2}\lra0$, while
$\|\xi_k\|_{\ups_k,p,1}=1$.
By Definition~\ref{obs_setup_dfn1}, a subsequence of
$\{\tilde{\Ga}_-(\ups_k)\}$ converges to a subspace~$V\subset\Ga(b)$.
On the other hand, a subsequence  of $\{\xi_k\}$
$C^0$-converges to a nonzero element $\xi^*\!\in\!\Ga_-(b^*)$,
which must be orthogonal to $V$ by Lemma~\ref{C0L2}.
This contradicts the second part of (1) of Definition~\ref{obs_setup_dfn1}.

\begin{dfn}
\label{obs_setup_dfn2}
Suppose $\Om$ is an open subset of $F^{(\eset)}{\cal T}$
such that $b(\ups)$ is defined for all $\ups\!\in\!\Om$.
A $\big(G_{\cal T}\!\times\!\hbox{Aut}({\cal T})\big)$-invariant
 smooth complex subbundle
$\Ga_-^{0,1}(\ups)\!\lra\!\Om$ of the Banach bundle $L^p\!\lra\!\Om$
with the same rank as \hbox{$\Ga_-^{0,1}\!\lra\!{\cal M}_{\cal T}^{(0)}$}
is an \under{obstruction bundle} if\\
(1) there exists $C\!\in\! C^{\i}({\cal M}_{\cal T}^{(0)};\Bbb{R})$
such that 
$$\|\eta\|_{\ups,p}\le C(b_{\ups})\|\eta\|_2
\quad\hbox{and}\quad
\|D_{\ups}^*\eta\|_{\ups,1}\le 
C(b_{\ups})|\ups|^{\frac{1}{p}}
\qquad\forall \ups\!\in\!\Om,~\eta\!\in\!\Ga_-^{0,1}(\ups);$$
(2) if $\pi_{\ups,-}^{0,1}\!:L^p(\ups)\!\lra\!\Ga_-^{0,1}(\ups)$
is the $(L^2,\ups)$-orthogonal projection,
there exists $\de\!\in\! C^{\i}({\cal M}_{\cal T}^{(0)};\Bbb{R}^+)$
such that for all $\ups\!\in\!\Om_{\de}$
and all $\eta\!\in\!\Ga^{0,1}(u_{\ups})$,\\
(2a)
$\big\|S_{\vp}\pi_{\vp,-}^{0,1}R_{\vp}\eta-
\pi_{\ups,-}^{0,1}\eta\big\|_{\ups,2}\le
C(b_{\ups})\|\vp\|_{\ups}\|\xi\|_{\ups,p}$
for all $\vp\!\in\! T_{\ups}F^{(\eset)}{\cal T}_{\de(b_{\nu})}$;\\
(2b) $\big\|S_{\vp}'\pi_{\vp,-}^{0,1}R_{\vp}'\eta-
\pi_{\ups,-}^{0,1}\eta\big\|_{\ups,2}\le
C(b_{\ups})\|\vp\|\|\xi\|_{\ups,p}$
for all $\vp\!\in\! {\cal K}_{b_{\ups}}{\cal T}\!\subset\! 
             T_{\ups}F^{(\eset)}{\cal T}_{\de(b_{\nu})}$.
\end{dfn}
 
\noindent
Such an obstruction bundle is related to 
the cokernel bundle $\Ga_-^{0,1}\!\lra\!{\cal M}_{\cal T}^{(0)}$.
However, if \hbox{$\hat{I}\neq\eset$}, 
the low-eigenspaces of $D_{\ups}D_{\ups}^*$
are too large to form an obstruction bundle; 
see Remark below.
Examples of bundles that satisfy Definition~\ref{obs_setup_dfn2}
can be found in~\cite{Z2}.
Given such an obstruction bundle, we denote by
$\pi_{\ups,-}^{0,1}$ and $\pi_{\ups,+}^{0,1}$
the $(L^2,\ups)$-orthogonal projection maps of $L^p(\ups)$
onto $\Ga_-^{0,1}(\ups)$ and its $(L^2,\ups)$-orthogonal 
complement~$\Ga_+^{0,1}(\ups)$, respectively.
The following lemma is clear from (1) of Definition~\ref{obs_setup_dfn2}.

\begin{lmm}
\label{obs_lmm1}
If $\Ga_-^{0,1}\lra\Om$ is an obstruction bundle,
there exists $C\!\in\! C^{\i}({\cal M}_{\cal T};\Bbb{R})$
such~that 
$$\|\pi^{0,1}_{\ups,\pm}\eta\|_{\ups,p}\le C(b_{\ups})\|\eta\|_{\ups,p}
\qquad\forall \ups\!\in\!\Om,~
\eta\!\in\!\Ga^{0,1}(u_{\ups}).$$
\end{lmm}

\begin{dfn}
\label{obs_setup_main_dfn}
If ${\cal T}$ is a semiregular bubble type, 
an \under{obstruction bundle setup} for $(V,J,{\cal T})$
is a tuple $(\de,\tilde{\Ga}_-,\Ga_-^{0,1},R)$, where\\
(1) $\de\!\in\!C^{\i}({\cal M}_{\cal T}^{(0)};\Bbb{R}^+)$ is
$\big(G_{\cal T}\!\times\!\hbox{Aut}({\cal T})\big)$-invariant 
and
$b(\ups)$ is defined for all $\ups\!\in\! F^{(0)}{\cal T}_{\de}$;\\
(2) $\tilde{\Ga}_-\!\!\lra\! F^{(\eset)}{\cal T}_{\de}$
and $\Ga_-^{0,1}\!\!\lra\! F^{(\eset)}{\cal T}_{\de}$
are a tangent-bundle model and an obstruction bundle,
respectively;\\
(3) $R\!:\pi^*\Ga^{0,1}_-\!\!\lra\!\Ga^{0,1}_-$ is
a smooth oriented 
$\big(G_{\cal T}\!\times\!\hbox{Aut}({\cal T})\big)$-equivarient 
bundle isomorphism over~$F^{(0)}{\cal T}_{\de}$,
where \hbox{$\pi\!:F^{(0)}{\cal T}_{\de}\lra{\cal M}_{\cal T}^{(0)}$}
is the bundle projection map.
\end{dfn}

\noindent
For the rest of the paper, we fix such an obstruction bundle
setup.
However, whenever we refer to
\hbox{$\de\!\in\! C^{\i}({\cal M}_{\cal T};\Bbb{R}^+)$},
we will mean any function smaller than the function~$\de$
in Definition~\ref{obs_setup_main_dfn}.
The following lemma states some of the consequences of our setup
that are crucial for the construction of the next subsection.
If ${\cal T}$ is a regular bubble type,
we take $\tilde{\Ga}_-(\ups)$ and $\Ga_-^{0,1}(\ups)$  to be
$\Ga_-(\ups)$ and~$\{0\}$, respectively,
and define the other bundles and the projection maps in the same~way.

\begin{lmm}
\label{approx_maps2}
If ${\cal T}$ is a simple bubble type,
there exist $\de,C\!\in\! C^{\i}({\cal M}_{\cal T}^{(0)};\Bbb{R}^+)$ 
such that for any \hbox{$\ups\!\in\! F^{(0)}{\cal T}_{\de}$} 
if ${\cal T}$ is regular
and any $\ups\!\in\! F^{(\eset)}{\cal T}_{\de}$ if ${\cal T}$ 
is semiregular,\\
(1) $\|\xi\|_{\ups,p,1}\le C(b_{\ups})\|D_{\ups}\xi\|_{\ups,p}$
for all $\xi\!\in\!\Ga_+(\ups)$ and
all $\xi\!\in\!\tilde{\Ga}_+(\ups)$;\\
(2) $\|\pi^{0,1}_{\ups,-}\eta\|_{\ups,p}\le C(b_{\ups})|\ups|^{\frac{1}{p}}
        \|\eta\|_{\ups,p}$ for all $\eta\!\in\!\tilde{\Ga}_+^{0,1}(\ups)$;\\
(3) $\pi^{0,1}_{\ups,+}\!: \tilde{\Ga}_+^{0,1}(\ups)\lra\Ga^{0,1}_+(\ups)$
is an isomorphism with the norm of the inverse bounded~by~$C(b_{\ups})$.
\end{lmm}

\noindent
{\it Proof:} (1) The first statement of the lemma is proved in the Appendix;
see Proposition~\ref{inverse_p}.
It is consequence of (2) and (4) of Lemma~\ref{approx_maps}
and of (1) of Definition~\ref{obs_setup_dfn1}.
The second claim is immediate from (1) of Definition~\ref{obs_setup_dfn2}
and the first claim.\\
(2) Let $W$ be the $(L^2,g_{\ups})$-orthogonal complement of
$\pi_{\ups,+}^{0,1}(\tilde{\Ga}_+^{0,1}(\ups))$ in $\Ga_+^{0,1}(\ups)$.
The second claim implies~that 
\begin{equation}
\label{decomp_e1}
L^p(\ups)=\big(\Ga_-^{0,1}(\ups)\oplus W\big)\oplus\tilde{\Ga}_+^{0,1}(\ups).
\end{equation}
Since $\tilde{\Ga}_+^{0,1}(\ups)$ is the image of $\tilde{\Ga}_+(\ups)$ 
under $D_{\ups}$, with respect to 
the decompositions \e_ref{decomp_e1} and
\hbox{$L^p_1(\ups)=\Ga_-(\ups)\oplus\tilde{\Ga}_+(\ups)$,} 
$$D_{\ups}=\l|\begin{array}{cc}
D_{\ups}^{(--)}& 0\\
D_{\ups}^{(+-)}& D_{\ups}^{(++)}
\end{array}\r|.$$
Since $D_{\ups}^{(++)}$ is an isomorphism by (1) of the lemma,
\begin{equation}
\label{index_calc1}
\begin{split}
\ind~D_{\ups}&=\ind~D_{\ups}^{(--)}=
\dim~\Ga_-(\ups)-
   \big(\dim~\Ga_-^{0,1}(\ups)+\dim~W\big)\\
&=\big(\dim~\Ga_-(b_{\ups})-\big(\dim~\Ga_-^{0,1}(b_{\ups})\big)+\dim~W
=\ind~D_{b_{\ups}}+\dim~W.
\end{split}
\end{equation}
On the other hand, by the Index Theorem, with $n\!=\!\dim_{\Bbb{C}}V$,
\begin{equation}\begin{split}
\label{index_calc2}
\ind~D_{\ups}&=2
\Big(\sum_{h\in\hat{I}(\ups)}\big(
\lan c_1(V,J),\la_i(\ups)\ran-n(g(\Si_{{\cal T},i})-1)\big)
-n\big(|\hat{I}(\ups)|-1\big)\Big)\\
&=2\Big(
\sum_{h\in\hat{I}(\ups)} \lan c_1(V,J),\la_i\ran- 
n(g(S)-1)\big)=\ind~D_{b_{\ups}}.
\end{split}\end{equation}
By equations \e_ref{index_calc1} and \e_ref{index_calc2},
$W\!=\!\{0\}$, and the last claim of the lemma follows from the second~one.\\

\noindent
{\it Remark:} It is essential for claim (1) of Lemma~\ref{approx_maps2}
that $p\!>\!2$. 
The operator $D_{\ups}^*D_{\ups}$ has at least
$|\hat{I}|(\dim~V)$ eigenvalues that tend to $0$ as $|\ups|\lra 0$.
The corresponding eigenfunctions converge to vector fields on
the components of $\Si_b$ that do not agree at the nodes.
If ${\cal T}$ is semiregular, the operator $D_{b_{\ups}}$ 
has cokernel~$\Ga_-^{0,1}(b)$.
In such a case, the number of low eigenvalues of 
$D_{\ups}^*D_{\ups}$, including~$0$,
is \hbox{$(\dim~\Ga_-^{0,1}(b))+|\hat{I}|(\dim~V)$}.\\

\noindent
Let $\tilde{\pi}_{\ups,+}^{0,1}\!: \Ga^{0,1}_+(\ups)
      \lra\tilde{\Ga}^{0,1}_+(\ups)$
denote the inverse of 
$\pi_{\ups,+}^{0,1}\!:\tilde{\Ga}^{0,1}_+(\ups)\lra\Ga^{0,1}_+(\ups)$.
We extend $\tilde{\pi}_{\ups,+}^{0,1}$ to all of 
$L^p(\ups)$ by taking it to be
$\tilde{\pi}_{\ups,+}^{0,1}\circ\pi_{\ups,+}^{0,1}$.
If $\eta\!\in\!\tilde{\Ga}^{0,1}_+(\ups)$, let 
 $P_{\ups}\eta\!\in\!\tilde{\Ga}_+(\ups)$ 
be the unique element such that $D_{\ups}P_{\ups}\eta=\eta$.
We extend $P_{\ups}$ to all of $L^p(\ups)$ by taking it to be
$P_{\ups}\circ\tilde{\pi}_{\ups,+}^{0,1}$.  
From Lemma~\ref{approx_maps2}, we immediately obtain

\begin{crl}
\label{approx_maps3}
If ${\cal T}$ is a simple bubble type,
there exist $\de,C\!\in\! C^{\i}({\cal M}_{\cal T}^{(0)};\Bbb{R}^+)$ 
such that for all  \hbox{$\ups\!\in\! F^{(0)}{\cal T}_{\de}$} 
if ${\cal T}$ is regular
and \hbox{$\ups\!\in\! F^{(\eset)}{\cal T}_{\de}$} 
if ${\cal T}$ is semiregular,\\
(1) $\|\tilde{\pi}_{\ups,+}^{0,1}\eta\|_{\ups,p}\le C(b_{\ups})\|\eta\|_{\ups,p}$
for all $\eta\!\in\!\Ga^{0,1}(\ups)$;\\
(2) $\|P_{\ups}\eta\|_{\ups,p,1}\le C(b_{\ups})\|\eta\|_{\ups,p}$
for all $\eta\!\in\!\Ga^{0,1}(\ups)$.
\end{crl}

\subsection{The Gluing Map}
\label{gluing_map}

\noindent
In this subsection, we look for small vector fields 
$\xi\!\in\!\tilde{\Ga}_+(\ups)$ such that $\exp_{b_{\ups},u_{\ups}}\xi$
is holomorphic if ${\cal T}$ is regular
and lies in ${\cal M}_{\Si,t\nu,\la^*}$ if ${\cal T}$ is semiregular.
In Subsection~\ref{surj_c}, we show that all holomorphic maps if  
${\cal T}$ is regular
and all maps in ${\cal M}_{\Si,t\nu,\la^*}\!\times\Si^L$ if 
${\cal T}$ is semiregular
that lie near ${\cal M}_{\cal T}$ with respect to the Gromov topology
can be obtained in this~way.\\

\noindent
If $\xi\!\in\!\Ga(u_{\ups})$, define $\exp_{\ups}\xi\!: \Si_{\ups}\lra V$ 
and $\bar{\partial}_{\ups}\xi\!\in\!\Ga^{0,1}(u_{\ups})$ by
$$\{\exp_{\ups}\xi\}(z)=\exp_{b_{\ups},u_{\ups}(z)}\xi(z),\quad
\{\bar{\partial}_{\ups}\xi\}|_z=
\Pi_{b_{\ups},\xi(z)}^{-1}\circ\l.\bar{\partial}\{\exp_{\ups}\xi\}\r|_z.$$
If $S=\Si$ and
$\nu\!\in\!\Ga(\Si;\La^{0,1}\pi_{\Si}^*T^*\Si\otimes\pi_V^*TV)$,
let $\nu_{\ups,\xi}\in\Ga^{0,1}(u_{\ups})$ be given by 
$$\nu_{\ups,\xi}|_z=\Pi_{b_{\ups},\xi(z)}^{-1}\circ
\nu|_{(z,\{\exp_{\ups}\xi\}z)}.$$
Then, 
\begin{equation}
\label{pert_e1}
\bar{\partial}\{\exp_{\ups}\xi\}(\cdot)=
t\nu|_{(\cdot,\{\exp_{\ups}\xi\}(\cdot))}
\Llra \bar{\partial}_{\ups}\xi=t\nu_{\ups,\xi}.
\end{equation}
Write
\begin{equation}
\label{expan_e}
\bar{\partial}_{\ups}\xi=\bar{\partial}u_{\ups}+D_{\ups}\xi+
N_{\ups}\xi \quad\hbox{and}\quad 
\nu_{\ups,\xi}|_z=\nu|_{(z,u_{\ups}(z))}+L_{\nu,\ups}\xi|_z.
\end{equation}
Then the second equation in \e_ref{pert_e1} is equivalent to
\begin{equation}
\label{pert_e2}
D_{\ups}\xi+N_{\ups,t\nu}\xi=t\nu-\bar{\partial}u_{\ups},
\end{equation}
and by Proposition~\ref{an-dbar_prp} in~\cite{Z1} 
and (1) of Lemma~\ref{approx_maps},
there exist 
$C_{\bar{\partial}},\de\!\in\! C^{\i}({\cal M}_{\cal T}^{(0)};\Bbb{R}^+)$
such that for any $\ups\!\in\! F^{(0)}{\cal T}_{\de}$
and $\xi_1,\xi_2\!\in\!\Ga(u_\ups)$,
\begin{equation}
\label{pert_e3}
\|N_{\ups,t\nu}\xi_1-N_{\ups,t\nu}\xi_2\|_{\ups,p}
\le C_{\bar{\partial}}(b_{\ups})
\l(\|\xi_1\|_{\ups,p,1}+\|\xi_2\|_{\ups,p,1}+t\r)\|\xi_1-\xi_2\|_{\ups,p,1}.
\end{equation}
If ${\cal T}$ is semiregular, the term $\nu$ will be fixed, and 
we will be looking for solutions of \e_ref{pert_e2}
with $t\!>\!0$ very small for $\ups\!\in\!  F^{(\eset)}{\cal T}_{\de}$.
If ${\cal T}$ is regular, we will consider \e_ref{pert_e2} with $t\!=\!0$ 
and $\ups\!\in\! F^{(0)}{\cal T}_{\de}$.
In both cases, we will consider only solutions $\xi$
of \e_ref{pert_e2} that lie in 
the subspace $\tilde{\Ga}_+(\ups)$ of $L^p_1(\ups)$, since the subspace
$\Ga_-(\ups)$ corresponds to moving along 
the image of the pregluing \hbox{map $\ups\!\lra\! b(\ups)$.}\\

\noindent
Vector field $\xi\!=\!P_{\ups}\eta$ with $\eta\!\in\!\Ga_+^{0,1}(\ups)$
solves equation~\e_ref{pert_e2} if and only if
\begin{align}
\label{pert_e4a}
&\eta+\pi_{\ups,+}^{0,1}N_{\ups,t\nu}P_{\ups}\eta=
\pi_{\ups,+}^{0,1}\l(t\nu-\bar{\partial}u_{\ups}\r)\\
\label{pert_e4b}
\hbox{and~~~~~~~} 
&\pi_{\ups,-}^{0,1}\Big(t\nu-\bar{\partial}u_{\ups}-
\tilde{\pi}_{\ups,+}^{0,1}\eta-N_{\ups,t\nu}P_{\ups}\eta\Big)=0.
\end{align}
Denote the map 
$\eta\!\lra\!\pi_{\ups,+}^{0,1}N_{\ups,t\nu}P_{\ups}\eta$
by $N_{\ups,t\nu}^+$.
By Corollary~\ref{approx_maps3} and equation~\e_ref{pert_e3}, 
there exist  \hbox{$\tilde{C}_{\bar{\partial}},\de
 \!\in\! C^{\i}({\cal M}_{\cal T}^{(0)};\Bbb{R}^+)$   such that for any 
$\ups\!\in\! F^{(0)}{\cal T}_{\de}$}  if ${\cal T}$ is regular and
$\ups\!\in\! F^{(\eset)}{\cal T}_{\de}$ if ${\cal T}$ is semiregular,
\begin{equation}
\label{quadratic_term}
\|N_{\ups,t\nu}^+\eta_1-N_{\ups,t\nu}^+\eta_2\|_{\ups,p}\le 
\tilde{C}_{\bar{\partial}}(b_{\ups})
\l(\|\eta_1\|_{\ups,p}+\|\eta_2\|_{\ups,p}+t\r)\|\eta_1-\eta_2\|_{\ups,p}
\qquad\forall\eta_1,\eta_2\!\in\!\Ga_+^{0,1}(\ups).
\end{equation}

\begin{lmm}
\label{fixed_pt}
There exist $\ep,\de\!\in\! C^{\i}({\cal M}_{\cal T}^{(0)};\Bbb{R}^+)$ 
such that for all  
$\ups\!\in\! F^{(0)}{\cal T}_{\de}$ and $t\!=\!0$ if ${\cal T}$ is regular,
$\ups\!\in\! F^{(\eset)}{\cal T}_{\de}$ and $t\!\in\![0;\de(b_{\ups})]$ 
if ${\cal T}$ is semiregular,
and
$\al\!\in\!\Ga^{0,1}_+(\ups)$ with $\|\al\|_{\ups,p}\!<\!\ep(b_{\ups})$, 
the~equation
$$\eta+N_{\ups,t\nu}^+\eta=\al$$
has a unique solution $\eta_{\al}$ in $\Ga^{0,1}_+(\ups)$
such that $\|\eta_{\al}\|_{\ups,p}\le 2\ep(b_{\ups})$. 
Furthermore, such a solution satisfies 
$\|\eta_{\al}\|_{\ups,p}\le 2\|\al\|_{\ups,p}$.
\end{lmm}

\noindent
{\it Proof:} 
Put $\ep(b)=(6\tilde{C}_{\bar{\partial}}(b))^{-1}$, 
where $\tilde{C}_{\bar{\partial}}$ is as in \e_ref{quadratic_term}.
Define
$$\Psi_{\al}\!:
\{\eta\!\in\!\Ga^{0,1}_+(\ups)\!:\|\eta\|_{\ups,p}\le 2\|\al\|_{\ups,p}\}
\lra\Ga^{0,1}_+(\ups)$$
by $\Psi_{\al}(\eta)=\al-N_{\ups,t\nu}^+\eta$. 
By equation~\e_ref{quadratic_term},
\begin{gather*}
\|\Psi_{\al}(\eta)\|_{\ups,p}\le \|\al\|_{\ups,p}+
\tilde{C}_{\bar{\partial}}(b_{\ups})
\big(\|\eta\|_{\ups,p}+t\big)\|\eta\|_{\ups,p} 
\le 2\|\al\|_{\ups,p};\\
\|\Psi_{\al}(\eta_1)-\Psi_{\al}(\eta_2)\|_{\ups,p}\le 
\tilde{C}_{\bar{\partial}}(b_{\ups})
  (\|\eta_1\|_{\ups,p}+\|\eta_2\|_{\ups,p}+t)\|\eta_1-\eta_2\|_{\ups,p}
 \le \frac{5}{6}\|\eta_1-\eta_2\|_{\ups,p}.
\end{gather*}
It follows that $\Psi_{\al}$ is a contracting operator, and thus
has a unique fixed point $\eta_{\al}$,~i.e.
$$\eta_{\al}+N_{\ups,t\nu}^+\eta_{\al}=\al, \quad\hbox{and}\quad
\|\eta_{\al}\|_{\ups,p}\le 2\|\al\|_{\ups,p}.$$
The uniqueness claim follows immediately by taking
the difference of the corresponding equations.

\begin{crl}
\label{fixed_pt_crl}
If ${\cal T}$ is a simple bubble type, there exist
$\de,\ep,C\!\in \!
   C^{\i}({\cal M}_{\cal T}^{(0)};\Bbb{R}^+)$ 
such that for all 
$\ups\!\in\! F^{(0)}{\cal T}_{\de}$ and $t\!=\!0$ if ${\cal T}$ is regular and
$\ups\!\in\! F^{(\eset)}{\cal T}_{\de}$ and $t\!\in\![0;\de(b_{\ups})]$ 
if ${\cal T}$ is semiregular,
there exists a unique $\eta_{\ups,t\nu}\!\in\!\Ga^{0,1}(\ups)$
such that $\eta_{\ups,t\nu}$ satisfies equation~\e_ref{pert_e4a}
and $\|\eta_{\ups,t\nu}\|_{\ups,p}\le\ep(b_{\ups})$.
Furthermore, 
$$\|\eta_{\ups,t\nu}\|_{\ups,p}\le 
C(b_{\ups})\big(t+|\ups|^{\frac{1}{p}}\big).$$ 
\end{crl}

\noindent
{\it Proof:} This corollary follows from Lemmas~\ref{fixed_pt} 
and~\ref{approx_maps}.\\
 
\noindent
We now put $\xi_{\ups,t\nu}=P_{\ups}\eta_{\ups,t\nu}$ and
$\tilde{u}_{\ups,t\nu}=\exp_{\ups}\xi_{\ups,t\nu}$.
Replacing $u_{\ups}$ in $b(\ups)$ by $\tilde{u}_{\ups,t\nu}$,
we obtain a~new bubble map that will be called $\tilde{b}_{t\nu}(\ups)$.
If ${\cal T}$ is regular 
(and thus $t\!=\!0$), we will write $\tilde{u}_{\ups}$
and $\tilde{b}(\ups)$ for $\tilde{u}_{\ups,0}$
and $\tilde{b}_0(\ups)$, respectively.
We can assume that the functions $\de,\ep$ and $C$ of 
Corollary~\ref{fixed_pt_crl} are 
$\big(G_{\cal T}\!\times\!\hbox{Aut}({\cal T})\big)$-invariant.
For ${\cal T}$ regular, we have thus constructed a {\it gluing} map
$$\tilde{\ga}_{\cal T}^{(0)}\!:  F^{(0)}{\cal T}_{\de}
\lra \bar{\cal M}_{\lr{\cal T}},\quad
\ups\lra\tilde{b}(\ups).$$
Since this map is 
$\big(G_{\cal T}\!\times\!\hbox{Aut}({\cal T})\big)$-invariant, 
as can be seen from the
construction, $\tilde{\ga}_{\cal T}^{(0)}$ induces a map
on the quotient
\begin{equation}
\label{gluing_map2a}
\tilde{\ga}_{\cal T}\!: F{\cal T}_{\de}
\lra \bar{\cal M}_{\lr{\cal T}}.
\end{equation}
By the  smooth dependence of solutions of \e_ref{pert_e4a},
the restrictions
$$\tilde{\ga}_{\cal T}^{(0)}\!: F^{(H)}{\cal T}_{\de}\lra
 {\cal M}_{{\cal T}(H)}^{(0)}$$
are smooth. 
However, continuity of $\tilde{\ga}_{\cal T}$ on all of 
$F{\cal T}_{\de}$ is not immediate.
In the next section, we show the map 
$\tilde{\ga}_{\cal T}$ is a homeomorphism onto a neighborhood
of ${\cal M}_{\cal T}$ in~$\bar{\cal M}_{\lr{\cal T}}$.
If ${\cal T}$ is semiregular and $t\!>\!0$, we have constructed a map
$$\tilde{\ga}_{{\cal T},t\nu}^{(0)}:
F^{(\eset)}{\cal T}_{\de}\Big|
\ep^{-1}(-t,t)\lra C^{\i}_{(\la^*;M)}(\Si;V),$$
which again is $\big(G_{\cal T}\!\times\!\hbox{Aut}({\cal T})\big)$-invariant
 and thus descends to a map
\begin{equation}
\label{gluing_map2b}
\tilde\ga_{{\cal T},t\nu}:
F^{\eset}{\cal T}_{\de}\Big|
\big(\ep^{-1}(-t,t)/G_{\cal T}\big)\lra C^{\i}_{(\la^*;M)}(\Si;V).
\end{equation}
The map $u_{\tilde{b}_{t\nu}(\ups)}$ lies in ${\cal M}_{\Si,t\nu,\la^*}$ 
if and only if equation~\e_ref{pert_e4b} is satisfied, i.e.
\begin{equation}
\label{psi_dfn}
R_{\ups}\psi_{{\cal T},t\nu}(\ups)\equiv
t\nu-\bar{\partial}u_{\ups}-\tilde{\pi}_{\ups,+}^{0,1}\eta_{\ups,t\nu}-
 N_{\ups,t\nu}P_{\ups}\eta_{\ups,t\nu}=0\in\Ga_-^{0,1}(\ups),
\end{equation}
since $\eta_{\ups,t\nu}$ satisfies equation~\e_ref{pert_e4a}.

\subsection{An Implicit Function Theorem}
\label{impl}

\noindent
In this subsection, we prove a refined version of 
the Implicit Function Theorem.
It will be used in the rest of this section to modify
the gluing maps of Subsection~\ref{gluing_map}
for the spaces ${\cal M}_{\cal T}(\mu)$, ${\cal U}_{\cal T}(\mu)$,~etc.\\

\noindent
Let ${\cal S}$ be a smooth oriented manifold, and  
${\cal NS}$, ${\cal N}^{\mu}$, and $F$ oriented Riemannian vector 
bundles over~${\cal S}$.
We denote by $b$, $(b,\vec{n})$, $(b,\si)$, and $(b,v)$ 
general elements of ${\cal S}$, ${\cal NS}$, ${\cal N}^{\mu}$,
and~$F$, respectively.
If $\Om$ is any subset of $F$ and $\de\!>\!0$, let
$$\Om(\de)=\big\{(b,\vec{n},v)\!\in\!{\cal NS}\oplus F\!:
(b,v)\!\in\!\Om;~|\vec{n}|,|v|\!<\!\de\big\}.$$
Let $U$ be an open neighborhood of ${\cal S}$ in 
${\cal NS}\oplus{\cal N}^{\mu}\oplus F$ and
$h\!: U\lra\Bbb{R}^k$ a smooth map such~that
$$h(b,\vec{n},\si,v)=h(b,\vec{n},\si,0),\quad
h|{\cal S}=0,\quad\hbox{and}\quad
d(h\!:{\cal N}^{\mu}_b\lra\Bbb{R}^n)_{(b,0)}\!: 
 {\cal N}^{\mu}_b\lra\Bbb{R}^k$$
is an orientation-preserving isomorphism for all $b\!\in\!{\cal S}$.
Let $\tilde{U}$ be a subset of $U$ such that $\tilde{U}$
is the fiber product along ${\cal S}$ of an open neighborhood
of ${\cal S}$ in ${\cal NS}\oplus {\cal N}^{\mu}$ and an open
subset $\Om$ of~$F$.
Suppose $\de_{\cal S}\!>\!0$, $C\!\in\! C^{\i}({\cal S};\Bbb{R}^+)$,  
and 
\hbox{$\tilde{h}_t\!: \tilde{U}\!\lra\Bbb{R}^k$} is a family of smooth
functions with $t\!\in\![0,\de_{\cal S}]$ such~that
$$\big|\tilde{h}_t-h\big|_{(b,\vec{n},\si,v)},
\Big|\frac{\partial\tilde{h}_t}{\partial\si}-
\frac{\partial h}{\partial\si}\Big|_{(b,\vec{n},\si,v)}
\le C(b)\big(|v|^{\frac{1}{p}}+t\big)\quad
\forall t\!\in\!(0,\de_{\cal S}),~ (b,\vec{n},\si,v)\!\in\!\tilde{U},$$
where $\frac{\partial h}{\partial\si}$ denotes the differential
of $h$ along the fibers of ${\cal N}^{\mu}$.

\begin{lmm}
\label{impl_l1}
Let $B$ be an open ball about $0\!\in\!\Bbb{R}^k$.
If $f\!: B\lra\Bbb{R}^k$ is a smooth function and 
$$k\big|Df|_z-Df|_0\big|<\big|(Df|_0)^{-1}\big|^{-1}\quad\forall 
z\!\in\! B,$$
then $f$ is injective on $B$.
\end{lmm}

\noindent
{\it Proof:} Let $f_i$ denote the $i$th component of $f$.
By the Mean Value Theorem, for all $x,y\!\in\! B$, there exists
$z_i(x,y)\!\in\! B$ such that
$$\big|f_i(x)-f_i(y)\big|=\big|Df_i|_{z_i(x,y)}\big||x-y|.$$
Adding up these equations over all $i$, we obtain
\begin{equation*}\begin{split}
\sum_{i=1}^{i=k}\big|f_i(x)-f_i(y)\big|
&\ge\sum_{i=1}^{i=k}\big|Df_i|_0\big||x-y|-
k\sup_{z\in B}\big|Df|_z-Df|_0\big||x-y|\\
&\ge\Big(\big|(Df|_0)^{-1}\big|^{-1}-
k\sup_{z\in B}\big|Df|_z-Df|_0\big|\Big)|x-y|.
\end{split}\end{equation*}

\begin{lmm}
\label{impl_l2}
For every precompact subset $K$ of ${\cal S}$,
there exists $\ep\!>\!0$ such that for all $t\!\in\!(0,\ep)$ and
\hbox{$(b,\vec{n},v)\!\in\!\Om(\ep)|K$}, the map
$$\big\{(b,\si)\!\in\!{\cal N}^{\mu}\!: |\si|<\ep\big\}
               \lra\tilde{h}_t(b,\vec{n},\si,v)$$
is defined and injective, and its differential defines an 
orientation-preserving isomorphism between ${\cal N}_b^{\mu}$
and~$\Bbb{R}^k$.
\end{lmm}

\noindent
{\it Proof:} The map above is defined as long as
$$\l\{ (b,\vec{n},\si,v)\!\in\! {\cal NS}\oplus{\cal N}^{\mu}_b\oplus F\!:
b\!\in\! K, (b,\vec{n},v)\!\in\!\Om(\ep), |\si|<\ep\r\}\subset\tilde{U}.$$ 
Since $K$ is precompact, existence of $\de>0$ such that
the last inclusion holds is trivial.
The other two statements follow from the third property of $h$ and
the second property of $\tilde{h}_t$ (see above);
Lemma~\ref{impl_l1} is needed to prove the injectivity.
Note that
the variation of $\frac{\partial\tilde{h}_t}{\partial\si}$
over $K$ can be bounded from the variation 
$\frac{\partial h}{\partial\si}$ and the second property of~$\tilde{h}_t$.

\begin{lmm}
\label{impl_l3}
For every precompact subset $K$ of ${\cal S}$ and $\ep>0$ sufficiently
small, there exists $\de>0$ such that for all $t\!\in\!(0,\de)$ and
$(b,\vec{n},v)\!\in\!\Om(\de)|K$, the image of the map
$$\{(b,\si)\!\in\!{\cal N}^{\mu}\!:|\si|<\ep\} 
                       \lra\tilde{h}_t(b,\vec{n},\si,v)$$
contains $0\in\Bbb{R}^k$.
\end{lmm}

\noindent
{\it Proof:} We assume $\ep>0$ does not exceed the number
provided by Lemma~\ref{impl_l2}.
Then by precompactness of $K$ and the proof of Lemma~\ref{impl_l2}, 
\begin{equation}
\label{impl_l3_e1}
\varepsilon\equiv\min
\big\{|h(b,\vec{n},\si,v)|\!: (b,\vec{n},v)\!\in\!\Om_{\ep}|K, 
(b,\si)\!\in\!{\cal N}^{\mu}, |\si|=\frac{1}{2}\ep\big\}>0.
\end{equation}
Since for each $(b,\vec{n},v)\in\Om(\ep)|K$,
the image of the map
$$\{(b,\si)\!\in\!{\cal N}^{\mu}\!:|\si|<\ep\} 
                       \lra h(b,\vec{n},\si,v)$$
contains a neighborhood of $0$ in $\Bbb{R}^k$
and $\tilde{h}_t$ is continuous, 
the claim follows from the first property of $\tilde{h}_t$
along with equation~\e_ref{impl_l3_e1}.

\begin{crl}
\label{impl_c4}
For every precompact open subset $K$ of ${\cal S}$, 
there exist $\de,C>0$ with the following property.
For all $t\!\in\!(0,\de)$,
there exists a smooth section
$$\varphi_t\in\Ga\big(\Om(\de)|K;\pi^*{\cal N}^{\mu}\big),$$
where $\pi:\Om(\de)|K\lra K$ is the bundle projection map, such that
$$\Om(\de)|K\lra\tilde{h}_t^{-1}(0),\quad
(b,\vec{n},v)\lra\big(b,\vec{n},\varphi_t(b,\vec{n},v),v\big),$$
is an orientation-preserving diffeomorphism.
Furthermore,
$$\big|\varphi_t(b,\vec{n},v)\big|\le 
C\big(|v|^{\frac{1}{p}}+t+|\vec{n}|\big)\quad
\forall (b,\vec{n},v)\!\in\!\Om(\de)|K.$$
Finally, if $G$ is a group that acts on the space ${\cal S}$
and bundles ${\cal NS}$, ${\cal N}^{\mu}$, and $F$,  and preserves 
$h$, $\tilde{h}_t$, $\Om$, and $K$, then $\varphi_t$ is $G$-equivarent.
\end{crl}

\noindent
{\it Proof:} With $\ep$ as provided by Lemma~\ref{impl_l2},
let $\de>0$ be as provided by Lemma~\ref{impl_l3}. 
Then
$$F_t: \l\{(b,\vec{n},\si,v)\!: (b,\vec{n},v)\!\in\!\Om(\de)|K,
|\si|<\ep\r\}\lra  \Om(\de)\times\Bbb{R}^k,~~
F_t(b,\vec{n},\si,v)=\big(b,\vec{n},v,\tilde{h}_t(b,\vec{n},v)\big)$$
is a diffeomorphism onto an open subset $W$ of the target space.
The inverse of~$F_t$ must have the~form
$$F_t^{-1}(b,\vec{n},v,\tilde{\si})=
\big(b,\vec{n},\phi_t(b,\vec{n},v,\tilde{\si}),v\big)$$
for some smooth function $\phi_t$.
By Lemma~\ref{impl_l3}, \hbox{$(\Om(\de)|K)\times\{0\}\subset W$.}
Thus,
$$\varphi_t\!\in\!\Ga\big(\Om(\de)|K;\pi^*{\cal N}^{\mu}\big),\quad
\varphi_t(b,\vec{n},v)=\phi_t(b,\vec{n},v,0),$$
is a well-defined section, and by definition of $\phi_t$,
$$\Om(\de)|K\lra\tilde{h}_t^{-1}(0),\quad
(b,\vec{n},v)\lra\big(b,\vec{n},\varphi_t(b,\vec{n},v),v\big),$$
is a diffeomorphism.
It is orientation-preserving by Lemma~\ref{impl_l2}.
The estimate on $\varphi_t$ follows from
the three properties of~$h$, the first property of~$\tilde{h}_t$,
and the proof of Lemma~\ref{impl_l1}.
The final statement of the lemma is clear, since 
our construction commutes with the \hbox{$G$-action}.

\subsection{The Orientation of ${\cal M}_{\Si,t\nu,\la^*}(\mu)$
and the Gluing Map}
\label{orient2_sec}

\noindent
At this point, our treatments of regular and semiregular cases diverge.
In this subsection, we assume that 
${\cal T}\!=\!(\Si,[N],I;j,\la)$ is a semiregular bubble type and
$\mu$ is an $N$-tuple of constraints in general position
as defined below.
We recall how each element of 
${\cal M}_{\Si,t\nu,\la^*}(\mu)$ is assigned a sign
and then specialize to the elements 
$\tilde{b}_{t\nu}(\ups)\!\in\!{\cal M}_{\Si,t\nu,\la^*}(\mu)$.
We conclude this subsection with Theorem~\ref{si_str} that 
describes the elements of ${\cal M}_{\Si,t\nu,\la^*}(\mu)$ lying near
the space~${\cal M}_{\cal T}(\mu)$.

\begin{dfn}
\label{regul_dfn}
(1) Section 
$\nu\!\in\!\Ga^{0,1}(\Si\times V;\La^{0,1}_{J,j}
                          \pi_{\Si}^*T^*\Si\otimes\pi_V^*TV)$
is \under{$\la^*$\!-regular} if for all \hbox{$t\!\in\!(0,1)$} and 
\hbox{$u\!\in\!{\cal M}_{\Si,t\nu,\la^*}$}, the operator
$D_{V,u}\!:\Ga(u)\lra\Ga^{0,1}(u)$ is surjective.\\
(2) If $\nu$ is $\la^*$-regular,
tuple $\mu$ of oriented submanifolds of $V$ is \under{$\nu$-regular} 
if for all $t\!\in\!(0,1)$,
\begin{gather*}
\bigoplus_{l\in[N]}T_{u_{\hat{0}}(y_l)}V\!=\!
\hbox{Im}~d\ev_{[N]}\big|_b +
\bigoplus_{l\in [N]}T_{u_{\hat{0}}(y_l)}\mu_l
\qquad\forall~b\!=\!(\Si,[N],\{\hat{0}\};,(\hat{0},y),u_{\hat{0}})\!\in\! 
                         {\cal M}_{\Si,t\nu,\la^*}(\mu),\\
\hbox{where}\qquad
d\ev_{[N]}\big|_b\!: 
\ker D_{V,u_{\hat{0}}}\oplus\bigoplus_{l\in[N]}T_{y_l}\Si
\lra  \bigoplus_{l\in[N]}T_{u_{\hat{0}}(y_l)}V,\quad
d\ev_l\big|_b\big(\xi,w_{[N]}\big)=\xi(y_l)+du_{\hat{0}}\big|_{y_l}w_l.
\end{gather*}
(3) If ${\cal T}$ is a $(V,J)$-semiregular bubble type,
tuple $\mu$ of oriented submanifolds of $V$ is \under{${\cal T}$-regular}
if 
\begin{gather*}
\bigoplus_{l\in {[N]}}T_{u_{j_l}(y_l)}V=
 \hbox{Im}~d\ev_{[N]}\big|_b
+\bigoplus_{l\in[N]}T_{u_{j_l}(y_l)}\mu_l
\qquad\forall~ b\!=\!(\Si,[N],I;x,(j,y),u)
                  \!\in\! {\cal H}_{\cal T}(\mu);\\
\hbox{where}\qquad
d\ev_{[N]}\big|_b\!:{\cal K}_b{\cal T}\lra
\bigoplus_{l\in {[N]}}T_{u_{j_l}(y_l)}V,\quad
d\ev_l\big|_b\big(\xi_I,w_{\hat{I}+[N]}\big)=
       \xi_{j_l}(y_l)+du_{j_l}\big|_{y_l}w_l.
\end{gather*}
(4) If ${\cal T}$ is a $(V,J)$-semiregular bubble type,
${\cal S}\!\subset\!{\cal M}_{\cal T}$ is a smooth submanifold,
and $\tilde{\cal S}\!\subset\!{\cal M}_{\cal T}^{(0)}$
is the preimage of ${\cal S}$ under the quotient projection map,
tuple $\mu$ of oriented submanifolds of $V$ is 
\under{${\cal S}$-regular}~if
$$\bigoplus_{l\in {[N]}}T_{u_{j_l}(y_l)}V=
d\ev_{[N]}\big|_b\big({\cal K}_b{\cal T}\cap T_b\tilde{\cal S}\big)
+\bigoplus_{l\in[N]}T_{u_{j_l}(y_l)}\mu_l
\qquad\forall~ b\!\in\!\tilde{\cal S}(\mu)\equiv 
\tilde{\cal S}\cap {\cal M}_{\cal T}^{(0)}(\mu).$$
\end{dfn}

\noindent
Note that all four definitions above are independent
of the choice of metrics on~$V$.
Throughout this subsection, we assume that 
$\nu$ is $\la^*$-regular, ${\cal T}$ is semiregular, 
and $\mu$ is $\nu$- and ${\cal T}$-regular.\\

\noindent
The space ${\cal M}_{\Si,t\nu,\la^*}$ consists of the maps 
$u\!:\Si\!\lra\! V$
such that $\bar{\partial}u|_z\!=\!t\nu(z,u(z))$ for all \hbox{$z\!\in\!\Si$.}
Thus, the tangent space at $u$ can be described as 
$$T_u{\cal M}_{\Si,t\nu,\la^*}=\big\{\xi\!\in\!\Ga(\Si; u^*TV)\!: 
 D_{V,u}\xi\!-\!tL_{\nu,u}\xi=0\big\},$$
where $L_{\nu,u}\xi$ is defined by
$$\{L_{\nu,u}\xi\}(z)=\na_{\xi(z)}^V\nu\big|_{(z,u(z))}.$$
The operator $D_{V,u}\!-\!tL_{\nu,u}$ is 
independent of the choice of the connection along~${\cal M}_{\Si,t\nu,\la^*}$
and by assumption has no cokernel if \hbox{$t\!\in\!(0,1)$}.
An~orientation on~${\cal M}_{\Si,t\nu,\la^*}$ is determined by an orientation
of the bundle $\La^{top}_{\Bbb{R}}T{\cal M}_{\Si,t\nu,\la^*}$ 
over~${\cal M}_{\Si,t\nu,\la^*}$, which is the determinant line bundle 
of the elliptic operator~$D_{V,u}\!-\!tL_{\nu,u}$.
Since $L_{\nu,u}$ has order zero,
the operator $D_{V,u}\!-\!tL_{\nu,u}$ 
is homotopic through elliptic operators
to the operator~$D_{V,u}$.
Thus, $\La^{top}_{\Bbb{R}}T{\cal M}_{\Si,t\nu,\la^*}$ is homotopic
to 
$$\hbox{det}\big(D_{V,u}\big)=
\big(\La_{\Bbb{R}}^{top}(\ker~D_{V,u})\big)\otimes 
\big(\La_{\Bbb{R}}^{top}(\coker~D_{V,u})\big);$$
see \cite{LM}.
Since $D_{V,u}$ commutes with $J$, 
$\ker D_{V,u}$ and $\coker~D_{V,u}$
are both complex vector spaces and thus have natural orientations,
which induce an orientation on the determinant line bundle of~$D_{V,u}$
and via a homotopy of operators on the determinant bundle
of~$D_{V,u}\!-\!tL_{\nu,u}$.
It follows that ${\cal M}_{\Si,t\nu,\la^*}\!\times\!\Si^N$ 
is naturally oriented.
If $\mu$ is a $\nu$-regular tuple of submanifold of~$V$ of total codimension
$$\codim~\mu=\dim~{\cal M}_{\Si,t\nu,\la^*}\!\times\!\Si^N
=\ind~D_{V,u}+2|N|
=2\Big(\lan c_1(V,J),\la^*\ran+(\dim~V)\big(1-g(\Si)\big)+|N|\Big),$$
the differential of the map
$$\ev_{[N]}\!:  {\cal M}_{\Si,t\nu,\la^*}\times\Si^N\lra\prod_{l\in [N]}V
\quad
(\Si,[N],\{\hat{0}\};,(\hat{0},y),u_{\hat{0}})\lra 
\big(u_{\hat{0}}(y_l)\big)_{l\in[N]},$$
i.e.~$d\ev_{[N]}$ as defined in (2) of Definition~\ref{regul_dfn},
induces an isomorphism between $T{\cal M}_{\Si,t\nu,\la^*}\oplus T\Si^N$
and the normal bundle of $\mu$
in $V^N$ at each point of~${\cal M}_{\Si,t\nu,\la^*}(\mu)$.
Here we identify the $N$-tuple~$\mu$ with the submanifold
$$\prod_{l\in [N]}\mu_l\subset\prod_{l\in [N]}V\equiv V^N.$$
Since the normal bundle of~$\mu$ is oriented, the evaluation map also induces
an orientation on $T{\cal M}_{\Si,t\nu,\la^*}\oplus  T\Si^N$
along~${\cal M}_{\Si,t\nu,\la^*}(\mu)$.
Each element $b\!\in\!{\cal M}_{\Si,t\nu,\la^*}(\mu)$
is assigned a plus sign or is positively oriented 
if the two orientations agree, and a minus sign otherwise.\\

\noindent
For any $\ups\!\in\! F^{(\eset)}{\cal T}$ such that $q_{\ups}$ is defined,
let $L_{\nu,\ups}\!: \Ga(u_{\ups})\!\lra\!\Ga^{0,1}(u_{\ups})$
be given~by
$$\{L_{\nu,\ups}\xi\}(z)=\na_{\xi(z)}^{b_{\ups}}\nu\big|_{(z,f(z))}.$$
Denote by $\Ga_{t,+}^{0,1}(\ups)$ the image of 
$\Ga_+(\ups)$ under the map \hbox{$D_{\ups}\!-\!tL_{\nu,\ups}$}.

\begin{lmm}
\label{orient_l1}
For any compact subset $K$ of ${\cal M}_{\cal T}^{(0)}$,
there exist $\de,C>0$
such that for all \hbox{$\ups\!\in\! F^{(\eset)}{\cal T}_{\de}|K$}
and $t\!\in\!(0,\de)$,\\
(1) $\|\xi\|_{\ups,p,1}\le C\|D_{\ups}\xi-tL_{\nu,\ups}\xi\|_{\ups,p}$
for all $\xi\!\in\!\Ga_+(\ups)$;\\
(2) $L^p(\ups)=\Ga_{t,+}^{0,1}(\ups)\oplus\Ga^{0,1}_-(\ups)$;\\
(3) if $D_{\ups,t}^{--}$ and $L_{\nu,\ups,t}^{--}$ are 
the $(-,-)$-components of $D_{\ups}$ and 
$L_{\nu,\ups}$ with respect to the decompositions
$L^p_1(\ups)=\Ga_+(\ups)\oplus\Ga_-(\ups)$ and
$L^p(\ups)=\Ga_{t,+}^{0,1}(\ups)\oplus\Ga^{0,1}_-(\ups)$, then
$$\pi_{\ups,-}\!:
\ker\big\{ D_{\ups}-tL_{\nu,\ups}\!: L^p_1(\ups)\!\lra\! L^p(\ups)\big\}\lra
\ker\big\{ D_{\ups,t}^{--}-tL_{\nu,\ups,t}^{--}\!: 
 \Ga_-(\ups)\!\lra\!\Ga^{0,1}_-(\ups)\big\}$$
is an orientation-preserving isomorphism,
provided one of the two operators is surjective.
\end{lmm}

\noindent
{\it Proof:} (1) The first claim is immediate from (1) of 
Lemma~\ref{approx_maps2} and (2) of Lemma~\ref{approx_maps}.
The second is obtained by the same argument as in the proof 
of (3) of Lemma~\ref{approx_maps2}.\\
(2) By construction, $\pi_{\ups,-}$ is an isomorphism of 
the two kernels of the lemma.
In particular, \hbox{$D_{\ups}-tL_{\nu,\ups}$} is surjective if 
and only if $D_{\ups,t}^{--}-tL_{\nu,\ups,t}^{--}$ is.
Define 
\begin{gather*}
\Phi_{\tau}\!: L^p_1(\ups)\oplus\Ga^{0,1}_-(\ups)\lra L^p(\ups)
\quad\hbox{and}\quad
\Psi_{\tau}\!: \Ga_-(\ups)\oplus\Ga^{0,1}_-(\ups)\lra\Ga^{0,1}_-(\ups)
\qquad\hbox{by}\\
\Phi_{\tau}(\xi,\eta)= 
D_{\ups}\xi+\tau tL_{\nu,\ups}\xi+\eta
\quad\hbox{and}\quad
\Psi_{\tau}(\xi,\eta)= 
\tau \big(D_{\ups,t}^{--}+tL_{\nu,\ups,t}^{--}\big)\xi+\eta.
\end{gather*}
The first map is surjective for all $\tau\!\in\![0,1]$
by (2) of the lemma,
while the surjectivity of the second map is immediate from the definition.
Furthermore, the maps
$$\phi_{\tau}\!:\ker \Phi_{\tau}\lra\Ga_-(\ups),~~
\phi_{\tau}(\xi,\eta)=\pi_{\ups,-}\xi,\quad
\quad\hbox{and}\quad
\psi_{\tau}\!: \ker \Psi_{\tau}\lra\Ga_-(\ups),~~
\psi_{\tau}(\xi,\eta)=\xi$$
are isomorphisms such that 
$$\psi_1^{-1}\phi_1(\xi,0)=\pi_{\ups,-}\xi;\quad
\hbox{if~}\phi_1(\xi,\eta)=\psi_1(\xi',\eta'),
~\eta=\eta';\quad
\psi_0^{-1}\phi_0J=J\psi_0^{-1}\phi_0.$$
It follows that $\pi_{\ups,-}$ is an orientation-preserving map
between the two kernels of the lemma.\\

\noindent
If $K$ is a precompact open subset of ${\cal M}_{\cal T}$
and $\de\!>\!0$ is such that $\tilde{b}_{t\nu}(\ups)$ is defined
for all \hbox{$\ups\!\in\! F^{\eset}{\cal T}_{\de}|K$} and 
\hbox{$t\!\in\!(0,\de)$}, let
${\cal M}(K,\de)$ and $\tilde{\cal M}_{t\nu}(K,\de)$
denote the images of $F^{\eset}{\cal T}_{\de}|K$
under the maps $\ga_{\cal T}$ and $\tilde{\ga}_{{\cal T},t\nu}$, respectively.
Both maps are continuous and injective;
see Subsection~\ref{local_inject}.
The smooth structure of $F{\cal T}$ induces smooth structures
on ${\cal M}(K,\de)$ and $\tilde{\cal M}_{t\nu}(K,\de)$, 
with the tangent bundles described by
\begin{gather*}
T_{b(\ups)}{\cal M}(K,\de)=
  \big\{\ze_{\vp}'\!=\!\frac{d}{d\tau}\ze_{\tau\vp}\big|_{\tau=0}\!:
  \vp\!\in\! T_{\ups}F{\cal T}\big\}
  \oplus\bigoplus_{l\in[N]}T_{y_l(\ups)}\Si;\\
T_{\tilde{b}_{t\nu}(\ups)}\tilde{\cal M}_{t\nu}(K,\de)=
 \big\{\tilde{\ze}_{\vp}'\!=\!
\frac{d}{d\tau}\tilde{\ze}_{\tau\vp}\big|_{\tau=0}\!:
    \vp\!\in\! T_{\ups}F{\cal T}\big\}
 \oplus\bigoplus_{l\in[N]}T_{y_l(\ups)}\Si,
\end{gather*}
where $T_{\ups}F{\cal T}$ denotes $T_{\ups}F^{\eset}{\cal T}$;
see Subsection~\ref{scale_sec}.
It is easy to see that $\vp\lra\ze_{\vp}'$
is nearly complex linear and $\pi_{\ups,-}$ is almost the identity on
the first component of  $T_{b(\ups)}{\cal M}(K,\de)$;
both error terms are bounded by~$C_K|\ups|$.
Furthermore, by (1) of Lemma~\ref{inj_l1} and Corollary~\ref{inj_c1}, 
\hbox{$\vp\!\lra\!\tilde{\ze}_{\vp}'$}
also nearly computes with the complex structures and 
$\Pi_{b_{\ups},\xi_{\ups,t\nu}}\pi_{\ups,-}
       \Pi_{b_{\ups},\xi_{\ups,t\nu}}^{-1}$
is almost the identity on the first component of  
$T_{\tilde{b}_{t\nu}(\ups)}\tilde{\cal M}_{t\nu}(K,\de)$;
in the given case, the error terms are bounded 
by~$C_K\big(t\!+\!|\ups|^{\frac{1}{p}}\big)$.
Thus, the orientations of ${\cal M}(K,\de)$ and 
$\tilde{\cal M}_{t\nu}(K,\de)$ 
induced by the natural orientation of $F{\cal T}$
agree with the orientations induced from the natural orientation on
$\Ga_-(\ups)\oplus\bigoplus\limits_{l\in[N]}T_{y_l(\ups)}\Si$
via the maps
\hbox{$\pi_{\ups,-}\oplus id$} and
$\pi_{\ups,-}\Pi_{b_{\ups},\xi_{\ups,t\nu}}^{-1}\oplus id$,
respectively.\\

\noindent
By construction in Subsection~\ref{gluing_map}, 
$$\tilde{\psi}_{t\nu}\!: 
\tilde{\cal M}_{t\nu}(K,\de)\lra\Ga^{0,1},\quad
\ups\lra t\nu|_{\tilde{u}_{\ups,t\nu}}-
\bar{\partial}\tilde{u}_{\ups,t\nu}\in\Ga^{0,1}(\tilde{u}_{\ups,t\nu}),$$
determines a section of the bundle $\Pi\Ga_-^{0,1}$ over 
$\tilde{\cal M}_{t\nu}(K,\de)$, given by
$$\Pi\Ga_-^{0,1}(\tilde{b}_{t\nu}(\ups))=
\Pi_{b_{\ups},\xi_{\ups,t\nu}}\Ga_-^{0,1}(\ups).$$
Note that the zero set of this section is precisely the space
$({\cal M}_{\Si,t\nu,\la}\times\Si^N)\cap\tilde{\cal M}_{t\nu}(K,\de)$.
A linearization of this section is given by
\begin{equation*}\begin{split}
\na_{\tilde{\ze}_{\vp}'}
(t\nu-\bar{\partial}\tilde{u}_{\ups,t\nu})
&\equiv \Pi_{b_{\ups},\xi_{\ups,t\nu}}
\pi^{0,1}_{\ups,t,-}\na_{\pi_{\ups,-}
      \Pi_{b_{\ups},\xi_{\ups,t\nu}}^{-1}\tilde{\ze}_{\vp}'}^{\ups}
\Pi_{b_{\ups},\xi_{\ups,t\nu}}^{-1}(t\nu-\bar{\partial}\tilde{u}_{\ups,t\nu})
\\
&=-\Pi_{b_{\ups},\xi_{\ups,t\nu}}\l(
 D_{\ups,t}^{--}-tL_{\ups,\nu,t}^{--}\r)
      \pi_{\ups,-}\Pi_{b_{\ups},\xi_{\ups,t\nu}}^{-1}\tilde{\ze}_{\vp}',
\end{split}\end{equation*}
where $\pi^{0,1}_{\ups,t,-}\!: 
L^p(\ups)=\Ga^{0,1}_{+,t}(\ups)\oplus{\Ga}^{0,1}_-(\ups)\lra
{\Ga}^{0,1}_-(\ups)$ is the projection map.

\begin{crl}
\label{orient_c2}
For any compact subset $K$ of ${\cal M}_{\cal T}^{(0)}$,
there exists \hbox{$\de>0$}
such that for all $t\!\in\!(0,\de)$, the orientation of 
$({\cal M}_{\Si,t\nu,\la}\times\Si^N)\cap\tilde{\cal M}_{t\nu}(K,\de)$
as the zero set of the section $\tilde{\psi}_{t\nu}$ agrees
with its natural orientation.
\end{crl}

\noindent
{\it Proof:} Suppose $\tilde{b}_{t\nu}(\ups)\!\in\!
     ({\cal M}_{\Si,t\nu,\la^*}\times\Si^N)\cap\tilde{\cal M}_{t\nu}(K,\de)$.
Since we can use any connection in
$\tilde{u}_{\ups,t\nu}^*TV$ to define the natural orientation
on $T_{\tilde{u}_{\ups,t\nu}}{\cal M}_{\Si,t\nu,\la^*}$,
we can write
$$\big\{D_{\tilde{u}_{\ups,t\nu}}-tL_{\nu,\tilde{u}_{\ups,t\nu}}\big\}\xi
=\Pi_{b_{\ups},\xi_{\ups,t\nu}}\big\{D_{\ups}-tL_{\nu,\ups}\big\}
\Pi_{b_{\ups},\xi_{\ups,t\nu}}^{-1}\xi
\quad\forall\xi\in\Ga(\tilde{u}_{\ups,t}).$$
Thus, by Lemma~\ref{orient_l1}, 
$\pi_{\ups,-}\circ\Pi_{b_{\ups},\xi_{\ups,t\nu}}^{-1}\oplus id$
induces an orientation-preserving isomorphism between
$T_{\tilde{b}_{t\nu}(\ups)}{\cal M}_{\Si,t\nu,\la^*}
\oplus\bigoplus\limits_{l\in[N]}T_{y_l(\ups)}\Si$
and $\ker(D_{\ups,t}^{--}-L_{\nu,\ups,t}^{--})
\oplus\bigoplus\limits_{l\in[N]}T_{y_l(\ups)}\Si$
with their natural orientations. 
By the preceding paragraph, the same is true
for the zero set of~$\tilde{\psi}_{t\nu}$.\\

\noindent
If $\mu$ is an $N$-tuple of constraints as above, let 
\begin{gather*}
{\cal M}(K,\de;\mu)=
\big\{b(\ups)\!\in\!{\cal M}(K,\de)\!: \ev_{[N]}(b(\ups))\!\in\!\mu\big\},\\
\tilde{\cal M}_{t\nu}(K,\de;\mu)=
\big\{\tilde{b}_{t\nu}(\ups)\!\in\!\tilde{\cal M}_{t\nu}(K,\de)\!: 
      \ev_{[N]}(\tilde{b}_{t\nu}(\ups))\!\in\!\mu\big\}.
\end{gather*}
Then ${\cal M}(K,\de;\mu)$ and $\tilde{\cal M}_{t\nu}(K,\de;\mu)$
are smooth manifolds. 
In fact, the smoothness of  ${\cal M}(K,\de;\mu)$
is immediate from the smoothness of
$F{\cal T}|{\cal M}_{\cal T}(\mu)$, 
which is a consequence of ${\cal T}$-regularity of~$\mu$,
while 
the smoothness of $\tilde{\cal M}_{t\nu}(K,\de;\mu)$ follows from 
Lemma~\ref{orient_l3} below.
Furthermore, since $\mu$ is $\nu$-regular,
the section $\tilde{\psi}_{t\nu}$ is transversal to zero 
in $\Pi\Ga_-^{0,1}$ over $\tilde{\cal M}_{t\nu}(K,\de;\mu)$. 
By Corollary~\ref{orient_c2}, the sign of 
\hbox{$\tilde{b}_{t\nu}(\ups)\!\in\!{\cal M}_{\Si,t\nu,\la^*}(\mu)$}
defined at the beginning of this subsection is its sign
as an element of the zero set  
of the section $\tilde{\psi}_{t\nu}$ of 
  $\Pi\Ga_-^{0,1}$ over~$\tilde{\cal M}_{t\nu}(K,\de;\mu)$.\\

\noindent
If  $b\!=\!\big(\Si,[N],I;x,(j,y),u\big)
 \!\in\!{\cal M}_{\cal T}^{(0)}(\mu)$, let
$$ {\cal K}^{\mu}_b{\cal T}=\big\{
\big(\xi,w_{\hat{I}+[N]}\big)\!\in\!{\cal K}_b{\cal T}\!:
\xi_{j_l}(y_l)+du_{j_l}|_{y_l}w_l\!\in\! T_{u_{j_l}(y_l)}\mu_l
~\forall l\!\in\![N]\big\}.$$
Denote by ${\cal N}^{\mu}_b{\cal T}$
the $(L^2,b)$-orthogonal complement
of ${\cal K}^{\mu}_b{\cal T}$ in~${\cal K}_b{\cal T}$.
Note that by (3) of Definition~\ref{regul_dfn}, 
$$\bigoplus_{l\in {[N]}}T_{u_{j_l}(y_l)}V=
d\ev_{[N]}\big|_b \big({\cal N}^{\mu}_b{\cal T}\big)
\oplus \bigoplus_{l\in[N]}T_{u_{j_l}(y_l)}\mu_l.$$
We denote by $\tilde{\cal N}^{\mu}{\cal T}$ the bundle over 
${\cal M}_{\cal T}^{(0)}(\mu)$ with fibers ${\cal N}^{\mu}_b{\cal T}$
and by ${\cal N}^{\mu}{\cal T}\lra{\cal M}_{\cal T}(\mu)$
its quotient by the natural $G_{\cal T}$-action.\\

\noindent
Suppose ${\cal S}\!\subset\!{\cal M}_{\cal T}$ is 
a smooth oriented submanifold such that
$\mu$ is ${\cal S}$-regular.
Denote by $\tilde{\cal S}\!\subset\!{\cal M}_{\cal T}^{(0)}$
the preimage of ${\cal S}$ under the quotient projection map.
Let ${\cal NS}\!\lra\!{\cal S}$ and 
${\cal N}\tilde{\cal S}\!\lra\!\tilde{\cal S}$ be the normal bundles.
Choose a $\big(G_{\cal T}\!\times\!\hbox{Aut}({\cal T})\big)$-equivarent 
orientation-preserving identification
\hbox{$\tilde{\phi}_{\cal S}\!:{\cal N}\tilde{\cal S}_{\de}
              \lra{\cal M}_{\cal T}^{(0)}$}
of neighborhoods of $\tilde{\cal S}$ in
${\cal N}\tilde{\cal S}$ and~${\cal M}_{\cal T}^{(0)}$.
Let \hbox{$\tilde{\Phi}_{\cal S}\!:
\pi_{{\cal N}\tilde{\cal S}}^*F^{(0)}{\cal T}\!\lra\!
F^{(0)}{\cal T}$}
be a $\big(G_{\cal T}\!\times\!\hbox{Aut}({\cal T})\big)$-equivarent 
vector-bundle isomorphism
covering $\tilde{\phi}_{\cal S}$ such that
$\tilde{\Phi}_{\cal S}$ is the identity
on~$\tilde{\cal S}$.
Let 
$$\phi_{\cal S}\!:{\cal NS}_{\de}\lra{\cal M}_{\cal T}
\quad\hbox{and}\quad
\Phi_{\cal S}\!:
\pi_{\cal NS}^*F{\cal T}\lra F{\cal T}$$
be the maps induced by $\tilde{\phi}_{\cal S}$
and $\tilde{\Phi}_{\cal S}$, respectively.
Put 
$${\cal S}(\mu)={\cal S}\cap{\cal M}_{\cal T}(\mu),\quad
\tilde{\cal S}(\mu)=\tilde{\cal S}\cap{\cal M}_{\cal T}^{(0)}(\mu).$$
Since $\mu$ is ${\cal S}$-regular, we can choose
a $\big(G_{\cal T}\!\times\!\hbox{Aut}({\cal T})\big)$-equivarent 
orientation-preserving identification\\
\hbox{$\tilde{\phi}_{\cal S}^{\mu}\!:
\tilde{\cal N}^{\mu}{\cal T}_{\de}|\tilde{\cal S}(\mu)\lra\tilde{\cal S}$}.
Let
$$\tilde{\Phi}_{\cal S}^{\mu}\!:
\pi_{\tilde{\cal N}^{\mu}{\cal T}}^*
\big( {\cal N}\tilde{\cal S}\oplus F^{(0)}{\cal T}\big)
\lra{\cal N}\tilde{\cal S}\oplus F^{(0)}{\cal T}$$
be a $\big(G_{\cal T}\!\times\!\hbox{Aut}({\cal T})\big)$-equivarent 
splitting-preserving vector-bundle
isomorphism covering $\tilde{\phi}_{\cal S}^{\mu}$ such that
$\tilde{\Phi}_{\cal S}^{\mu}$ is the identity
on~$\tilde{\cal S}(\mu)$.
Denote by
$$\phi_{\cal S}^{\mu}\!:
{\cal N}^{\mu}{\cal T}_{\de}|{\cal S}(\mu)\lra\tilde{\cal S}
\quad\hbox{and}\quad
\Phi_{\cal S}^{\mu}\!:
\pi_{{\cal N}^{\mu}{\cal T}}^*
\big( {\cal NS}\oplus F{\cal T}\big)
\lra{\cal NS}\oplus F{\cal T}$$
the maps induced by $\tilde{\phi}_{\cal S}^{\mu}$
and $\tilde{\Phi}_{\cal S}^{\mu}$, respectively.

\begin{dfn}
\label{regulariz_dfn}
With notation as above and in Subsection~\ref{scale_sec}, tuple 
$(\Phi_{\cal S},\Phi_{\cal S}^{\mu})$ is a 
\under{regularization} of ${\cal S}(\mu)$ if for all
\hbox{$b\!\in\!\tilde{\cal S}(\mu)$}, 
$\vec{n}\!\in\!{\cal N}_b\tilde{\cal S}_{\de(b)}$, and 
$\si\!\in\!\tilde{\cal N}_b^{\mu}{\cal T}_{\de(b)}$,
there exists 
$\vp(\vec{n},\si)\!\in\!{\cal K}_{\tilde{\phi}_{\cal S}(b,\vec{n})}{\cal T}$
such~that 
$$ \tilde{\Phi}_{\cal S}\tilde{\Phi}_{\cal S}^{\mu}
(b,\si;\vec{n},v_{\hat{I}})
=\big\{\tilde{\Phi}_{\cal S}(b,\vec{n};v)\big\}
    \big(\vp(\vec{n},\si)\big)
\qquad\forall~v\!\in\! F_b^{(0)}{\cal T}.$$
\end{dfn}

\noindent
Note that if $\mu$ is ${\cal S}$-regular,
${\cal S}(\mu)$ admits a normalization.
In fact, we can start with any choice of 
$\Phi_{\cal S}$ and 
$\Phi_{\cal S}^{\mu}\big|
\pi_{{\cal N}^{\mu}{\cal T}}^*{\cal NS}$
as in the preceding paragraph,
and then choose $\Phi_{\cal S}^{\mu}\big|
\pi_{{\cal N}^{\mu}{\cal T}}^*F{\cal T}$ 
so that the triple satisfies the requirements of the definition.
In applications of Theorem~\ref{si_str} in~\cite{Z2}, 
the exact choice of $\Phi_{\cal S}^{\mu}$ does not matter,
but that of $\Phi_{\cal S}$ does play a role.\\

\noindent
For the purposes of Theorem~\ref{si_str},
we assume that $\Phi_{\cal S}$ and $\Phi_{\cal S}^{\mu}$
also encode the lifts of $\phi_{\cal S}$ and $\phi_{\cal S}^{\mu}$
to the bundles \hbox{$\pi_{\cal NS}^*\Ga_-^{0,1}\!\lra\!{\cal S}$}
and \hbox{$\pi_{{\cal N}^{\mu}{\cal T}}^*\Ga_-^{0,1}\!\lra\!{\cal S}(\mu)$}, 
respectively.
Put
\begin{gather*}
F^{(0)}{\cal S}={\cal N}\tilde{\cal S}\oplus F^{(0)}{\cal T},\qquad
F^{(\eset)}{\cal S}=
\big\{(b,\vec{n},v_{\hat{I}})\!\in\! F^{(0)}{\cal S}\!:
v_{\hat{I}}\!\in\! F^{(\eset)}_b{\cal T}\big\};\\
F{\cal S}={\cal NS}\oplus F{\cal T},\qquad
F^{\eset}{\cal S}=
\big\{[b,\vec{n},v_{\hat{I}}]\!\in\! F^{(0)}{\cal S}\!:
[b,v_{\hat{I}}]\!\in\! F^{\eset}_b{\cal T}\big\}.
\end{gather*}

\begin{lmm}
\label{orient_l3}
For any $\big(G_{\cal T}\!\times\!\hbox{Aut}({\cal T})\big)$-invariant 
precompact open subset $K$ of $\tilde{\cal S}(\mu)$,
there exist an open neighborhood $U_K$ of $K$ in 
${\cal M}_{\cal T}^{(0)}$ and $\de,C\!>\!0$ 
with the following property.
If $t\!\in\!(0,\de)$,  there exists a smooth 
$\big(G_{\cal T}\!\times\!\hbox{Aut}({\cal T})\big)$-equivarent section 
$$\tilde{\varphi}_{{\cal S},t\nu}^{\mu}\in
\Ga\big( F^{(\eset)}{\cal S}_{\de}|K;
\pi_{F^{(0)}{\cal S}}^*\tilde{\cal N}^{\mu}{\cal T}\big),$$
such that
$\big\|\tilde{\varphi}_{{\cal S},t\nu}^{\mu}(\ups)\|_{b_{\ups},C^0}\le
 C\!\big(t+ |\ups|^{\frac{1}{p}}\big)$
for all $\ups\!\in\! F^{(\eset)}{\cal S}_{\de}|K$
and
$$ F^{\eset}{\cal S}_{\de}|K
\lra\tilde{\cal M}_{t\nu}(U_K,\de;\mu),\quad
[b,\vec{n},v_{\hat{I}}]\lra 
\tilde{\ga}_{{\cal T},t\nu}\Big( 
\tilde{\Phi}_{\cal S}\big( 
\tilde{\Phi}_{\cal S}^{\mu}
  \tilde{\varphi}_{{\cal S},t\nu}^{\mu}
      \big(b,\vec{n},v_{\hat{I}})\big)\Big),$$
is an orientation-preserving diffeomorphism.
\end{lmm}

\noindent
{\it Proof:} 
Since $\mu$ is a regular value of $\ev_{[N]}|{\cal S}$ and $K$ is precompact,
there exists $\de\!>\!0$ \hbox{such that the~map}
\begin{gather*}
\big\{(b,\vec{n},v,\si)\!\in\! F^{(0)}{\cal S}\!\oplus\!
 \tilde{\cal N}^{\mu}{\cal T}\big|K\!:
|\vec{n}|,\|\si\|_{b,C^0}\!<\!\de\big\}
\lra F^{(0)}{\cal T},~~
(b,\vec{n},v,\si)\lra 
\tilde{\Phi}_{\cal S}\tilde{\Phi}_{\cal S}^{\mu}
(b,\si;\vec{n},v),
\end{gather*}
is a $G_{\cal T}$-equivarent orientation-preserving diffeomorphism 
onto its image.
Thus, if $\de\!>\!0$ is sufficiently small,
there exists $C\!>\!0$ such that, with notation as 
in Definition~\ref{regulariz_dfn},
$$C^{-1}\|\si-\si'\|_{b,C^0}\le
\big\|\vp(\vec{n},\si)-\vp(\vec{n},\si')\big\|_{b,C^0}
\le C\|\si-\si'\|_{b,C^0}
\quad \forall~b\!\in\!\tilde{\cal S}(\mu),~
\vec{n}\!\in\!{\cal N}_b\tilde{\cal S}_{\de},~
\si,\si'\!\in\!\tilde{\cal N}_b^{\mu}{\cal T}.$$
Then by Corollary~\ref{inj_c1c} and definition of 
$S_{\vp}'$ in Subsection~\ref{scale_sec},
\begin{gather*}
\Big| 
d_V\big( 
 \tilde{\phi}_{\cal S}\tilde{\Phi}_{\cal S}^{\mu}(b,\si;\vec{n}),
\tilde{\ga}_{{\cal T},t\nu}\big(
\tilde{\Phi}_{\cal S}\tilde{\Phi}_{\cal S}^{\mu}(b,\si;\vec{n},v)
\big)\big)-
d_V\big( 
 \tilde{\phi}_{\cal S}\tilde{\Phi}_{\cal S}^{\mu}(b,\si';\vec{n}),
\tilde{\ga}_{{\cal T},t\nu}\big(
\tilde{\Phi}_{\cal S}\tilde{\Phi}_{\cal S}^{\mu}(b,\si';\vec{n},v)
\big)\big)\Big|\\
\le C\big(t+|v|^{\frac{1}{p}}\big) \|\si-\si'\|_{b,C^0}
\qquad 
\forall~t\!\in\!(0,\de),~
b\!\in\!\tilde{\cal S}(\mu),~
\vec{n}\!\in\!{\cal N}_b\tilde{\cal S}_{\de},~
\si,\si'\!\in\!\tilde{\cal N}_b^{\mu}{\cal T},
v\!\in\! F^{(\eset)}_b{\cal T}_{\de}.
\end{gather*}
On a neighborhood of $\ev_{[N]}(b)\!\in\!\mu$, 
we can identify the normal bundle of~$\mu$ in~$V^N$
$g_V$-isometrically with the trivial hermitian bundle of the same rank.
Let $\pi$ denote the projection onto the fiber.
Since $\mu$ is ${\cal S}$-regular,
$$\|\si-\si'\|_{b,C^0}\le
C\big|\pi\ev_{[N]}\big( 
 \tilde{\phi}_{\cal S}\tilde{\Phi}_{\cal S}^{\mu}(b,\si;\vec{n})\big)
-\pi\ev_{[N]}\big( 
 \tilde{\phi}_{\cal S}\tilde{\Phi}_{\cal S}^{\mu}(b,\si';\vec{n})\big)\big|
\quad
\forall~b\!\in\!\tilde{\cal S}(\mu),~
\vec{n}\!\in\!{\cal N}_b\tilde{\cal S}_{\de},~
\si,\si'\!\in\!\tilde{\cal N}_b^{\mu}{\cal T}.$$
Thus, we can apply Corollary~\ref{impl_c4} to
$$h=\pi\circ\ev_{[N]}\circ\tilde{\phi}_{\cal S}\circ
                             \tilde{\Phi}_{\cal S}^{\mu}
\quad\hbox{and}\quad
\tilde{h}_t=\pi\circ\ev_{[N]}\circ\tilde{\ga}_{{\cal T},t\nu}\circ
\tilde{\Phi}_{\cal S}\circ \tilde{\Phi}_{\cal S}^{\mu}.$$
We obtain $\de,\ep\!>\!0$ and for each $t\!\in\!(0,\de)$
a section $\tilde{\varphi}_{{\cal S},t\nu}^{\mu}$
with the claimed bound such that the~map
\begin{gather*}
F^{(\eset)}{\cal S}_{\de}|K\lra
\Big\{(b,\vec{n},v_{\hat{I}},\si)\!: \|\si\|_{b,C^0}<\ep,
\ev_{[N]}\tilde{\ga}_{{\cal T},t\nu}\big(
\tilde{\Phi}_{\cal S}
\big( \tilde{\Phi}_{\cal S}^{\mu}(b,\vec{n};\si,v)\big)\big)
    \!\in\!\mu\Big\},\\
\big(b,\vec{n},v\big)\lra
\big(b,\vec{n},v,\tilde{\varphi}_{{\cal S},t\nu}^{\mu}(b,\vec{n},v)\big)
\end{gather*}
is an orientation-preserving diffeomorphism.
Since 
\begin{gather*}
\big\{[b,\vec{n},v,\si]\!: 
[b,\vec{n},v]\!\in\! F^{\eset}{\cal S}_{\de}, 
\|\si\|_{b,C^0}<\ep\big\}\lra 
\tilde{\cal M}_{t\nu}(U_K,\de),\\
[b,\vec{n},\si,v]\lra
\tilde{\ga}_{{\cal T},t\nu}\big(
\tilde{\Phi}_{\cal S} 
\big(\tilde{\Phi}_{\cal S}^{\mu}(b,\vec{n};\si,v)\big)\big),
\end{gather*}
is orientation-preserving by the discussion above and 
our assumptions on $\tilde{\phi}_{\cal S}$, the claim follows.
Above
$$U_K=\tilde{\phi}_{\cal S}\Big( 
\tilde{\Phi}_{\cal S}^{\mu}\big(\big\{(b,\vec{n},\si)\!\in\! 
{\cal N}\tilde{\cal S}\oplus\tilde{\cal N}^{\mu}{\cal T}|K\!:
\|\vec{n}\|_{b,C^0}<\ep,
\|\si\|_{b,C^0}<\ep\big\}\big)\Big).$$

\begin{thm}
\label{si_str}
Suppose $\la^*\!\in\!H_2(V;\Bbb{Z})$,
${\cal T}\!=\!\big(\Si,[N],I;j,\la\big)$
is a \hbox{$(V,J)$-semiregular} bubble type, 
with \hbox{$\sum\limits_{i\in I}\la_i\!=\!\la^*$} 
and cokernel bundle \hbox{$\Ga_-^{0,1}\!\lra{\cal M}_{\cal T}$}, and
$(\tilde{\Ga}_-,\Ga_-^{0,1},R)$ is an obstruction bundle setup.
Let
${\cal S}\!\subset\!{\cal M}_{\cal T}$ be a smooth oriented submanifold,
$$\nu\!\in\!\Ga^{0,1}(\Si\times V;\La^{0,1}_{J,j}
                          \pi_{\Si}^*T^*\Si\otimes\pi_V^*TV)$$
a $\la^*$-regular section, 
$\mu$ a $\nu$-, ${\cal T}$-, and ${\cal S}$-regular
$N$-tuple of submanifolds of~$V$ of total codimension
$$\codim~\mu=
2\big(\lan c_1(V,J),\la^*\ran+(\dim~V)\big(1-g(S)\big)+|N|\big),$$
and $(\Phi_{\cal S},\Phi_{\cal S}^{\mu})$
is regularization of~${\cal S}(\mu)$.
Then for every precompact open subset $K$ of ${\cal S}(\mu)$,
there exist a neighborhood $U_K$ of $K$ in 
\hbox{$\bar{C}^{\i}_{(\la^*;N)}(\Si;\mu)$} and $\de,\ep,C\!>\!0$ 
with the following property.
For every $t\!\in\!(0,\ep)$, there exist  a section 
$$\varphi_{{\cal S},t\nu}^{\mu}\in\Ga\big(F^{\eset}{\cal S}_{\de}|K;
\pi_{F{\cal S}}^*{\cal N}^{\mu}{\cal T}\big),
\quad\hbox{with}\quad
\big\|\varphi_{{\cal S},t\nu}^{\mu}(\ups)\big\|_{b_{\ups},C^0}\le
C\big(t+|\ups|^{\frac{1}{p}}\big),$$
and a sign-preserving bijection between
${\cal M}_{\Si,t\nu,\la^*}(\mu)\cap U_K$ and
the zero set of the~section $\psi_{{\cal S},t\nu}^{\mu}$
defined~by
\begin{alignat*}{2}
\psi_{{\cal S},t\nu}^{\mu}\!\in\!\Ga\big(F^{\eset}{\cal S}_{\de}|K;
\pi_{F{\cal S}}^*\Ga_-^{0,1}\big),&&\quad&
\Phi_{\cal S}^{\mu}\big(\varphi_{{\cal S},t\nu}^{\mu}(\ups);
\psi_{{\cal S},t\nu}^{\mu}(\ups)\big)=
\psi_{{\cal S},t\nu}\big(
\Phi_{\cal S}^{\mu}(\varphi_{{\cal S},t\nu}^{\mu}(\ups))\big);\\
\psi_{{\cal S},t\nu}\!\in\!\Ga
\big(F^{\eset}{\cal S}_{\de}\big|({\cal S}\cap U_K);
\pi_{F{\cal S}}^*\Ga_-^{0,1}\big),&&\quad&
\Phi_{\cal S}\big(\ups;\psi_{{\cal S},t\nu}(\ups)\big)=
\psi_{{\cal T},t\nu}\big(\Phi_{\cal S}(\ups)\big);\\
\psi_{{\cal T},t\nu}\!\in\!
     \Ga\big(F^{\eset}{\cal T}_{\de}\big|({\cal M}_{\cal T}\cap U_K);
\pi_{F{\cal T}}^*\Ga_-^{0,1}\big),&&\quad&
R_{\ups}\psi_{{\cal T},t\nu}(\ups)= \pi_{\ups,-}^{0,1}\!
\big(
t\nu_{\ups,t}\!-\!\bar{\partial}(u_{b_{\ups}}\circ q_{\ups})\!-\!
              D_{\ups}\xi_{\ups,t\nu}
           -\tilde{\eta}_{\ups,t\nu}\!\big)
\end{alignat*}
for some 
$\nu_{\ups,t},\tilde{\eta}_{\ups,t\nu}
\!\in\!\Ga^{0,1}(u_{\ups})$
and $\xi_{\ups,t\nu}\!\in\!\tilde{\Ga}_+(\ups)$, 
dependent smoothly on $\ups$, such~that
$$\big\|\nu_{\ups,t}-\nu\big\|_{\ups,2}\le
C\big(t+|\ups|^{\frac{1}{p}}\big),\quad
\big\|\xi_{\ups,t\nu}\big\|_{\ups,p,1}\le 
C\big(t+|\ups|^{\frac{1}{p}}\big),\quad
\big\|\tilde{\eta}_{\ups,t\nu}\big\|_{\ups,p}\le 
C\big(t+|\ups|^{\frac{1}{p}}\big)^2.$$
Furthermore, if $z\!\in\!\Si$ and 
$(B_{b_{\ups}}(u_{\ups}(z),C\de),J,g_{V,b_{\ups}})$
is isometric to a ball in $\Bbb{C}^n$, then 
\hbox{$\tilde{\eta}_{\ups,t\nu}(z)=0$}.
\end{thm}

\noindent
{\it Remark:}
In specific applications, 
the main goal would be to express
the number of zeros of $\psi_{{\cal S},t\nu}^{\mu}$
in terms of the cohomology ring of 
a closure of~${\cal S}_{\cal T}(\mu)$.
One of the significant intermediate steps 
is to extract the leading-order terms from 
the section~$\psi_{{\cal S},t\nu}^{\mu}$.
If $\la_{\hat{0}}\!=\!0$, 
the estimate on $\nu_{\ups,t}$ given above
easily leads to  a sufficiently good estimate on 
$\pi_{\ups,-}^{0,1}\nu_{\ups,t}$; see~\cite{I} and~\cite{Z2}.
In such a case, one can also extract the first-order
term from $\pi_{\ups,-}^{0,1}\bar{\partial}u_{\ups}$,
which suffices for the computation in~\cite{I}.
A  power-series expansion for
$\pi_{\ups,-}^{0,1}\bar{\partial}u_{\ups}$ is given in~\cite{Z2},
where terms of up to third degree are used.
With the choice of metrics in~\cite{Z2},
the term $\pi_{\ups,-}^{0,1}\tilde{\eta}_{\ups,t\nu}$ vanishes.
The remaining term is shown to be secondary
for a good choice of the obstruction bundle setup.\\

\noindent
{\it Proof:} Let $\de,\ep\!>\!0$ be  as in Lemma~\ref{orient_l3}
and its proof. 
We take $\varphi_{{\cal S},t\nu}^{\mu }$ to be the section
descendent from the $G_{\cal T}$-equivarent section 
$\tilde{\varphi}_{{\cal S},t\nu}^{\mu}$.
Denote by $U_K'$ the open set $U_K$ of Lemma~\ref{orient_l3}.
By Corollary~\ref{surject_crl}, there exists a neighborhood  
$U_K$ of $K$ in \hbox{$\bar{C}^{\i}_{(\la^*;N)}(\Si;\mu)$}
such that ${\cal M}_{\Si,t\nu,\la^*}(\mu)\cap U_K$
is contained in 
$\tilde{\cal M}_{t\nu}(U_K',\de;\mu)$.
The neighborhood $U_K$ can always be chosen to contain 
all the zeros of the section $\tilde{\psi}_{t\nu}$
of the bundle $\Pi\Ga_-^{0,1}$ over~${\cal M}_{\Si,t\nu,\la^*}(\mu)\cap U_K$.
By Corollary~\ref{orient_c2}, ${\cal M}_{\Si,t\nu,\la^*}(\mu)\cap U_K$
is precisely the oriented zero set of the section~$\tilde{\psi}_{t\nu}$.
Since the map
$$F^{\eset}{\cal S}_{\de}|K
\lra\tilde{\cal M}_{t\nu}(U_K',\de;\mu),\quad
\ups\lra 
\tilde{\ga}_{{\cal T},t\nu}\big(\Phi_{\cal S}\Phi_{\cal S}^{\mu}
\big(\varphi_{{\cal S},t\nu}^{\mu}(\ups)\big)\big),$$
is an orientation-preserving diffeomorphism by Lemma~\ref{orient_l3},
it induces a sign-preserving bijection between the zero set of
$\tilde{\psi}_{t\nu}$ on $\tilde{\cal M}_{t\nu}(U_K',\de;\mu)$,
and the zero set of the section
$$(\tilde{\ga}_{{\cal T},t\nu}\Phi_{\cal S}
\Phi_{\cal S}^{\mu}\varphi_{{\cal S},t\nu}^{\mu})^*
\tilde{\psi}_{t\nu}\in
\Ga\big(F^{\eset}{\cal S}_{\de}|K;
(\tilde{\ga}_{{\cal T},t\nu}\Phi_{\cal S}
\Phi_{\cal S}^{\mu}\varphi_{{\cal S},t\nu}^{\mu})^*\Pi\Ga_-^{0,1}\big).$$
By equation~\e_ref{psi_dfn}, under the canonical identification
$$(\tilde{\ga}_{{\cal T},t\nu}\ga_{\cal T}^{-1})^*\Pi\Ga_-^{0,1}=\Ga_-^{0,1}
\lra {\cal M}(U_K',\de),$$
the section 
$(\tilde{\ga}_{{\cal T},t\nu}\ga_{\cal T}^{-1})^*\tilde{\psi}_{t\nu}$
corresponds to the section $\psi_{t\nu}$ given by
\begin{equation}\label{str_e3}\begin{split}
\psi_{t\nu}\big(b(\ups)\big)&=
t\nu|_{u_{\ups}}-\bar{\partial}u_{\ups}-
 \tilde{\pi}^{0,1}_{\ups,+}\eta_{\ups,t\nu}
-N_{\ups,t\nu}P_{\ups}\eta_{\ups,t\nu}\\
&=\pi_{\ups,-}^{0,1}\big(t\nu|_{u_{\ups}}-\bar{\partial}u_{\ups}-
 \tilde{\pi}^{0,1}_{\ups,+}\eta_{\ups,t\nu}
-N_{\ups,t\nu}P_{\ups}\eta_{\ups,t\nu}\big)\\
&=\pi_{\ups,-}^{0,1}\big(t\nu_{\ups,\xi_{\ups,t\nu}}-\bar{\partial}u_{\ups}-
 \tilde{\pi}^{0,1}_{\ups,+}\eta_{\ups,t\nu}
-N_{\ups}P_{\ups}\eta_{\ups,t\nu}\big).
\end{split}\end{equation}
The second equality above is automatic, 
since $\psi_{t\nu}\big(b(\ups)\big)\!\in\!\Ga_-^{0,1}(\ups)$;
the third follows from the definition of 
$N_{\ups,t\nu}$ in Subsection~\ref{gluing_map}.
The bounds on the terms 
$\nu_{\ups,\xi_{\ups,t\nu}}$,
$\eta_{\ups,t\nu}$,
$N_{\ups}P_{\ups}\eta_{\ups,t\nu}$ also 
follow from Subsection~\ref{gluing_map}.
By definition of~$\Ga_-^{0,1}$ in Subsection~\ref{scale_sec}
and equation~\e_ref{str_e3}, 
under the canonical identification
$$\ga_{\cal T}^*\Ga_-^{0,1}=
\pi_{F{\cal T}}^*\Ga_-^{0,1}
\lra F^{\eset}{\cal T}_{\de}\big|U_K',$$
the section $\ga_{\cal T}^*\psi_{t\nu}$ corresponds
to the section~$\psi_{{\cal T},t\nu}$, described in the statement
of the theorem.\\

\noindent
The next proposition describes a special case of 
the above theorem.
It is obtained by fixing a metric~$g$ on~$\Si$
and going through the construction analogous to that 
in Subsection~\ref{gluing_map} and
then modification for constraints as above.
The sign statement below follows
from the fact that the $(L^2,g,g_V)$-projection
\hbox{$\ker(D_{V,u_b}\!-\!\tau tL_{\nu,u_b})\!\lra\! \ker D_{V,u_b}$}
is an isomorphism for all $\tau\!\in\![0,1]$,
$t$~sufficiently small, and 
$b\!\in\!{\cal M}_{\Si,t\nu,\la^*}(\mu)$ sufficiently close to~$K$.

\begin{prp}
\label{holom_prp}
Suppose $\la^*\!\in\!H_2(V;\Bbb{Z})$,
${\cal T}\!=\!\big(\Si,[N],\{\hat{0}\};\hat{0},\la^*\big)$
is a \hbox{$(V,J)$-regular} bubble type,
$$\nu\!\in\!\Ga^{0,1}(\Si\times V;\La^{0,1}_{J,j}
                          \pi_{\Si}^*T^*\Si\otimes\pi_V^*TV)$$
is any section, and 
$\mu$ is a ${\cal T}$-regular
$N$-tuple of submanifolds of~$V$ of total codimension
$$\codim~\mu=
2\big(\lan c_1(V,J),\la^*\ran+(\dim~V)\big(1-g(S)\big)+|N|\big).$$
Then ${\cal M}_{\cal T}(\mu)$ is a discrete set
and for every finite subset $K$ of ${\cal S}(\mu)$,
there exist a neighborhood $U_K$ of $K$ in 
\hbox{$\bar{C}^{\i}_{(\la^*;N)}(\Si;\mu)$}, $\ep>0$,
and for each $t\!\in\!(0,\ep)$
a sign-preserving bijection between $K$
and \hbox{${\cal M}_{\Si,t\nu,\la^*}(\mu)\cap U_K$}.
\end{prp}

\subsection{Gluing Maps for Spaces $\bar{\cal U}_{\cal T}^{(0)}(\mu)$
and Orientations}
\label{orient1_sec}

\noindent
We now return to spaces the case
${\cal T}\!=\!(S,M,I;j,\la)$ is a regular bubble type.
Our primary interest is in the \hbox{case $S\!=\!S^2$},
so we assume that this is the case.
However, most of the analysis in this subsection applies to
any regular bubble type~${\cal T}$.
Let $\mu$ be a generic $\tilde{M}$-tuple of submanifolds in~$V$,
as defined below.
If $I\!=\!\bigsqcup\limits_{k\in K}\!I_k$ is the decomposition
of $I$ into rooted trees and $\{{\cal T}_k\}$ are the corresponding
simple types derived from~${\cal T}$, the product gluing map,
$$\big(\tilde{\ga}_{{\cal T}_k}\big)_{k\in K}\!:
\prod_{k\in K}F{\cal T}_{k,\de_k}\lra
\prod_{k\in K}\bar{\cal M}_{\lr{{\cal T}_k}},$$
may not map the total space of the bundle over 
${\cal U}_{\cal T}^{(0)}(\mu)$ 
into~$\bar{\cal U}_{\lr{\cal T}}^{(0)}(\mu)$.
In this subsection, we remedy this deficiency of the product gluing map.
We also show that all the spaces 
$\bar{\cal U}_{\lr{\cal T}}^{(0)}(\mu)$ and 
$\bar{\cal U}_{\lr{\cal T}}(\mu)$
are naturally oriented topological orbifolds and 
the gluing maps defined below preserve orientations.

\begin{dfn}
\label{regul_dfn2}
If ${\cal T}$ is a $(V,J)$-regular bubble type,
tuple $\mu$ of oriented submanifolds of $V$ is \under{${\cal T}$-regular}
if 
\begin{gather*}
\bigoplus_{l\in M}T_{u_{j_l}(y_l)}V=
 \hbox{Im}~d\ev_{[N]}\big|_b
+\bigoplus_{l\in M}T_{u_{j_l}(y_l)}\mu_l
\qquad\forall~ b=(S^2,M,I;x,(j,y),u)
                  \!\in\! {\cal B}_{\cal T}(\mu);\\
\hbox{where}\qquad
d\ev_{[N]}\big|_b\!:\tilde{\cal K}_b{\cal T}\lra
\bigoplus_{l\in M}T_{u_{j_l}(y_l)}V,\quad
d\ev_l\big|_b\big(\xi_I,w_{\hat{I}+[N]}\big)=
       \xi_{j_l}(y_l)+du_{j_l}\big|_{y_l}w_l,
\end{gather*}
\end{dfn}

\noindent
Let ${\cal T}$, ${\cal T}_k$, $K$, $\mu$, and $b$ be as above.
Denote by 
$b_k\!=\!\big(S^2,M_k,I_k;x|\hat{I}_k,(j,y)|M_k,u|I_k\big)$
the corresponding  ${\cal T}_k$-bubble map;
see Subsection~\ref{stratification}.
Let ${\cal N}_b^{\mu}{\cal T}$ be the 
$(L^2,b)$-orthogonal complement~of
\begin{equation*}
\tilde{\cal K}_b^{\mu}{\cal T}=
\big\{ (\xi,w_{\hat{I}+M})\!\in\!\tilde{\cal K}_b{\cal T}\!:
\xi_{j_l}(y_l)+du_{j_l}|_{y_l}w_l\!\in\! T_{u_{j_l}(y_l)}\mu_l
~\forall l\!\in\! M\cap\tilde{M}\big\}\subset 
\tilde{\cal K}_b{\cal T}
\end{equation*}
in $\bigoplus\limits_{k\in K}\!{\cal K}_{b_k}{\cal T}_k$.
Denote by 
$\tilde{\cal N}^{\mu}{\cal T}\!\!\lra\!{\cal B}_{\cal T}(\mu)$
and  ${\cal N}^{\mu}{\cal T}\!\!\lra\!{\cal U}_{\cal T}^{(0)}(\mu)$
the corresponding vector bundles.~Let 
\begin{gather*}
{\cal N}^{\mu}{\cal B}_{\cal T}=
\tilde{\cal N}^{\mu}{\cal T}\oplus(\Bbb{C}\oplus\Bbb{R})^K
\lra {\cal B}_{\cal T}(\mu),\quad
{\cal N}^{\mu}{\cal U}_{\cal T}^{(0)}=
{\cal N}^{\mu}{\cal T}\oplus(\Bbb{C}\oplus\Bbb{R})^K\lra
{\cal U}_{\cal T}^{(0)}(\mu);\\
F^{(0)}{\cal T}=\bigoplus_{k\in K}F^{(0)}{\cal T}_k
\lra {\cal B}_{\cal T},\quad
F{\cal T}=\bigoplus_{k\in K}F{\cal T}_k\lra
{\cal U}_{\cal T}^{(0)}.
\end{gather*}
The last two vector bundles carry norms induced from
the norms on $F^{(0)}{\cal T}_k$, while we 
define norms on the first two by
$$\l|\l(b,\si,(c,r)\r)\r|=\|\si\|_{b,C^0}+|(c,r)|,$$
if $\si\!\in\!{\cal N}_b^{\mu}{\cal T}\subset
\bigoplus\limits_{k\in K}\tilde{\cal K}_{b_k}{\cal T}_k$ and
$(c,r)\!\in\!(\Bbb{C}\oplus\Bbb{R})^K$.
If $\de$ is sufficiently small, define
$$\tilde{\phi}_{\cal T}^{\mu}\!: 
{\cal N}^{\mu}{\cal B}_{{\cal T},\de}\lra\!
\prod_{k\in K}{\cal M}_{{\cal T}_k}^{(0)}
\quad\hbox{by}\quad
\tilde{\phi}_{\cal T}^{\mu}\big(\si,(c,r)\big)=
\big((c_k,r_k)\cdot H_{{\cal T}_k}^{(0)}(\si_k)\big)_{k\in K}
\!\in\! \prod_{k\in K}{\cal M}_{{\cal T}_k}^{(0)},$$
where $H_{{\cal T}_k}^{(0)}$ is as 
at the end of Subsection~\ref{str_sub1}
and $(c_k,r_k)\cdot$ denotes the action~of 
a neighborhood of 
$$0\in\Bbb{C}\times\Bbb{R}=\Bbb{C}\times\Bbb{R}\times\{0\}\subset
\Bbb{C}\times\Bbb{R}\times\Bbb{R}$$
described in Subsection~\ref{summary3}.
Since $\tilde{\phi}_{\cal T}^{\mu}$ is 
$\big(\tilde{G}_{\cal T}\!\times\!\hbox{Aut}({\cal T})\big)$-equivarent, 
it descends to a $G_{\cal T}^*$-equivarent map
$$\phi_{\cal T}^{\mu}\!:
{\cal N}^{\mu}{\cal U}_{{\cal T},\de}^{(0)}
\lra \prod_{k\in K}{\cal M}_{{\cal T}_k}.$$
Let $\Phi_{\cal T}^{\mu}\!: 
     \pi_{{\cal N}^{\mu}{\cal U}_{\cal T}}^*F{\cal T}\!\lra\! F{\cal T}$
be a $G_{\cal T}^*$-equivarent vector-bundle map covering
the map~$\phi_{\cal T}^{\mu}$ such that $\Phi_{\cal T}^{\mu}$ is 
the identity over~${\cal U}_{\cal T}^{(0)}(\mu)$.
Denote by $\tilde{\Phi}_{\cal T}^{\mu}$ the lift of $\Phi_{\cal T}^{\mu}$ 
to~${\cal N}^{\mu}{\cal B}_{{\cal T},\de}$.
Let  $\Phi_{{\cal T},k}^{\mu}$ and $\tilde{\Phi}_{{\cal T},k}^{\mu}$
be $k$th components of $\Phi_{\cal T}^{\mu}$ and 
$\tilde{\Phi}_{\cal T}^{\mu}$, respectively.

\begin{lmm}
\label{rat_l1}
With notation as above, there exist 
$\big(\tilde{G}_{\cal T}\!\times\!\hbox{Aut}({\cal T})\big)$-invariant 
functions 
$\de,C\!\in\! C^{\i}({\cal B}_{\cal T}(\mu);\Bbb{R}^+)$ and 
a $\big(\tilde{G}_{\cal T}\!\times\!\hbox{Aut}({\cal T})\big)$-equivarent 
section
$$\tilde{\varphi}_{\cal T}^{\mu}\!\in\!
\Ga\big( F^{(0)}{\cal T}_{\de};
\pi_{F^{(0)}{\cal T}}^*{\cal N}^{\mu}{\cal B}_{{\cal T},\de}\big),$$
such that
$\big|\tilde{\varphi}_{\cal T}^{\mu}(\ups)\big|\le 
C(b_{\ups})|\ups|^{\frac{1}{p}}$ and
$$ F{\cal T}_{\de}\lra
 \bar{\cal U}_{\lr{\cal T}}^{(0)}(\mu),\quad
\ups\lra \big(\tilde{\ga}_{{\cal T}_k}\big(\Phi_{{\cal T},k}^{\mu}
\varphi_{\cal T}^{\mu}(\ups)\big)\big)_{k\in K},$$
is a homeomorphism onto an open neighborhood of 
 ${\cal U}_{\cal T}^{(0)}(\mu)$ in
$\bar{\cal U}_{\lr{\cal T}}^{(0)}(\mu)$.
Furthermore, the restriction of this map to  $F^{\eset}{\cal T}_{\de}$
is an orientation-preserving diffeomorphism onto an open subset 
of~${\cal U}_{\lr{\cal T}}^{(0)}(\mu)$.
\end{lmm}

\noindent
{\it Proof:} Denote by ${\cal N}_{\cal T}\mu$ the normal bundle of
$$X_{\cal T}(\mu)\equiv
\big\{x_K\!\in\! V^K\!\!: x_{k_1}\!=\!x_{k_2}~\forall k_1,k_2\!\in\! K\big\}
\times\prod_{l\in M\cap\tilde{M}}\mu_l\subset
V^K\times V^{M\cap\tilde{M}}.$$
Let $\tilde{\cal N}_{\cal T}\mu=
{\cal N}_{\cal T}\mu\oplus(\Bbb{C}\oplus\Bbb{R})^K$.
Since the $\tilde{G}_{\cal T}$-action does not change any evaluation maps
and the constraints are in general position,
the differential of the map
\begin{gather*}
\Psi_{{\cal T},\tilde{M}}\!: 
\prod_{k\in K}{\cal M}_{{\cal T}_k}^{(0)}\lra
V^K\times V^{M\cap\tilde{M}}\times(\Bbb{C}\times\Bbb{R})^K,\\
\Psi_{{\cal T},\tilde{M}}=
\l((\ev_{{\cal T}_k}(b_k))_{k\in K};
(\ev_l(b))_{l\in M\cap\tilde{M}};
(\Psi_{{\cal T}_k,\hat{0}_{I_k}}(b_k))_{k\in K}\r),
\end{gather*}
where $\Psi_{{\cal T}_k,\hat{0}_{I_k}}(b_k)\!\in\!\Bbb{C}\times\Bbb{R}$ is
as in Subsection~\ref{str_sub1}, induces an isomorphism
between ${\cal N}_b^{\mu}{\cal B}_{\cal T}$ 
and~$\tilde{\cal N}_{\cal T}\mu$.
This isomorphism is orientation-preserving by definition of 
orientations.
Thus,
$$\tilde{\phi}_{\cal T}^{\mu}\!:
{\cal N}^{\mu}{\cal B}_{{\cal T},\de}\lra
\prod_{k\in K}{\cal M}_{{\cal T}_k}^{(0)}$$
is an orientation-preserving diffeomorphism onto an open neighborhood
of ${\cal B}_{\cal T}(\mu)$ 
in~$\prod\limits_{k\in K}{\cal M}_{{\cal T}_k}^{(0)}$,
provided $\de\!\in\! C^{\i}({\cal B}_{\cal T}(\mu);\Bbb{R}^+)$
is sufficiently small.
By the same argument as in Subsection~\ref{orient2_sec},
for any simple bubble type ${\cal T}'$, the map
$$\tilde{\ga}_{{\cal T}'}\!:F^{\eset}{\cal T}'_{\de}\lra
{\cal M}_{\lr{{\cal T}'}}= {\cal H}_{\lr{{\cal T}'}}$$
is an orientation-preserving diffeomorphism onto an open subset
of ${\cal M}_{\lr{{\cal T}'}}$ provided 
$\de\!\in\! C^{\i}({\cal M}_{{\cal T}'};\Bbb{R}^+)$ is sufficiently small.
Along with Corollary~\ref{surject_crl2}, this implies that
the product map
$$\prod_{k\in K}\tilde{\ga}_{{\cal T}_k}\!: \prod_{k\in K}F{\cal T}_{k,\de}
\lra \prod_{k\in K}\bar{\cal M}_{\lr{{\cal T}_k}}$$
is a homeomorphism onto an open neighborhood of 
$\prod\limits_{k\in K}{\cal M}_{{\cal T}_k}$ in
$\prod\limits_{k\in K}\bar{\cal M}_{\lr{{\cal T}_k}}$
and its restriction to the preimage of 
$\prod\limits_{k\in K}{\cal M}_{\lr{{\cal T}_k}}$ is 
an orientation-preserving diffeomorphism.
The lemma now follows by applying an argument similar to
the proof of Lemma~\ref{orient_l3} to the functions
$$h(\ups)=\Psi_{{\cal T},\tilde{M}}\big(
       \tilde{\phi}_{\cal T}^{\mu}\big(\si,(c,r)\big)\big),\quad
\tilde{h}(\ups)=\tilde{h}_0(\ups)=
\Psi_{{\cal T},\tilde{M}}\big(\big(\tilde{\ga}_k
\Phi_{{\cal T},k}^{\mu}(\ups)\big)_{k\in K}\big),$$
where we write
$\ups\!=\!\big(\si,(c,r),v\big)$,
with $\big(\si,(c,r)\big)\!\in\!{\cal N}^{\mu}{\cal B}_{\cal T}$
and $v\!\in\! F^{(0)}{\cal T}$.
Since ${\cal B}_{\cal T}(\mu)$ is generally not precompact
in $\prod\limits_{k\in K}{\cal M}_{{\cal T}_k}^{(0)}$,
we end up with \hbox{$\de,C\!\in\! C^{\i}({\cal B}_{\cal T}(\mu);\Bbb{R}^+)$},
instead \hbox{of $\de,C\!\in\!\Bbb{R}^+$}.
Another difference is that $\tilde{h}$ is not necessarily smooth
with respect to the standard smooth structure 
\hbox{on~${\cal N}^{\mu}{\cal B}_{\cal T}\oplus F^{(0)}{\cal T}$}.
However, we can put a smooth structure on the total space 
such that the composite maps
$${\cal N}^{\mu}{\cal B}_{{\cal T},\de}\oplus F^{(0)}{\cal T}
\lra F^{(0)}{\cal T} \lra\Bbb{R},\quad
\ups\lra\tilde{\Phi}_{\cal T}^{\mu}(\ups),~
v_h\lra|v_h|^{\frac{1}{3p}},\quad h\!\in\!\hat{I}_k,~k\!\in\! K,$$
are smooth, whenever $\de\!\in\! C^{\i}({\cal B}_{\cal T}(\mu);\Bbb{R}^+)$
is sufficiently small.
Then by Corollary~\ref{cont_c4}, $\tilde{h}$ is $C^2$,
which is sufficient for the arguments of Subsection~\ref{impl}.
Finally, in the given case $\tilde{h}$ is defined on 
all of 
$\big({\cal N}^{\mu}{\cal B}_{\cal T}\!\oplus\! F^{(0)}{\cal T}\big)_{\de}$
and thus the second condition on $\tilde{h}_t$ in 
Subsection~\ref{impl} is redundant.\\

\noindent
Suppose ${\cal T}$ is a bubble type and 
$\mu$ is an $\tilde{M}$-tuple of constraints in general position.
By Lemma~\ref{rat_l1}, there exist 
$G_{\cal T}^*$-invariant functions
$\de,C\!\in\! C^{\i}({\cal U}^{(0)}_{\cal T}(\mu);\Bbb{R}^+)$
and a $G_{\cal T}^*$-equivarent section 
$$\varphi_{\cal T}^{\mu}\in
\big(F{\cal T}_{\mu};
\pi_{F{\cal T}}^*{\cal N}^{\mu}{\cal U}_{{\cal T},\de}^{(0)}\big)$$
such~that
$\big|\varphi_{\cal T}^{\mu}(\ups)\big|\le
C(b_{\ups})|\ups|^{\frac{1}{p}}$ and
$$\tilde{\ga}_{\cal T}^{\mu}\!: 
F{\cal T}_{\de}\lra \bar{\cal U}_{\lr{\cal T}}^{(0)}(\mu),\quad
\tilde{\ga}_{\cal T}^{\mu}(\ups)=
\big(\tilde{\ga}_{{\cal T}_k}\big( 
 \Phi_{{\cal T},k}^{\mu}\varphi_{\cal T}^{\mu}(\ups)\big)\big)_{k\in K},$$
is a homeomorphism onto a neighborhood of 
${\cal U}_{\cal T}^{(0)}(\mu)$ in
$\bar{\cal U}_{\lr{\cal T}}^{(0)}(\mu)$, which is an orientation-preserving
diffeomorphism on a dense open subset of the domain.
If ${\cal T}'$ is another regular bubble type such that
\hbox{$\lr{\cal T}=\lr{{\cal T}'}$}
and $\mu$ is \hbox{${\cal T}'$-regular},
 it follows~that
$$\tilde{\ga}_{\cal T}^{\mu~-1}\tilde{\ga}_{{\cal T}'}^{\mu}\!:
\tilde{\ga}_{{\cal T}'}^{\mu~-1}\big(
\tilde{\ga}_{\cal T}^{\mu}(F{\cal T}_{\de})\big)\lra
\tilde{\ga}_{\cal T}^{\mu~-1}\big(
\tilde{\ga}_{{\cal T}'}^{\mu}(F{\cal T}'_{\de'})\big)$$
is a homeomorphism provided 
$\de'\!\in\! C^{\i}({\cal U}^{(0)}_{{\cal T}'}(\mu);\Bbb{R}^+)$
is sufficiently small.
Furthermore, by the above it is orientation-preserving
on a dense open subset of its domain.
It follows that 
$\tilde{\ga}_{\cal T}^{\mu-1}\tilde{\ga}_{{\cal T}'}^{\mu}$
is an orientation-preserving homeomorphism everywhere.
We thus obtain

\begin{thm}
\label{str_global}
Let ${\cal T}^*\!=\!(S^2,M,I^*;j,\la^*)$ 
be a basic bubble type and $\mu$ an $M$-tuple of constraints 
such that $\mu$ is ${\cal T}$-regular for every bubble 
\hbox{type ${\cal T}\le{\cal T}^*$}.\\
(1) The spaces  $\bar{\cal U}_{{\cal T}^*}^{(0)}(\mu)$ and 
$\bar{\cal U}_{{\cal T}^*}(\mu)$ are oriented topological orbifolds.\\
(2) Suppose ${\cal T}\!<\!{\cal T}^*$, 
$\phi_{\cal T}^{\mu}\!:
{\cal N}^{\mu}{\cal T}_{\tilde{\de}}\lra{\cal U}_{\cal T}^{(0)} $
is a $G_{{\cal T}^*}$-equivarent identification
of neighborhoods of ${\cal U}_{\cal T}^{(0)}(\mu)$ 
in ${\cal N}^{\mu}{\cal T}$ and in ${\cal U}_{\cal T}^{(0)}$, 
and $\Phi_{\cal T}^{\mu}\!:\pi_{{\cal N}^{\mu}{\cal T}}^*F{\cal T}
\!\lra\! F{\cal T}$ is a lift of $\phi_{\cal T}^{\mu}$
such that \hbox{$\Phi_{\cal T}^{\mu}|{\cal U}_{\cal T}^{(0)}(\mu)\!=\!1$}.
Then there exist \hbox{$G_{{\cal T}^*}$-invariant} functions
$\de,C\!\in\! C^{\i}\big({\cal U}_{\cal T}^{(0)}(\mu);\Bbb{R}^+\big)$
and a $G_{{\cal T}^*}$-equivarent continuous orientation-preserving
identification 
$$\tilde{\ga}_{\cal T}^{\mu}:  F{\cal T}_{\de}
       \lra \bar{\cal U}_{{\cal T}^*}^{(0)}(\mu),$$
of neighborhoods of ${\cal U}_{\cal T}^{(0)}(\mu)$
in $F{\cal T}$ and in $\bar{\cal U}_{{\cal T}^*}^{(0)}(\mu)$,
which is smooth on 
$F^{\eset}{\cal T}_{\de}\!\lra\!{\cal U}_{{\cal T}^*}(\mu)$.
Furthermore, for every $\ups\!\in\! F{\cal T}_{\de}$,
there exists $\si(\ups)\!\in\!{\cal N}^{\mu}{\cal T}$ such that
$$\|\si(\ups)\|_{b^*}\le C(b^*)|\ups|^{\frac{1}{p}}
\quad\hbox{and}\quad
u_{\tilde{\ga}_{\cal T}^{\mu}(\ups)}=
\exp_{V,u_{b'}\circ q_{\Phi_{\cal T}^{\mu}(\si(\ups))}}\xi_{\ups},
\quad\hbox{where}\quad
\Phi_{\cal T}^{\mu}(\si(\ups))\equiv[b',v'],$$
for some $\xi_{\ups}\!\in\!
 \Ga\big(u_{b'}\circ q_{\Phi_{\cal T}^{\mu}(\si(\ups))}\big)$
with $\|\xi_{\ups}\|_{C^0}\le C(b_{\ups})|\ups|^{\frac{1}{p}}$.
\end{thm}

\noindent
{\it Remark:} The descriptive statement (2) of Theorem~\ref{str_global}
is used in~\cite{Z2} for local excess-intersection type of
computations on the spaces~$\bar{\cal U}_{{\cal T}^*}(\mu)$.

\section{Technical Issues}
\label{tech_sec}

\subsection{Continuity of the Gluing Map}
\label{map_cont}

\noindent
Let ${\cal T}\!=\!(S,M,I;j,\la)$ be a simple regular bubble
type and $H$ a nonempty subset of $\hat{I}$.
Suppose $\ups_k\!\in\! F^{(\eset)}{\cal T}_{\de}$ and the sequence
$\ups_k$ converges to $\ups^*\!\in\! F^{(H)}{\cal T}_{\de}$.
In this section, we show that $\tilde{\ga}_{\cal T}(\ups_k)$
converges to $\tilde{\ga}_{\cal T}(\ups^*)$ in the Gromov topology.
Our main interest is the \hbox{case $S=S^2$}.\\

\noindent
It is sufficient to assume that $\pi_h(\ups_k)\!=\!\pi_h(\ups^*)$ if 
$h\!\not\in\! H$.
In particular, $b\equiv b_{\ups^*}=b_{\ups_k}$.
Denote by \hbox{$\tilde{\cal T}={\cal T}(H)$} the bubble type of $b(\ups^*)$.
For each $k$, define
$$\tilde{\ups}_k=\big(\tilde{b}(\ups^*),(\tilde{v}_k)_H\big)
\in F_{b(\ups^*)}^{(0)}\tilde{\cal T}$$
as  follows. 
If $h\!\in\! H$, put
$$i_Hh=\min\l\{i<h\!: \hbox{if~}h'\!\in\! I~\&~ i< h'<h, 
       h'\!\not\in\! H\r\}.$$
Since $I$ is a rooted $i_Hh$ is well-defined.
Let
$$\tilde{v}_{k,h}=\prod_{i_Hh< h'\le h}v_{k,h'}.$$
Since $\ups_k\!\lra\!\ups^*$, $\tilde{\ups}_{k,h}\!\lra\!0$ 
for all $h\!\in\! H$.
Furthermore, by construction $\Si_{\ups_k}=\Si_{\tilde{\ups}_k}$ and 
$q_{\ups_k}\!=\!q_{\ups^*}\!\circ\! q_{\tilde{\ups}_k}$.
In particular, $u_{\ups_k}\!=\!u_{\ups^*}\!\circ\! q_{\tilde{\ups}_k}$.\\

\noindent
For any $h\!\in\! H$ and $\de>0$, let
\begin{gather*}
A_{h,\de,k}=q_{\ups_k}^{-1}
\l(\l\{(\io_h,z)\!: r_{b,h}(z)\le\de\r\}\cup
\l\{(h,z)\!: |q_S^{-1}(z)|\le\de\r\}\r)\subset\Si_{\ups_k},\\
A_{h,\de}^*=q_{\ups^*}^{-1}
\l(\l\{(\io_h,z)\!: r_{b,h}(z)\le\de\r\}\cup
\l\{(h,z)\!: |q_S^{-1}(z)|\le\de\r\}\r)\subset\Si_{\ups^*},\quad
\Si_{\de}^*=\Si_{\ups^*}-\!\bigcup_{h\in H}A_{h,\de}^*.
\end{gather*}
It is convenient to make the following definitions.
If $\eta_k\!\in\! L^p(\ups_k)$, the sequence $\{\eta_k\}$ {\it converges to}
\hbox{$\eta^*\!\in\! L^p(\ups^*)$} if $q_{\tilde{\ups}_k}^{-1*}\eta_k$
converges to $\eta^*$ in the $L^p$-norm on all precompact open
subsets of $\Si_{\ups^*}^*$ and
\begin{equation}\label{map_cont_e3a}
\lim_{\de\lra0}\lim_{k\lra 0}\|\eta_k\|_{\ups_k,L^p(A_{h,\de,k})}=0
\qquad\forall h\!\in\! H.
\end{equation}
If $\xi_k\!\in\! L^p_1(\ups_k)$, the sequence $\{\xi_k\}$ {\it converges to}
$\xi^*\!\in\! L^p_1(\ups^*)$ if $\xi_k\circ q_{\tilde{\ups}_k}^{-1}$
converges to $\xi^*$ in the $L^p_1$-norm on all precompact open
subsets of $\Si_{\ups^*}^*$ and
\begin{equation}\label{map_cont_e3b}
\lim_{\de\lra0}\lim_{k\lra 0}\|\xi_k\|_{\ups_k,L^p_1(A_{h,\de,k})}=0
\qquad\forall h\!\in\! H.
\end{equation}
In \e_ref{map_cont_e3a} and \e_ref{map_cont_e3b},
we use the modified Sobolev norms.

\begin{lmm}
\label{cont_l1a}
There exist $C,\de\!\in\! C^{\i}\big({\cal M}_{\cal T}^{(0)};\Bbb{R}^+\big)$
such that for any sequence 
\hbox{$\{\ups_k\!\in\! F^{(\eset)}{\cal T}_{\de}\}$} converging to
$\ups^*$ as above and  $\xi\!\in\!\Ga_+(\ups^*)$
$$\l\|\pi_{\ups_k,-}(\xi\circ q_{\tilde{\ups}_k})\r\|_{\ups_k,p,1}
\le C(b)|\ups_k-\ups^*|\|\xi\|_{\ups^*,p,1}.$$
\end{lmm}

\noindent
{\it Proof:} Note that $\Ga_-(\ups_k)=
 \big\{\xi^-\circ q_{\tilde{\ups}_k}\!: \xi^-\!\in\!\Ga(\ups^*)\big\}$.
Thus, the difference 
$q_{\tilde{\ups}_k}^*\pi_{\ups^*,-}\!-\!\pi_{\ups_k,-}q_{\tilde{\ups}_k}^*$
arises entirely from the difference between the metrics
$q_{\tilde{\ups}_k}^*g_{\ups^*}$ and $g_{\ups_k}$.
By construction, the two metrics differ only on
the annuli $A_{h,2|v_{k,h}|^{\frac{1}{2}},k}$ for $h\!\in\! H$.
Thus, the claim follows from (2) of Lemma~\ref{approx_maps}.

\begin{lmm}
\label{cont_l1}
If $\eta_k$ converges to $\eta^*$, then $P_{\ups_k}\eta_k$ converges
to $P_{\ups^*}\eta^*$.
\end{lmm}

\noindent
{\it Proof:} (1) Let $\{\ep_k\},\{\de_k\}\subset\Bbb{R}^+$ be sequences 
converging to zero such that
\begin{equation}
\label{cont_l1_e1a}
\l\|\eta^*\r\|_{\ups^*,L^p(A_{h,\de_k}^*)}\le\ep_k;\quad
\l\|P_{\ups^*}\eta^*\r\|_{\ups^*,L^p_1(A_{h,\de_k}^*)}\le\ep_k,\quad
\lim_{k\lra\i}\l\|\eta_k\r\|_{\ups_k,L^p(A_{h,\de_{k^*},k})}<\ep_k
~~\forall h\!\in\! H.
\end{equation}
For every $k^*>0$, choose $N_{k^*}$ such that for all $k>N_{k^*}$
\begin{gather}
\label{cont_l1_e1b}
\big\|q_{\tilde{\ups}_k}^{-1*}\eta_k-\eta^*
 \big\|_{\ups^*,L^p(\Si_{\de_{k^*}}^*)}\le \ep_{k^*} 
\quad\hbox{and}\quad
\l\|\eta_k\r\|_{\ups_k,L^p(A_{h,\de_{k^*},k})}\le\ep_{k^*}~~
\forall h\!\in\! H.
\end{gather}
It can be assumed that $2|\ups_k-\ups^*|^{\frac{1}{2}}
 \le\de_{k^*},\ep_{k^*}$ whenever $k>N_{k^*}$.
For any $k>N_{k^*}$, let $\tilde{\eta}_{k^*,k}\!\in\! L^p(\ups^*)$ be given~by
$$\tilde{\eta}_{k^*,k}=
\begin{cases}
q_{\tilde{\ups}_k}^{-1*}\tilde{\eta}_k,
&\hbox{on~}\Si_{\de_{k^*}}^*;\\
0,&\hbox{outside of~}\Si_{\de_{k^*}}^*.\\
\end{cases}$$
Then
$\|\tilde{\eta}_{k^*,k}\|_{\ups^*,p}\le \|\eta_k\|_{\ups_k,p}$.
Let
$$\tilde{P}_{k^*,k}\eta_k=	 
q_{\tilde{\ups}_k}^*P_{\ups^*}\tilde{\eta}_{k^*,k}\in L^p_1(\ups_k).$$
Then by Lemma~\ref{approx_maps2} and the first assumptions
of \e_ref{cont_l1_e1a} and \e_ref{cont_l1_e1b},
\begin{equation}\label{cont_l1_e2a}\begin{split}
\l\|q_{\tilde{\ups}_k}^{-1*}\tilde{P}_{k^*,k}\eta_k-
P_{\ups^*}\eta^*\r\|_{\ups^*,L^p_1(\Si_{\de_{k^*}}^*)}
&\le \l\|P_{\ups^*}\tilde{\eta}_{k^*,k}-P_{\ups^*}\eta^*\r\|_{\ups^*,p,1}\\
&\le C(b)\l\|\tilde{\eta}_{k^*,k}-\eta^*\r\|_{\ups^*,p}
\le 2C(b)\ep_{k^*}.
\end{split}\end{equation}
Since $\|dq_{\tilde{\ups}_k}\|_{C^0}\le C(b)$,
by \e_ref{cont_l1_e1a} and the first assumption of \e_ref{cont_l1_e1b}
for all $h\!\in\! H$,
\begin{equation}\label{cont_l1_e2b}\begin{split}
&\big\|\tilde{P}_{k^*,k}\eta_k\big\|_{\ups_k,L^p_1(A_{h,\de_{k^*},k})}
\le C(b)
\l\|P_{\ups^*}\tilde{\eta}_{k^*,k}\r\|_{\ups^*,L^p_1(A_{h,\de_{k^*}}^*)}\\
&\qquad\qquad\qquad\le  C(b)\Big( 
\big\|P_{\ups^*}\eta^*\big\|_{\ups^*,L^p_1(A_{h,\de_{k^*}}^*)}
+\l\|P_{\ups^*}\tilde{\eta}_{k^*,k}-P_{\ups^*}\eta^*\r\|_{\ups^*,p,1}\Big)
\le C'(b)\ep_{k^*}.
\end{split}\end{equation}
(2) We now show that $\tilde{P}_{k^*,k}\eta_k$ is close to $P_{\ups_k}\eta_k$.
By Lemmas~\ref{approx_maps2} and~\ref{cont_l1a},
\begin{equation}
\label{cont_l1_e3}\begin{split}
\big\|\tilde{P}_{k^*,k}\eta_k-P_{\ups_k}\eta_k\big\|_{\ups_k,p,1}
&\le C(b)\Big( \big\|D_{\ups_k}\tilde{P}_{k^*,k}\eta_k-\eta_k\big\|_{\ups_k,p}
+\big\|\pi_{\ups_k,-}\tilde{P}_{k^*,k}\eta_k\big\|_{\ups_k,p,1}\Big)\\
&\le C(b)\Big(
\big\|D_{\ups_k}\tilde{P}_{k^*,k}\eta_k-\eta_k\big\|_{\ups_k,p}+
|\ups_k-\ups^*|\|\eta_k\|_{\ups_k,p}\Big).
\end{split}
\end{equation}
Since $D_{\ups^*}P_{\ups^*}\tilde{\eta}_{k^*,k}=\tilde{\eta}_{k^*,k}$ and
$q_{\tilde{\ups}_k}$ is holomorphic outside of the annuli
$A_{h,\de_{k^*},k}$,
\begin{equation}
\label{cont_l1_e4a}
D_{\ups_k}\tilde{P}_{k^*,k}\eta_k=\eta_k
\quad\hbox{on}~~
 \Si_{\ups_k}-\bigcup\limits_{h\in H}A_{h,\de_{k^*},k};
\end{equation}
By equation~\e_ref{cont_l1_e2b},
\begin{equation}\label{cont_l1_e4b}
\big\|D_{\ups_k}\tilde{P}_{k^*,k}\eta_k\big\|_{\ups_k,L^p(A_{h,\de_{k^*},k})}
\le C(b)\big\|\tilde{P}_{k^*,k}\eta_k\big\|_{\ups_k,L^p_1(A_{h,\de_{k^*},k})}
\le C'(b)\ep_{k^*}.
\end{equation}
Thus, from equations~\e_ref{cont_l1_e3}-\e_ref{cont_l1_e4b}
and the second assumption of~\e_ref{cont_l1_e1b},
we conclude that  for all $k>N_{k^*}$
\begin{equation}\label{cont_l1_e5}
\big\|\tilde{P}_{k^*,k}\eta_k-P_{\ups_k}\eta_k\big\|_{\ups_k,p,1}
\le C(b)\ep_{k^*}\l(1+\|\eta^*\|_{\ups^*,p}\r).
\end{equation}
Since $\|dq_{\tilde{\ups}_k}^{-1}\|_{C^0}\le C(b)$ on $\Si_{\de_{k^*}}^*$,
by equations~\e_ref{cont_l1_e2a}, \e_ref{cont_l1_e2b},
and \e_ref{cont_l1_e5},
\begin{gather}\label{cont_l1_e6a}
\big\|q_{\tilde{\ups}_k}^{-1*}P_{\ups_k}\eta_k-P_{\ups^*}\eta^*
\big\|_{\ups^*,L^p_1(\Si_{\de_{k^*}}^*)}\le
C(b)\ep_{k^*}\big(1+\|\eta^*\|_{\ups^*,p}\big);\\
\label{cont_l1_e6b}
\big\|P_{\ups_k}\eta_k\big\|_{\ups_k,L^p_1(A_{h,\de_{k^*},k})}
\le C(b)\ep_{k^*}\big(1+\|\eta^*\|_{\ups^*,p}\big)\quad\forall h\!\in\! H.
\end{gather}
By equations~\e_ref{cont_l1_e6a} and \e_ref{cont_l1_e6b},
$P_{\ups_k}\eta_k$ converges to $P_{\ups^*}\eta^*$.

\begin{lmm}
\label{cont_l2}
There exist $\tilde{C},\de\!\in\! C^{\i}({\cal M}_{\cal T}^{(0)};\Bbb{R}^+)$
such that for all $\ups^*\!\in\! F^{(H)}{\cal T}_{\de}$ and $h\!\in\! H$,
$$\big\|P_{\ups^*}\bar{\partial}u_{\ups^*}\big\|_{g_{\ups^*},
C^1(A_{h,\de(b_{\ups^*})}^*)}\le \tilde{C}(b_{\ups^*}).$$
\end{lmm}

\noindent
{\it Proof:} For each $h\!\in\! H$, this lemma is obtained 
by pasting 
$(P_{\ups^*}\bar{\partial}u_{\ups^*})|(A_{h,\de(b_{\ups^*}^*)}\cap 
 \Si_{\ups^*,\io_h})$
and $(P_{\ups^*}\bar{\partial}u_{\ups^*})|(A_{h,\de(b_{\ups^*}^*)}\cap 
 \Si_{\ups^*,h})$
onto $\Si_{b_{\ups^*},\io_h}$ and $\Si_{b_{\ups^*},h}$
via a cutoff function.
We then use the usual elliptic estimates and Sobolev inequalities
on $\Si_{b_{\ups^*},\io_h}$ and $\Si_{b_{\ups^*},h}$
along with 
$$\l\|P_{\ups^*}\bar{\partial}u_{\ups^*}\r\|_{\ups^*,p,1}\le 
C(b_{\ups^*})|\ups^*|^{\frac{1}{p}}.$$
The bound obtained in this way is actually 
$C(b_{\ups^*})|\ups^*|^{\frac{1}{p}}$.

\begin{crl}
\label{cont_c3}
There exist $C,\de\!\in\! C^{\i}({\cal M}_{\cal T}^{(0)};\Bbb{R}^+)$
such that for any sequence $\ups_k\!\in\! F^{(\eset)}{\cal T}_{\de}$
converging to $\ups^*\!\in\! F^{(H)}{\cal T}_{\de}$ as above,
\begin{gather*}
\big\|q_{\tilde{\ups}_k}^{-1*}\eta_{\ups_k}-\eta_{\ups^*}
\big\|_{\ups^*,L^p(\Si_{2|\ups_k-\ups^*|^{1/2}}^*)}\le
C(b)|\ups_k-\ups^*|^{\frac{1}{p}};\\
\l\|\eta_{\ups_k}\r\|_{\ups_k,L^p(A_{h,2|\ups_k-\ups^*|^{1/2},k})}
\le C(b)|\ups_k-\ups^*|^{\frac{1}{p}}\quad\forall h\!\in\! H.
\end{gather*}
\end{crl}

\noindent
{\it Proof:} Let $\de_k=2|\ups_k-\ups^*|^{\frac{1}{2}}$ and 
$\ep_k=(2\|\be'\|_{C^0}+\tilde{C}(b))|\ups_k-\ups^*|^{\frac{1}{p}}$,
where $\tilde{C}$ is the function given by Lemma~\ref{cont_l2}.
Put
\begin{gather*}
\eta^{(0)}=-\bar{\partial}u_{\ups^*},\qquad
\eta^{(m+1)}=-\bar{\partial}u_{\ups^*}-N_{\ups^*}P_{\ups^*}\eta^{(m)}
\quad m\ge 0;\\
\eta_k^{(0)}=-\bar{\partial}u_{\ups_k},\qquad
\eta_k^{(m+1)}=-\bar{\partial}u_{\ups_k}-N_{\ups_k}P_{\ups_k}\eta_k^{(m)}
\quad m\ge 0.
\end{gather*}
By Lemma~\ref{cont_l2} and the explicit description of 
$\bar{\partial}q_{\ups_k}$ in Lemma~\ref{basic_gluing_lmm},
$\ep_k$, $\de_k$, $\eta^{(0)}$, and 
$\eta_k^{(0)}$ satisfy \e_ref{cont_l1_e1a} and~\e_ref{cont_l1_e1b}.
Suppose $\ep_k^{(m)}$ is such that 
$\ep_k^{(m)}$, $\de_k$, $\eta^{(m)}$, and $\eta_k^{(m)}$
satisfy \e_ref{cont_l1_e1a} and~\e_ref{cont_l1_e1b}.
Since the map $q_{\tilde{\ups}_k}$ is holomorphic on 
$q_{\tilde{\ups}_k}^{-1}(\Si_{\de_k}^*)$,
by \e_ref{cont_l1_e6a}, \e_ref{cont_l1_e6b}, 
the estimates in the proof of Lemma~\ref{fixed_pt}
and the derivation of equation~\e_ref{pert_e3} in~\cite{Z1},
\begin{gather*}
\label{cont_c3_e2a}\begin{split}
&\big\|q_{\ups_k}^{-1*}N_{\ups_k}P_{\ups_k}\eta_k^{(m)}-
N_{\ups^*}P_{\ups^*}\eta^{(m)}\big\|_{\ups^*,L^p(\Si_{\de_k}^*)}
=\big\|N_{\ups^*}q_{\ups_k}^{-1*}P_{\ups_k}\eta_k^{(m)}-
N_{\ups^*}P_{\ups^*}\eta^{(m)}\big\|_{\ups^*,L^p(\Si_{\de_k}^*)}\\
&\le C(b)\Big( 
\big\|q_{\ups_k}^{-1*}P_{\ups_k}\eta_k^{(m)}
              \big\|_{\ups^*,L^p(\Si_{\de_k}^*)}+
\big\|P_{\ups_k}\eta_k^{(m)}\big\|_{\ups^*,L^p(\Si_{\de_k}^*)}\Big)
\big\|q_{\ups_k}^{-1*}P_{\ups_k}\eta_k^{(m)}-
P_{\ups^*}\eta^{(m)}\big\|_{\ups^*,L^p(\Si_{\de_k}^*)}\\
&\le C'(b)\big(\ep_{k^*}^{(m)}+|\ups_k|^{\frac{1}{p}}\big)\ep_{k^*}^{(m)};
\end{split}\\
\label{cont_c3_e2b}\begin{split}
\big\|N_{\ups_k}P_{\ups_k}\eta_k^{(m)}\big\|_{\ups_k,L^p(A_{h,\de_k,k})}
\le C(b)|\ups_k|^{\frac{1}{p}}
 \big\|P_{\ups_k}\eta_k^{(m)}\big\|_{\ups_k,L_1^p(A_{h,\de_k,k})}
\le C'(b)|\ups_k|^{\frac{1}{p}}\ep_{k^*}^{(m)}.
\end{split}
\end{gather*}
Thus, we can take $\ep_k^{(m+1)}\!=\!\ep_k^{(m)}\!+\!
C'(b)\big(\ep_{k^*}^{(m)}\!+\!|\ups_k|^{\frac{1}{p}}\big)\ep_{k^*}^{(m)}$.
This sequence is bounded as long as $|\ups_k|^{\frac{1}{p}}$
is sufficiently small (depending only on $b$).
Since $\eta_{\ups^*}$ is the limit in the $(\ups^*,p)$-norm of the sequence
$\eta^{(m)}$ and 
$\eta_{\ups_k}$ is the limit in the $(\ups_k,p)$-norm of the sequence
$\eta_k^{(m)}$, the claim follows.

\begin{crl}
\label{cont_c4}
If ${\cal T}$ is a simple regular bubble type, 
there exist $\de,C\!\in\! C^{\i}({\cal M}_{\cal T};\Bbb{R}^+)$
such that for any sequence $\{\ups_k\!\in\! F^{\eset}{\cal T}_{\de}\}$
converging to $\ups^*\!\in\! F^H{\cal T}_{\de}$, 
$\tilde{\ga}(\ups_k)$ converges to $\tilde{\ga}(\ups^*)$
with respect to the Gromov topology.
Furthermore, 
\begin{gather*}
d_V\big(\ev(\tilde{\ga}(\ups^*)),\ev(\tilde{\ga}(\ups_k))\big)\le
C(b_{\ups^*})\l|\ups_k-\ups^*\r|^{\frac{1}{p}}
\qquad\hbox{if~}S=S^2;\\
d_V\big(\ev_l(\tilde{\ga}(\ups^*)),\ev_l(\tilde{\ga}(\ups_k))\big)\le
C(b_{\ups^*})\l|\ups_k-\ups^*\r|^{\frac{1}{p}}
\quad\forall l\!\in\! M;\\
\Big| \Psi_{\lr{\cal T},\hat{0}}\big(\tilde{\ga}(\ups_k)\big)-
\Psi_{{\cal T}(\ups^*),\hat{0}}\big(\tilde{\ga}(\ups^*)\big)\Big|
\le C(b_{\ups^*})\l|\ups_k-\ups^*\r|^{\frac{1}{p}}
\qquad\hbox{if~}S=S^2.
\end{gather*}
\end{crl}

\noindent
{\it Proof:} It is sufficient to consider the case 
$\pi_h(\ups_k)\!=\!\pi_h(\ups^*)$ for all $h\!\not\in\! H$
if $\ups^*\!\in\! F^{(H)}{\cal T}_{\de}$.
In such a case, $q_{\tilde{\ups}_k}$ maps the marked points of
$\Si_{\ups_k}$ to the marked points of $\Si_{\ups^*}$
and $u_{\ups_k}\!=\!u_{\ups^*}\!\circ\! q_{\ups_k}$.
By construction,
$$\tilde{u}_{\ups_k}=\exp_{b_{\ups^*},u_{\ups_k}}P_{\ups_k}\eta_{\ups_k},
\qquad
\tilde{u}_{\ups_k}=\exp_{b_{\ups^*},u_{\ups^*}}P_{\ups^*}\eta_{\ups^*}.$$
By Corollary~\ref{cont_c4} and the proof of Lemma~\ref{cont_l1},
\begin{gather}
\label{cont_c4_e2a}
\big\|q_{\tilde{\ups}_k}^{-1*}P_{\ups_k}\eta_{\ups_k}-P_{\ups^*}\eta_{\ups^*}
\big\|_{\ups^*,L^p_1(\Si_{2|\ups_k-\ups^*|^{1/2}}^*)}
\le C(b_{\ups^*})|\ups_k-\ups^*|^{\frac{1}{p}};\\
\label{cont_c4_e2b}
\l\|P_{\ups_k}\eta_{\ups_k}
\r\|_{\ups_k,L^p_1(A_{h,2|\ups_k-\ups^*|^{1/2},k})}
\le C(b_{\ups^*})|\ups_k-\ups^*|^{\frac{1}{p}}
\quad\forall h\!\in\! H.
\end{gather}
Let $\ze_k\!\in\!\Ga(\Si^*;\tilde{u}_{\ups^*})$ be given by
$$\exp_{b_{\ups^*},\tilde{u}_{\ups^*}}\ze_k=\tilde{u}_{\ups_k}
\circ q_{\tilde{\ups}_k},\qquad
\|\ze_k\|_{C^0}<\inj~g_{V,b_{\ups^*}}.$$
By equation~\e_ref{cont_c4_e2a} and the proof of 
(2) of Lemma~\ref{approx_maps},
\begin{equation}
\label{cont_c4_e3a}
\|\ze_k\|_{C^0(\Si_{2|\ups_k-\ups^*|^{1/2}}^*)}\le
C(b_{\ups^*})\|\ze_k\|_{\ups^*,L^p_1(\Si_{2|\ups_k-\ups^*|^{1/2}}^*)}
\le C'(b_{\ups^*})|\ups_k-\ups^*|^{\frac{1}{p}}.
\end{equation}
On the other hand, by~\e_ref{cont_c4_e2b} and 
by the same argument as in (3) of the proof of Lemma~\ref{inver1},
the variations of 
$P_{\ups^*}\eta_{\ups^*}$ on  $A_{h,2|\ups_k-\ups^*|^{1/2}}^*$ and 
$P_{\ups_k}\eta_{\ups_k}$ on $A_{h,2|\ups_k-\ups^*|^{1/2},k}^*$
are bounded $C(b_{\ups})|\ups_k-\ups^*|^{\frac{1}{p}}$.
This can be seen from equation~\e_ref{inver1_e7};
observe that an argument similar to the proof Lemma~\ref{cont_l2} shows 
that we can take $\de$ to be 
any small number bigger than $2|\ups_k-\ups^*|^{\frac{1}{2}}$.
Equation~\e_ref{cont_c4_e3a} and the small variation on the annuli
imply that 
$$\sup_{z\in\Si_{\ups_k}} 
  d_V\big(\tilde{u}_{\ups^*}(q_{\tilde{\ups}_k}(z)),\tilde{u}_{\ups_k}(z)\big)
\le C(b_{\ups^*})|\ups_k-\ups^*|^{\frac{1}{p}}.$$
It follows that $\tilde{\ga}_{\cal T}(\ups_k)$ converges
to $\tilde{\ga}_{\cal T}(\ups^*)$ in the Gromov topology.
The estimate on the evaluation maps is immediate from the above bound.
The last estimate follows from equations~\e_ref{cont_c4_e2a}
and~\e_ref{cont_c4_e2b}, along with a Sobolev estimate 
on a neighborhood of $\i\!\in\!\Si_{\ups^*,\hat{0}}$
which implies that the $C^1$-norm of $\ze_k$ there
is bounded by $C(b_{\ups^*})|\ups_k-\ups^*|^{\frac{1}{p}}$.

\subsection{Injectivity of the Gluing Map}
\label{local_inject}

\noindent
The goal of this subsection is to prove that the gluing maps of 
\e_ref{gluing_map2a} and \e_ref{gluing_map2b} are injective, 
as long as $\de\!\in\! C^{\i}({\cal M}_{\cal T}^{(0)};\Bbb{R}^+)$ is 
sufficiently small.
We start by showing local injectivity on the subspaces
$F^H{\cal T}_{\de}$ of $F{\cal T}_{\de}$,
where $H$ is a subset of $\hat{I}$.\\

\noindent
If ${\cal T}$ is regular, we are only interested in the case $t\!=\!0$.
If ${\cal T}$ is semiregular, we only consider the case $H\!=\!\eset$.
We use the same notation as in Subsection~\ref{scale_sec}.
If $\|\vp\|_{\ups}$ is sufficiently small, define
\hbox{$\tilde{\ze}_{\vp,t\nu}\!\in\!\Ga(\tilde{u}_{\ups,t\nu})$ by}
$$\exp_{b_{\ups},\tilde{u}_{\ups,t\nu}}\tilde{\ze}_{\vp,t\nu}=u_{\vp,t\nu},
\quad
\|\tilde{\ze}_{\vp,t\nu}\|_{b_{\ups},C^0}<\inj~g_{V,b}.$$

\begin{lmm}
\label{inj_l1b}
There exist
$\de,C\!\in\! C^{\i}({\cal M}_{\cal T}^{(0)};\Bbb{R}^+)$ such that
for all $\ups\!\in\! F^{(H)}{\cal T}_{\de}$,
where $H\!=\!\eset$ if ${\cal T}$ is semiregular,
and $\vp\!\in\!\tilde{T}_{\ups}F^H{\cal T}_{\de(b_{\ups})}$,\\
(1) $\big\|S_{\vp}N_{\vp,t\nu}R_{\vp}\xi\!-\!N_{\ups,t\nu}\xi
\big\|_{\ups,p}\!\le\! C(b_{\ups})\|\vp\|_{\ups}\|\xi\|_{\ups,p,1}^2$
for all $\xi\!\in\!\Ga(u_{\ups})$ with $\|\xi\|_{\ups,p,1}\le\de(b_{\ups})$
\hbox{and $t\!\in\![0,1]$};\\
(2) $\big\|S_{\vp}\tilde{\pi}_{\vp,\pm}R_{\vp}\xi-\tilde{\pi}_{\ups,\pm}\xi
\big\|_{\ups,p,1}\le C(b_{\ups})\|\vp\|_{\ups}\|\xi\|_{\ups,p,1}$
for all $\xi\!\in\!\Ga(u_{\ups})$;\\
(3) $\big\|S_{\vp}P_{\vp}R_{\vp}\eta-P_{\ups}\eta\big\|_{\ups,p,1}
\le C(b_{\ups})\|\vp\|_{\ups}\|\eta\|_{\ups,p}$
for all $\xi\!\in\!\Ga^{0,1}(u_{\ups})$.
\end{lmm}

\noindent
{\it Proof:} 
Claim (1) follows from (2) of Lemma~\ref{inj_l1} and
Riemannian geometry estimates such as in~\cite{Z1}.
Claim (2) is a consequence of (5) of Lemma~\ref{inj_l1}
and (2) of Definition~\ref{obs_setup_dfn1}.
Finally, (3) is obtained from (1), (2), 
Definitions~\ref{obs_setup_dfn1} and~\ref{obs_setup_dfn2}, and
Lemmas~\ref{inj_l1} and ~\ref{approx_maps2} as follows.
Writing $\lap_{\vp}P$ for $S_{\vp}P_{\vp}R_{\vp}-P_{\ups}$, etc.,
\begin{equation*}\begin{split}
\lap_{\vp}P
&=P_{\ups}\pi_{\ups,+}^{0,1}D_{\ups}\tilde{\pi}_{\ups,+}\!\lap_{\vp}\!P
\!+\!\tilde{\pi}_{\ups,-}S_{\vp}P_{\vp}R_{\vp}
=P_{\ups}\pi_{\ups,+}^{0,1}D_{\ups}\!\lap_{\vp}\!P
\!-\!\big(P_{\ups}\pi_{\ups,+}^{0,1}D_{\ups}\!-\!1\big)
\!\lap_{\vp}\!\tilde{\pi}_{\cdot,+}S_{\vp}P_{\vp}R_{\vp}\\
&=P_{\ups}\pi_{\ups,+}^{0,1}\!\lap_{\vp}\!\tilde{\pi}_{\cdot,+}^{0,1}
-\Big(P_{\ups}\pi_{\ups,+}^{0,1}\!\lap_{\vp}\!D+
\big(P_{\ups}\pi_{\ups,+}^{0,1}D_{\ups}-1\big)
\!\lap_{\vp}\!\tilde{\pi}_{\cdot,+}\Big)S_{\vp}P_{\vp}R_{\vp}\\
&=-P_{\ups}\!\lap_{\vp}\!\pi_{\cdot,+}^{0,1}
S_{\vp}\tilde{\pi}_{\vp,+}^{0,1}R_{\vp}
-\Big(P_{\ups}\pi_{\ups,+}^{0,1}\!\lap_{\vp}\!D+
\big(P_{\ups}\pi_{\ups,+}^{0,1}D_{\ups}-1\big)
\!\lap_{\vp}\!\tilde{\pi}_{\vp,+}S_{\vp}P_{\vp}R_{\vp}\Big).
\end{split}\end{equation*}

\begin{crl}
\label{inj_c1}
There exist 
$\de,C\!\in\! C^{\i}({\cal M}_{\cal T}^{(0)};\Bbb{R}^+)$
such that for all $t\!\in\![0;\de(b_{\ups})]$,
$\ups\!\in\! F^{(H)}{\cal T}_{\de}$,
where $H=\eset$ if ${\cal T}$ is semiregular, 
and \hbox{$\vp\!\in\! \tilde{T}_{\ups}F^H{\cal T}_{\de(b_{\ups})}$}, 
$$C(b_{\ups})^{-1}\|\vp\|_{\ups}\le
   \|\tilde{\ze}_{\vp,t\nu}\|_{\ups,p,1}+\sum_{h\in H}|w_h(\vp)|_{g_{\ups}}
      +\sum_{l\in M}|w_l(\vp)|_{g_{\ups}}
\le C(b_{\ups})\|\vp\|_{\ups}.$$
Furthermore, 
$\big\|S_{\vp}\xi_{\vp,t\nu}-\xi_{\ups,t\nu}\big\|_{\ups,p,1}
\le C(b_{\ups})\big(t+|\ups|^{\frac{1}{p}}\big)\|\vp\|_{\ups}$.
\end{crl}

\noindent
{\it Proof:} 
The first claim of the lemma is immediate from the second 
and (1) of Lemma~\ref{inj_l1}.
On the other hand, by construction in Subsection~\ref{gluing_map},
$$\xi_{\vp,t\nu}
=tP_{\vp}\nu-P_{\vp}\bar{\partial}u_{\vp}
-P_{\vp}N_{\vp,t\nu}\xi_{\vp,t\nu}.$$
Thus, if $t$ and $|\ups|$ are sufficiently small (depending on $b_{\ups}$), 
the second claim follows from
Lemmas~\ref{inj_l1},~\ref{inj_l1b}, Corollary~\ref{fixed_pt_crl},
and equation~\e_ref{pert_e3}.

\begin{crl}
\label{inj_main}
If ${\cal T}$ is a simple bubble type and
$K$ is an open subset of ${\cal M}_{\cal T}$
with compact closure,
there exists $\de>0$ such that for any $t\!\in\![0,\de]$,
the map
$$\tilde{\ga}_{{\cal T},t\nu}\!:
F^{\eset}{\cal T}_{\de}|K\lra  C^{\i}_{(\la^*;L)}(S;V),\qquad
\tilde{\ga}_{{\cal T},t\nu}(\ups)=\tilde{b}_{t\nu}(\ups),$$
is a differentiable embedding.
\end{crl}

\noindent
{\it Proof:} We first deduce from Corollary~\ref{inj_c1} that 
$\tilde{\ga}_{{\cal T},t\nu}$ is injective if $\de$ is sufficiently small.
Suppose not, i.e. there exist sequences 
$\ups_k,\ups_{k'}\!\in\! F^{\eset}{\cal T}_{\de}|K$ such that
$$\ups_k\!\lra\! b\!\in\!\bar{K},\quad
\ups_k'\lra b'\!\in\!\bar{K},~~\hbox{and}~~
\tilde{b}_{t\nu}(\ups_k)\!=\!\tilde{b}_{t\nu}(\ups_k').$$
It follows that $b\!=\!b'$,
after possibly modifying the sequence $\{\ups_{k'}\}$
by actions of~$G_{\cal T}$.
If for some $k$, $\ups_k'\!=\!\ups_k(\vp_k)$ with $\|\vp_k\|_{\ups_k}$
sufficiently small, then by Corollary~\ref{inj_c1}, $\ups_k'=\ups_k$.
Otherwise, the difference between $q_{\ups_k}$ and $q_{\ups_k'}$
is uniformly bounded below outside of the preimage of the zeroth
component and the necks $A_{\ups_k,h}$.
Thus, the bubble maps $b(\ups_k)$ and $b(\ups_k')$ are far apart unless
$b$ has an automorphism.
In the latter case, $\ups_k'$ can be replaced by an equivalent element
of~$F^{\eset}{\cal T}_{\de}$.
In the former case,  $\tilde{u}_{\ups_k}$ and $\tilde{u}_{\ups_k'}$
cannot be the same because 
$$\big\|P_{\ups_k}\eta_{\ups_k,t\nu}\big\|_{C^0}\le 
 C\big(t+|\ups_k|^{\frac{1}{p}}\big)\le C'\de^{\frac{1}{p}}
\quad\hbox{and}\quad
\big\|P_{\ups_k'}\eta_{\ups_k',t\nu}\big\|_{C^0}\le 
 C\big(t+|\ups_k'|^{\frac{1}{p}}\big)\de^{\frac{1}{p}}.$$
Thus $\tilde{\ga}_{{\cal T},t\nu}$ is injective on
$F^{\eset}{\cal T}_{\de}|K$ provided $\de$ is sufficiently small
(depending on $K$).
The smoothness of $\tilde{\ga}_{{\cal T},t\nu}$ follows from
the smooth dependence of solutions of equation~\e_ref{pert_e4a}
on the parameters.
Finally, the differential of 
$\tilde{\ga}_{{\cal T},t\nu}$ is nondegenerate by Corollary~\ref{inj_c1}.

\begin{crl}
\label{inj_main2}
If ${\cal T}$ is regular, there exists 
$\de\!\in\! C^{\i}\l({\cal M}_{\cal T};\Bbb{R}^+\r)$ 
such that for all $m$, the map 
$$\tilde{\ga}_{\cal T}: 
\bigcup_{|H|=m}\!\! F^H{\cal T}_{\de}\lra
\bigcup_{|H|=m}\!\!{\cal M}_{{\cal T}(H)},
\qquad
\tilde{\ga}_{\cal T}(\ups)=\tilde{b}(\ups),$$
is injective. 
\end{crl}

\noindent
{\it Proof:} The same argument as in the proof of Corollary~\ref{inj_main}
shows that map
$$\tilde{\ga}_{\cal T}\!: F^H{\cal T}_{\de}\lra 
 {\cal H}_{{\cal T}(H)}$$
is an embedding if $\de$ is sufficiently small.
It remains to see that 	
$\tilde{\ga}_{\cal T}^{(0)}(\ups)\!\neq\! g\!\cdot\!
              \tilde{\ga}_{\cal T}^{(0)}(\ups')$
for any \hbox{$g\!\in\! G_{{\cal T}(H)}$} whenever $[\ups]\!\neq\![\ups']$.
For each $\ups\!\in\! F^{(H)}{\cal T}_{\de}$  and $i\!\in\! H$, 
we construct \hbox{$\l(c_i(\ups),r_i(\ups)\r)\!\in\!\Bbb{C}\!\times\!\Bbb{R}$}
such~that 
$$\l(c(\ups),r(\ups)\r)\cdot\tilde{\ga}_{\cal T}^{(0)}(\ups)\in 
{\cal M}_{{\cal T}(H)}^{(0)}.$$
We define $c_i(\ups)\!\in\!\Bbb{C}$ and $r_i(\ups)\!\in\!\Bbb{R}$ by
\begin{gather*}
\tilde{\Psi}\big((c_i(\ups),0)\cdot\tilde{u}_{\ups,i}\big)+
\sum_{\io_h(\ups)=i}d_h\big({\cal T}(H)\big)\big(x_h(\ups)+c_i(\ups)\big)
+\sum_{j_l(\ups)=i}\big(y_l(\ups)+c_i(\ups)\big)=0;\\
\begin{split}
\Psi^{(3)}\big((c_i(\ups),r_i(\ups))\cdot\tilde{u}_{\ups,i}\big)+
\sum_{\io_h(\ups)=i} d_h\big({\cal T}(H)\big)
 \be\big((1+ r_i(\ups))|x_h(\ups)+c_i(\ups)|\big)&\\
+\sum_{j_l(\ups)=i}\be\big((1+r_i(\ups))|y_l(\ups)+ c_i(\ups)|\big)
&=\frac{1}{2}.
\end{split}\end{gather*}
Since the metric  $g_{\ups,i}$ for $i>0$ agrees with the standard
metric on $S^2$ on a neighborhood of the south pole
and $\Psi_{{\cal T},i}(b_{\ups})=0$, by 
Corollary~\ref{inj_c1} for any $\vp\!\in\!T_{\ups}F^H{\cal T}$
with  $\|\vp\|_{\ups}$ sufficiently small,
\begin{equation}
\label{inj_main2_e1}
\big|c_i(\vp)-c_i(\ups)\big|
\le C(b_{\ups})|\ups|^{\frac{1}{p}}\|\vp\|_{\ups},\qquad
\big|r_i(\vp)-r_i(\ups)\big|
\le C(b_{\ups})|\ups|^{\frac{1}{p}}\|\vp\|_{\ups}.
\end{equation}
Let $\bar{b}(\ups)=(c(\ups),r(\ups))\cdot\tilde{b}(\ups)$.
Write
$$\bar{b}(\ups)
=\big(S,M,H\cup\{\hat{0}\};\bar{x}(\ups),(j(\ups),\bar{y}),\bar{u}_{\ups}
\big).$$
If $\|\vp\|_{\ups}$ is sufficiently small, define
$\bar{\ze}_{\vp}\!\in\!\Ga(\bar{u}_{\ups})$ by
$$\exp_{b_{\ups},\bar{u}_{\ups}}\bar{\ze}_{\vp}=u_{\vp},\quad
\|\bar{\ze}_{\vp}\|_{b_{\ups},C^0}<\inj~g_{V,b_{\ups}}.$$
Similarly, for $h\!\in\! H$ and $l\!\in\! M$, define 
$\bar{w}_h(\vp)\!\in\! T_{\bar{x}_h(\ups)}\Si_{\ups,\io_h(\ups)}$
and 
$\bar{w}_l(\vp)\!\in\! T_{\bar{y}_l}\Si_{\ups,j_l(\ups)}$ by
\begin{gather*}
\exp_{g_{\ups},\bar{x}_h(\ups)}\bar{w}_h(\vp)=\bar{x}_h(\vp),\quad
|\bar{w}_h(\vp)|\equiv|\bar{w}_h(\vp)|_{g_{\ups}}
<\hbox{inj}(g_{\ups},x_h(\ups));\\
\exp_{g_{\ups},\bar{y}_l(\ups)}\bar{w}_l(\vp)=\bar{y}_l(\vp),\quad
|\bar{w}_l(\vp)|\equiv|\bar{w}_l(\vp)|_{g_{\ups}}
<\hbox{inj}(g_{\ups},y_l(\ups)).
\end{gather*}
Then by equation~\e_ref{inj_main2_e1} and Corollary~\ref{inj_c1},
\begin{equation}
\label{inj_est2}
C''(b_{\ups})^{-1}\|\vp\|_{\ups}
\le\|\bar{\ze}_{\vp}\|_{\ups,p,1}+\sum_{h\in H}|\bar{w}_h(\vp)|+
\sum_{l\in M}|\bar{w}_l(\vp)|
\le C(b_{\ups})\|\vp\|_{\ups}.
\end{equation}
It follows that the map 
$$F^H{\cal T}_{\de}\lra{\cal M}_{{\cal T}(H)}^{(0)},\qquad
\ups\lra \bar{b}(\ups),$$
is a local embedding.
By the same argument as in the proof of Lemma~\ref{inj_main},
we can conclude that this map is injective as long as 
$\de\!\in\! C^{\i}({\cal M}_{\cal T}^{(0)};\Bbb{R}^+)$
is sufficiently small.
Since this map is $G_{{\cal T}(H)}$-equivarent by construction,
it follows that the induced map on the quotient, i.e.
the map of Corollary~\ref{inj_main2}, is injective.

\begin{crl}
\label{inj_main2c}
If $S=S^2$, there exists 
$\de\!\in\! C^{\i}\l({\cal M}_{\cal T};\Bbb{R}^+\r)$ 
such that the map 
$$\tilde{\ga}_{\cal T}\!: 
F{\cal T}_{\de}\big| {\cal M}_{\cal T}\lra\bar{\cal M}_{\lr{\cal T}}$$
is injective. Furthermore, the restriction
$$\tilde{\ga}_{\cal T}\!: F^{\eset}{\cal T}_{\de}\big| 
              {\cal M}_{\cal T}\lra{\cal M}_{\lr{\cal T}}$$
is a differentiable embedding.
\end{crl}

\noindent
In order to adjust the gluing procedure in the presence
of constraints, below we state the analogue of 
Corollary~\ref{inj_c1} for 
\hbox{$\vp\!\in\!{\cal K}_{b_{\ups}}{\cal T}
           \!\subset\! T_{\ups}F^{\eset}{\cal T}$}.
It is obtained in the same way as Corollary~\ref{inj_c1},
except the analogue of Lemma~\ref{inj_l1b} 
would make use of Lemma~\ref{inj_l1c}, instead of  Lemma~\ref{inj_l1},
and of (2b), instead of (2a), of Definitions~\ref{obs_setup_dfn1}
and~\ref{obs_setup_dfn2}.
We also use (3) of Lemma~\ref{approx_maps}.

\begin{crl}
\label{inj_c1c}
There exist 
$\de,C\!\in\! C^{\i}({\cal M}_{\cal T}^{(0)};\Bbb{R}^+)$
such that for all $t\!\in\![0;\de(b_{\ups})]$,
$\ups\!\in\! F^{(\eset)}{\cal T}_{\de}$, 
and \hbox{$\vp\!\in\! {\cal K}_{b_{\ups}}{\cal T}_{\de(b_{\ups})}$}, 
$$C(b_{\ups})^{-1}\|\vp\|\le
   \|\tilde{\ze}_{\vp,t\nu}\|_{\ups,p,1}+
\sum_{l\in M}|w_l(\vp)|_{g_{\ups}}
\le C(b_{\ups})\|\vp\|.$$
Furthermore, 
$\big\|S_{\vp}'\xi_{\vp,t\nu}-\xi_{\ups,t\nu}\big\|_{\ups,C^0}
\le C(b_{\ups})\big(t+|\ups|^{\frac{1}{p}}\big)\|\vp\|$.
\end{crl}

\subsection{The Basic Gluing Map and the Space of Balanced Maps}
\label{surject_b}

\noindent
Our next goal is to show that the gluing map of Subsection~\ref{gluing_map}
is surjective in the appropriate sense.
More precisely, if ${\cal T}$ is a regular bubble type,
we show that the image of 
$\tilde{\ga}_{\cal T}$ contains a neighborhood of 
${\cal M}_{\cal T}$ in~$\bar{\cal M}_{\lr{\cal T}}$. 
If ${\cal T}$ is a semiregular, 
we show that all elements in ${\cal M}_{\Si,t\nu,\la^*}$
that are close to any given compact subset of ${\cal M}_{\cal T}$ are 
in the image of the gluing map $\tilde{\ga}_{{\cal T},t\nu}$
if $t$ is sufficiently small.
The major difficulty in doing this is the following.
If $\ups\!\in\! F{\cal T}$, 
a small change in the singular points of $b_{\ups}$
may lead to a very large change in the map~$u_{\ups}$.
This is precisely the reason we used the norm $\|\vp\|_{\ups}$
on $T_{\ups}F^H{\cal T}$ instead of just $\|\vp\|$ 
in Subsection~\ref{scale_sec}.
In order to deal with this issue, we need Corollary~\ref{surjb_c2},
which is proved in this subsection.
We continue to assume that ${\cal T}$ is a simple bubble~type.
\\

\noindent 
Recall that ${\cal H}_{\cal T}$ is the set of tuples 
$b\!=\!\big(S,M,I;x,(j,y),u\big)$ such that 
$u_{\io_h}(x_h)\!=\!u_h(\i)$ for all $h\!\in\!\hat{I}$ and 
$\bar{\partial}u_i\!=\!0$ for \hbox{all $i\!\in\! I$}.
Furthermore, ${\cal M}_{\cal T}^{(0)}$
is the subset of ${\cal H}_{\cal T}$ consisting of 
the tuples~$b$ such that $\Psi_{{\cal T},h}(b)=\!0\!$ for 
\hbox{all $h\!\in\!\hat{I}$}.
It is convenient to make the following definitions.
If $H$ is a subset of $\hat{I}$ and $\ep\!\ge\! 0$, let
\begin{equation*}\begin{split}
{\cal M}_{{\cal T},\ep}^{(H)}\!=\!\big\{
b\!=\!\big(S,M,I;x,(j,y),u\big)\!:
\bar{\partial}u_i\!=\!0~\forall i\!\in\! I;~
d_V\big(u_{\io_h}\!(x_h),u_h(\i)\big)\!\le\!\ep~\forall h\!\in\!\hat{I};
\quad&\\
\big|\Psi_{{\cal T},h}(b)\big|\!\le\!\ep~\forall h\!\in\!\hat{I}\!-\! H
\big\}.&
\end{split}\end{equation*}

\begin{lmm}
\label{surjb_l1}
There exist $\de,C\!\in\! C^{\i}\big({\cal M}_{\cal T}^{(0)};\Bbb{R}^+\big)$
with the following property.
Suppose \hbox{$b^*\!\in\!{\cal M}_{\cal T}^{(0)}$}, $\ep<\de(b^*)$,
$b\!\in\!{\cal M}_{{\cal T},\ep}^{(H)}$ is such that $d(b^*,b)\le\de(b^*)$, 
and $\ups\!=\!(b,v_{\hat{I}})\!\in\! F^{(H)}_b{\cal T}_{\de(b^*)}$.
Then there exist 
$$\tilde{b}\in{\cal M}_{{\cal T},\ep^2}^{(H)}
\quad\hbox{and}\quad
\tilde{\ups}=\big(\tilde{b},\tilde{v}_{\hat{I}}\big)\in 
F^{(H)}_{\tilde{b}}{\cal T}$$
such that\\
(1) $d(b,\tilde{b})\le C (b^*)\ep$ and 
$\l|\tilde{v}_h-v_h\r|_b\le C (b^*)\ep|v_h|_b$ for all $h\!\in\!\hat{I}$;\\
(2) if $q_{\ups}(z)\!\in\!\Si_{{\cal T},i}$, 
$\l|r_{b,h}\l(q_{\ups}(z)\r)\r|\!\ge\! 2|v_h|^{\frac{1}{2}}$
for all $h\!\in\!\hat{I}\!-\! H$ such that $\io_h\!=\! i$ and 
\hbox{$\l|q_S^{-1}\l(q_{\ups}(z)\r)\r|\!\ge\! 2|v_i|^{\frac{1}{2}}$} if 
$i\in\hat{I}-H$, then
$d_b\l(q_{\ups}(z),q_{\tilde{\ups}}(z)\r)\le\de(b^*)\ep.$
\end{lmm}

\noindent
{\it Proof:} (1) Let $b\!=\!\big(S,M,I;x,(j,y),u\big)$.
If $\de$ is sufficiently small, by Proposition~\ref{balanced_case},
we can choose \hbox{$\xi_i\!\in\!\Ga(u_i)$}
such that $\|\xi_i\|_{g_{b,i},C^1}\le C(b^*)\ep$ and
$$b'\equiv\big(S,M,I;x,(j,y),u'\big)\in{\cal H}_{\cal T},$$
where $u'_i\!=\!\exp_{u_i}\xi_i$.
The $C^1$-bound on $\xi_i$ and the assumption 
$b\!\in\!{\cal M}_{{\cal T},\ep}^{(H)}$ imply that
\hbox{$\l|\Psi_{{\cal T},h}(b')\r|\!\le\! C'(b^*)\ep$} for all 
$h\!\in\!\hat{I}\!-\!H$.\\
(2) We now define 
$\tilde{b}'\!\equiv\!\big(S,M,I;\tilde{x}',(j,\tilde{y}'),
 \tilde{u}'\big)\!\in\!{\cal H}_{\cal T}$
and $\tilde{\ups}'\!=\!(\tilde{b}',\tilde{v}_{\hat{I}}')$ as follows.
Suppose \hbox{$i^*\!\in\! I$} and for all $i\!\in\!\hat{I}$ with $i\!>\!i^*$,
$h\!\in\!\hat{I}$ with $\io_h\!=\!i$, and $l\!\in\! M$ with $j_l\!=\!i$,
we have constructed \\
(i) $(c_i,r_i)\!\in\!\Bbb{C}\!\times\!\Bbb{R}$ such
that $|(c_i,r_i)|\le C(b^*)\ep$;\\
(ii) $\tilde{x}_h',\tilde{y}_l'\!\in\!\Si_{{\cal T},i}$
such that $|r_{b,h}(\tilde{x}_h')|\le C(b^*)\ep$
and $|\phi_{y_l}\tilde{y}_l'|\le C(b^*)\ep$;\\
(iii) $\tilde{v}_h'\!\in\!\Bbb{C}$ such that 
$\l|\tilde{v}_h'-v_h\r|\le C(b^*)\ep|v_h|$;\\
(iv) if $x_i\!\in\! S^2$, $\bar{x}_i\!\in\! S^2$, such that
$|r_{b,i}(\bar{x}_i)|\le C(b^*)\ep|v_i|$;\\
(v) if $x_i\!\in\! S^2$, $\bar{v}_i\!\in\!\Bbb{C}$
such that $\l|\bar{v}_i-v_i\r|\le C(b^*)\ep|v_i|$,\\
such that\\
(I1) if $i\!\not\in\! H$, $\Psi_{{\cal T},i}(\tilde{b}')=0$
where $\tilde{u}_i'=(c_i,r_i)\cdot u_i'$;\\
(I2) if $\Si_{{\cal T},\io_i}=S^2$, $z\!\in\!\Si_{{\cal T},\io_i}$,  and 
$\big|\phi_{\bar{x}_h}q_{i,(x_i,v_i)}(z)\big|
                  \le \frac{2}{3}|\bar{v}_h|^{\frac{1}{2}}$
for some $h\!\in\!\hat{I}\!-\!H$, then
\begin{gather*}
\l|\phi_{\tilde{x}_h'}q_{i,(\bar{x}_i,\bar{v}_i)}(z)\r|
\le |\tilde{v}_h'|^{\frac{1}{2}} \hbox{~~and~~}
q_{h,(\tilde{x}_h',\tilde{v}_h')}
 \l(q_{i,(\bar{x}_i,\bar{v}_i)}(z)\r)
=q_{h,(\bar{x}_h,\bar{v}_h)}\l(q_{i,(x_i,v_i)}(z)\r),
\end{gather*}
where $q_{i,(x_i,v_i)}$, etc., are the maps defined in 
Section~\ref{basic_gluing}.\\
Note that while we have not defined $\tilde{b}'$ completely yet,
(I1) is still a well-defined statement.
The function $\Psi_{{\cal T},i}$ depends only the $i$th
bubble component of $\tilde{b}'$, 
which has already been constructed by the induction assumptions.\\

\noindent
If $i^*\!\in\! H$, we take $c_{i^*}\!=\!0$ and $r_{i^*}\!=\!0$.
If $i^*\!\in\!\hat{I}\!-\!H$,
let $(c_{i^*},r_{i^*})\!\in\!\Bbb{C}\!\times\!\Bbb{R}$ be given by
\begin{gather*}
\tilde{\Psi}\l((c_{i^*},0)u_{i^*}'\r)
+\sum_{\io_h=i^*}d_h({\cal T})\l(\bar{x}_h+c_{i^*}\r)
+\sum_{j_l=i^*}\l(y_l+c_{i^*}\r)=0;\\
\Psi^{(3)}\l((c_{i^*},r_{i^*})u_{i^*}'\r)
+\sum_{\io_h=i^*}d_h({\cal T})\be
            \big((1+r_{i^*})|\bar{x}_h+c_{i^*}|\big)
+\sum_{j_l=i^*}\be\big((1+r_{i^*})|y_l+c_{i^*}|\big)
   =\frac{1}{2}.
\end{gather*}
If $\ep$ is sufficiently small, by the proof of Lemma~\ref{balanced_case}
such $(c_{i^*},r_{i^*})\!\in\!\Bbb{C}\!\times\!\Bbb{R}$ exists and satisfies
\hbox{$|c_{i^*}|,|r_{i^*}|\le C(b^*)\ep$}. 
For all $h\!\in\!\hat{I}$ with $\io_h=i^*$ and $l\!\in\! M$ with $j_l=i^*$, 
put
\begin{gather*}
\tilde{x}_h'=(1+r_{i^*})\l(\bar{x}_h+c_{i^*}\r),\quad
\tilde{v}_h'=(1+r_{i^*})\bar{v}_h,\quad
\tilde{y}_l'=(1+r_{i^*})\l(y_l+c_{i^*}\r);\\
\bar{x}_{i^*}=x_{i^*}-c_{i^*}v_i,\quad
\bar{v}_{i^*}=(1+r_{i^*})^{-1}v_i  \quad\hbox{if}\quad x_{i^*}\!\in\! S^2.
\end{gather*}
Continuing in this way, for all $i\!\in\!\hat{I}$, 
$h\!\in\!\hat{I}$ with $\io_h=i$, and $l\!\in\! M$ with $j_l=i$,
we obtain elements (i)-(v) satisfying (I1),(I2).
Let $\tilde{u}_{\hat{0}}'=u_{\hat{0}}'$.
If $l\!\in\! M$ and $j_l={\hat{0}}$, take $\tilde{y}_l'=y_l$.\\
(3) If $S\!=\!S^2$, let
$(\tilde{x}_h',\tilde{v}_h')\!=\!(\bar{x}_h,\bar{v}_h)$ if $\io_h\!=\!0$, 
$\tilde{b}\!=\!\tilde{b}'$, and $\tilde{\ups}\!=\!\tilde{\ups}'$.
By the inductive construction, $\tilde{b}$ and $\tilde{\ups}$
satisfy the requirements of the lemma.
In fact, $\tilde{b}\!\in\!{\cal M}_{{\cal T},0}^{(H)}$.
If $S\!=\!\Si$, we could extend the above construction to the principal
component $\Si$ as we did for $S\!=\!S^2$ if
$q_{\tilde{\ups}'}$  were defined using the metric 
$g_{b,\hat{0}}$ on $\Si$, which may differ slightly from 
$g_{\tilde{b}',\hat{0}}$.
This problem is fixed in below.\\
(4) If $l\!\in\! M$ and $j_l\!=\!\hat{0}$, we take 
$\tilde{y}_l\!=\!y_l$ as before.
For all $h\!\in\!\hat{I}$ with $\io_h\!=\!\hat{0}$, 
let $\tilde{x}_h\!\in\!\Si$, $\tilde{v}_h\in T_{\tilde{x}_h}\Si$, and 
$\Th_h\!: B_{2|v_h|_b^{-\frac{1}{2}}}(0;\Bbb{C})\!\lra\!\Bbb{C}$
be such~that\\
($\Si_{\hat{0}}$1) $d_b\l(x_h,\tilde{x}_h\r)\le C(b^*)\ep |v_h|$,
$\l||\tilde{v}_h|_{\tilde{b}}-|v_h|_b\r|\le C(b^*)\ep|v_h|_b$;\\
($\Si_{\hat{0}}$2) for all
$z\!\in\! B_b\big(x_h,2|v_h|_b^{\frac{1}{2}}\big)$,
$$\frac{\phi_{\tilde{b},h}z}{\tilde{v}_h}= \big(1+r_h\big)
 \Big\{c_h+\frac{\phi_{b,h}z}{v_h}+
  \Th_h\Big(\frac{\phi_{b,h}z}{v_h}\Big)\Big\};$$
($\Si_{\hat{0}}$3) $\Th_h$ is holomorphic, $\Th_h(0)=0$, $\Th_h'(0)=0$, and 
$\l\|\Th_h''\r\|_{C^0}\le C(b^*)|\ups|^2\ep$.\\
Note that even though we have not defined $\tilde{b}$ completely yet,
($\Si_{\hat{0}}$1) and ($\Si_{\hat{0}}$2) are still well-defined statements,
since the metric $g_{\tilde{b},\hat{0}}$ on $\Si$ depends only
on the singular points $\{\tilde{x}_h\!: \io_h\!=\!\hat{0}\}$ on~$\Si$.
Existence of such $\tilde{x}_h$, $\tilde{v}_h$, and 
$\Th_h$ follows from Corollary~\ref{si_metrics_c4},
provided $\de$ is sufficiently~small.\\

\noindent
If $\io_i\!=\!\hat{0}$ and $j_l\!=\!i$, let 
$(i,\tilde{y}_l)= q_{\tilde{\ups},i}q_{\ups,i}^{-1}(i,y_l).$
The map $q_{\tilde{\ups},i}$ is well-defined even though 
$\tilde{\ups}$ has not been defined completely yet.
By ($\Si_{\hat{0}}$2),
\begin{equation}\label{surj_si10b}
\tilde{y}_l
=q_{\tilde{\ups},i}q_{\ups,i}^{-1}(y_l)
=\frac{\phi_{\tilde{b},\tilde{x}_i}\phi_{b,x_i}^{-1}(y_lv_i)}{\tilde{v}_i}
=\big(1+r_i\big) \big\{c_i+y_l+\Th_i(y_l)\big\}.
\end{equation}
Since $\tilde{y}_l'=\l(1+r_i\r)\l(y_l+c_i\r)$, 
$\l|\tilde{y}_l-\tilde{y}_l'\r|\le C(b^*)|\ups|^2\ep$ by~($\Si_{\hat{0}}$3).\\

\noindent
Suppose $h\!\in\!\hat{I}$, $\io_h\!\in\!\hat{I}$, and 
for every $i\!\in\!\hat{I}$ with $i<h$
and $j\!\in\! M$ with $j_l=i$, we have defined\\
\begin{gather*}
\tilde{x}_i\in\Si_{{\cal T},\io_i},\quad
\tilde{y}_l\in\Si_{{\cal T},i},\quad
\tilde{v}_i\in 
\begin{cases}
T_{\tilde{x}_i}\Si,&\hbox{if~}\io_i=0;\\
\Bbb{C},&\hbox{if~}\io_i>0;
\end{cases}\quad
\tilde{c}_i\in\Bbb{C},\quad
\Th_i: B_{2|v_i|_b^{-\frac{1}{2}}}(0;\Bbb{C}) \lra\Bbb{C}
\end{gather*}
such that\\
($\Si$1) 
$|{\phi}_{\tilde{b}',i}\tilde{x}_i|_{\tilde{b}'}\le C(b^*)|\ups|^2\ep$
if $\io_i>0$ and 
$|\phi_{\tilde{b}',\tilde{y}_l'}\tilde{y}_l|_{\tilde{b}'}\
\le C(b^*)|\ups|^2\ep$;\\
($\Si$2) $\l||\tilde{v}_i|_{\tilde{b}}-|v_i|_b\r|\le C(b^*)\ep|v_i|_b$;\\
($\Si$3) $\l|\tilde{c}_i-c_i\r|\le C(b^*)|\ups|^2\ep$;\\
($\Si$4) for all $z\!\in\!\Si$ such that 
$r_{b,i}q_{\ups,\io_i}(z)\le 2|v_i|_b^{\frac{1}{2}}$,
$$\frac{\phi_{\tilde{b},i}
q_{\tilde{\ups},\io_i}z}{\tilde{v}_i}
=\big(1+r_i\big)\Big\{\tilde{c}_i
+\frac{\phi_{b,i}q_{\ups,\io_i}z}{v_i}
+\Th_i\Big(\frac{\phi_{b,i}q_{\ups,\io_i}z}{v_i}\Big)\Big\}$$
($\Si$5) $\Th_i$ is holomorphic, $\Th_i(0)=0$, $\Th_i'(0)=0$, and 
$\l\|\Th_i''\r\|_{C^0}\le C(b^*)|\ups|^2\ep$.\\
If $h\!\in\! H$, we take $\tilde{x}_h=\tilde{x}_h'$, 
$\tilde{v}_h=\tilde{v}_h'=0$, $\tilde{y}_l=\tilde{y}_l'$ if $j_l=h$,
$\tilde{c}_h=c_h=0$, and $\Phi_h(z)=0$.
If $h\!\not\in\! H$, let
$$(\io_h,\tilde{x}_h)=q_{\tilde{\ups},\io_h}
q_{\ups,\io_h}^{-1}(h,\bar{x}_h).$$
By an argument similar to \e_ref{surj_si10b}, from ($\Si$4) we obtain 
\begin{equation}\label{surj_si12b}
\tilde{x}_h=\big(1+r_{\io_h}\big) \big\{ \tilde{c}_{\io_h}+
 \bar{x}_h+\Th_{\io_h} ( \bar{x}_h)\big\}.
\end{equation}
Since $\tilde{x}_h'=(1+r_{\io_h})(\bar{x}_h+c_{\io_h})$,
\e_ref{surj_si12b}, ($\Si$3), and ($\Si$5) imply the first part of 
($\Si$1) with $i=h$.
Furthermore, by assumption ($\Si$4), 
\begin{equation}\begin{split}\label{surj_si12c}
\phi_{\tilde{b},h}q_{\tilde{\ups},\io_h}(z)
&=q_{\tilde{\ups},\io_h}(z)-\tilde{x}_h
=\frac{\phi_{\tilde{b},\io_h}q_{\tilde{\ups},\io_{\io_h}}(z)}
{\tilde{v}_{\io_h}}- \tilde{x}_h\\
&=\big(1+r_{\io_h}\big) \bigg\{\Big( 
\frac{\phi_{b,\io_h}q_{\ups,\io_{\io_h}}(z)}{v_{\io_h}}-\bar{x}_h\Big) 
+ \Big(\Th_{\io_h} 
   \Big(\frac{\phi_{b,\io_h}q_{\ups,\io_{\io_h}}(z)}{v_{\io_h}}\Big)-
\Th_{\io_h}\big(\bar{x}_h\big)\Big)\bigg\}.
\end{split}\end{equation}
Since $\Th_{\io_h}$ is holomorphic, and 
$$\frac{\phi_{b,\io_h}q_{\ups,\io_{\io_h}}(z)}{v_{\io_h}}-\bar{x}_h
=\phi_{b,h}q_{\ups,\io_h}(z)+c_hv_h,$$
we can rewrite \e_ref{surj_si12c} as 
\begin{equation}\label{surj_si12d}
\phi_{\tilde{b},h}q_{\tilde{\ups},\io_h}(z)=
\big(1+r_{\io_h}\big)\big(1+a_h\big)v_h \Big\{\tilde{c}_h+
\frac{\phi_{b,h}q_{\ups,\io_h}(z)}{v_h}
+\Th_h\Big( \frac{\phi_{b,h}q_{\ups,\io_h}(z)}{v_h}\Big)\Big\},
\end{equation}
where the complex numbers $a_h,\tilde{c}_h\!\in\!\Bbb{C}$ and  
the holomorphic function 
$\Th_h\!: B_{2|v_h|^{-\frac{1}{2}}}(0,\Bbb{C})\lra\Bbb{C}$
are given by
\begin{gather}
\label{surj_si15a}
a_h=\frac{d}{dz}\Th_{\io_h}(z)\Big|_{z=x_h},\qquad
\big(1+a_h\big) \tilde{c}_h= c_h+ 
\frac{\Th_{\io_h}\l(x_h\r)-\Th_{\io_h}\l(x_h-c_hv_h\r)}{v_h},\\
\label{surj_si15c}
\Th_h(z)=\frac{\Th_{\io_h}\l(v_hz+x_h\r)-v_hz\Th_{\io_h}'\l(x_h\r)
           -\Th_{\io_h}\l(x_h\r)}{(1+a_h)v_h}.
\end{gather}
By \e_ref{surj_si15c}, $\Th_h(0)=0$ and $\Th_h'(0)=0$.
By our assumptions on $\Th_{\io_h}$ and 
\e_ref{surj_si15a},\e_ref{surj_si15c},
\begin{gather}
\label{surj_si16a}
\l| a_h \r| \le C(b^*)|\ups|^2\ep |x_{\io_h}|
\le C'(b^*)|\ups|^2\ep;\\
\l| \tilde{c}_h- c_h\r|\le 
C(b^*)\l( \ep |a_h|+ |v_h|^{-1}|\ups|^2\ep |c_hv_h| \r)
\le C'(b^*)|\ups|^2\ep    ,\\
\label{surj_si16c}
\l\| \Th_h''\r\|_{C^0} \le C(b^*)|v_h|^{-1}|\ups|^2\ep|v_h|^2
\le C'(b^*)|\ups|^2\ep.
\end{gather}
We now take
$$\tilde{v}_h=\l(1+a_h\r)\tilde{v}_h'=
\l(1+a_h\r)\l(1+r_{\io_h}\r)\l(1+r_h\r)^{-1}v_h.$$
It follows from \e_ref{surj_si16a}-\e_ref{surj_si16c} that
the induction hypotheses ($\Si$2)-($\Si$5) with $i\!=\!h$ are satisfied.
If $j_l\!=\!h$, let 
$(h,\tilde{y}_l)\!=\!q_{\tilde{\ups},h}q_{\ups,h}^{-1}(h,y_l)$.
By the same argument as in the case $\io_h\!=\!\hat{0}$ above,
($\Si$3)-($\Si$5) of the $i=\io_h$ case imply that
the second part of ($\Si$1) with $i\!=\!h$ is satisfied.
Continuing in this way, we obtain tuples
$$\tilde{b}=\big(\Si,M,I;\tilde{x},(j,\tilde{y}), \tilde{u}'\big),
\quad \tilde{c}=c_{\hat{I}},\quad
\tilde{\ups}=\big(\tilde{b},\tilde{v}_{\hat{I}}\big),$$
satisfying ($\Si$1)-($\Si$5).
Since $\tilde{b}'\!\in\!{\cal M}_{{\cal T},0}^{(H)}$, 
by ($\Si$1) $\tilde{b}\!\in\!{\cal M}_{{\cal T},\ep}^{(H)}$
if $\de$ is sufficiently small.
Finally, \hbox{($\Si$1)-($\Si$5)} along with (I1) and (I2) 
show that $\tilde{b}$ and $\tilde{\ups}$ 
satisfy the two requirements of the~lemma.\\

\begin{crl}
\label{surjb_c2}
If ${\cal T}$ is a simple bubble type,
there exist $\de,C\!\in\! C^{\i}({\cal M}_{\cal T}^{(0)};\Bbb{R}^+)$
with the following property.
Suppose $b^*\!\!\in\!{\cal M}_{\cal T}^{(0)}$, $\ep\!<\!\de(b^*)$,
$b\!\in\!{\cal M}_{{\cal T},\ep}^{(H)}$ is such that 
$d(b^*,b)\!\le\!\de(b^*)$, 
and \hbox{$\ups\!=\!(b,v_{\hat{I}})\!\in\! F_b^{(H)}{\cal T}_{\de(b^*)}$}.
Then there exist 
$\tilde{b}\in{\cal M}_{{\cal T},0}^{(H)}$ and
$
\tilde{\ups}=(\tilde{b},\tilde{v}_{\hat{I}})\in F^{(H)}_{\tilde{b}}{\cal T}$ 
such that\\
(1) $d(b,\tilde{b})\le C(b^*)\ep$ and 
$\l|\tilde{v}_h-v_h\r|\le C (b^*)\ep |v_h|$ for all $h\!\in\!\hat{I}$;\\
(2) if $q_{\ups}(z)\!\in\!\Si_{{\cal T},i}$, 
$\l|r_{b,h}\l(q_{\ups}(z)\r)\r|\!\ge\! 3|v_h|^{\frac{1}{2}}$
for all $h\!\in\!\hat{I}\!-\! H$ such that $\io_h\!=\! i$ and 
\hbox{$\l|q_S^{-1}\l(q_{\ups}(z)\r)\r|\!\ge\! 3|v_i|^{\frac{1}{2}}$} if 
$i\!\in\!\hat{I}\!-\!H$, then
$d_b\big(q_{\ups}(z),q_{\tilde{\ups}}(z)\big)\le\ep$.
\end{crl}

\noindent
{\it Proof:} If $S=S^2$, the tuples $\tilde{b}$ and $\tilde{\ups}$
constructed in the first half of the proof of Lemma~\ref{surjb_l1}
satisfy the requirements of the corollary.
In fact, $d_b\l(q_{\ups}(z),q_{\tilde{\ups}}(z)\r)=0$ if $z$ is as in 
(2)~above.
If $S=\Si$,~let 
$$\tilde{\ep}=\ep^2\prod\limits_{h\in[I]-H}|v_h|_b^2>0.$$ 
If $C(b^*)\de(b^*)$ is sufficiently small, by repeated 
applications of Lemma~\ref{surjb_l1}, we can replace the tuples 
$b$ and $\ups$ by $b'\!\in\!{\cal M}_{{\cal T},\tilde{\ep}}^{(H)}$
and $\ups'\!=\!(b',v'_{\hat{I}})\!\in\! F^{(H)}{\cal T}$ such that\\
(1) $d(b,b')\le C'(b^*)\ep$ and 
$\l|v_h'-v_h\r|\le C'(b^*)\ep |v_h|_b$ for all $h\!\in\!\hat{I}$;\\
(2) if $q_{\ups}(z)\!\in\!\Si_{{\cal T},i}$, 
$\l|r_{b,h}\l(q_{\ups}(z)\r)\r|\!\ge\! \frac{5}{2}|v_h|^{\frac{1}{2}}$
for all $h\!\in\!\hat{I}\!-\! H$ such that $\io_h\!=\! i$ and 
\hbox{$\l|q_S^{-1}\l(q_{\ups}(z)\r)\r|\!\ge\! 
                          \frac{5}{2}|v_i|^{\frac{1}{2}}$} 
if $i\!\in\!\hat{I}\!-\!H$, then
$d_b\l(q_{\ups}(z),q_{\tilde{\ups}}(z)\r)\le 2\de(b)\ep$.\\
Applying the construction of the first half of the proof of 
Lemma~\ref{surjb_l1} to the tuples $b'$ and $\ups'$, we obtain
tuples $\tilde{b}\!\in\!{\cal M}_{\cal T}^{(0)}$ and 
$\tilde{\ups}\!=\!(\tilde{b},\tilde{v}_{\hat{I}})
             \!\in\! F^{(H)}_{\tilde{b}}{\cal T}$ such that
$$d(b',\tilde{b})\le C(b^*)\tilde{\ep}
\quad\hbox{and}\quad
\l|\tilde{v}_h-v_h'\r|\le C(b^*)\tilde{\ep}|v_h'|_{b'} 
~~\forall h\in\hat{I}.$$
Then if $z$ is as in the requirement (2) of the corollary,
$$d_{b'}\l(q_{\ups'}(z),q_{\tilde{\ups}}(z)\r)\le 
 C(b^*)\tilde{\ep}\Big(\prod_{h\in\hat{I}-H}\!\!|v_h'|\Big)^{-1}
 \le \ep^2$$
if $\de$ is sufficiently small. Thus, the tuples $\tilde{b}$ and
$\tilde{\ups}$ satisfy both requirements of the corollary.

\subsection{Gromov Convergence and the $L^p$-norm of the Differential}
\label{surject_a}

\noindent
Let $b_k\!=\!\big(S,M,I;x,(j,y_k),u_k\big)$ 
be a sequence of smooth maps converging to 
$$b^*=\big(S,M,I^*;x^*,(j^*,y^*),u^*\big)
\in{\cal M}_{{\cal T}^*}^{(0)}$$
with respect to the Gromov topology such that 
$\bar{\partial}u_{k,\hat{0}}\!=\!t_k\nu$ with $t_k\!\lra\! 0$
and $\bar{\partial}u_{k,h}\!=\!0$ \hbox{if $h\!\in\!\hat{I}$.}
We assume that ${\cal T}^*$ is a simple bubble type.
In the next subsection, it is proved
that $b_k$ lies in the image of the gluing map 
$\tilde{\ga}_{{\cal T},t_k\nu}$ for some $k$. 
In this subsection, we show the differentials of $du_{k,i}$
satisfy a certain condition which holds for all bubble
maps in the image of~$\tilde{\ga}_{{\cal T},t_k\nu}$.\\

\noindent
By definition of convergence, for all $k$ sufficiently we can choose\\
(a) curves ${\cal C}_k=\big(S,M,I^*;x_k',(j^*,y^*)\big)$ 
with  $\lim\limits_{k\lra\i} x_{k,h}'\lra x_h^*$ 
for all $h\!\in\!\hat{I}^*$, and\\
(b) vectors $(v_k)_{\hat{I}^*}\!\in\! F_{{\cal C}_k}^{(0)}$
with $16|v_k|_{g_b}\le r_{{\cal C}_k}g_b$,\\
such that $\lim\limits_{k\lra\i}|\ups_k|=0$,
${\cal C}(\ups_k)=\big(S,M,I;x_k,(j_k,y(\ups_k))\big)$, and
$$\lim_{k\lra\i}\sup_{z\in\Si_{{\cal C}(\ups_k)}}
   d_V(u_{b^*}(q_{\ups_k}(z)),u_{b_k}(z))=0,\qquad
\lim_{k\lra\i}q_{\ups_k}(j_{k,l},y_{k,l})=(j_l^*,y_l^*)
~~\forall l\!\in\! M,$$
where $\ups_k=({\cal C}_k,(v_k)_{\hat{I}^*})$
and $g_b$ denotes the standard metric on $\Si_{{\cal C}_k}$ if $S=S^2$.
Let
$$\phi_{k,h}=
\begin{cases}
\phi_{x_{k,h}'},&\hbox{if~}x_{k,h}'\!\in\! S^2;\\
\phi_{g_{b,\hat{0}},x_{k,h}'},&\hbox{if~}x_{k,h}'\!\in\!\Si;
\end{cases}~~~
r_{k,h}=
\begin{cases}
r_{x_{k,h}'},&\hbox{if~}x_{k,h}'\!\in\! S^2;\\
r_{g_{b,\hat{0}},x_{k,h}'},&\hbox{if~}x_{k,h}'\!\in\!\Si.
\end{cases}$$
Let $g_{\ups_k}$ be the metric 
on $\Si_{b_k}\!=\!\Si_{\ups_k}$ defined as in Section~\ref{approx_maps_sct1},
using the metric $g_{b,\hat{0}}$ on $\Si$ if $S\!=\!\Si$.\\

\noindent
For any element in the image $\tilde{\ga}_{{\cal T},t\nu}$
that lies near $b^*$,
the modified $(L^p,g_{\ups_k})$-norm of $d\tilde{u}_{\ups}$ is bounded
above by a constant dependent only on $b^*$.
Furthermore, as $\ups\!\lra\! 0$ and the size of the necks is reduced,
the modified $(L^p,g_{\ups})$-norm of $d\tilde{u}_{\ups}$ on such necks tends 
to zero.
The modified $(L^p,g_{\ups})$-norm  is bounded above by
the usual $(L^{2p},g_{\ups})$-norm times some 
constant dependent only on~$b^*$.
In this subsection, we show that the $(L^{2p},g_{\ups_k})$-norm
of $du_{b_k}$ is uniformly bounded and tends to zero
on the ``necks.''
Instead of using our usual cutoff function $\be$,
we will use the family of cutoff functions provided by the following
lemma. The proof can be found in \cite[p166]{MS}.
The statement below is somewhat sharper than in \cite{MS},
but the proof in \cite{MS} suffices.

\begin{lmm}
\label{cutoff_ms}
For every $\ep>0$, there exists a smooth function 
$\tilde{\be}_{\ep}\!: \Bbb{R}\lra[0,1]$ such that 
$$\int_{\Bbb{C}}|\tilde{\be}_{\ep}'(r)|^2rdrd\th\le 8\ep,
\hbox{~~~and~~~}
\tilde{\be}_{\ep}(r)=
\begin{cases}
1,&\hbox{if~}r\ge 1;\\
0,&\hbox{if~}r\le e^{-1/\ep}.\\
\end{cases}$$
\end{lmm}

\noindent
Similarly to the above, given $r\!>\!0$, we denote by $\tilde{\be}_{\ep,r}$
the cutoff function defined by 
$\tilde{\be}_{\ep,r}(t)\!=\!\tilde{\be}_{\ep}(r^{-\frac{1}{2}}t)$.\\

\noindent
We now define nearly holomorphic maps 
$f_{k,i}\!\in\! C^{\i}(\Si_{{\cal C}_k,i};V)$.
In order to simplify computations, we fix a finite family of 
$J$-invariant metrics on~$V$ such that for some fixed $\ve>0$ and
for every \hbox{$q\!\in\! V$}
there exists a metric $g_{V,q}$ in this family
such that $\big(B_{g_{V,q}}(q,\ve),J,g_{V,q}\big)$
is isomorphic to a ball in~$\Bbb{C}^n$.
Since $V$ is compact and the family of metrics is finite, 
all estimates below that depend on a particular metric~$g_{V,q}$
will involve bounds dependent only on~$V$.
We denote by $\exp_q$ the exponential map of
(the Levi-Civita connection of) the metric~$g_{V,q}$
and the $g_{V,q}$-geodesic about $q$ or radius $\ep$ by~$B_q(\ep)$.
If $\de>0$ and $h\!\in\! I^*\!-\!I$, let
\begin{equation}\label{surj_a2}\begin{split}
&B_{h,k}^+(\de)=\big\{ (\io_h^*,z)\!\in\!\Si_{{\cal C}_k,\io_h^*}\!: 
 r_{k,h}(\io_h^*,z)\le\de \big\},\\
&B_{h,k}^-(\de)=\big\{ (h,z)\!\in\!\Si_{{\cal C}_k,h}\!: 
 |q_S^{-1}(z)|\le\de\big\}.\\
\end{split}\end{equation}
Choose a sequence $\ep_k\!\in\!\Bbb{R}^+$ converging to zero.
Let $r_k=\Big(2\sum\limits_{i\in I^*}\|du_i^*\|_{b^*,C^2}\Big)^{-1}\!\!\ep_k$.
By taking a subsequence if necessary, it can be assumed that
\begin{equation}\label{surj_a3a}
|t_k|\le\ep_k,\quad
d_V\big(u_{b^*}(q_{\ups_k}(z)),u_{b_k}(z)\big)\le\ep_k,\quad
r_{b^*,h}(\io_h^*,x_{k,h}')\le r_k,\quad
e^{\frac{2p}{\ep_k}}|v_{k,h}|_{b^*}^{\frac{1}{2}}\le r_k.
\end{equation}
Let $q_h=u_h^*(\i)$ and
$$\tilde{A}_{h,k}^{\pm}=
B_{h,k}^{\pm}\Big(|v_{k,h}|_{b^*}^{\frac{1}{2}}\Big)-
 B_{h,k}^{\pm}
   \Big(e^{-\frac{1}{\ep_k}}|v_{k,h}|_{b^*}^{\frac{1}{2}}\Big).$$
By \e_ref{surj_a3a},
$u_{b_k}\!\Big(\! q_{\ups_k}^{-1}\big(\! B_{h,k}^{\pm}(e^{\frac{1}{\ep_k}}
     |v_{k,h}|_b^{\frac{1}{2}})\!\big)\!\Big)
\!\subset\! B_{q_h}\!\l(C(b^*)\ep_k\r)$.
Thus, we can define $\xi_{k,h}^{\pm}\!\in\! 
 C^{\i}\!\big(\tilde{A}_{h,k}^{\pm};T_{q_h}V\big)$~by
\begin{equation}\label{surj_a4}
\begin{array}{ll}
\exp_{q_h,q_h}
 \xi_{k,h}^+(z)=u_{b_k}
 \big(q_{\ups_k,\io_h^*}^{-1}(\io_h^*,z)\big),&\quad
|\xi_{k,h}^+(z)|_{g_{V,q_h}}<\ve;\\
\exp_{q_h,q_h}
 \xi_{k,h}^-(z)=u_{b_k}\big(q_{\ups_k,\io_h^*}^{-1}
     \big(\io_h^*,\phi_{k,h}^{-1}(zv_{k,h})\big)\big),&\quad
|\xi_{k,h}^-(z)|_{g_{V,q_h}}<\ve,
\end{array}\end{equation}
provided $k$ is sufficiently large (depending on $b^*$).
Let $\bar{\xi}_{k,h}^{\pm}\in T_{q_h}V$ be given~by
\begin{equation}\label{surj_a5}
\bar{\xi}_{k,h}^{\pm}=\frac{1}{\text{Area}(\tilde{A}_{k,h}^{\pm})}
\int_{\tilde{A}_{k,h}^{\pm}}\xi_{k,h}^{\pm},
\end{equation}
where the area and the integral are computed using the metric
$g_{b^*,\io_h^*}$ on $\Si_{b^*,\io_h^*}$ and $g_{b^*,h}$ on $\Si_{b^*,h}$.
Define $f_{k,i}\!\in\! C^{\i}(\Si_{b^*,i};V)$ by
\begin{equation*}\label{surj_a6}
f_{k,i}(z)=
\begin{cases}
\exp_{q_h,q_h}
\Big\{\bar{\xi}_{k,h}^+ +\tilde{\be}_{\ep_k,|v_{k,h}|_{b^*}}\l(r_{k,h}(z)\r)
\big(\xi_{k,h}^+(z)-\bar{\xi}_{k,h}^+\big)\Big\},&
\hbox{if~}r_{k,h}(z)\le |v_{k,h}|_{b^*}^{\frac{1}{2}};\\
\exp_{q_i,q_i}
\Big\{\bar{\xi}_{k,i}^- +\tilde{\be}_{\ep_k,|v_{k,i}|_{b^*}}
 \l(|q_S^{-1}(z)|\r)\big(\xi_{k,i}^-(z)-\bar{\xi}_{k,i}^-\big)\Big\},&
\hbox{if~}i\!\in\! I^*\!-\!I~\&~
|q_S^{-1}(z)|\le |v_{k,i}|_{b^*}^{\frac{1}{2}};\\
u_{b_k}\big(q_{\ups_k}^{-1}(i,z)\big),&\hbox{otherwise}.
\end{cases}
\end{equation*}
Let $\ze_{k,i}'\!\in\!\Ga(u_i^*)$ be given by
$$\exp_{b^*,u_i^*}\ze_{k,i}'=f_{k,i},\quad
\|\ze_{k,i}'\|_{b^*,C^0}<\inj~g_{V,b^*}.$$

\begin{lmm}
\label{surj_al1}
There exists $C>0$ such that for all $k$ sufficiently large 
and $i\!\in\! I^*$,
$$\|\ze_{k,i}'\|_{b^*,C^0}\le C\ep_k,\qquad
   \|\bar{\partial}f_{k,i}\|_{g_{b^*,i},2p}\le 
           C\ep_k^{\frac{1}{2p}}\big(\|df_{k,i}\|_{g_{b^*,i},2p}+1\big).$$
\end{lmm}

\noindent
{\it Proof:} The first statement is clear from \e_ref{surj_a3a}
 and the construction of $f_{k,i}$ above.
Suppose \hbox{$z\!\in\!\Si_{b^*,i}$}.
If $z\!\not\in\! B_{h,k}^+\big(|v_{k,h}|_{b^*}^{\frac{1}{2}}\big)$ for all 
$h\!\in\! I^*\!-\!I$
and $z\!\not\in\! B_{i,k}^-\big(|v_{k,i}|_{b^*}^{\frac{1}{2}}\big)$ if 
$i\!\in\! I^*\!-\!I$, then
\begin{equation}
\label{surj_al1_e1}
|\bar{\partial}f_{k,i}|_{g_{b^*,i},z}\le C_{\nu}t_k\le C_{\nu}\ep_k.
\end{equation}
Suppose $z\!\in\!\tilde{A}_{h,k}^+$ with $h\!\in\! I^*\!-\!I$. 
Since the metric $g_{\P,q_h}$ is flat near $q_h$,
\begin{equation}
\label{surj_al1_e2}
\bar{\partial}f_{k,i}\big|_z=
d\exp_{q_h,q_h}\bar{\partial}
\Big\{\bar{\xi}_{k,h}^+ 
+\tilde{\be}_{\ep_k,|v_{k,h}|_{b^*}}\l(r_{k,h}(\cdot)\r)
\big(\xi_{k,h}^+-\bar{\xi}_{k,h}^+\big)\Big\}_z.
\end{equation}
It follows from \e_ref{surj_al1_e2} that
\begin{equation}\label{surj_al1_e3}\begin{split}
|\bar{\partial}f_{k,i}|_{g_{b^*,i},z} &\le 
 C\Big(|v_{k,h}|_{b^*}^{-\frac{1}{2}}
           |d\tilde{\be}_{\ep}|_{|v_{k,h}|_{b^*}^{-\frac{1}{2}}r_{k,h}(z)}
\big|\xi_{k,h}^+-\bar{\xi}_{k,h}^+\big|_z+
\big|\bar{\partial}\xi_{k,h}^+\big|_z\Big)
\end{split}\end{equation}
By Poincare Lemma, 
see Lemma~\ref{an-poin_4} in~\cite{Z1} 
applied with $r=|v_{h,k}|_{b^*}^{-\frac{1}{2}}$
and $2p$ instead of~$p$,
\begin{equation*}\label{surj_al1_e4}\begin{split}
\Big\| |v_{k,h}|_{b^*}^{-\frac{1}{2}}
   |d\tilde{\be}_{\ep_k}|_{|v_{k,h}|_{b^*}^{-\frac{1}{2}}r_{k,h}(\cdot)}
\big|\xi_{k,h}^+-\bar{\xi}_{k,h}^+\big|
 \Big\|_{g_{b^*,i},L^{2p}(\tilde{A}_{h,k}^+)}  
 \le C|v_{k,h}|_{b^*}^{-\frac{p-1}{2p}}
\|d\tilde{\be}_{\ep_k}\|_2^{\frac{1}{p}}
\big\|\xi_{k,h}^+-\bar{\xi}_{k,h}^+\big\|_{b^*,C^0}&\\
\le C'\ep_k^{\frac{1}{2p}}
      \big\|d\xi_{k,h}^+\big\|_{g_{b^*,i},L^{2p}(\tilde{A}_{h,k}^+)}
\le C'\ep_k^{\frac{1}{2p}}\l\|df_{k,h}\r\|_{g_{b^*,i},2p}.&
\end{split}\end{equation*}
The last two equations give 
\begin{equation}\label{surj_al1_e5} 
\l\|\bar{\partial}f_{k,i}\r\|_{g_{b^*,i},L^{2p}(\tilde{A}_{h,k}^+)}
\le C\big(\ep_k^{\frac{1}{2p}}\l\|df_{k,h}^+\r\|_{g_{b^*,i},2p}+\ep_k\big).
\end{equation}
The same estimate applies to 
$\l\|\bar{\partial}f_{k,i}\r\|_{g_{b^*,i},L^{2p}(\tilde{A}_{i,k}^-)}$
if $i\!\in\! I^*\!-\!I$.
Here the exponent of $\frac{2p}{\ep_k}$ in \e_ref{surj_a3a} is crucial:
\begin{equation}\label{surj_al1_e6}\begin{split}
\|\bar{\partial}\xi_{k,i}^-\|_{g_{b^*,i},L^{2p}(\tilde{A}_{i,k}^-)}^{2p}
&\le \int\limits_{|v_{k,i}|_{b^*}^{-\frac{1}{2}}\le r\le 
                        e^{\frac{1}{\ep_k}}|v_i|_{b^*}^{-\frac{1}{2}}}
t_k^{2p}|\nu\circ dq_{\ups_k,\io_i^*}^{-1}|_{g_{b^*,\io_i^*}}^{2p}
|v_{k,i}|_{b^*}^{2p}(1+r^2)^{2p-2}rdrd\th\\
&\le C t_k^{2p}|v_{k,i}|_{b^*}^{2p}
\l(|v_{k,i}|^{-\frac{1}{2}}e^{\frac{1}{\ep_k}}\r)^{4p-2}
\le C t_k^{2p}.
\end{split}\end{equation}
Since $f_{k,i}$ is constant on 
$B_{h,k}^+\big(e^{-\frac{1}{\ep_k}}|v_{k,h}|_{b^*}^{\frac{1}{2}}\big)$
for $h\!\in\! I^*\!-\!I$ with $\io_h^*=i$ and on 
$B_{i,k}^-\big(e^{-\frac{1}{\ep_k}}|v_{k,i}|_{b^*}^{\frac{1}{2}}\big)$
if $i\!\in\! I^*\!-\!I$, the second claim is proved.

\begin{crl}
\label{surj_ac2}
There exists $C>0$ such that for all $k$ sufficiently large,
$$\|df_{k,i}\|_{g_{b^*,i},2p}\le C \quad\hbox{and}\quad
 \|\ze_{k,i}'\|_{g_{b^*,i},2p,1}\le C\ep_k^{\frac{1}{2p}}.$$
\end{crl}

\noindent
{\it Proof:} By the quadratic expansion of $\bar{\partial}_{u_i^*}\ze_{k,i}'$
as in Subsection~\ref{gluing_map},
\begin{equation}\label{surj_ac2_e1}
D_{b^*,u_i^*}\ze_{k,i}'+N_{\bar{\partial},u_i^*}\ze_{k,i}'
                                   =\bar{\partial}_{u_i^*}\ze_{k,i}',
\end{equation}
where 
\begin{equation}\label{surj_ac2_e2}
\|\bar{\partial}_{u_i^*}\ze_{k,i}'\|_{g_{b^*,i},2p}\le 
C\ep_k^{\frac{1}{2p}}\big( \|df_{k,i}\|_{g_{b^*,i},2p}+1\big)
\end{equation}
by Lemma~\ref{surj_al1} and 
\begin{equation}\label{surj_ac2_e3}
\|N_{\bar{\partial},u_i^*}\ze_{k,i}'\|_{g_{b^*,i},2p}
\le C\|\ze_{k,i}'\|_{C^0} \|\ze_{k,i}'\|_{g_{b^*,i},2p,1}
\le C\ep_k \|\ze_{k,i}'\|_{g_{b^*,i},2p,1},
\end{equation}
by Proposition~\ref{an-dbar_prp} in~\cite{Z1} and Lemma~\ref{surj_al1}.
Thus, by standard elliptic estimates for $u_{b^*}$ and 
\e_ref{surj_ac2_e1}-\e_ref{surj_ac2_e3}, 
\begin{equation}\label{surj_ac2_e4}\begin{split}
\|\ze_{k,i}'\|_{g_{b^*,i},2p,1} &\le 
C\big( \|D_{b^*,u_i}\ze_{k,i}'\|_{g_{b^*,i},2p}+ 
                       \|\ze_{k,i}'\|_{g_{b^*,i},2p} \big)\\
&\le C'\ep_k^{\frac{1}{2p}}\big(\|\ze_{k,i}'\|_{g_{b^*,i},2p,1}+
\|df_{k,i}\|_{g_{b^*,i},2p}+1\big).
\end{split}\end{equation}
On the other hand, since $f_{k,i}=\exp_{b^*,u_i^*}\ze_{k,i}'$,
\begin{equation}\label{surj_ac2_e5}
\|df_{k,i}\|_{g_{b^*,i},2p}\le C\big(\|du_i^*\|_{g_{b^*,i},2p}+
\|\ze'_{k,i}\|_{g_{b^*,i},2p,1}\big).
\end{equation}
If $\ep_k$ is sufficiently small, the claim follows
from equations~\e_ref{surj_ac2_e4} and~\e_ref{surj_ac2_e5}.

\begin{crl}
\label{surj_ac3}
There exists $C>0$ such that for all $k$ sufficiently large, 
$h\!\in\!\hat{I}^*$, and $\de>0$,
$$\l\|du_{b_k}\r\|_{g_{\ups_k},L^{2p}
\l(q_{\ups_k}^{-1}(B_{h,k}^{\pm}(\de))\r)}
\le C\big(\ep_k^{\frac{1}{2p}}+\de^{\frac{1}{p}}\big).$$
\end{crl}

\noindent
{\it Proof:} If $h\!\in\!\hat{I}$, the statement is immediate 
from Corollary~\ref{surj_ac2}; so we assume $h\!\in\! I^*\!-\!I$.
The metric~$g_{\ups}$ on $q_{\ups_k}^{-1}(B_{h,k}^+(\de))$
differs by a bounded factor from the metric $q_{\ups_k,\io_h^*}^*g_{b^*,i}$.
Thus, 
\begin{equation}\label{surj_ac3_e1}\begin{split}
\l\|du_{b_k}\r\|_{g_{\ups_k},L^{2p}\l(q_{\ups_k}^{-1}(B_{h,k}^+(\de))\r)}
&\le C\big\|d(f_k\circ q_{\ups_k,\io_h^*}^{-1})
\big\|_{g_{b^*,\io_h^*},L^{2p}\l(B_{h,k}^+(\de)-
       B_{h,k}^+(|v_{k,h}|_{b^*}^{\frac{1}{2}})\r)}\\
&=C\big\|df_{k,\io_h^*}
\big\|_{g_{b^*,\io_h^*},L^{2p}\l(B_{h,k}^+(\de)-
          B_{h,k}^+(|v_{k,h}|_{b^*}^{\frac{1}{2}})\r)}\\
&\le C\big\|df_{k,\io_h^*}
\big\|_{g_{b^*,\io_h^*},L^{2p}\l(B_{h,k}^+(\de)\r)}.
\end{split}\end{equation}
Since $f_{k,\io_h^*}=\exp_{b^*,u_{\io_h^*}}\ze_{k,\io_h^*}'$,
 by Corollary~\ref{surj_ac2},
\begin{equation}\label{surj_ac3_e2}
\big\|df_{k,\io_h^*}\big\|_{g_{b^*,\io_h^*},L^{2p}\l(B_{h,k}^+(\de)\r)}
\!\le\!
C\Big(\big\|du_{\io_h^*}^*\big\|_{g_{b^*,\io_h^*},L^{2p}\l(B_{h,k}^+(\de)\r)}
+\big\|\ze_{k,\io_h^*}'\big\|_{g_{b^*,\io_h^*},2p,1}\Big)
\!\le\! C'\big(\de^{\frac{1}{p}}+\ep^{\frac{1}{2p}}\big).
\end{equation}
The claim for $B_{h,k}^+(\de)$ follows from \e_ref{surj_ac3_e1} and
\e_ref{surj_ac3_e2}.
The metric $g_{\ups}$ on $q_{\ups_k}^{-1}(B_{h,k}^-(\de))$
differs by a bounded factor from the metric 
which is the pullback of the metric
$g_{b^*,h}$ by the map
$$z\lra  q_N\Big(\frac{\phi_{k,h}q_{\ups_k,\io_h^*}(z)}{v_{k,h}}\Big).$$
Thus, similarly to the above,
\begin{gather}
\label{surj_ac3_e3}
\l\|du_{b_k}\r\|_{g_{\ups_k},L^{2p}\l(q_{\ups_k}^{-1}(B_{h,k}^-(\de))\r)}
\le C\l\|df_{k,h}\r\|_{g_{b^*,h},L^{2p}\l(B_{h,k}^-(\de)\r)};\\
\label{surj_ac3_e4}
 \l\|df_{k,h}\r\|_{g_{b^*,h},L^{2p}\l(B_{h,k}^-(\de)\r)}
\le C\big(\de^{\frac{1}{p}}+\ep^{\frac{1}{2p}}\big).
\end{gather}
The claim for $B_{h,k}^-(\de)$ follows from \e_ref{surj_ac3_e3} and
\e_ref{surj_ac3_e4}.

\subsection{Surjectivity of the Gluing Map}
\label{surj_c}

\noindent
We continue with the notation of Subsection~\ref{surject_a}.
In this subsection, for $k$ sufficiently, we use Corollary~\ref{surjb_c2}
to construct 
$$\tilde{\ups}_k=\big(\tilde{b}_k,(\tilde{v}_k)_{\hat{I}^*}\big)\in
F^{(0)}{\cal T}_{\de}$$ 
and $\tilde{\ze}_k\!\in\!\Ga(u_{\tilde{\ups}_k})$
such that $\tilde{b}_k$ is very close to $b$ in 
${\cal M}_{\cal T}^{(0)}$, $\|\tilde{\ze}_k\|_{\tilde{\ups}_k,p,1}$
is small, and $u_{b_k}\!=\!\exp_{\tilde{\ups}_k}\tilde{\ze}_k$.
We then look at the elements of $F^{(0)}{\cal T}_{\de}$
near $\tilde{\ups}_k$ to find $\tilde{\ups}_k'$ and 
$\tilde{\ze}_k\!\in\!\tilde{\Ga}_+(\tilde{\ups}_k')$
such that $u_{b_k}\!=\!\exp_{\tilde{\ups}_k'}\tilde{\ze}_k'$.
If ${\cal T}$ is semiregular, we consider only  the case $I\!=\!\hat{0}$;
if ${\cal T}$ is regular, we \hbox{assume $t\!=\!0$}.\\

\noindent
Let $H=\hat{I}\subset\hat{I}^*$.
If $\de>0$ and $i\!\in\! I^*$, put
$$\Si_{i,\de}=\big\{(i,z)\!\in\!\Si_{b^*,i}\!: r_{b^*,h}(i,z)\le\de~\forall
h\!\in\! \hat{I}\!-\!H\hbox{~s.t.~}\io_h=i,~
|q_S^{-1}(z)|\ge\de\hbox{~if~}i\!\in\! \hat{I}\!-\!H\big\}.$$
In addition to \e_ref{surj_a3a}, we can assume that our sequence satisfies
\begin{equation}\label{surj_a3b}
\l\|\ze_{k,i}\r\|_{g_{b^*,i},C^2(\Si_{i,r_k})} \le\ep_k.
\end{equation}
Let $b_k'\!=\!\big(S,M,I^*;x_k',(j^*,y^*),u^*\big)$.
By the second assumption in \e_ref{surj_a3a},
$$d(b^*,b_k')\le C\ep_k\Lra b_k'\in{\cal M}_{{\cal T}^*,C\ep_k}^{(H)},$$
since $b^*\!\in\!{\cal M}_{{\cal T}^*}^{(0)}$, 
where $C>0$ depends only on $b^*$.
By the last assumption of \e_ref{surj_a3a}, $|\ups_k|_{b_k'}\le C\ep_k$.
Thus, if $\ep_k>0$ is sufficiently small, by Corollary~\ref{surjb_c2}, 
there exist 
$$\tilde{b}_k\in{\cal M}_{{\cal T}^*}^{(H)} 
\quad\hbox{and}\quad
\tilde{\ups}_k=\big(\tilde{b}_k,(\tilde{v}_k)_{\hat{I}^*}\big)
\in F^{(H)}{\cal T}^*$$
such that\\
(1) $d(b,\tilde{b})\le C'\ep_k$ and 
$\l|\tilde{v}_{k,h}-v_{k,h}\r|\le C'\ep_k |v_{k,h}|_b$ 
for all $h\!\in\!\hat{I}^*$;\\
(2) if $q_{\ups_k}(z)\!\in\!\Si_{{\cal T}^*,i}$, 
$\l|r_{b^*,h}\l(q_{\ups_k}(z)\r)\r|\!\ge\! 3|v_{k,h}|^{\frac{1}{2}}$
for all $h\!\in\!\hat{I}^*\!-\! H$ such that $\io_h^*\!=\! i$ and 
\hbox{$\l|q_S^{-1}\l(q_{\ups_k}(z)\r)\r|\!\ge\! 3|v_{k,i}|^{\frac{1}{2}}$} if 
$i\!\in\! I^*\!-\!H$, then
$d_b\l(q_{\ups_k}(z),q_{\tilde{\ups}_k}(z)\r)\le\ep_k$.\\
It then follows from the second and third assumptions of \e_ref{surj_a3a}
that there exist $\tilde{\ze}_k\!\in\!\Ga(u_{\tilde{\ups}_k})$,
\hbox{$\tilde{w}_{k,h}\!\in\! 
        T_{x_h(\tilde{\ups}_k)}\Si_{\tilde{\ups}_k,\io_h}$} 
for $h\!\in\! H$,
and $\tilde{w}_{k,l}\!\in\! T_{y_l(\tilde{\ups}_k)}\Si_{\tilde{\ups}_k,j_l}$ 
for $l\!\in\! M$
such~that
\begin{gather*}
\exp_{\tilde{\ups}_k}\tilde{\ze}_k=u_{b_k},\quad
\exp_{g_{\tilde{\ups}_k},x_{k,h}(\tilde{\ups}_k)}\tilde{w}_{k,h}
                                                           =x_{k,h},\quad
\exp_{g_{\tilde{\ups}_k},y_{k,l}(\tilde{\ups}_k)}\tilde{w}_{k,l}
                                                            =y_{k,l};\\
\|\tilde{\ze}_k\|_{b^*,C^0},
|\tilde{w}_{k,h}|_{g_{\tilde{\ups}_k},x_{k,h}(\tilde{\ups}_k)}
|\tilde{w}_{k,l}|_{g_{\tilde{\ups}_k},y_{k,l}(\tilde{\ups}_k)}\le C'\ep_k.
\end{gather*}

\begin{lmm}
\label{surj_bl3}
There exists $C>0$ such that for all $k$,
$$\|\tilde{\ze}_k\|_{\tilde{\ups}_k,p,1}\le C\ep_k^{\frac{1}{2p}}.$$
\end{lmm}

\noindent
{\it Proof:} By \e_ref{surj_a3b}, (1), and (2), 
$\|\tilde{\ze}_k\|_{g_{\ups},C^1}\le C\ep_k$
outside of the necks 
$$\tilde{A}_{k,h}=q_{\ups_k}^{-1}
\Big(B_{k,h}^+(r_k)\cup B_{k,h}^-(r_k)\Big).$$
On the other hand, $\|du_{\tilde{\ups}_k}\|_{\tilde{\ups}_k,C^0}\le C$
by Lemma~\ref{approx_maps} and 
$$\|du_{\tilde{\ups}_k}\|_{\tilde{\ups}_k,L^p(\tilde{A}_{k,h})}\le 
C(\ep_k^{\frac{1}{2p}}+r_k^{\frac{1}{p}})\le C'\ep_k^{\frac{1}{2p}}$$
by Corollary~\ref{surj_ac3}. 
The three estimates imply the claim.\\

\noindent
Suppose $\hat{I}\!=\!\eset$ and thus $H\!=\!\eset$.
If $k$ is sufficiently large and 
$\vp\!\in\! T_{\tilde{\ups}_k}F^{\eset}{\cal T}^*$ 
is such that \hbox{$2\|\vp\|_{\tilde{\ups}_k}\!<\!\de(b^*)$}, 
where $\de$ is as 
in Lemmas~\ref{inj_l1} and~\ref{inj_l1b},~let
$$\tilde{b}_k(\vp)=\tilde{b}_{t\nu}\big(\tilde{\ups}_k(\vp)\big)=
\big(S,M,\{\hat{0}\};,(\hat{0},\bar{y}(\vp)),\tilde{u}_{\vp,t\nu}\big)$$
be the tuple defined as in Subsections~\ref{scale_sec}
and~\ref{gluing_map}.
Let $\tilde{\ze}_k(\vp)\!\in\!\Ga(\tilde{u}_{\vp,t\nu})$ 
and $\tilde{w}_{k,l}(\vp)\!\in\! T_{y_l(\vp)}\Si_{\vp,j_l}$ 
for $l\!\in\! M$ be given by
$$\exp_{\tilde{\ups}_k(\vp)}\tilde{\ze}_k(\vp)=u_{b_k},~~
\exp_{g_{\tilde{\ups}_k(\vp)},y_{k,l}(\vp)}\tilde{w}_{k,l}(\vp)=y_{k,l};
\quad
\|\tilde{\ze}_k(\vp)\|_{b^*,C^0},
|\tilde{w}_{k,l}(\vp)|_{g_{\tilde{\ups}_k},y_{k,l}(\vp)}\le 2C'\ep_k.$$ 
We need to find $\vp$ such that 
$\tilde{\pi}_{\vp,-}\tilde{\ze}_k(\vp)=0$ 
and  $y_l(\vp)=y_{k,l}$, or equivalently
\begin{equation}\label{surj_e20}
S_{\vp}\tilde{\pi}_{\vp,-}\tilde{\ze}_k(\vp)=0
\quad\hbox{and}\quad
S_{\vp}\tilde{w}_{k,l}(\vp)=0,
\end{equation}
where $S_{\vp}\tilde{w}_{k,l}(\vp)$
denotes the parallel transform of $\tilde{w}_{k,l}(\vp)$
back to $y_l(\tilde{\ups}_k)$ along the $g_{\ups}$-geodesic
$$s\!\lra\!\exp_{y_l(\tilde{\ups}_k)}sw_l(\vp).$$

\begin{lmm}
\label{surj_cl1}
There exists $C>0$ such that for all $k$ sufficiently large
and $\vp,\vp'\!\in\! T_{\tilde{\ups}_k}F^{(\eset)}{\cal T}^*$ with 
$2\|\vp\|_{\tilde{\ups}_k}<\de(b^*)$,
\begin{equation*}\begin{split}
&S_{\vp}\tilde{\pi}_{\vp,-}\tilde{\ze}_k(\vp)=
\tilde{\pi}_{\ups,-}\tilde{\ze}_k+
\tilde{N}^{(0)}(\tilde{\ze}_k,\vp)-
\tilde{\pi}_{\ups,-}\ze_{\vp}+N^{(0)}(\vp),\\
&S_{\vp}\tilde{w}_{k,h}(\vp)=\tilde{w}_{k,h}
+\tilde{N}^{(h)}(\tilde{w}_{k,h},\vp)
-w_h(\vp)+N^{(h)}(\vp)\quad\forall h\!\in\!M,
\end{split}\end{equation*}
where $\ze_{\vp}$ is as in Section~\ref{scale_sec} and
$\tilde{N}^{(l)}$ and $N^{(l)}$ satisfy
\begin{equation}\label{surj_cl1_e1}\begin{split}
&\big\|\tilde{N}^{(0)}(\tilde{\ze}_k,\vp)-
      \tilde{N}^{(0)}(\tilde{\ze}_k,\vp')\big\|_{\tilde{\ups}_k,2}\le 
C\|\tilde{\ze}_k\|_{\tilde{\ups}_k,2}\|\vp-\vp'\|_{\tilde{\ups}_k};\\
&\big|\tilde{N}^{(l)}(\vp,\tilde{w}_{k,l})-
 \tilde{N}^{(l)}(\vp',\tilde{w}_{k,l})
     \big|_{g_{\tilde{\ups}_k},y_l(\tilde{\ups}_k)}
\le C\|\tilde{w}_{k,l}\|_{g_{{\ups}_k},y_l(\tilde{\ups}_k)}
      \|\vp-\vp'\|_{\tilde{\ups}_k} \quad \forall l\!\in\! M;\\
&\|N^{(0)}(\vp)-N^{(0)}(\vp')\|_{\tilde{\ups}_k,2}\le 
C\l(\|\vp\|_{\tilde{\ups}_k}+\|\vp'\|_{\tilde{\ups}_k}\r)
\|\vp-\vp'\|_{\tilde{\ups}_k};\\
&\big|N^{(l)}(\vp)-N^{(l)}(\vp')
\big|_{g_{\tilde{\ups}_k},y_l(\tilde{\ups}_k)}\le 
C\l(\|\vp\|_{\tilde{\ups}_k}+\|\vp'\|_{\tilde{\ups}_k}\r)
\|\vp-\vp'\|_{\tilde{\ups}_k}\quad \forall l\!\in\! M.
\end{split}\end{equation}
\end{lmm}

\noindent
{\it Proof:} 
This lemma follows from a pointwise Riemannian geometry
estimate on 
$S_{\vp}\tilde{\ze}_k(\vp)-(\tilde{\ze}_k-\tilde{\ze}_{\vp})$
and the fact that all statements in 
Lemmas~\ref{inj_l1} and~\ref{inj_l1b}
can be written in a form similar to \e_ref{surj_cl1_e1}, 
e.g.~for all $\xi\in\Ga(\tilde{\ups}_k)$
$$\big\|S_{\vp}\tilde{\pi}_{\vp,-}R_{\vp}\xi-
S_{\vp'}\tilde{\pi}_{\vp',-}R_{\vp'}\xi
\big\|_{\tilde{\ups}_k,2}
\le C\|\vp-\vp'\|_{\tilde{\ups}_k,2}\|\xi\|_{\tilde{\ups}_k,2}.$$
The latter fact can be seen from the two lemmas and
the definitions of $R_{\vp}$ and $S_{\vp}$ in Subsection~\ref{scale_sec}.

\begin{lmm}
\label{surj_l15}
There exist $C,\de\!\in\! C^{\i}({\cal M}_{{\cal T}^*}^{(0)},\Bbb{R}^+)$
such that for all $\ups\!\in\! F^{(H)}{\cal T}_{\de}^*$ and 
$\vp\!\in\! T_{\ups}F^{(H)}{\cal T}_{\de}^*$ with 
$\|\vp\|_{\tilde{\ups}_k}\le\de(b)$,
$$\|\ze_{\vp}\|_{\ups,2}\le C(b)\|\tilde{\pi}_{\ups,-}\ze_{\vp}\|_{\ups,2}.$$
\end{lmm}

\noindent
{\it Proof:} It can be seen directly from the definitions that
$$\|\ze_{\vp}\|_{\ups,2}\le \big(1+C(b_{\ups}\big)|\ups|\big)
          \|\pi_{\ups,-}\ze_{\vp}\|_{\ups,2}.$$
The claim then follows from the proof of (2b) 
of Lemma~\ref{tangent_l1}.

\begin{crl}
\label{surl_c17}
There exist a neighborhood $U$ of $b^*$ in 
${\cal M}_{{\cal T}^*}^{(0)}$ and $\de,\ep\!>\!0$
such that for all \hbox{$\ups\!\in\! F^{(\eset)}{\cal T}_{\de}^*|U$},
$\xi\!\in\!\tilde{\Ga}_-(\ups)$ with $\|\xi\|_{\ups,2}<\de$,
$w_h\!\in\! T_{x_h(\ups)}
 \Si_{\ups,\io_h}$ with $|w_h|_{g_{\ups},x_h(\ups)}\!<\!\de$
and \hbox{$w_l\!\in\! T_{y_l(\ups)}\Si_{\ups,j_l}$} for $l\!\in\! M$ with 
$|w_l|_{g_{\ups},y_l(\ups)}<\de$, the system of equations
$$\tilde{\pi}_{\ups,-}\ze_{\vp}-N^{(0)}(\vp)=\xi,\qquad
w_h(\vp)-N^{(l)}(\vp)=w_h~~\forall h\!\in\! H,M$$
has  a (unique) solution $\vp\!\in\! T_{\ups}F^{(\eset)}{\cal T}$
with $\|\vp\|_{\ups}<\ep$.
\end{crl}

\noindent
{\it Proof:} By Lemmas~\ref{inj_l1} and \ref{surj_l15}, 
$$ C^{-1}\|\vp\|_{\ups}\le
 \|\tilde{\pi}_{\ups,-}\xi\|_{\ups,2}+
\sum\limits_{h\in H}|w_h(\vp)|_{g_{\ups},x_h(\ups)}+
\sum\limits_{l\in M}|w_l(\vp)|_{g_{\ups},y_l(\ups)}\le C\|\vp\|_{\ups}.$$
whenever $b_{\ups}$ lies near $b^*$.
Thus, the claim follows from \e_ref{surj_cl1_e1}
by the usual contraction-principle argument.

\begin{crl}
\label{surject_crl}
Let ${\cal T}^*=(S,M,\{\hat{0}\};j^*,\la^*)$ 
be a simple bubble type.
If ${\cal T}^*$ is regular, the map 
$$\tilde{\ga}_{{\cal T}^*}\!: 
   F{\cal T}_{\de}^*\lra\bar{\cal M}_{\lr{{\cal T}^*}}$$
contains a neighborhood of ${\cal M}_{{\cal T}^*}$ 
in~$\bar{\cal M}_{\lr{{\cal T}^*}}$.
If ${\cal T}^*$ is semiregular, 
$H\!=\!\eset$, and $k$ is sufficiently large, there exists 
$\tilde{\ups}_k\!\in\! F^{(\eset)}{\cal T}_{\de}^*$ such that
$b_k=\tilde{\ga}_{{\cal T}^*,t_k\nu}(\tilde{\ups}_k)$.
\end{crl}

\noindent
{\it Proof:} The second statement is immediate from Lemmas~\ref{surj_bl3} 
and \ref{surj_cl1} and Corollary~\ref{surl_c17}.
If ${\cal T}^*$ is regular, 
what we have shown is that the image of 
$\tilde{\ga}_{{\cal T}^*}$ contains a neighborhood of 
${\cal M}_{{\cal T}^*}$ in 
${\cal M}_{\lr{{\cal T}^*}}\cup{\cal M}_{{\cal T}^*}$.
Furthermore, there exists a sequence of neighborhoods
$U_1\supset U_2\supset\ldots$ of $b^*$ in 
$\bar{\cal M}_{\lr{{\cal T}^*}}$
such that $\bigcap\limits_{k}U_k=\{[b^*]\}$.
If $[b_k]\!\in\!{\cal M}_{\cal T}$ is a sequence of bubble maps
converging to \hbox{$[b^*]\!\in\!{\cal M}_{{\cal T}^*}$},
it can be assumed that $[b_k]\!\in\! U_k$.
By the above statement applied to~${\cal T}$,
we can choose sequences 
$$\{[b_{kr}]\}\subset{\cal M}_{\lr{{\cal T}^*}}=
{\cal M}_{\lr{\cal T}}$$
such that for each fixed $k$ the sequence $\{[b_{kr}]\}$ converges
to~$[b_k]$.
Since $U_k$ is an open neighborhood of~$[b_k]$, it can be assumed
that $[b_{kr}]\!\in\! U_k$ for all~$r$.
By the above, the image of 
$\tilde{\ga}_{{\cal T}^*}|F{\cal T}^*_{\frac{1}{2}\de}$
contains $U_k\cap {\cal M}_{\lr{{\cal T}^*}}$ if $k$ is sufficiently large.
Thus, for all $r$ there exists 
$\ups_{kr}\!\in\! F{\cal T}^*_{\frac{1}{2}\de}$
such that $\tilde{\ga}_{{\cal T}^*}(\ups_{kr})=[b_{kr}]$.
Let $\tilde{\ups}_k\!\in\! F{\cal T}^*_{\de}$ be the limit of 
the sequence $\ups_{kr}$ with $k$ fixed.
Then, by continuity of the map $\tilde{\ga}_{{\cal T}^*}$, 
see Corollary~\ref{cont_c4},
$$\tilde{\ga}_{{\cal T}^*}(\tilde{\ups}_k)=
\lim_{r\lra\i}\tilde{\ga}_{{\cal T}^*}(\tilde{\ups}_{kr})=
\lim_{r\lra\i}[b_{kr}]=[b_k].$$
Thus, the image of  $\tilde{\ga}_{{\cal T}^*}$ contains
a neighborhood of 
${\cal M}_{{\cal T}^*}$ in~$\bar{\cal M}_{\lr{{\cal T}^*}}$.

\begin{crl}
\label{surject_crl2}
If ${\cal T}^*\!=\!(S,M,I^*;j^*,\la^*)$ is 
a simple regular bubble type, the map 
$$\tilde{\ga}_{{\cal T}^*}\!: 
   F{\cal T}_{\de}^*\lra\bar{\cal M}_{\lr{{\cal T}^*}}$$
is a homeomorphism onto an open neighborhood of 
${\cal M}_{{\cal T}^*}$ in $\bar{\cal M}_{\lr{{\cal T}^*}}$
provided $\de\!\in\! C^{\i}({\cal M}_{{\cal T}^*};\Bbb{R}^+)$
is sufficiently small.
\end{crl}

\noindent
{\it Proof:}
By Corollaries~\ref{cont_c4}, \ref{inj_main2c}, \ref{surject_crl},
the map
$\tilde{\ga}_{{\cal T}^*}\!: 
F{\cal T}_{\de}^*\lra\bar{\cal M}_{\lr{{\cal T}^*}}$
is a continuous bijection onto a neighborhood of 
${\cal M}_{{\cal T}^*}$ in $\bar{\cal M}_{\lr{{\cal T}^*}}$.
In addition, the proof of Corollary~\ref{surject_crl} shows that 
$\tilde{\ga}_{{\cal T}^*}$ is an open~map.

\section{Appendix}

\subsection{Properties of Smooth Families of Metrics on $\Si$}
\label{si_metrics}

\noindent
Let $m$ be a positive integer and
$$\aleph=\l\{x=x_{[m]}\!: x_h\!\in\!\Si, x_h\!\neq\! x_l
\hbox{~if~}h\!\neq\! l\r\}.$$
Suppose $\{g_x\!: x\!\in\! \aleph\}$ is a smooth family of metrics
on $\Si$ such that for any $x=x_{[m]}\in\aleph$
the metric $g_x$ is flat on a neighborhood of $x_h$ in $\Si$
for all $h\!\in\![m]$.
If $x=x_{[m]}\in\aleph$ and $v\!\in\! T_y\Si$,~let
$$T_x\aleph=\bigoplus_{h\in[m]}T_{x_h}\Si,~~~~
|v|_x=|v|_{g_x,y}.$$
If $w=w_{[m]}\in T_x\aleph$, let $|w|$ denote $\sum\limits_{h\in[m]}|w_h|_x$.
Define $x(w)\!\in\!\Si^m$ by
$$x(w)=\big(x_1(w),\ldots,x_m(w)\big)=
\big(\exp_{g_x,x_1}w_1,\ldots,\exp_{g_x,x_m}w_m\big).$$
We denote by $\phi_{x,y}$ the map $\phi_{g_x,y}$
and by $B_x(y,\de)$ the set $B_{g_x}(y,\de)$
described in Subsection~\ref{notation}.
If \hbox{$\de\!: \aleph\lra\Bbb{R}$}, let
$$T\aleph_{\de}=\l\{(x,w)\!: x\!\in\!\aleph;~ w\!\in\! T_x\aleph,~ 
|w|_x<\de(x)\r\}.$$

\begin{lmm}
\label{si_metrics_l1a}
There exist $\de\!\in\! C^{\i}(\aleph;\Bbb{R}^+)$
and a smooth families of holomorphic maps
$$\big\{\tilde{p}_{h,(x,w)}\!:
\{z\in B_x(x_h,\de(x))\}\lra\Si~
|~ (x,w)\!\in\! T\aleph_{\de}\big\},$$
such that each map $\tilde{p}_{h,(x,w)}$ is a $(g_x,g_{x(w)})$-isometry,
\begin{gather}
\label{si_metrics_l1a_e}
d\phi_{x,x_h}\big|_{x_h(w)}\phi_{x(w),x_h(w)}
\tilde{p}_{h,(x,w)}(z)=
\phi_{x,x_h}(z),\\
\hbox{and}\qquad
d_{g_x}\big(z,\tilde{p}_{h,(x,w)}(z)\big)\le 2|w|_x
\qquad\forall z\!\in\! B_x\big(x_h,\de(x)\big).\notag
\end{gather}
In particular, both sides of \e_ref{si_metrics_l1a_e} are defined.
\end{lmm}

\noindent
{\it Proof:}
We choose $\de$ such that if $w\!\in\! T_x\aleph$
and $|w|\le 4\de(x)$, then $x(w)\!\in\!\aleph$ and 
the metric $g_{x(w)}$ is flat on $B_x\big(x_h,2\de(x)\big)$.
This choice of $\de$ insures that both sides 
of \e_ref{si_metrics_l1a_e} are defined.
Equation~\e_ref{si_metrics_l1a_e} is equivalent to
\begin{equation}
\label{si_metrics_l1a_e1}
\phi_{x(w),x_h(w)}\tilde{p}_{h,(x,w)}(z)
=d\exp_{g_x,x_h}\big|_{w_h}\phi_{x,x_h}(z)
=\phi_{x,x_h(w)}z+d\exp_{g_x,x_h}\big|_{w_h}w_h,
\end{equation}
since the metric $g_x$ is flat on $B_x\big(x_h,2\de(x)\big)$.
This equation defines the required map~$\tilde{p}_{h,(x,w)}$.
Since the metrics $g_x$ and $g_{x(w)}$ are flat on $B_x\big(x_h,2\de(x)\big)$,
the maps $\phi_{x,x_h(w)}z$ and $\phi_{x(w),x_h(w)}$
are holomorphic, and thus $\tilde{p}_{h,(x,w)}$ is holomorphic.
Taking the differential of~\e_ref{si_metrics_l1a_e1},
we obtain
\begin{equation}
\label{si_metrics_l1a_e2}
d\phi_{x(w),x_h(w)}\big|_{\tilde{p}_{h,(x,w)}(z)}
\circ d\tilde{p}_{h,(x,w)}\big|_z=
d\phi_{x,x_h(w)}\big|_z.
\end{equation}
Since $\phi_{x(w),x_h(w)}$ and $\phi_{x,x_h(w)}$
are $(g_{x(w)},g_{x(w)})$- and $(g_x,g_x)$-isometries,
respectively, on $B_x\big(x_h,2\de(x)\big)$,
it follows that $\tilde{p}_{h,(x,w)}$ is 
a $(g_x,g_{x(w)})$-isometry on $B_x\big(x_h,2\de(x)\big)$.
By \e_ref{si_metrics_l1a_e1},
\begin{equation}
\label{si_metrics_l1a_e3}
d_{g_x}\big(z,\tilde{p}_{h,(x,w)}(z)\big)
\le |w_h|_x+
\big|\big(\phi_{x(w),x_h(w)}-\phi_{x,x_h(w)}\big)
\tilde{p}_{h,(x,w)}(z)\big|_x
\le |w_h|_x+C(x)|w|\de(x),
\end{equation}
since the family of metrics is smooth.
If $C(x)\de(x)\!<\!1$, the remaining claim of the lemma 
follows from~\e_ref{si_metrics_l1a_e3}.

\begin{lmm}
\label{si_metrics_l1}
There exist $\de,C_k\!\in\! C^{\i}(\aleph;\Bbb{R}^+)$,
where $k$ is a positive integer, 
$\al_h\!\in\! C^{\i}(T\aleph_{\de};\Bbb{C})$,
and smooth families of maps
$$\big\{\Th_{w,h}\!:\{v\!\in\! T_{x_h}\S\!i: |v|_x\!<\!\de(x)\}
\lra T_{x_h}\Si~ |~ (x,w)\!\in\! T\aleph_{\de}\big\}$$
such that each map $\Th_{w,h}$ is holomorphic,
$\Th_{w,h}(0)=0$, $\Th_{w,h}'(0)=0$,
$\|\Th_{w,h}^{\lr{k}}\|_{C^0}\le C_k(x)|w|$, 
\hbox{$|\al_h(w)|\!\le\! C_0(x)|w|$}, and 
\begin{equation}\label{si_metrics_l1e1}
d\phi_{x,x_h}\big|_{x_h(w)}  d\phi_{x(w),x_h(w)}\big|_{x_h}
\big(\phi_{x(w),x_h}z\big)=
\big(1+\al_h(w)\big)\phi_{x,x_h}z+\Th_{w,h}\big(\phi_{x,x_h}z\big).
\end{equation}
for all  $z\!\in\! B_x\big(x_h,\de(x)\big)$.
In particular, both sides of \e_ref{si_metrics_l1e1} are defined.
\end{lmm}

\noindent
{\it Proof:} We choose $\de$ such that if $w\!\in\! T_x\aleph$
and $|w|\!\le\! 4\de(x)$, then $x(w)\!\in\!\aleph$ and 
the metric $g_{x(w)}$ is flat on~$B_x\big(x_h,4\de(x)\big)$.
This choice of $\de$ insures that both side of \e_ref{si_metrics_l1e1} 
are defined.
If $w$ and $z$ are as in the statement of the lemma,
by the flatness of the metric $g_{x(w)}$ near $x_h$,
$\Bbb{C}$-linearity of the differential of the exponential map near zero,
and the smoothness of the family of the metrics 
\begin{equation}\label{si_metrics_l1e2}
 d\phi_{x,x_h}\big|_{x_h(w)} d\phi_{x(w),x_h(w)}\big|_{x_h}
\big(\phi_{x(w),x_h}z\big)=\big(1+a_h(w)\big)\big(\phi_{x(w),x_h}z\big),
\end{equation}
for some $a_h\!\in\! C^{\i}(T\aleph_{\de};\Bbb{C})$ such that $a_h(0)=0$.
Note that if $g_{x(w)}=g_x$, $a_h(w)=0$, since the metric~$g_x$
is flat on~$B_x\big(x_h,|w|\big)$.
The map
$$\{ v\!\in\! T_{x_h}\Si\!: |v|_x<2\de(x)\}\lra T_{x_h}\Si,\quad
v\lra \phi_{x(w),x_h}\phi_{x,x_h}^{-1}v-v$$
is holomorphic since $\phi_{x(w),x_h}$ and $\phi_{x,x_h}$ are,
and vanishes at~$0$.
Thus,
\begin{equation}\label{si_metrics_l1e3}
\phi_{x(w),x_h}\phi_{x,x_h}^{-1}v=
\big(1+b_h(w)\big)v+\Th_{w,h}(v),
\end{equation}
for some $b_h(w)\!\in\!\Bbb{C}$ and holomorphic function
$\Th_{w,h}$ such that $\Th_{w,h}(0),\Th_{w,h}'(0)\!=\!0$.
Equation~\e_ref{si_metrics_l1e1} follows from 
\e_ref{si_metrics_l1e2} and \e_ref{si_metrics_l1e3}.
Smoothness of $b_h$ and $\Th_{\cdot,h}$ in $w$
follows from the smoothness of the family of the metrics.
The bounds on $\al_h$ and the derivatives of $\Th_{w,h}$
follow from their smoothness and compactness of the fibers of 
$$\big\{w\!\in\! T_x\aleph\!: |w|\le\de(x)\big\}\lra\aleph.$$

\begin{lmm}
\label{si_metrics_l2}
There exist $\de,C\!\in\! C^{\i}(\aleph;\Bbb{R}^+)$ 
and smooth families of maps
$$N_h\!: \l\{(x,w)\!: x\!\in\!\aleph;~ 
 (x,w)\!\in\! T\aleph_{\de}\r\}\lra T\aleph$$
such that $|N_h(w,w_h)|_x\le C(x)|w||w_h|$ and 
\begin{equation}\label{si_metrics_l2e1}
 d\phi_{x,x_h}\big|_{x_h(w)} 
\big(\phi_{x(w),x_h(w)}x_h\big)=-w_h+N_h(w,w_h).
\end{equation}
In particular, the left-hand side of \e_ref{si_metrics_l2e1} is defined.
\end{lmm}

\noindent
{\it Proof:} We take $\de$ as in Lemma~\ref{si_metrics_l1}.
Then, the left-hand side of \e_ref{si_metrics_l2e1} is defined and 
\begin{equation}\label{si_metrics_l2e2}\begin{split}
& d\phi_{x,x_h}\big|_{x_h(w)} \big(\phi_{x(w),x_h(w)}x_h\big)\\
&\qquad= d\phi_{x,x_h}\big|_{x_h(w)}\big(\phi_{x,x_h(w)}x_h\big)
+ d\phi_{x,x_h}\big|_{x_h(w)}\l\{
\l(\phi_{x(w),x_h(w)}x_h\r)-\l(\phi_{x,x_h(w)}x_h\r)\r\}\\
&\qquad=-w_h+N_h(w,w_h),
\end{split}\end{equation}
where $N(\cdot,\cdot)$ is some smooth function of both
variables, that vanishes if either input is zero.
Equation~\e_ref{si_metrics_l2e1} is thus proved,
while the bound on $N_h$ is obtained from
its smoothness and compactness of the fibers as in the 
proof of Lemma~\ref{si_metrics_l1}.

\begin{lmm}
\label{si_metrics_l3}
There exists $\de\!\in\! C^{\i}(\aleph;\Bbb{R}^+)$ such that
for any $x\!\in\!\aleph$, $v\!\in\! T_x\aleph$ with $|v|\!<\!\de(x)$,
and \hbox{$c\!=\!c_{[m]}\!\in\!\Bbb{C}^m$} with $|c||v|\!<\!\de(x)$,
there exists $w\!\in\! T\aleph$ with $|w_h|_x\!<\!2|c_h||v_h|_x$
such that for \hbox{$z\!\in\! B_x\big(x_h,4\de(x)^{\frac{1}{2}}\big)$}
\begin{equation}\label{si_metrics_l3e1}
 d\phi_{x,x_h}\big|_{x_h(w)}\big(\phi_{x(w),x_h(w)}z\big)=
\l(1+\al_h(w)\r)\l(c_hv_h+\phi_{x,x_h}z\r)
+\Th_{w,h}\l(\phi_{x,x_h}z\r),
\end{equation}
where $\al_h(w)$ and $\Th_{w,h}$ are as in Lemma~\ref{si_metrics_l1}.
In particular, both sides of \e_ref{si_metrics_l3e1} are defined.
\end{lmm}

\noindent
{\it Proof:}
We start by choosing $\de$ so that $8\de^{\frac{1}{2}}$ is smaller
that the function $\de$ of Lemmas~\ref{si_metrics_l1}
and~\ref{si_metrics_l2}.
By flatness of the metric $g_{x(w)}$ on 
$B\big(x_h,8\de(x)^{\frac{1}{2}}\big)$
for $w\!\in\! T_x\aleph$
with $|w|<\de(x)$
\begin{equation}\label{si_metrics_l3e2}
\phi_{x(w),x_h(w)}z=d\phi_{x(w),x_h(w)}\big|_{x_h}\phi_{x(w),x_h}z+
\phi_{x(w),x_h(w)}x_h
\end{equation}
for any $z\!\in\! B\big(x_h,4\de(x)^{\frac{1}{2}}\big)$.
Taking $d\phi_{x,x_h}\big|_{x_h(w)}$ of both sides of 
\e_ref{si_metrics_l3e2} and applying Lemmas~\ref{si_metrics_l1}
and \ref{si_metrics_l2}, we obtain
\begin{equation}\label{si_metrics_l3e3}
 d\phi_{x,x_h}\big|_{x_h(w)}\big(\phi_{x(w),x_h(w)}z\big)=
\big(1+\al_h(w)\big)\phi_{x,x_h}z+\Th_{w,h}\l(\phi_{x,x_h}z\r)
-w_h+N_h(w,w_h).
\end{equation}
Thus, we need to solve the equations
\begin{equation}\label{si_metrics_l3e4}
-w_h+N_h(w,w_h)=\l(1+\al_h(w)\r)c_hv_h.
\end{equation}
Let $\Psi_h(w)=N_h(w,w_h)-\l(1+\al_h(w)\r)c_hv_h$.
If $|w|\le 2|c_h||v_h|$, then by Lemmas~\ref{si_metrics_l1}
and \ref{si_metrics_l2},
\begin{equation}\label{si_metrics_l3e5}\begin{split}
&|\Psi(w)|\le C(x)|c||v|\l(2|c||v|+1\r)
\le 2|c||v|,\\
&|\Psi(w)-\Psi(w')|\le C(x)|c||v||w-w'|\le \frac{1}{2}|w-w'|,
\end{split}\end{equation}
provided $4C(x)\de(x)<1$.
In such a case, $\Psi$ is a contracting operator, and thus
\e_ref{si_metrics_l3e4} has a unique solution $w\!\in\! T_x\aleph$
with $|w|<2|c||v|$.
The estimate $|w_h|<2|c_h||v_h|$ follows directly 
from~\e_ref{si_metrics_l3e4} if $\de(x)$ is sufficiently small.

\begin{crl}
\label{si_metrics_c4}
There exist $\de,C_k\!\in\! C^{\i}(\aleph;\Bbb{R}^+)$,
where $k$ is a positive integer, such that
for any \hbox{$x\!\in\!\aleph$}, $v\!\in\! T_x\aleph$ with $|v|<\de(x)$,
$c\!=\!c_{[m]}\!\in\!\Bbb{C}^m$ with $|c|<\de(x)$, and 
$r\!=\!r_{[m]}\!\in\!\Bbb{R}^m$ with $|r|<\frac{1}{2}$,
there exists $\tilde{x}\!\in\!\aleph$ and 
$\tilde{v}\!\in\! T_{\tilde{x}}\aleph$
such that\\
(1) $\tilde{x}_h\!\in\! B_x\l(x_h,2|c_h||v_h|\r)$,
$\big|\frac{g_{\tilde{x}}}{g_x}-1\big|\le C_1(x)|c||v|$,
$\big||\tilde{v}_h|_{\tilde{x}}-|v_h|_x\big|\le 
C_1(x)\l(|c||v|+|r_h|\r)|v_h|$;\\
(2) for any $z\!\in\! B_x\l(x_h,4\de(x)^{1/2}\r)$,
\begin{equation}\label{si_metrics_c4e1}
\frac{\phi_{\tilde{x},\tilde{x}_h}z}{\tilde{v}_h}=
\l(1+r_h\r)\Big\{c_h+\frac{\phi_{x,x_h}z}{v_h}+
\Th_{v,c,r,h}\Big(\frac{\phi_{x,x_h}z}{v_h}\Big)\Big\},
\end{equation}
where $\Th_{v,c,r,h}$ is a holomorphic function,
varying smoothly with the parameters, such that
$$\Th_{v,c,r,h}(0)=0,~~~\Th_{v,c,r,h}'(0)=0,~~~
\big\|\Th_{v,c,r,h}^{\lr{k}}\big\|_{C^0}\le C_k(x)|c||v||v_h|^{k-1}.$$
\end{crl}

\noindent
{\it Proof:}
Let $\de$ be as in Lemma~\ref{si_metrics_l3}.
Given $v$ and $c$ as in the statement of the lemma,
let $w\!\in\! T_x\aleph$ be the element provided by Lemma~\ref{si_metrics_l3}.
Take
$$\tilde{x}_h=x_h(w),\qquad
\tilde{v}_h=(1+r_h)^{-1}\big(1+a_h(w)\big)
d\phi_{x,x_h}^{-1}|_{w_h}v_h.$$
The estimates in (1) are immediate from Lemma~\ref{si_metrics_l3},
provided $\de$ is sufficiently small.
The inequalities in (2) arise from the smooth dependence of $w$
on $x$, $v$, and $c$ in Lemma~\ref{si_metrics_l3},
and the fact that $w$ is zero if either $v=0$ or $c=0$.

\subsection{Sobolev Inequalities for the Metrics $g_{\ups}$}
\label{sob_main}

\noindent
In this subsection, we prove (3) of Lemma~\ref{approx_maps}.
The reason this estimate holds is that $(\Si_{\ups},g_{\ups})$
can be written as a union of the surfaces
$(\Si_{b_{\ups},i},g_{\ups})$ with small disks missing 
and annuli $(\tilde{A}_{\ups,h}^{\pm},g_{\ups})$
that are uniformly equivalent to annuli in~$\Bbb{R}^2$
with the smaller radius less than half of the~larger.\\

\noindent
Suppose ${\cal T}\!=\!(S,M,I;j,\la)$ 
is a bubble type and 
$$\ups=\big(b,v_{\hat{I}}\big)=
\big((S,M,I;x,(j,y),u),v_{\hat{I}}\big)
\in F^{(0)}{\cal T}_{\de}.$$
For any $h\!\in\!\hat{I}$ and $i\!\in\!I$, put
\begin{gather*}
\tilde{A}_{\ups,h}^-=q_{\ups,i_h}^{-1}\Big(\big\{
(\io_h,z)\!\in\!\Si_{b_{\ups},\io_h}\!: 
(2\de_{\cal T}(b_{\ups}))^{-1}|v_h|_b\le r_{b,h}(z)\le 
        |v_h|_b^{\frac{1}{2}}\big\}\Big);\\
\tilde{A}_{\ups,h}^+=q_{\ups,i_h}^{-1}\Big(\big\{
(\io_h,z)\!\in\!\Si_{b_{\ups},\io_h}\!: 
|v_h|_b^{\frac{1}{2}}\le r_{b,h}(z)\le 2\de_{\cal T}(b_{\ups})\big\}\Big);\\
S_{\ups,i}=q_{\ups,i}^{-1}\Big(\big\{(i,z)\!\in\! S_{b_{\ups},i}\!: 
    r_{b,h}(z)\ge\de_{\cal T}(b_{\ups}) \hbox{~if~}\io_h=i;~
   |q_S^{-1}(z)|\ge\de_{\cal T}(b_{\ups}) \hbox{~if~}i>0\big\}\Big).
\end{gather*}
Let $\tilde{A}_{\ups,h}$ denote $\tilde{A}_{\ups,h}^-\cup\tilde{A}_{\ups,h}^+$.

\begin{lmm}
\label{c0_bound_l} For any $p\!>\!2$,
there exists $C_p\!\in\! C^{\i}({\cal M}^{(0)}_{\cal T};\Bbb{R})$ 
such that for any $\ups\!\in\! F^{(0)}{\cal T}_{\de}$ 
\hbox{and $h\!\in\!\hat{I}$},
$$\xi\in\tilde{\Ga}_c(\tilde{A}_{\ups,h};u_{\ups}^*TV)\Lra
\|\xi\|_{C^0}\le C_p(b_{\ups})\|\xi\|_{g_{\ups},p,1}.$$
\end{lmm}

\noindent
{\it Proof:} By construction of the metric $g_{\ups}$,
$g_{\ups}|\tilde{A}_{\ups,h}$ is the pullback of the metric $g_{\ups,\io_h}$
on $q_{\ups,\io_h}(\tilde{A}_{\ups,h})$ by the map~$q_{\ups,\io_h}$.
Furthermore, the metric~$g_{\ups,\io_h}$ on $q_{\ups}(\tilde{A}_{\ups,h}^{\pm})$
differs from the standard metric on
the annulus $B_{2\de_{\cal T}(b_{\ups}),|v_h|^{\frac{1}{2}}}\subset\Bbb{R}^2$ 
by factors bounded by~$C(b_{\ups})$.
Since $\|du_{\ups}\|_{g_{\ups},p}\le C_p(b_{\ups})$ by 
(1) of Lemma~\ref{approx_maps}, 
the claim follows from Proposition~\ref{an-c0bound_prp}
in~\cite{Z1}.

\begin{prp}
\label{c0_bound}
For any $p>2$,
there exists $C_p\!\in\! C^{\i}({\cal M}^{(0)}_{\cal T};\Bbb{R})$ 
such that for all 
\hbox{$\ups\!\in\! F^{(0)}{\cal T}_{\de}$},
$$\|\xi\|_{C^0}\le C_p(b_{\ups})\|\xi\|_{g_{\ups},p,1}
\qquad\hbox{for all~~}\xi\in\Ga(u_\ups).$$
\end{prp}

\noindent
{\it Proof:} (1) Note that $g_{\ups}|S_{\ups,i}$ is the pull-back of
the metric $g_{b_{\ups},i}$ on $q_{\ups,i}(S_{\ups,i})$ by the map
$q_{\ups,i}$. Thus, by Proposition~\ref{an-c0bound_prp} in~\cite{Z1}, if
$\xi\!\in\!\Ga_c(S_{\ups,i};u_{\ups}^*TV)$,
$$\|\xi\|_{C^0}=\|\xi\circ q_{\ups,i}\|_{C^0}
\le C_p(\|du_{\ups}\circ q_{\ups,i}\|_{g_{b_{\ups},i},p})
\|\xi\circ q_{\ups,i}\|_{g_{b_{\ups},i},p,1}
= C_p(b_{\ups})\|\xi\|_{g_{\ups},p,1},$$
since $\xi$ vanishes outside of $S_{\ups,i}$.\\
(2)  We now define a partition of unity subordinate to
$\{S_{\ups,i},\tilde{A}_{\ups,h}\!: i\!\in\! I, h\!\in\!\hat{I}\}$. 
Put
\begin{gather*}
\eta_{\ups,h}^+(z)=
\begin{cases}
1-\be_{\de_{\cal T}^2(b_{\ups})}\big(r_{b_{\ups},h}(q_{\ups,\io_h}(z))\big),&
\hbox{if~}q_{\ups,\io_h}(z)\in\Si_{b_{\ups},\io_h};\\
1,&\hbox{otherwise};
\end{cases}\\
\eta_{\ups,h}^-(z)=
\begin{cases} 
1-\be_{\de_{\cal T}^2(b_{\ups})}\big(\big|q_S^{-1}q_{\ups}(z)\big|\big),&
\hbox{if~}q_{\ups,h}(z)\in\Si_{b_{\ups},h};\\
1,&\hbox{otherwise};
\end{cases}\quad
\tilde{\eta}_{\ups}(z)=1-\prod_{h\in\hat{I}}\eta_{\ups,h}^-(z)\eta_{\ups,h}^+(z).
\end{gather*}
Note that $d\eta_{\ups,h}^{\pm}$ is supported in $\tilde{A}_{\ups,h}^{\pm}$.
It follows from the definition of $g_{\ups}$ that
$$\|d\eta_{\ups,h}^{\pm}\|_{g_{\ups},C^1}=
\|d(\eta_{\ups,h}^{\pm}\circ q_{\ups,i_h}^{-1})\|_{g_{\ups,i_h},C^1}
\le C(b_{\ups}).$$
Thus, if $\xi\!\in\!\Ga(u_{\ups})$ by (1) and Lemma~\ref{c0_bound_l},
\begin{equation*}\begin{split}
\|\xi\|_{C^0}&\le\sum_{i\in I}\|\tilde{\eta}\xi\|_{C^0(S_{\ups,i})}+
\sum_{h\in\hat{I}}\|\eta_{\ups,h}^-\eta_{\ups,h}^+\xi\|_{C^0}
\le C_p(b_{\ups})\Big(\|\tilde{\eta}\xi\|_{g_{\ups},p,1}+
\sum_{h\in\hat{I}}\|\eta_{\ups,h}^-\eta_{\ups,h}^+\xi\|_{g_{\ups},p,1}\Big)\\
&\le C_p(b_{\ups})\Big(|I|\|\xi\|_{g_{\ups},p,1}+
2\sum_{h\in\hat{I}}\|\eta_{\ups,h}^-\eta_{\ups,h}^+\|_{g_{\ups},C^1}
\|\xi\|_{g_{\ups},p}\Big)
\le C'_p(b_{\ups})\|\xi\|_{g_{\ups},p,1}.
\end{split}\end{equation*}

\subsection{Elliptic Estimates for the Metrics $g_{\ups}$}

\noindent
This subsection contains the proof of (4) of Lemma~\ref{approx_maps},
the main elliptic estimate  
for the operators~$D_{\ups}$ and the modified Sobolev norms.
This estimate does {\it not} hold for the standard Sobolev norms.
The argument is essentially the same as in~\cite{LT},
but we do include all of the details, based on~\cite{Z1},
and state a sharper estimate.\\

\noindent
Let ${\cal T}$, $\ups$, 
$\tilde{A}_{\ups,h}=\tilde{A}_{\ups,h}^-\cup\tilde{A}_{\ups,h}^+$,
and $S_{\ups,i}$ be as in Subsection~\ref{sob_main}.
If $\io_h\!=\!0$, the metric $g_{b_{\ups},\hat{0}}$ is flat on 
$B_{b_{\ups},h}(\de_{\cal T}(b_{\ups})^{\frac{1}{2}})$.
Thus, for any $h\!\in\!\hat{I}$, we can choose conformal
polar coordinates $(r,\th)$ on $\tilde{A}_{\ups,h}$ such that
\hbox{$r(z)=r_{b,\ups}\big(q_{\ups,\io_h}(z)\big)$}.
Since $g_{\ups}|\tilde{A}_{\ups,h}$ is the pullback of the metric 
$g_{\ups,\io_h}$ on $q_{\ups,\io_h}(\tilde{A}_{\ups,h})$ by 
the map~$q_{\ups,\io_h}$, 
\begin{equation}
\label{elliptic_main_e1}
g_{\ups}=\Big\{  
 \big(1-\be_{|v_h|}(2r)\big)\frac{2C(b_{\ups})}{|v_h|+|v_h|^{-1}r^2}
             +\be_{|v_h|}(r)\Big\}\big(dr^2+r^2d\th^2\big)
\qquad\hbox{on}~~\tilde{A}_{\ups,h}.
\end{equation}
Similarly, since $\rho_{\ups}=\rho_{\ups,\io_h}\circ q_{\ups,\io_h}$
\begin{equation}
\label{elliptic_main_e2}
\rho_{\ups}=r^2+\frac{|v_h|^2}{r^2}
\qquad\hbox{on}~~\tilde{A}_{\ups,h}.
\end{equation}

\begin{lmm}
\label{elli_main1}
For all $p>1$, there exists 
$C_p\!\in\! C^{\i}({\cal M}^{(0)}_{\cal T};\Bbb{R})$ 
such that for all $\ups\!\in\! F^{(0)}{\cal T}_{\de}$ 
\hbox{and $h\!\in\!\hat{I}$}, 
$$\xi\in\Ga_c(\tilde{A}_{\ups,h};u_{\ups}^*TV)\Lra
\Big(\int_{\tilde{A}_{\ups,h}}\rho_{\ups}^{-\frac{p-2}{p}}|\na^{b_{\ups}}\xi|^2
                    \Big)^{\frac{1}{2}} 
\le C_p(b_{\ups})\big(\|D_{\ups}\xi\|_{\ups,p}+\|\xi\|_{\ups,p}\big).$$
\end{lmm}

\noindent
{\it Proof:} (1)  Let $\ep_1$ and $\ep_2$ denote 
$\l(2\de_{\cal T}(b_{\ups})\r)^{-1}|v_h|$ and 
$2\de_{\cal T}(b_{\ups})$, respectively.
Note that the integral on the left-hand side in the statement
of the lemma is conformally invariant. 
With respect to the metric $dr^2+r^2d\th^2$,
$$\big|D_{\ups}\xi\big|_{(r,\th)}=
 \Big|\frac{D}{dr}\xi+Jr^{-1}\frac{D}{d\th}\xi\Big|_{(r,\th)},$$
where $\frac{D}{dr}$ and $\frac{D}{d\th}$ denote covariant 
differentiation with respect to the connection~$\na^{b_{\ups}}$
and the norms are taken with respect to the metric~$g_{V,b}$
on~$V$.
Thus,
\begin{equation}\label{elliptic_main_e3}
a_h^2\equiv
\int_{\tilde{A}_{\ups,h}}\rho_{\ups}^{-\frac{p-2}{p}}|\na^{b_{\ups}}\xi|^2
\le \|D_{\ups}\xi\|_{\ups,p}^2
-2\int_0^{2\pi}\!\!\!\int_{\ep_1}^{\ep_2}\rho_{\ups}^{-\frac{p-2}{p}}
\Big\lan\frac{D}{dr}\xi,J\frac{D}{d\th}\xi\Big\ran  drd\th.
\end{equation}
Since $\na^{b_{\ups}} J=0$, using integration by parts twice, we obtain
\begin{gather}
\int_0^{2\pi}\!\!\!\int_{\ep_1}^{\ep_2}\rho_{\ups}^{-\frac{p-2}{p}} 
        \Big\lan\frac{D}{dr}\xi,J\frac{D}{d\th}\xi\Big\ran drd\th
= -\int_0^{2\pi}\!\!\!\int_{\ep_1}^{\ep_2}\rho_{\ups}^{-\frac{p-2}{p}}
\Big\lan\frac{D}{d\th}\frac{D}{dr}\xi,J\xi\Big\ran  drd\th\notag\\
\label{elliptic_main_e3a}
=-\int_0^{2\pi}\!\!\!\int_{\ep_1}^{\ep_2}\rho_{\ups}^{-\frac{p-2}{p}}
\Big(\Big\lan\frac{D}{dr}\frac{D}{d\th}\xi,J\xi\Big\ran-
   \big\lan {\cal R}(u_r,u_{\th})\xi,J\xi\big\ran\Big)drd\th\\
=\int_0^{2\pi}\!\!\!\int_{\ep_1}^{\ep_2}\rho_{\ups}^{-\frac{p-2}{p}}
\Big(\Big\lan\frac{D}{d\th}\xi,J\frac{D}{dr}\xi\Big\ran-
\frac{(p-2)}{p}\frac{\rho_{\ups}'(r)}{\rho_{\ups}(r)}
       \Big\lan\frac{D}{d\th}\xi,J\xi\Big\ran
+\big\lan {\cal R}_{b_{\ups}}(u_r,u_{\th})\xi,J\xi\big\ran\Big)drd\th\notag,
\end{gather}
where $u_r$ and $u_{\th}$ denote $\frac{d}{dr}u_{\ups}$
and $\frac{d}{d\th}u_{\ups}$, respectively,
and ${\cal R}_{b_{\ups}}$ is the curvature tensor of the 
connection~$\na^{b_{\ups}}$.
Since \hbox{$\l\lan\frac{D}{d\th}\xi,J\frac{D}{dr}\xi\r\ran=-
\l\lan\frac{D}{dr}\xi,J\frac{D}{d\th}\xi\r\ran$},
by \e_ref{elliptic_main_e3a}
and (1) of Lemma~\ref{approx_maps},
\begin{equation}\label{elliptic_main_e4}\begin{split}
&\Big|\int_0^{2\pi}\!\!\int_{\ep_1}^{\ep_2}\rho_{\ups}^{-\frac{p-2}{p}}
\Big\lan\frac{D}{dr}\xi,J\frac{D}{d\th}\xi\Big\ran drdr\th\Big|\\ 
&\qquad\qquad
\le\frac{|p-2|}{2p}\Big|\int_{\ep_1}^{\ep_2}\rho_{\ups}^{-\frac{p-2}{p}}
\frac{\rho_{\ups}'(r)}{\rho_{\ups}(r)}
\int_0^{2\pi}\Big\lan\frac{D}{d\th}\xi,J\xi\Big\ran d\th dr\Big|+
C(b_{\ups})\|\xi\|_{\ups,p}^2.
\end{split}\end{equation}
(2) By Poincare Lemma, see Proposition~\ref{an-poincare} in~\cite{Z1}, 
for every circle with $r$ fixed,
\begin{equation}\label{elliptic_main_e5}\begin{split}
\Big|\int_0^{2\pi} \Big\lan\frac{D}{d\th}\xi,J\xi\Big\ran d\th\Big|
\le& \int_0^{2\pi}\big|\frac{D}{d\th}\xi\big|^2 d\th
+C(g_{b_{\ups}})\bigg\{    \Big(\int_0^{2\pi}|u_{\th}|^2d\th\Big)
   \Big(\int_0^{2\pi}|\xi|^2d\th\Big)\\
&\qquad\quad+\Big(\int_0^{2\pi}|u_{\th}|d\th\Big)
     \Big(\int_0^{2\pi}|\xi|^2d\th\Big)^{\frac{1}{2}}
\Big(\int_0^{2\pi}\Big|\frac{D}{d\th}\xi\Big|^2 
                    d\th\Big)^{\frac{1}{2}}\bigg\}.
\end{split}\end{equation}
Since $\l|\frac{\rho_{\ups}'(r)}{\rho_{\ups}(r)}\r|\le 2r^{-1}$ on 
$\tilde{A}_{\ups,h}$,
by Holder's inequality and the first part of Lemma~\ref{approx_maps},
\begin{equation}\label{elliptic_main_e6}\begin{split}
\frac{1}{2}\int_{\ep_1}^{\ep_2}\rho_{\ups}^{-\frac{p-2}{p}}
\Big|\frac{r\rho_{\ups}'(r)}{\rho_{\ups}(r)}\Big|
\Big(\int_0^{2\pi}\!r^{-1}|u_{\th}|d\th\Big)
  \Big(\int_0^{2\pi}\!|\xi|^2d\th\Big)^{\frac{1}{2}}
\Big(\int_0^{2\pi}\!
  r^{-2}\Big|\frac{D}{d\th}\xi\Big|^2 d\th\Big)^{\frac{1}{2}}
     rdr&\\
\le C\|\xi\|_{\ups,p}\Big(\int_{\tilde{A}_{\ups,h}}\rho_{\ups}^{-\frac{p-2}{p}}
r^{-2}\Big|\frac{D}{d\th}\xi\Big|^2\Big)^{\frac{1}{2}}.\quad&
\end{split}\end{equation}
Similarly,
\begin{equation}\label{elliptic_main_e7}\begin{split}
\frac{1}{2}\int_{\ep_1}^{\ep_2}
\rho_{\ups}^{-\frac{p-2}{p}}\l|\frac{r\rho_{\ups}'(r)}{\rho_{\ups}(r)}\r|
\l(\int_0^{2\pi}\!r^{-2}|u_{\th}|^2d\th\r)
   \l(\int_0^{2\pi}\!|\xi|^2d\th\r)rdrd\th
\le C(b_{\ups})\|\xi\|_{\ups,p}^2.
\end{split}\end{equation}
Combining equations~\e_ref{elliptic_main_e5}-\e_ref{elliptic_main_e7},
we obtain
\begin{equation}\label{elliptic_main_e8}\begin{split}
\frac{1}{2}&\Big|\int_{\ep_1}^{\ep_2}\rho_{\ups}^{-\frac{p-2}{p}}
\frac{\rho_{\ups}'(r)}{\rho_{\ups}(r)}
\int_0^{2\pi}\Big\lan\frac{D}{d\th}\xi,J\xi\Big\ran d\th dr\Big|\\
&\qquad\qquad\qquad\quad
\le \int_0^{2\pi}\!\!\!\int_{\ep_1}^{\ep_2}\rho_{\ups}^{-\frac{p-2}{p}}
r^{-2}\Big|\frac{D}{d\th}\xi\Big|^2rdrd\th
+C(b_{\ups})\l(\|\xi\|_{\ups,p}^2+\|\xi\|_pa_h\r).
\end{split}\end{equation}
Note that
\begin{equation}\label{elliptic_main_e9}\begin{split}
&\int_0^{2\pi}\!\!\!\int_{\ep_1}^{\ep_2}\rho_{\ups}^{-\frac{p-2}{p}}
r^{-2}\Big|\frac{D}{d\th}\xi\Big|^2rdrd\th\\
&\quad=\frac{1}{2}\int_{\tilde{A}_{\ups,h}}\rho_{\ups}^{-\frac{p-2}{p}}
\Big(r^{-2}\Big|\frac{D}{d\th}\xi\Big|^2+
\Big|\Big(\frac{D}{dr}\xi+Jr^{-1}\frac{D}{d\th}\xi\Big)-
               \frac{D}{dr}\xi\Big|^2\Big)
\le \frac{1+\ep}{2}a_h^2+C_{\ep}\|D_{\ups}\xi\|_{\ups,p}^2
\end{split}\end{equation}
for any $\ep>0$.
Combining equations~\e_ref{elliptic_main_e3},\e_ref{elliptic_main_e4},
\e_ref{elliptic_main_e8} and \e_ref{elliptic_main_e9}, we obtain
$$a_h^2\le\frac{|p-2|}{p}(1+\ep)a_h^2+  
\big(C(b_{\ups})+C_{\ep}\big)\big(
\|D_{\ups}\xi\|_{\ups,p}^2+\|\xi\|_{\ups,p}^2+\|\xi\|_{\ups,p}a_h\big).$$
Since $\frac{|p-2|}{p}<1$, the claim follows by choosing $\ep$
sufficiently small.

\begin{lmm}
\label{elli_main2}
For all $p\ge 1$, there exists 
$C_p\!\in\! C^{\i}({\cal M}^{(0)}_{\cal T};\Bbb{R})$ 
such that for all $\ups\!\in\! F^{(0)}{\cal T}_{\de}$ 
\hbox{and $h\!\in\!\hat{I}$}, 
$$\xi\in\Ga_c(\tilde{A}_{\ups,h}|u_{\ups}^*TV)\Lra
\|\na^{b_{\ups}}\xi\|_{g_{\ups},p} 
\le C_p(b_{\ups})\Big( \|D_{\ups}\xi\|_{\ups,p}+\|\xi\|_{\ups,p}+
\Big(\int_{\tilde{A}_{\ups,h}}\rho_{\ups}^{-\frac{p-2}{2p}}|\na\xi|^2
            \Big)^{\frac{1}{2}} \Big).$$
\end{lmm}

\noindent
{\it Proof:} Choose a sequence 
$$\de_0>\ldots>\de_{N+1}>0\hbox{~~~such that~~~}
\de_0=\ep_2,~\de_{N+1}=\ep_1,~~
\frac{1}{3}\le\frac{\de_{l+1}}{\de_l}\le\frac{1}{2}.$$
For each $l=1,\ldots,N-1$, let $g_l$ denote the metric 
$$g_l=\Big(\de_l^2+\frac{|v_h|^2}{\de_l^2}\Big)^{-1}g_{\ups}
\hbox{~~on~~} 
\tilde{A}_l=\big\{(r,\th)\!\in\!\tilde{A}_{\ups,h}\!: 
       \de_{l+2}\le r\le\de_{l-1} \big\}.$$
Let $\rho_l=\de_l^2+|v_h|^2\de_l^{-2}$ and denote by $A_l$ the annulus 
$\{(r,\th)\!\in\!\tilde{A}_{\ups,h}\!: \de_{l+1}\le r\le\de_l\}$.
The pullback of the metric $g_l$ on $\tilde{A}_l$  to the annulus
$(\frac{\de_{l+2}}{\de_l},\frac{\de_{l-1}}{\de_l})\times S^1\subset\Bbb{R}^2$
by the map $(r,\th)\lra(\de_lr,\th)$ differs from the Eucledian
metric by a conformal factor bounded by $C(b_{\ups})$, since
$$\frac{1}{100}\le 
\Big\{\l(1-\be_{|v_h|}(\de_lr)\r)
\frac{2}{|v_h|+|v_h|^{-1}\de_l^2r^2}+\be_{|v_h|}(\de_lr) \Big\}
\Big(\de_l^2+\frac{|v_h|^2}{\de_l^2}\Big)^{-\frac{1}{2}}\de_l
\le 200,$$
whenever $r\!\in\!(\frac{1}{9},3)$ and $\de_l\!\in\!(|v_h|,1)$.
Thus, by Proposition~\ref{an-elli_prp1} in~\cite{Z1},  
\begin{equation}\label{elliptic_main_e10}
\|\na^{b_{\ups}}\xi\|_{g_l,L^p(A_l)}\le 
C\Big(\|D_{\ups}\xi\|_{g_l,L^p(\tilde{A}_l)}+
\|\na^{b_{\ups}}\xi\|_{g_l,L^p(\tilde{A}_l)}+
\|\xi du\|_{g_l,L^p(\tilde{A}_l)}\Big),
\end{equation}
or equivalently
\begin{equation}\label{elliptic_main_e11}
\|\na^{b_{\ups}}\xi\|_{g_{\ups},L^p(A_l)}\le 
C\Big(\|D_{\ups}\xi\|_{g_{\ups},L^p(\tilde{A}_l)}+
\|\rho_l^{-\frac{p-2}{2p}}\na^{b_{\ups}}\xi\|_{L^2(\tilde{A}_l)}+
\|\xi du\|_{g_{\ups},L^p(\tilde{A}_l)}\Big).
\end{equation}
Since $\frac{\rho_{\ups}(r)}{\rho_l}\in[\frac{1}{81},81]$ 
when $r\!\in\![\de_{l+2},\de_{l-1}]$,
\e_ref{elliptic_main_e11} is equivalent to
$$\Big(\int_{A_l}|\na^{b_{\ups}}\xi|^p\Big)^{\frac{1}{p}} \le C_p(b_{\ups})
 \Big( \Big(\int_{\tilde{A}_l}|D_{\ups}\xi|^p\Big)^{\frac{1}{p}}+
\Big(\int_{\tilde{A}_l}\rho_{\ups}^{-\frac{p-2}{p}}
           |\na^{b_{\ups}}\xi|^2\Big)^{\frac{1}{2}}
+\|\xi du_{\ups}\|_{g_{\ups},L^p(\tilde{A}_l)}\Big).$$
The claim follows by
summing up the last inequality over all~$l$
and using (1) of Lemma~\ref{approx_maps}.\\

\noindent
{\it Remark:} The above proof does not quite apply to 
the two outermost annuli~$A_1$ and~$A_N$.
However, since $\xi$ is compactly supported in $\tilde{A}_{\ups,h}$,
the proof of Proposition~\ref{an-elli_prp1} in~\cite{Z1} can be applied
to $A_1$ with $A_1\cup A_2$ replacing~$\tilde{A}_1$
to~\e_ref{elliptic_main_e10},
and similarly for~$A_N$.
Alternatively, for the purposes of proving Proposition~\ref{elliptic_bound}
below, it is sufficient to prove Lemma~\ref{elli_main2} and
Corollary~\ref{elli_main3} for $\xi$ that vanish on $A_1$ and~$A_N$.

\begin{crl}
\label{elli_main3}
For all $p> 1$, there exists 
$C_p\!\in\! C^{\i}( {\cal M}^{(0)}_{\cal T};\Bbb{R})$ 
such that for all $\ups\!\in\! F^{(0)}{\cal T}_{\de}$ \hbox{and $h\in\hat{I}$}
$$\xi\in\Ga_c(\tilde{A}_{\ups,h};u^*TV)\Lra
\|\xi\|_{\ups,p,1}\le C_p(b_{\ups})
\big(\|D_{\ups}\xi\|_{\ups,p}+\|\xi\|_{\ups,p}\big).$$
\end{crl}

\noindent
{\it Proof:} This corollary follows immediately from 
Lemmas~\ref{elli_main3} and \ref{elli_main3}.

\begin{prp}
\label{elliptic_bound}
For all $p>1$, there exists 
$C_p\!\in\! C^{\i}({\cal M}^{(0)}_{\cal T};\Bbb{R})$ 
such that for all $\ups\!\in\! F^{(0)}{\cal T}_{\de}$,
$$\|\xi\|_{\ups,p,1}\le C_p(b_{\ups})
\big(\|D_{\ups}\xi\|_{\ups,p}+\|\xi\|_{\ups,p}\big)
\qquad\hbox{for all}~~\xi\in\Ga(u_\ups).$$
\end{prp}

\noindent
{\it Proof:} This proposition follows from
Corollary~\ref{elli_main3}
and Proposition~\ref{an-elli_prp2} in~\cite{Z1}
by taking a partition of unity as in the proof of 
Proposition~\ref{c0_bound}.

\subsection{Fiber-Uniform Inverse for the Operator $D_{\ups}$}

\begin{lmm}
\label{inver1}
Let $\{\ups_k\}$ be a sequence in $F^{(0)}{\cal T}_{\de}$ that converges
to $b^*\!\in\!{\cal M}_{\cal T}^{(0)}$. 
Suppose $\xi_k\!\in\!\Ga(u_{\ups_k})$ is such that 
      $\|\xi_k\|_{\ups_k,p,1}\le 1$ for all $k$, 
while $\|D_{\ups_k}\xi_k\|_{\ups_k,p}\lra 0$ as $k\lra\i$
for some $p>2$.
Then a subsequence of $\{\xi_k\}$ $C^0$-converges $\xi^*\!\in\!\Ga_-(b^*)$.
Furthermore, $\|\xi_k\|_{\ups_k,C^0}$ converges to  $\|\xi^*\|_{b^*,C^0}$.
\end{lmm}

\noindent
 {\it Proof:} (1) Write 
$b^*=\big(S,M,I;x^*,(j,y^*),u_I\big)$ and 
$$\ups_k=(b_k,v_k)=
\big((S,M,I;x_k,(j,y_k),u_k),(v_k)_{\hat{I}}\big).$$
For each $i\!\in\! I$ and $\de>0$, put
$$S_{i,\de}^*=\{z\!\in\!\Si_{b,i}\!:  
r_{b^*,h}(z)\ge\de ~\forall h\in\hat{I};  
|q_S^{-1}(z)|\ge\de \hbox{~if~}i>0\}.$$
For $i\!\in\! I$ and all $k$ sufficiently large
(depending on $\de$), define 
$\ze_{k,i},\xi_{k,i}'\!\in\!\Ga(u_i^*|S_{i,\de}^*)$ by
$$\exp_{b^*,u_i^*(z)}\ze_{k,i}(z)=u_{\ups_k}(q_{\ups_k}^{-1}(z)),~
\|\ze_{k,i}\|_{b^*,C^0}<\inj~g_{V,b^*};\quad
\Pi_{b^*,\ze_{k,i}(z)}\xi_{k,i}'(z)=\xi_k(q_{\ups_k}^{-1}(z)).$$
Since $\|\na^{b^*}\ze_{k,i}\|_{b^*,C^0}\le C$
for $k$ sufficiently large, (1) of Lemma~\ref{approx_maps} and 
by Corollary~\ref{an-pt_l1c} in~\cite{Z1},
\begin{equation}\label{inver1_e1}\begin{split}
&\|\xi_{k,i}'\|_{b^*,p,1}\le (1+\ep_k)\|\xi_k\|_{\ups_k,p,1}
+\ep_k\|\xi_k\|_{\ups_k,C^0},\\
&\|D_{b,u_i}\xi_{k,i}'\|_{b^*,p}
\le(1+\ep_k)\|D_{\ups_k}\xi_k\|_{\ups_k,p}+
\ep_k\|\xi_k\|_{\ups_k,C^0},
\end{split}\end{equation}
where $\ep_k\!\lra\! 0$ as $k\!\lra\!\i$. 
Since $\|\xi_k\|_{\ups_k,p,1}\!\le\! 1$, (2) of Lemma~\ref{approx_maps}
applied to~\e_ref{inver1_e1},
\begin{equation}\label{inver1_e2}
\|\xi_{k,i}'\|_{b^*,p,1}\le (1+\tilde{\ep}_k)\|\xi_k\|_{\ups_k,p,1}
+\tilde{\ep}_k,~~~
\|D_{b,u_i}\xi_{k,i}'\|_{b^*,p}
\le(1+\tilde{\ep}_k)\|D_{\ups_k}\xi_k\|_{\ups_k,p}+\tilde{\ep}_k,
\end{equation}
where $\tilde{\ep}_k\!\lra\! 0$ as $k\!\lra\!\i$. 
Sobolev's embedding theorem then implies that $\xi_{k,i}'$ converges
to a vector field $\xi_i^*\!\in\!\Ga(u_i|\Si_{b^*,i}^*)$ in the $C^0$-norm 
on the compact subsets of $\Si_{b^*,i}^*$. 
Furthermore,  \hbox{$\|\xi_i^*\|_{b^*,C^0}<\i$}, since
$$\|\xi_{k,i}'\|_{b^*,C^0}\le (1+\ep_k)\|\xi_k\|_{\ups_k,C^0}\le C.$$
(2) We will now show that $D_{b^*,u_i^*}\xi_i^*\!=\!0$ weakly, 
i.e.~$\llan\xi_i^*,D_{b^*,u_i^*}^*\eta\rran_{b^*}\!=\!0$ 
for any $\eta\!\in\!\Ga^{0,1}(u_i^*)$.
We~have
\begin{equation}\label{inver1_e3}
\llan\xi_i^*,D_{b^*,u_i^*}^*\eta\rran_{b^*}
=\lim_{\de\lra 0}\int_{S_{i,\de}^*}
        \lr{\xi_i^*,D_{b^*,u_i^*}^*\eta}_{b^*}
=\lim_{\de\lra 0}\lim_{k\lra\i}
\int_{S_{i,\de}^*}\lr{\xi_{k,i}',D_{b^*,u_i^*}^*\eta}_{b^*},
\end{equation}
since $\xi_{k,i}'\lra\xi_i^*$ in the $C^0$-norm on $S_{i,\de}$.
By integration by parts,
\begin{equation}\label{inver1_e4}
\Big|\int_{S_{i,\de}^*}\lr{\xi_{k,i}',D_{b^*,u_i^*}^*\eta}_{b^*}-
\int_{S_{i,\de}^*}\lr{D_{b^*,u_i^*}\xi_{k,i}',\eta}_{b^*}\Big|\le
C\int_{\partial S_{i,\de}^*}|\xi_{k,i}'||\eta|\le
C'\|\xi_{k,i}'\|_{b^*,C^0}\|\eta\|_{b^*,C^0}\de.
\end{equation}
Since $\|D_{b^*,u_i*}\xi_{k,i}'\|_{b^*,p}\!\lra\! 0$ as
$k\!\lra\!\i$ on $S_{i,\de}^*$ and $\|\xi_{k,i}'\|_{b^*,C^0}\le C$, by
\e_ref{inver1_e3} and \e_ref{inver1_e4},
$$\llan\xi_i^*,D_{b^*,u_i^*}^*\eta\rran_{b^*}=0\qquad
    \forall\eta\!\in\!\Ga^{0,1}(u_i^*).$$
(3) Since $D_{b^*,u_i^*}\xi_i^*=0$ weakly on $S_i$ and
$D_{b^*,u_i^*}$ is an elliptic operator, 
it follows that $\xi_i^*$ is smooth and $D_{b^*,u_i^*}\xi_i^*=0$.
It will now be shown that $\xi_{\io_h}^*(x_h^*)=\xi_h^*(\i)$
for all $h\in\hat{I}$, i.e. 
\hbox{$\xi^*\!\equiv\! \xi^*_I\!\in\!\Ga(b^*)$}.
For each $h\!\in\!\hat{I}$,
let $A_{h,\de,k}\!\subset\! S$ denote the small cylinder connecting 
$q_{\ups_k}^{-1}(S_{h,\de}^*)$ and $q_{\ups_k}^{-1}(S_{\io_h,\de}^*)$.
Let \hbox{$\ep\!>\!0$} be any small number.
Choose small $\de>0$ such that $u_h(B_{b^*,h}(\i,2\de))$ and 
$u_{\io_h}^*(B_{b^*,\io_h^*}(x_h^*,2\de))$ lie in 
$B_{b^*}(u_h^*(\i),\ep)$.
Then we can write 
$$u_{b^*}^*(z)=\exp_{b^*,u_{b^*}^*(x_h^*)}\bar{u}_{b^*}(z),~~
|\bar{u}_{b^*}(z)|< \inj~g_{V,b^*};\quad
\bar{\xi}_k'(z)\equiv \Pi_{b^*,\bar{u}_b(z)}^{-1}\xi_k'(z)$$
for $z\!\in\! B_{b^*,h}(\i,\de)\cup B_{b^*,\io_h^*}(x_h^*,\de).$
Similarly, put
$$\bar{\xi}_h^*(z)=\Pi_{b^*,\bar{u}_{b^*}(z)}^{-1}\xi_h^*(z)
\hbox{~~and~~}
\bar{\xi}_{\io_h}^*(z)=\Pi_{b^*,\bar{u}_{b^*}(z)}^{-1}\xi_{\io_h}^*(z)$$
for $z$ in $B_{b^*,h}(\i,\de)$ and in $B_{b^*,i_h}(x_h^*,\de)$, respectively.
We can also assume that $\de$ is so small that
$\l|\bar{\xi}_h^*-\xi_h^*(\i)\r|_{b^*}$ and 
$\l|\bar{\xi}_{\io_h}^*-\xi_{i_h}^*(x_h^*)\r|_{b^*}$ do not exceed 
$\ep$ on $B_{b^*,h}(\i,\de)$ and on $B_{b^*,\io_h}(x_h^*,\de)$, respectively. 
Choose large $k^*$ such that all $k>k^*$
$$\|\xi^*-\xi_k'\|_{C^0(S_{h,\de}^*\cup S_{\io_h,\de}^*)}\le\ep.$$ 
It can be assumed that $u_k(A_{h,2\de,k})$ lies in $B_{b^*}(u^*(x_h^*);2\ep)$ 
for $k>k^*$.    Thus, we can write
$$u_k(z)=\exp_{b^*,u(x_h^*)}\bar{u}_k(z),~~
|\bar{u}_k(z)|_{b^*}< \inj~g_{V,b^*};\quad
\bar{\xi}_k(z)\equiv \Pi_{b^*,\bar{u}_k(z)}^{-1}\xi_k(z)$$
if $z\!\in\! A_{h,\de,k}.$
Pick points $z_1$ and $z_2$, one on each component of the boundary
of $A_{h,\de,k}$.    Then
\begin{equation}\label{inver1_e6}\begin{split}
\big|\xi_h^*(\i)-\xi_{\io_h}^*(x_h^*)\big|_{b^*}
&\le 2\big(\ep+ \big|\bar{\xi}_h^*(q_{\ups_k}(z_1))-
\bar{\xi}_{\io_h}^*(q_{\ups_k}(z_2))\big|_{b^*}\big)\\
&\le 4\big(\ep+ \big|\bar{\xi}_{k,h}'(q_{\ups_k}(z_1))-
\bar{\xi}_{k,i_h}'(q_{\ups_k}(z_2))\big|_{b^*}\big)\\
&\le C\big(\ep+\big|\bar{\xi}_k(z_1)-\bar{\xi}_k(z_2)\big|_{b^*}
+\|\ze_k\|_{b^*,C^0(S_{h,\de}^*\cup S_{i_h,\de}^*)}
\|\bar{\xi}_k\|_{b^*,C^0(A_{h,\de,k})}\big).
\end{split}\end{equation}
Since $A_{h,\de,k}$ is uniformly equivalent to the union of
two annuli with the larger radius bounded above by $\de$
and the smaller radius less than half of the larger, 
by Lemma~\ref{an-c0_plane_lmm} in~\cite{Z1} and
Holder's inequality,
\begin{equation}\label{inver1_e7}
\big|\bar{\xi}_k(z_1)-\bar{\xi}_k(z_2)\big|_{b^*}
\le C\big|\bar{\xi}_k(z_1)-\bar{\xi}_k(z_2)\big|_{b_k}
\le C'\de^{\frac{2(p-2)}{p}}\|d\bar{\xi}_k\|_{\ups_k,L^p(A_{h,\de,k})}.
\end{equation}
By Corollary~\ref{an-pt_l1c} in~\cite{Z1}
and Proposition~\ref{elliptic_bound},
\begin{equation}\label{inver1_e8}
\|d\bar{\xi}_k\|_{\ups_k,L^p(A_{h,\de,k})}\le 
\|\xi_k\|_{\ups_k,p,1}+\|du_{\ups_k}\|_{\ups_k,p}\|\xi_k\|_{\ups_k,C^0}
\le C.
\end{equation}
Combining equations \e_ref{inver1_e6}-\e_ref{inver1_e8}, we obtain
\begin{equation}\label{inver1_e9}
\big|\xi_h^*(\i)-\xi_{\io_h}^*(x_h^*)\big|_{b^*}\le 
C\big(\ep+\de^{\frac{2(p-2)}{p}}+
\|\ze_k\|_{b^*,C^0(S_{h,\de}^*\cup S_{\io_h,\de}^*)}\big).
\end{equation}
Since the last term in \e_ref{inver1_e9} tends to zero 
as $k\lra\i$ and $\ep$ and $\de$
can be chosen to be arbitrarily small, it follows 
$\xi_h^*(\i)=\xi_{\io_h}^*(x_h^*)$.

\begin{prp}
\label{inverse_p}
For any simple bubble type ${\cal T}$,
there exist $C,\de\!\in\! C^{\i}({\cal M}_{\cal T}^{(0)};\Bbb{R})$ 
such that for all $\ups\!\in\! F^{(0)}{\cal T}_{\de}$
if ${\cal T}$ is regular and $\ups\!\in\! F^{(\eset)}{\cal T}_{\de}$
if ${\cal T}$ is semiregular,
$$\|\xi\|_{\ups,p,1}\le C_p(b_{\ups})\|D_{\ups}\xi\|_{\ups,p}
\qquad
\forall\xi\!\in\!\Ga_+(\ups)
\quad\hbox{and}\quad\forall\xi\!\in\!\tilde{\Ga}_+(\ups).$$
\end{prp}

\noindent
{\it Proof:} If not, we can choose a sequence 
$\ups_k\!\in\! F^{(0)}{\cal T}_{\de}$, converging to some 
$b\!\in\! {\cal M}_{\cal T}^{(0)}$ and vector fields
\hbox{$\xi_k\!\in\!\Ga_+(\ups_k)$} 
(or $\xi_k\!\in\!\tilde{\Ga}_+(\ups_k)$) such that 
$\|\xi_k\|_{\ups_k,p,1}=1$, while $\|D_{\ups_k}\xi\|_{\ups_k,p}\lra 0$.
If $\xi_k\!\in\!\Ga_+(\ups_k)$, note that 
$\{\Ga_-(\ups_k)\}$ $C^0$-converges to \hbox{$V\!\equiv\! \Ga_-(b)$}.
If $\xi_k\!\in\!\tilde{\Ga}_+(\ups_k)$,
by Definition~\ref{obs_setup_dfn1}, 
a subsequence of $\{\tilde{\Ga}_-(\ups_k)\}$
$C^0$-converges to a subspace $V\!\subset\! L^p_1(b)$
such that $\pi_{b,-}\!:V\lra\Ga_-(b)$ is an isomorphism.
In either case, by the first statement of Lemma~\ref{inver1}, 
a subsequence of $\{\xi_k\}$ 
$C^0$-converges to a vector field $\xi^*\!\in\!\Ga_-(b)$.
By the second statement of Lemma~\ref{inver1}, 
$\xi^*$ must be orthogonal~$V$,
since $\xi_k\!\in\!\Ga_+(\ups_k)$ (or $\xi_k\!\in\!\tilde{\Ga}_+(\ups_k)$).
Thus, $\xi^*\!=\!0$. 
On the other hand, by Proposition~\ref{elliptic_bound},
there exists $\ep>0$ such that $\|\xi_k\|_{\ups_k,p}\ge\ep$ 
for all $k$ sufficiently large. 
However, by Lemma~\ref{inver1}, $\|\xi_k\|_{\ups_k,C^0}\lra 0$,
which is a contradiction.

\end{document}